\documentclass[11pt]{article}
\usepackage{authblk,amssymb,fullpage,latexsym,amsfonts,amsmath,amsthm}
\usepackage{color}
\usepackage{enumerate}
\usepackage{algorithm}
\usepackage{algpseudocode}
\usepackage{bm}
\usepackage{subfigure} 
\usepackage{mathtools}
\usepackage[titletoc,title]{appendix}
\usepackage[T1]{fontenc}
\usepackage{dsfont}
\usepackage{natbib}
\usepackage{url}
\usepackage{float}
\usepackage{placeins}

\usepackage[font=small,labelfont=bf]{caption}
\usepackage{authblk} 
\usepackage{hyperref} 
\usepackage{cleveref}
\hypersetup{
colorlinks = true,
linkcolor = blue, 
filecolor = magenta,
urlcolor = blue, 
citecolor = blue 
}

\newcommand{\mP}{\mathbb{P}}
\newcommand{\mE}{\mathbb{E}}
\newcommand{\mV}{\mathrm{Var}}

\newcommand{\bpi}{\boldsymbol{\pi}}
\newcommand{\blambda}{\boldsymbol{\lambda}}
\newcommand{\bmu}{\boldsymbol{\mu}}
\newcommand{\bx}{\boldsymbol{x}}

\newcommand{\given}{\,|\,} 
\newcommand{\iid}{\overset{\mathrm{i.i.d.}}{\sim}}
\newcommand{\dd}{\mathrm{d}}

\newcommand{\convP}{\overset{p}{\longrightarrow}}
\newcommand{\convD}{\overset{d}{\longrightarrow}}

\DeclarePairedDelimiter\ceil{\lceil}{\rceil}
\DeclarePairedDelimiter\floor{\lfloor}{\rfloor}
\DeclareMathOperator{\TV}{TV}
\DeclareMathOperator{\supp}{supp}

\newcommand{\mathttt}{\mathrm}

\newcommand{\blue}[1]{#1}

\newtheorem{definition}{Definition}
\newtheorem{theorem}{Theorem}
\newtheorem{lemma}{Lemma}

\newtheorem{proposition}{Proposition}

\newtheorem{remark}{Remark}
\newtheorem{example}{Example}

\newtheorem{procedure}{Procedure}

\let\hat\widehat

\begin{document}

\begingroup
\renewcommand{\thefootnote}{\fnsymbol{footnote}}
\begin{center}
	\textbf{\LARGE \blue{Differentially Private Permutation Tests}}\footnotemark[1]
	
	\vspace*{.2in}
	
	{\large Ilmun Kim$^{\dagger,\diamond}$  \quad Antonin Schrab$^{\ddagger,\diamond}$}
	
	\vspace*{.1in}
	
	{
	\begin{tabular}{c}
		$^{\dagger}$Department of Mathematical Sciences, Korea Advanced Institute of Science and Technology \\
		$^{\ddagger}$Department of Computer Science and Technology, University of Cambridge \\[.3em]
		$^\diamond$Both authors contributed equally to this work
	\end{tabular}
	}

	\vspace*{.1in}
	
	\today
	
	\vspace*{.1in}
\end{center}
\footnotetext[1]{This arXiv version incorporates selected revisions from the journal version, including the extension to R\'{e}nyi differential privacy and the expanded empirical study. The high-privacy minimax lower bound in \Cref{Theorem: High-privacy minimax separation over L2}, together with its proof, is new to this version and does not appear in the journal version.}
\endgroup

\begin{abstract}
	Recent years have witnessed growing concerns about the privacy of sensitive data. In response to these concerns, differential privacy has emerged as a rigorous framework for privacy protection, gaining widespread recognition in both academic and industrial circles. While substantial progress has been made in private data analysis, existing methods often suffer from impracticality or a significant loss of statistical efficiency. This paper aims to alleviate these concerns in the context of hypothesis testing by introducing differentially private permutation tests. The proposed framework extends classical non-private permutation tests to private settings, maintaining both finite-sample validity and differential privacy in a rigorous manner. The power of the proposed test depends on the choice of a test statistic, and we establish general conditions for consistency and non-asymptotic uniform power. To demonstrate the utility and practicality of our framework, we focus on reproducing kernel-based test statistics and introduce differentially private kernel tests for two-sample and independence testing:~$\mathttt{dpMMD}$ and $\mathttt{dpHSIC}$. The proposed kernel tests are straightforward to implement, applicable to various types of data, and attain minimax optimal power across different privacy regimes. Our empirical evaluations further highlight their competitive power under various synthetic and real-world scenarios, emphasizing their practical value. The code is publicly available to facilitate the implementation of our framework.
\end{abstract}

\section{Introduction}
Ensuring the privacy of sensitive data has become a critical concern in modern data analysis. As organizations collect and analyze vast amounts of personal information, safeguarding individual privacy has emerged as a crucial ethical and legal imperative. In response to these challenges, differential privacy (DP), introduced by \cite{dwork2006calibrating}, has emerged as a rigorous framework for addressing privacy concerns, and has gained widespread recognition not only in academia but also in industry. For instance, major industry players such as Apple~\citep{apple2017}, Google~\citep{erlingsson2014rappor} and Microsoft~\citep{ding2017collecting} have embraced differential privacy as a robust definition of privacy. This growing trend has sparked a recent surge of research in statistics and related fields aimed at integrating differential privacy and its variants~\citep{dwork2006our,bun2016concentrated,mironov2017renyi,dong2022gaussian} into data analysis and developing privacy-preserving methodologies. In this line of work, a major challenge is to strike a balance between privacy guarantees and statistical efficiency. Notably, a high privacy guarantee requires substantial data perturbation, which in turn degrades statistical performance. Conversely, releasing less perturbed data can improve statistical efficiency but at the expense of reduced privacy guarantees. Therefore, balancing this trade-off between privacy and efficiency has been a central topic in the existing literature~\citep[\emph{e.g.},][]{duchi2018minimax,cai2021cost,kamath2020primer}. 

Broadly, there are two major statistical problems tackled under privacy constraints:~estimation and hypothesis testing~\citep[][for a recent review]{kamath2020primer}. This paper focuses on the latter problem, which requires access to the null distribution of a test statistic in order to effectively calibrate a test statistic. Analyzing the distribution of a test statistic becomes particularly more challenging in private settings due to additional random sources arising from privacy mechanisms. As we review in \Cref{Section: Related Work}, substantial efforts have been made to address private testing problems. These efforts involve adapting classical hypothesis tests to private settings or developing new testing procedures that achieve an optimal balance between privacy and statistical power. 

Despite the significant progress made over the last decade, there are still several areas where further improvements can be made. One such area includes the reliance on asymptotic methods for determining the critical value of a test statistic. The practical quality of this asymptotic approach depends on the convergence rate of a privatized statistic to the limiting distribution. This convergence rate is often slow in private settings, and more importantly, the limiting distribution may vary depending on the delicate interplay between privacy and other parameters. This issue puts practitioners in a bind as it is unclear which limiting distribution should be considered a priori. All these concerns lead to the unreliability of asymptotic private tests in real-world applications.

Another area of concern is the limited practicality of existing methods. Many private statistical tests are designed specifically for discrete data, not directly applicable to handling continuous or mixed-type data. Moreover, existing methods often rely on unspecified constants and heuristics, making them less user-friendly and potentially undermining their reliability. We also point out that the majority of research has concentrated on theoretical aspects of private testing, and only a handful of papers are equipped with thorough empirical evaluations and open-source code. 

In this work, we aim to tackle the aforementioned concerns by introducing differentially private permutation tests. The primary goal is to extend classical non-private permutation tests to differentially private settings, applicable to any test statistic with finite global sensitivity. The proposed private permutation test inherits the finite-sample validity of the classical permutation test under the exchangeability condition, while ensuring differential privacy. The power of the proposed test depends on the choice of a test statistic, and we establish sufficient conditions for consistency and non-asymptotic uniform power. To demonstrate the effectiveness of our framework, we focus on the two-sample and independence testing problems and propose differentially private versions of the maximum mean discrepancy (MMD) and Hilbert--Schmidt independence criterion (HSIC), which we call ``dpMMD'' and ``dpHSIC'', respectively. On the theoretical side, we prove minimax optimality of the $\mathttt{dpMMD}$ and $\mathttt{dpHSIC}$ tests across all privacy regimes in terms of kernel metrics. On the empirical side, we showcase the competitive power of the proposed tests across various practical scenarios. The code that implements our methods is publicly available at \url{https://github.com/antoninschrab/dpkernel} to allow practitioners to build on our findings.

\subsection{Related Work} \label{Section: Related Work}

In recent years, there has been a growing body of research on hypothesis testing problems under privacy constraints. Since the early work by \cite{vu2009differential} and \cite{fienberg2011privacy}, numerous attempts have been made to extend classical non-private tests to their private counterparts. Examples include ANOVA~\citep{campbell2018differentially,swanberg2019improved}, likelihood ratio tests~\citep{canonne2019structure}, tests for regression coefficients~\citep{sheffet2017differentially,alabi2022hypothesis}, rank or sign-based nonparametric tests~\citep{task2016differentially,couch2019differentially}, conditional independence tests~\citep{kalemaj2023differentially} and $\chi^2$-tests~\citep{fienberg2011privacy,wang2015revisiting,gaboardi2016differentially,rogers2017new,kakizaki2017differentially,friedberg2023privacy}. While most of the aforementioned work focuses on asymptotic settings where the sample size goes to infinity, a recent line of work within computer science has placed greater emphasis on finite-sample analysis. In particular, \cite{cai2017priv} propose a two-step algorithm for identity testing for discrete distributions, and study the sample complexity under DP. The work of \cite{acharya2018differentially} explores both identity (goodness-of-fit) testing and closeness (two-sample) testing in finite-sample settings, and improves the upper bound result by \cite{cai2017priv} and \cite{aliakbarpour2017differentially} for identity testing. \cite{aliakbarpour2019private} privatize the non-private test proposed by \cite{diakonikolas2016new}, and investigate the sample complexity for closeness testing and independence testing. \blue{We also highlight the uniformly most powerful tests for binomial data developed by \citet{awan2020}.} In line with these advancements, our work develops private permutation tests, and studies their non-asymptotic performance under DP settings. 

Despite the extensive body of literature, the majority of research has focused on private tests designed for discrete or bounded data. There are a few notable exceptions that have explored other data types. For example, \cite{canonne2020private} and \cite{narayanan2022private} have investigated goodness-of-fit testing for high-dimensional Gaussian distributions. In addition, \cite{raj2020differentially} have proposed private two-sample tests based on finite dimensional approximations of kernel mean embeddings. The flexibility offered by kernel methods enables these tests to handle a wide variety of data types. However, their tests are asymptotic in nature, which may introduce reliability concerns when working with small sample sizes. Moreover, their analysis requires the number of features to be fixed. Such requirement potentially limits the power of the test when dealing with alternatives which are not well-represented by these fixed numbers of features.  

Another line of work aims to develop a generic way to create private tests from non-private ones. The subsample-and-aggregate idea~\citep{nissim2007smooth} has emerged as a useful tool for this purpose. In particular, it offers a strategy to convert non-private sample complexity results into private ones in a black-box manner as pointed out by \cite{cai2017priv,canonne2019structure,canonne2020private}. Recent studies by \cite{pena2022differentially} and \cite{kazan2023test} have focused specifically on the practical implementation of the subsample-and-aggregate approach. However, it is worth mentioning that this generality typically comes at the cost of suboptimal power, and often fails to recover the optimal sample complexity \citep{canonne2019structure,canonne2020private}. Moreover, the performance of this subsample-and-aggregate approach is sensitive to the number of subsamples, and determining the optimal value of this parameter remains an open problem.

Beyond global differential privacy, there has been a substantial amount of work on hypothesis testing under local differential privacy. Some of the notable works include  \cite{liao2017hypothesis,gaboardi2018local,sheffet2018locally,acharya2019test,berrett2020locally,dubois2021goodness,lam2022minimax} and see the references therein. Local differential privacy requires data perturbation at the individual level, proving particularly useful in settings where data providers lack trust in data analysts. This individual-wise approach demands different analyses than global differential privacy, and the results under global and local differential privacy are not directly comparable.

Our work is also related to recent advances in kernel-based minimax testing~\citep{li2019optimality,albert2019adaptive,schrab2021mmd,kim2020minimax}. Specifically, we extend the non-private minimax testing rates established in this line of work to private counterparts. To achieve this, we leverage the techniques therein, such as the two moments method and exponential inequalities for permuted statistics in \cite{kim2020minimax}, and adapt them to private settings.

\subsection{An Overview of Our Results}
The main contributions of this work are summarized below.
\begin{itemize}
	\item \textbf{DP Permutation Tests (\Cref{Section: Differentially Private Permutation Tests}).} We introduce differentially private permutation tests in Algorithm~\ref{Algorithm: DP permutation test}, and establish their theoretical properties. A naive way of extending the classical permutation tests to private settings is to first make the original test statistic and its permuted counterparts differentially private, and then carry out the permutation test based on these individually privatized statistics. However, this naive approach results in an unnecessary power loss by adding more noise as the number of permutations increases. 
	The proposed framework addresses this issue by utilizing the quantile representation of a permutation test. This strategy leads to a substantial power gain over the naive approach, while being finite-sample valid. We present sufficient conditions for pointwise consistency~(\Cref{Theorem: Pointwise Consistency}) and non-asymptotic uniform power~(\Cref{Theorem: General uniform power condition}) of the proposed tests. The latter uniform power condition can be regarded as an extension of the two moments method~\citep{kim2020minimax} to private settings. 	 
	\blue{We further introduce a Gaussian-mechanism variant satisfying R\'{e}nyi differential privacy and show that it retains finite-sample validity, privacy, consistency, and a general uniform power guarantee (\Cref{Section: Renyi Differentially Private Permutation Tests}).}
	\item \textbf{DP Kernel Tests (\Cref{Section: Applications to Kernel-Based Inference}).} We showcase the versatility of our framework by applying it to two specific tasks: differentially private two-sample and independence testing based on reproducing kernel-based test statistics. We consider the plug-in estimators of the MMD and HSIC, and privatize them through the proposed method employing the standard Laplace mechanism. 
	The practical performance of the resulting differentially private kernel tests heavily depends on the global sensitivity used in the Laplace mechanism. To boost the empirical performance, we put significant effort into establishing sharp upper bounds for the global sensitivity of the plug-in estimators of MMD and HSIC, as well as matching lower bounds for popular kernels (\Cref{Lemma: Sensitivity of MMD} and \Cref{Lemma: Sensitivity of HSIC}). We then establish key properties of the proposed kernel tests, including non-asymptotic validity and consistency against any fixed alternatives in \Cref{Theorem: Properties of dpMMD} and \Cref{Theorem: Properties of dpHSIC}.
\item \textbf{Uniform Power and Optimality (\Cref{Section: Uniform separation and optimality}).} We characterize the trade-off between differential privacy and statistical power through the lens of minimax analysis. To this end, we analyze the minimum separation required for the differentially private MMD test to achieve significant power in terms of the MMD metric. We derive an upper bound on this minimum separation in \Cref{Theorem: Uniform separation for MMD}, and a lower bound in \Cref{Theorem: Minimax separation in MMD}, which matches in all relevant parameters including testing error rates and privacy levels. Our minimax results suggest that there is an unavoidable loss of power when the privacy parameter is smaller than a certain threshold (\emph{i.e.},~high privacy). On the other hand, the privacy guarantee comes for free in terms of separation rate when the privacy parameter exceeds this threshold (\emph{i.e.},~low privacy). We also derive the minimum separation in terms of the $L_2$ metric in \Cref{Theorem: Minimax Separation over L2} that extends the prior work~\citep{li2019optimality,schrab2021mmd} on non-private minimax testing to private settings. \blue{Our new lower bound in \Cref{Theorem: High-privacy minimax separation over L2} matches the high-privacy upper bound for $s\geq d/2$, proving that $\mathttt{dpMMD}$ is minimax rate-optimal in that regime.} In \Cref{Section: Separation in HSIC metric} and \Cref{Section: Separation in L2 for dpHSIC}, we present analogous findings for the differentially private HSIC test, which closely resemble the results obtained for the MMD test.
	\item \textbf{Negative Results of U-statistics (\Cref{Section: Private Test based on the MMD U-statistic}).} We also derive perhaps unexpected negative results of U-statistics in private settings. U-statistics have played an important role, arguably more popular than V-statistics, in deriving non-private minimax rates for various testing problems~\citep{li2019optimality,albert2019adaptive,kim2020minimax,berrett2021optimal,schrab2021mmd}.
	Given this trend, it is natural to consider U-statistics as an initial building block for obtaining private minimax rates. However, it turns out that the U-statistics suffer from higher global sensitivity than the corresponding V-statistics, requiring a higher level of noise to effectively privatize the resulting procedure. We formalize this observation in the context of kernel testing and show that the private tests based on U-statistics have sub-optimal power in high privacy regimes. This negative result naturally justifies our approach based on the V-statistics (equivalently, plug-in estimators) of the MMD and HSIC. 
	\item \textbf{Empirical Validation (\Cref{Section: Simulations}).} A significant portion of the prior work on differentially private testing has focused on theoretical aspects, often lacking practical value. On the other hand, practical approaches to differentially private testing frequently rely on heuristics without proper theoretical validation. Our work helps bridge the gap between theory and practice by balancing both aspects. In particular, we highlight that our method is simple to use and comes with strong theoretical guarantees, as we demonstrate throughout the paper. The empirical results in \Cref{Section: Simulations} and \Cref{Section: Additional Simulations} also illustrate the competitive performance of the proposed method across diverse scenarios, underscoring its practical value. \blue{The expanded study includes speech measurements from patients with Parkinson's disease, robustness to unequal sample sizes, and the TOT--SARRM baseline.}  
\end{itemize}

We further make contributions by presenting asymptotic distributions of privatized kernel statistics~(\Cref{Section: Limiting null distributions}), general consistency results for resampling-based tests~(\Cref{Section: General Pointwise Consistency Result}) as well as other technical innovations~(\Cref{Section: Technical Lemmas}). We also introduce private kernel tests obtained through the subsample-and-aggregate idea~(\Cref{Section: Alternative private tests}). \blue{In addition, we discuss computationally efficient variants and the interaction between kernel choice and privacy (\Cref{Section: Discussion}).} Due to space constraints, we relegate these additional results to the appendix. 

\subsection{Organization}
The rest of this paper is organized as follows. We begin in \Cref{Section: Background} by providing a brief overview of the fundamental concepts of differential privacy. \Cref{Section: Differentially Private Permutation Tests} presents our main proposal, namely the differentially private permutation test, and investigates its finite-sample validity and consistency in power. Moving forward to \Cref{Section: Applications to Kernel-Based Inference}, we apply our proposed permutation framework to specific scenarios, focusing on differentially private kernel testing. In particular, we explore the privatization of kernel MMD and HSIC tests, and delve into minimum separation rates of the resulting tests in \Cref{Section: Uniform separation and optimality}. To validate our theoretical findings, \Cref{Section: Simulations} presents empirical evaluations of the proposed algorithms, comparing their performance with existing methods. Finally, we conclude in \Cref{Section: Discussion} with a discussion outlining potential directions for future work. All the proofs and additional results are deferred to the appendix.

\subsection{Notation} Given two datasets $\mathcal{X}_n \coloneqq (X_1,\ldots,X_{n})$ and $\tilde{\mathcal{X}}_n\coloneqq (\tilde X_1,\ldots,\tilde X_{n})$, we denote the Hamming distance between $\mathcal{X}_n$ and $\tilde{\mathcal{X}}_n$ by $d_{\mathrm{ham}}(\mathcal{X}_n, \tilde{\mathcal{X}}_n) := \sum_{i=1}^n \mathds{1}(X_i \neq \tilde X_i)$. 
For two sequences of real numbers $a_n,b_n$, we write $a_n \lesssim b_n$ (and similarly $a_n \gtrsim b_n$) if there exists some positive constant $C>0$ independent of $n$ such that $a_n \leq C b_n$ for all $n \geq 1$. We also write $a_n \asymp b_n$ if $a_n \lesssim b_n$ and $b_n \lesssim a_n$. 
For $x \in \mathbb{R}$, $\floor{x}$ denotes the largest integer smaller than or equal to $x$.
For a natural number $k \in \mathbb{N}$, we use $[k]$ to denote the set $\{1,\ldots,k\}$. We let $\boldsymbol{\Pi}_n$ denote the set of all permutations of $[n]$. For a continuous function $f : \mathbb{R}^d \mapsto \mathbb{R}$, the $L_2$ and $L_\infty$ norms of $f$ are given as $\|f\|_{L_2} = \{\int_{\mathbb{R}^d} f^2(\bx) \mathrm{d}\bx \}^{1/2}$ and $\|f\|_{L_\infty} = \sup_{\bx \in \mathbb{R}^d}|f(\bx)|$, respectively. We say $X \sim \mathsf{Laplace}(0,1)$ if $X$ follows a Laplace distribution with location and scale parameters $(0,1)$. We often denote 
\begin{align} \label{Eq: definition of xi}
	\xi_{\varepsilon,\delta} \coloneqq \varepsilon + \log\bigl(1 / (1-\delta) \bigr)
\end{align}
to simplify the notation in various places.

\section{Background:~Differential Privacy} \label{Section: Background}

This section presents a brief overview of the basic concepts and properties regarding differential privacy. For a comprehensive treatment, we refer the readers to \cite{dwork2014algorithmic}. In our work, we adhere to the definition of differential privacy~\citep[][page 25]{dwork2014algorithmic}, allowing for the inclusion of additional auxiliary variables. This extended definition requires that the standard differential privacy condition holds for every possible value of the auxiliary variable. In our permutation testing framework, we treat random permutations as the auxiliary variables independent of the dataset.

\begin{definition}[Differential Privacy]
	\label{Definition: DP}
	Consider a randomized algorithm $\mathcal{A}$, which takes as input a dataset $\mathcal{X}_n $ and an additional auxiliary variable $w \in \mathcal{W}$.
	For $\varepsilon>0$ and $\delta\in[0,1)$, the algorithm $\mathcal{A}$ is said to be $(\varepsilon,\delta)$-differentially private if for (i)~all $S \in \mathrm{range}(\mathcal{A})$, (ii)~all $w \in \mathcal{W}$ and (iii)~all two datasets $\mathcal{X}_n$ and $\tilde{\mathcal{X}}_n$ with $d_{\mathrm{ham}}(\mathcal{X}_n, \tilde{\mathcal{X}}_n) \leq 1$, the following inequality holds:
	\begin{align*}
		\mP \bigl( \mathcal{A}(\mathcal{X}_n;w) \in S \given \mathcal{X}_n, w \bigr) ~\leq~ e^{\varepsilon} \mP \bigl( \mathcal{A}(\tilde{\mathcal{X}}_n;w) \in S \given \tilde{\mathcal{X}}_n, w \bigr) + \delta.
	\end{align*}
\end{definition}
We note that $(\varepsilon,0)$-differential privacy is often simply referred to as $\varepsilon$-DP or pure-DP. On the other hand, $(\varepsilon,\delta)$-differential privacy with $\delta \in (0,1)$ is referred to as approximate-DP, considered as a relaxation of pure-DP. As mentioned by \citet[][page 18]{dwork2014algorithmic}, employing a large value of $\delta$ may lead to serious privacy breaches, potentially exposing the complete information of a small number of individuals with a non-trivial probability. Hence, it is generally desirable to choose small values of $\delta$ such as $\delta \lesssim \varepsilon^2 n^{-1}$. Nevertheless, our interest lies in exploring entire privacy regimes, and developing comprehensive results applicable in a variety of settings. Consequently, we do not place restrictions on privacy parameters other than $\varepsilon >0$ and $\delta \in [0,1)$.

We collect several fundamental properties of differential privacy that are useful in our contexts. The first property is called post-processing~\citep[][Proposition 2.1]{dwork2014algorithmic}, which asserts that any arbitrary post-processing applied to the outcome of a differentially private algorithm preserves the same level of privacy. 

\begin{lemma}[Post-Processing] \label{Lemma: Post-Processing}
	Suppose that an algorithm $\mathcal{A}$ is $(\varepsilon,\delta)$-differentially private. Then for an arbitrary randomized function $f$, the composition $f  \circ \mathcal{A}$ also preserves $(\varepsilon,\delta)$-differentially privacy. 
\end{lemma}

Another important property of differential privacy, called the composition theorem \citep[Theorem 3.16]{dwork2014algorithmic}, presents the overall privacy guarantee for a composition of multiple DP mechanisms. 

\begin{lemma}[Composition]
	\label{Lemma: Composition Theorem}
Suppose that each algorithm $\mathcal{A}_i$ is $(\varepsilon_i,\delta_i)$-differentially private for $i \in [m]$. Then, the composed algorithm $\mathcal{A}_{1:m}$ defined as $\mathcal{A}_{1:m} \coloneqq (\mathcal{A}_1,\dots,\mathcal{A}_m)$ is $(\sum_{i=1}^m \varepsilon_i,\sum_{i=1}^m\delta_i)$-differentially private.
\end{lemma}

The definition of $(\varepsilon,\delta)$-DP immediately leads to the following group property~\citep[][Lemma 19]{acharya2021differentially}, which plays an important role in constructing minimax lower bounds under DP.

\begin{lemma}[Group Privacy] \label{Lemma: Group privacy}
	Suppose that an algorithm $\mathcal{A}$ is $(\varepsilon,\delta)$-differentially private. Then for (i)~all $S \in \mathrm{range}(\mathcal{A})$, (ii)~all $w \in \mathcal{W}$, and (iii)~all two datasets $\mathcal{X}_n$ and $\tilde{\mathcal{X}}_n$ with $d_{\mathrm{ham}}(\mathcal{X}_n, \tilde{\mathcal{X}}_n) \leq m$, the following inequality holds:
	\begin{align} \label{Eq: implication of DP definition}
		\mP \bigl( \mathcal{A}(\mathcal{X}_n;w) \in S \given \mathcal{X}_n, w \bigr) ~\leq~ e^{m\varepsilon} \mP \bigl( \mathcal{A}(\tilde{\mathcal{X}}_n;w) \in S \given \tilde{\mathcal{X}}_n, w \bigr) + m e^{(m-1)\varepsilon}\delta.
	\end{align}
\end{lemma}

Several mechanisms have been developed to safeguard differential privacy, with the Laplace mechanism~\citep{dwork2006calibrating} standing out as one of the most commonly used approaches. To formally state the Laplace mechanism, we first describe the global sensitivity, which is a keystone in the differential privacy framework. We stress that the definition presented below allows us to take into account an additional auxiliary variable $w$, and thus it is more general than the definition commonly encountered in the DP literature, such as in \citet{dwork2006calibrating}.

\begin{definition}[Global $\ell_p$-Sensitivity]
	\label{Definition: lp sensitivity} 
	Consider a function $f$ taking as input a dataset $\mathcal{X}_n$ and an additional auxiliary variable $w\in\mathcal{W}$, and assume that the output of $f$ lies in $\mathbb{R}^r$. For $p \geq 1$, the global $\ell_p$-sensitivity of $f$ is defined as
	\begin{equation*}
		\Delta_{\!f}^p ~\coloneqq~ \sup_{w \in \mathcal W } \sup_{\substack{\mathcal{X}_n, \tilde{\mathcal{X}}_n:\\d_{\mathrm{ham}}(\mathcal{X}_n, \tilde{\mathcal{X}}_n) \leq 1}} \big\|f(\mathcal{X}_n; w) - f(\tilde{\mathcal{X}}_n; w)\big\|_p,
	\end{equation*}
	where $\|x\|_p$ denotes the $\ell_p$ norm of a vector $x$. 
\end{definition}
The Laplace mechanism works with the $\ell_1$-sensitivity, which determines the scaling factor of the Laplace noise injected into the outputs of the function. Formally, the Laplace mechanism is given as follows.

\begin{definition}[Laplace Mechanism]
	\label{Definition: Laplace Mechanism} 
	Consider a function $f$ with the $\ell_1$-sensitivity $\Delta_f^1$ described in Definition~\ref{Definition: lp sensitivity}. For a given privacy parameter $\xi >0$, the Laplace mechanism is defined as the random function:
	\begin{equation*}
		\mathcal M_{\!f}^{\xi}(\mathcal{X}_n;w) \coloneqq f(\mathcal{X}_n;w) + \frac{\Delta_{\!f}^1}{\xi}(\zeta_1,\dots,\zeta_r)^\top,
	\end{equation*}
	where $\zeta_1,\dots,\zeta_r \iid  \mathsf{Laplace}(0,1)$ generated independent of $\mathcal{X}_n$ and $w$.
\end{definition}
The privacy guarantee of the Laplace mechanism depends crucially on the choice of privacy parameter $\xi$. It is well-known that the Laplace mechanism is $(\varepsilon,0)$-DP when $\xi = \varepsilon$~\citep[][Theorem 3.6]{dwork2014algorithmic}. In general, \citet[Lemma 5]{acharya2018differentially} shows that any $(\varepsilon+\delta,0)$-DP algorithm is also $(\varepsilon,\delta)$-DP. While this strategy is effective for small values of $\delta$, it returns a suboptimal result when $\delta$ is close to one. Concretely, when $\delta$ approaches one, we are entering the non-private regime where adding noise is unnecessary. However, the Laplace mechanism with $\xi = \varepsilon+ \delta$ injects a non-negligible amount of noise to the algorithm. The refined calibration result proposed by \cite{holohan2014differential} avoids such issue, proving that $\mathcal M_{\!f}^{\xi_{\varepsilon,\delta}}$ with $\xi_{\varepsilon,\delta} =  \varepsilon+\log(1/(1-\delta))$ is also $(\varepsilon,\delta)$-DP. We record this guarantee in the following lemma. 

\begin{lemma}[Differential Privacy of Laplace Mechanism]
	\label{Lemma: Laplace Mechanism}
	Let $\varepsilon>0$ and $\delta\in[0,1)$. The Laplace mechanism $\mathcal M_{\!f}^{\xi_{\varepsilon,\delta}}$ in Definition~\ref{Definition: Laplace Mechanism} with $\xi_{\varepsilon,\delta} = \varepsilon + \log\bigl(1/(1-\delta)\bigr)$ is $(\varepsilon,\delta)$-differentially private. 
\end{lemma}
It is worth noting that \cite{holohan2014differential} consider typical differential privacy without considering an auxiliary variable $w$. Nevertheless, the same proof can be applied to differential privacy involving auxiliary variables, provided that we consider the global sensitivity holding uniformly over $w \in \mathcal{W}$ as in Definition~\ref{Definition: lp sensitivity}.

\begin{remark}[Gaussian Mechanism] \normalfont
	For $(\varepsilon,\delta)$-DP, one can also consider the Gaussian mechanism with the $\ell_2$-sensitivity, another common method for preserving privacy~\citep{dwork2006calibrating,dwork2014algorithmic}. The Gaussian mechanism can be beneficial over the Laplace mechanism when the $\ell_2$-sensitivity is significantly smaller than the $\ell_1$-sensitivity. However, such benefit is not immediately clear when the outcome of $f$ is one-dimensional where the $\ell_p$ sensitivity remains the same for any $p \geq 1$. As our framework is mainly concerned with one-dimensional numeric outcomes, we simply focus on the Laplace mechanism and refer to the \emph{global $\ell_1$-sensitivity} as the \emph{global sensitivity} whenever it is clear from the context, and simply denote it by $\Delta_f$. 
\end{remark}

\section{Differentially Private Permutation Tests} \label{Section: Differentially Private Permutation Tests}
In this section, we introduce a general framework for constructing a differentially private permutation test. To begin, consider a class of distributions $\mathcal{P}$, which is formed by the union of two disjoint subclasses:~$\mathcal{P}_0$ and $\mathcal{P}_1$. Suppose that we observe a random sample $\mathcal{X}_n$ of size $n$ drawn from $P \in \mathcal{P}$. Given $\mathcal{X}_n$, our ultimate goal is to test whether $H_0: P \in \mathcal{P}_0$ or $H_1: P \in \mathcal{P}_1$, while preserving differential privacy. Consider a test statistic $T : \mathcal{X}_n \mapsto \mathbb{R}$, which is assumed to take a large value under the alternative hypothesis $H_1$. To build on the permutation principle~\citep[][Chapter 15.2]{lehmann2005testing}, we make the assumption that $\mathcal{X}_n$ is exchangeable under the null $H_0$. That is, for any permutation $\bpi \coloneqq (\pi_1,\ldots, \pi_n) \in \boldsymbol{\Pi}_n$, the joint distribution of $\mathcal{X}_n$ is the same as that of $\mathcal{X}_n^{\bpi} \coloneqq  (X_{\pi_1},\ldots,X_{\pi_n})$. Under the exchangeability assumption, the permutation test rejects the null when $T$ is significantly larger than the permuted counterparts. More formally, let $\bpi_1,\ldots,\bpi_B$ be i.i.d.~random permutations of $[n]$, and denote by $ T(\mathcal{X}_n^{\bpi_1}) ,\ldots, T(\mathcal{X}_n^{\bpi_B})$, the test statistics computed based on $\mathcal{X}_n^{\bpi_1},\ldots,\mathcal{X}_n^{\bpi_B}$, respectively. The permutation test then rejects the null hypothesis when the permutation $p$-value is less than or equal to significance level $\alpha$, \emph{i.e.},
\begin{align} \label{Eq: permutation p-value representation}
	\hat{p} \coloneqq  \frac{1}{B+1} \bigg\{ \sum_{i=1}^B \mathds{1}\bigl(T(\mathcal{X}_n^{\bpi_i}) \geq T(\mathcal{X}_n) \bigr) +1 \bigg\} \leq \alpha.
\end{align}
It is well-known that $\hat{p}$ is super-uniform, \emph{i.e.},~$\mP(\hat{p} \leq t) \leq t$ for all $t \in [0,1]$, under exchangeability of $\mathcal{X}_n$ (\emph{e.g.},~\Cref{Lemma: permutation p-value}). Therefore the permutation test $\mathds{1}(\hat{p} \leq \alpha)$ controls the type I error for any finite sample size $n$. Our aim is to privatize the permutation test under the DP constraint, while maintaining finite-sample validity and achieving competitive (potentially optimal) power. 

For notational convenience, we often write $T_0 = T(\mathcal{X}_n)$ and $T_i = T(\mathcal{X}_n^{\bpi_i})$ for $i \in [B]$, and set $\bpi_0 = (1,2,\ldots,n)$ in what follows.

\subsection{Proposed Privatization Method}
To describe the proposed method, suppose that the test statistic $T$ has the global sensitivity (\Cref{Definition: lp sensitivity}) with the permutation $\bpi$ as an auxiliary variable:
\begin{align} \label{Eq: global sensitivity}
	 \Delta_T \coloneqq \sup_{\bpi \in \boldsymbol{\Pi}_n} \sup_{\substack{\mathcal{X}_n, \tilde{\mathcal{X}}_n:\\d_{\mathrm{ham}}(\mathcal{X}_n, \tilde{\mathcal{X}}_n) \leq 1}} \big|T(\mathcal{X}_n^{\bpi}) - T(\tilde{\mathcal{X}}_n^{\bpi}) \big|.
\end{align}
Our tailored definition of global sensitivity to permutation tests above is stronger than the usual one since it measures the sensitivity over all possible permutations. However, this additional requirement is not overly restrictive as we demonstrate below for integral probability metrics. 

\begin{example}[Sensitivity of Integral Probability Metric] \label{Example: Sensitivity of Integral Probability Metric} \normalfont
	Consider a two-sample setting where we observe random variables $\mathcal{Y}_n = \{Y_1,\ldots,Y_n\}$ and $\mathcal{Z}_m = \{Z_1,\ldots,Z_m\}$, each supported on $\mathbb{S}$. Let $\mathcal{F}$ be a class of real-valued functions on $\mathbb{S}$. A plug-in estimator of the corresponding integral probability metric (IPM) is given as
	\begin{align} \label{Eq: plug-in IPM}
		T = \sup_{f \in \mathcal{F}} \bigg| \frac{1}{n} \sum_{i=1}^n f(Y_i) - \frac{1}{m}\sum_{i=1}^m f(Z_i) \bigg|.
	\end{align}	
	As detailed in \Cref{Section: Details on IPM}, the global sensitivity of $T$ is precisely equal to 
	\begin{align*}
		\Delta_{T} = \frac{1}{\min\{n,m\}} \sup_{X,X' \in \mathbb{S}} \sup_{f \in \mathcal{F}}|f(X) - f(X')|.
	\end{align*}
	The IPM includes several metrics commonly used in the literature~\citep{Sriperumbudur2012} and their sensitivity can be analyzed as follows.
	\begin{enumerate}[(a)]
		\item \emph{Mean difference in $\ell_p$}: Let $\mathbb{S} \subset \mathbb{R}^d$, and set $1 \leq p,q \leq \infty$ such that $1/p + 1/q=1$. Choosing $\mathcal{F} = \{f: \mathbb{S} \mapsto \mathbb{R}  \given f(x) = a^\top x, \|a\|_q \leq 1\}$, the IPM becomes the $p$th norm of the sample mean difference between $\mathcal{Y}_n$ and $\mathcal{Z}_m$. In this case, we have $\sup_{X,X' \in \mathbb{S}} \sup_{f \in \mathcal{F}}|f(X) - f(X')| = \sup_{X,X' \in \mathbb{S}}\|X - X'\|_p$.
        \item \emph{Wasserstein distance}: When $\mathcal{F} =\{f : \mathbb{S} \mapsto \mathbb{R}  \given  \|f\|_{\mathrm{Lip}} \leq 1\}$ where $\|f\|_{\mathrm{Lip}}$ denotes the minimal Lipschitz constant for $f$ on a metric space $(\mathbb{S},\|\cdot\|)$, the IPM corresponds to the Wasserstein 1-distance. By the Lipschitz property of $f$, we have $\sup_{X,X' \in \mathbb{S}} \sup_{f \in \mathcal{F}}|f(X) - f(X')| \leq \sup_{X,X' \in \mathbb{S}}\|X-X'\|$.
		\item \emph{Total variation distance}:  Let $\mathcal{F}=\{f : \mathbb{S} \mapsto \mathbb{R}  \given \sup_{x \in \mathbb{S}} |f(x)| \leq 1\}$. The corresponding IPM is the total variation distance for which we have $\sup_{X,X' \in \mathbb{S}} \sup_{f \in \mathcal{F}}|f(X) - f(X')| \leq 2$. 
		\item \emph{Kolmogorov distance}: When $\mathcal{F} =\{\mathds{1}(-\infty,x]:x \in\mathbb{R}^d \}$, the IPM is called the Kolmogorov distance. In this case, we have $\sup_{X,X' \in \mathbb{S}} \sup_{f \in \mathcal{F}}|f(X) - f(X')| \leq 1$ with $\mathbb{S} = \mathbb{R}^d$. 
		\item \emph{Maximum mean discrepancy}: Let $\|f\|_{\mathcal{H}_k}$ be the norm of a function $f$ in a reproducing kernel Hilbert space $\mathcal{H}_k$ equipped with kernel $k$. When $\mathcal{F} = \{f : \mathbb{S} \mapsto \mathbb{R}  \given \|f\|_{\mathcal{H}_k} \leq 1\}$, the IPM corresponds to the maximum mean discrepancy and it satisfies that  $\sup_{X,X' \in \mathbb{S}} \sup_{f \in \mathcal{F}}|f(X) - f(X')| \leq \sqrt{2K}$ where $K$ is the maximum value of a non-negative kernel $k$. See \Cref{Lemma: Sensitivity of MMD} for details.
	\end{enumerate}
\end{example}

In general, obtaining the exact value of the global sensitivity $\Delta_T$ can be challenging, and thus we often work with an upper bound for $\Delta_T$. We remark that the differential privacy guarantee of the Laplace mechanism in \Cref{Lemma: Laplace Mechanism} remains valid when we replace the sensitivity in the Laplace mechanism with any upper bound. With an abuse of notation, we also use $\Delta_T$ to denote an upper bound for the global sensitivity, when the exact value of the global sensitivity is not available.

\paragraph{Naive Approach.} Given the global sensitivity of $T$, one naive attempt to privatize the permutation test is to apply the basic composition theorem~(\Cref{Lemma: Composition Theorem}). To depict the idea, let $\{\zeta_i\}_{i=0}^B$ be a sequence of i.i.d.~$\mathsf{Laplace}(0,1)$ random variables, and define 
\begin{align*}
	\widetilde{M}_i \coloneqq  T_i + \frac{\Delta_T}{\varepsilon(B+1)^{-1} + \log \bigl(1 / \{1 - \delta(B+1)^{-1}\} \bigr)}\zeta_i, \quad \text{for $i \in \{0\} \cup [B]$.}
\end{align*}
By the Laplace mechanism, each $\widetilde{M}_i$ is $\big(\varepsilon/(B+1),\delta/(B+1)\big)$-DP and the composition theorem in \Cref{Lemma: Composition Theorem} then ensures that the permutation $p$-value given as
\begin{align} \label{Eq: naive permutation p-value}
	\hat{p}_{\mathrm{dp}}^{\,\mathrm{naive}} \coloneqq  \frac{1}{B+1} \bigg\{ \sum_{i=1}^B \mathds{1}\bigl(\widetilde{M}_i \geq \widetilde{M}_0 \bigr) +1 \bigg\}
\end{align}
is $(\varepsilon,\delta)$-DP. Moreover, $\{\widetilde{M}_i\}_{i=0}^B$ are exchangeable under the null, which in turn yields that $\hat{p}_{\mathrm{dp}}^{\,\mathrm{naive}}$ is a valid $p$-value by \Cref{Lemma: permutation p-value}. While this naive approach returns rigorous guarantees on both privacy and validity, there is room for improvement in terms of power performance. Specifically, observe that the noise level grows linearly in the number of permutations $B$. Consequently, when $B$ becomes significantly large, the Laplace noise overwhelms the signal, leading to a loss of power. It is also worth noting that the permutation $p$-value is lower bounded by $(B+1)^{-1}$. This implies that in order to have non-zero power, the number of permutations $B$ must exceed $\alpha^{-1} - 1$. Hence, one cannot take $B$ to be arbitrarily small. This issue serves as the motivation for our proposal, which is described below.

\paragraph{Refined Approach.} The factor of $B+1$ arises from an application of the composition theorem (\Cref{Lemma: Composition Theorem}), which cannot be improved in general~\citep[\emph{e.g.}, Section 2.1 of][]{steinke2022composition}. As one of our key contributions, we remove this unpleasant dependence on $B$ via the quantile representation of the permutation test~(\Cref{Lemma: Quantile representation}). To describe our proposal, define
\begin{align*}
	M_i \coloneqq  T_i + \frac{2 \Delta_T}{\xi_{\varepsilon,\delta}}\zeta_i, 
\end{align*}
for $i \in \{0\} \cup [B]$, where $\xi_{\varepsilon,\delta}$ can be recalled in \eqref{Eq: definition of xi}. Notably, the noise level $2 \Delta_T \xi_{\varepsilon,\delta}^{-1}$ is independent of $B$ and strictly smaller than that of the naive approach for any $B > 1$. Given $\{M_i\}_{i=0}^B$, we define the private permutation $p$-value as 
\begin{align} \label{Eq: permutation p-value}
	\hat{p}_{\mathrm{dp}} \coloneqq  \frac{1}{B+1} \bigg\{ \sum_{i=1}^B \mathds{1}\bigl(M_i \geq M_0 \bigr) +1 \bigg\}, 
\end{align}
and reject the null when $\hat{p}_{\mathrm{dp}} \leq \alpha$. We summarize the proposed method in Algorithm~\ref{Algorithm: DP permutation test}. 

\begin{algorithm}[t]\raggedright \caption{Differentially Private Permutation Test} \label{Algorithm: DP permutation test}
	\textbf{Input}: Data $\mathcal{X}_n$, significance level $\alpha\in(0,1)$, privacy parameters $\varepsilon>0$ and $\delta\in[0,1)$, test statistic $T$, global sensitivity (or its upper bound) $\Delta_T$, number of permutations $B\in\mathbb N$.\\
	\vskip .3em
	\textbf{For} $i \in [B]$ \textbf{do} \\
	\begin{algorithmic}
		\State Generate a random permutation $\bpi_i$ of $[n]$. 
		\State Generate $\zeta_i \sim \mathsf{Laplace}(0,1)$. 
		\State Set $M_i \leftarrow T(\mathcal{X}_n^{\bpi_i}) + 2 \Delta_T \xi_{\varepsilon,\delta}^{-1} \zeta_i$ where $\xi_{\varepsilon,\delta} \coloneqq \varepsilon + \log(1 / (1-\delta))$. 
	\end{algorithmic}
	\textbf{End For} 
	\begin{algorithmic}
		\State Generate $\zeta_0 \sim \mathsf{Laplace}(0,1)$ and set $M_0 \leftarrow T(\mathcal{X}_n) + 2 \Delta_T \xi_{\varepsilon,\delta}^{-1} \zeta_0$.
		\State Compute the permutation $p$-value $\hat{p}_{\mathrm{dp}}$ as in \eqref{Eq: permutation p-value}.
	\end{algorithmic}
	\textbf{Output}: Reject $H_0$ if $\hat{p}_{\mathrm{dp}} \leq \alpha$.
\end{algorithm}

\subsection{Validity and Privacy Guarantee}
Having introduced our method, we next investigate its theoretical guarantees and provide intuition behind our proposal. We start with the validity of the private test, which follows immediately from \Cref{Lemma: permutation p-value}. 

\begin{theorem}[Validity Guarantee] \label{Theorem: validity of private permutation tests}
	Suppose that $\mathcal{X}_n$ are exchangeable under the null $H_0: P \in \mathcal{P}_0$. Then for any $\alpha \in (0,1)$ and $B,n \geq 1$, the type I error of the test $\mathds{1}(\hat{p}_{\mathrm{dp}} \leq \alpha)$ from \Cref{Algorithm: DP permutation test} satisfies 
	\begin{align*} 
		\sup_{P \in \mathcal{P}_0} \mP_{P} (\hat{p}_{\mathrm{dp}} \leq \alpha) = \frac{\floor{(B+1)\alpha}}{B+1} \leq \alpha.
	\end{align*}
\end{theorem}

It is worth emphasizing that type I error control of the proposed test is both non-asymptotic and uniform over the entire class of null distributions $\mathcal{P}_0$. Another distinct feature is that the type I error is equal to ${\floor{(B+1)\alpha}}/{(B+1)}$, which can be strictly smaller than $\alpha$. If this small gap is a concern, one can make the type I error exactly equal to $\alpha$ through randomization~(\Cref{Lemma: Randomized tests}). We also remark that even if we replace the global sensitivity~$\Delta_T$ in the procedure with any other value, type I error control remains valid. In other words, the validity of the proposed test is not affected by the noise level of the Laplace mechanism. 

Next we turn to the privacy guarantee of the proposed test and show that it is $(\varepsilon,\delta)$-DP. 

\begin{theorem}[Privacy Guarantee] \label{Theorem: Differential privacy of permutation tests}
	For any $\alpha \in (0,1)$, the permutation test $\mathds{1}(\hat{p}_{\mathrm{dp}} \leq \alpha)$ from \Cref{Algorithm: DP permutation test} is $(\varepsilon,\delta)$-differentially private. 
\end{theorem}
It is worth highlighting that the privacy guarantee does not require the exchangeability of $\mathcal{X}_n$. Hence the proposed test is $(\varepsilon,\delta)$-DP under both the null and the alternative. As mentioned before, we prove the privacy guarantee of the proposed test via the quantile representation of the permutation test~(\Cref{Lemma: Quantile representation}). That is, rejecting the null when $\hat{p} \leq \alpha$ where $\hat{p}$ is given in \eqref{Eq: permutation p-value representation} is equivalent to rejecting the null when $T_0 > Q_{1-\alpha}$ where $Q_{1-\alpha}$ is the $1-\alpha$ quantile of $\{T_i\}_{i=0}^B$. Roughly speaking, our proof proceeds by privatizing $T_0$ and $Q_{1-\alpha}$, separately, which raises the factor of $2$ in $2 \Delta_T \xi_{\varepsilon,\delta}^{-1}$. However, a direct application of the Laplace mechanism to $T_0$ and $Q_{1-\alpha}$ destroys the exchangeability of random variables, thereby type I error control is no longer guaranteed. Making both $T_0$ and $Q_{1-\alpha}$ private while ensuring the finite-sample validity of the resulting test is non-trivial, and thus we highlight it as our main contribution. Along the way, we develop a general sensitivity result of quantiles in \Cref{Lemma: Sensitivity of quantiles}, which may be of independent interest. We also point out that the factor of $2$ in the noise level is a price to pay for not knowing the null distribution of $T$. When $T$ is distribution-free under the null, then it is possible to sharpen the constant factor from two to one. 

\subsection{Power Analysis}
Moving our focus to the power property, we aim to provide tractable conditions for pointwise consistency and non-asymptotic uniform power. Starting with pointwise consistency, the following result provides conditions under which the power converges to one as the sample size increases against a fixed alternative. Below, we add the subscript $n$ to $B_n$ to indicate that the number of permutations can vary with the sample size. 

\begin{theorem}[Pointwise Consistency] \label{Theorem: Pointwise Consistency}
	Let $\alpha \in (0,1)$ be a fixed constant. For a given alternative distribution $P$, suppose that $\lim_{n \rightarrow \infty} \mP_P(M_0 \leq M_1) = 0$. Then for any positive sequence of $B_n$ such that $\min_{n \geq 1} B_n > \alpha^{-1} - 1$, the differentially private permutation test is consistent in power as $\lim_{n \rightarrow \infty} \mP_P(\hat{p}_{\mathrm{dp}} \leq \alpha) = 1$. 
\end{theorem}

In view of the above result, proving consistency of the permutation test essentially boils down to verifying the condition $M_0 > M_1$, \emph{i.e.}, the original statistic is greater than a permuted statistic, with probability approaching one. In \Cref{Section: Applications to Kernel-Based Inference}, we showcase the consistency results based on kernel-based methods for two-sample and independence testing. We note that \Cref{Theorem: Pointwise Consistency} can be proven in a straightforward manner via a union bound when $B_n$ is fixed.  A similar result for fixed $B_n$ can be found in \citet[][Lemma 5.2]{dobriban2022consistency} and \citet[][Theorem 6]{rindt2020consistency}. Extending this result to any arbitrary sequence of $B_n$ requires a different technique that exploits the conditional i.i.d.~structure of given variables. To broaden the scope of our paper, we develop a consistency result for general resampling-based tests in \Cref{Lemma: General conditions for consistency} of \Cref{Section: General Pointwise Consistency Result}, from which we can derive \Cref{Theorem: Pointwise Consistency} as a corollary.

While pointwise consistency is a useful property, it is often regarded as a relatively weak guarantee. We now shift our focus to the second result of this subsection, providing a non-asymptotic, uniform guarantee on the power under stronger assumptions. In particular, we identify the moment conditions under which the proposed test has significant power. These conditions can be regarded as an extension of \citet[][Lemma 3.1]{kim2020minimax} to differentially private settings, recovering their conditions when $\xi_{\varepsilon,\delta} = \infty$. Below, the symbols $\mE_{P,\boldsymbol{\pi}}$ and $\mV_{P,\boldsymbol{\pi}}$ denote the expectation and the variance, respectively, taken over both $\mathcal{X}_n$ and $\bpi$. 
\begin{theorem}[Uniform Power] \label{Theorem: General uniform power condition}
	For $\alpha \in (0,1)$, $\beta \in (0,1-\alpha)$ and $\xi_{\varepsilon,\delta} > 0$, assume that $B \geq 6 \alpha^{-1} \log(2\beta^{-1})$ and for any $P \in \mathcal{P}_1$,
	\begin{equation}
	\begin{aligned}  \label{Eq: Uniform Power Condition}
		\mE_{P}[T(\mathcal{X}_n)] - \mE_{P,\boldsymbol{\pi}}[T(\mathcal{X}_n^{\boldsymbol{\pi}})] ~ \geq ~ &C_1 \sqrt{\frac{\mV_{P}[T(\mathcal{X}_n)] + \mV_{P,\boldsymbol{\pi}}[T(\mathcal{X}_n^{\boldsymbol{\pi}})]}{\alpha \beta}} \\[.5em]
		+ ~ & C_2 \frac{\Delta_T}{\xi_{\varepsilon,\delta}} \max\bigg\{ \! \log \biggl(\frac{1}{\alpha}\biggr),  \, \log \biggl(\frac{1}{\beta}\biggr) \bigg\},
	\end{aligned}
	\end{equation}
	where $C_1$ and $C_2$ are universal constants. Then the uniform power of the private permutation test is bounded below by $1- \beta$ as
	\begin{align*}
		\inf_{P \in \mathcal{P}_1} \mP_{P} (\hat{p}_{\mathrm{dp}} \leq \alpha) \geq 1- \beta.
	\end{align*}
\end{theorem}
A few remarks are in order.
\begin{itemize}
\item The above theorem ensures that the private permutation test has significant power as long as the signal of the problem, namely the difference between the expected values of the original test statistic and of the permuted test statistic, is larger than the noise of the problem, namely the square root of the variances and the noise level of the Laplace mechanism.
\item The proof of \Cref{Theorem: General uniform power condition}, given in \Cref{Section: Proof Theorem: General uniform power condition}, builds on the proof of \citet[][Lemma 3.1]{kim2020minimax} where the key idea is to replace the random permutation threshold with a deterministic one using concentration inequalities. The main distinction from \cite{kim2020minimax} is the incorporation of Laplace noises in the analysis, which results in the second line of the condition~\eqref{Eq: Uniform Power Condition}. We also note that \Cref{Theorem: General uniform power condition} concerns the Monte Carlo permutation test, which is computationally more efficient than the full permutation test analyzed in \citet[][Lemma 3.1]{kim2020minimax}. 
\item Notably, the condition on $B$ is independent of the sample size. Similar conditions can be found in \cite{schrab2021mmd} and \cite{schrab2022efficient}, while we improve the dependence on $\alpha$ from $\alpha^{-2}$ to $\alpha^{-1}$ using a multiplicative Chernoff bound~(\Cref{Lemma: Quantile Approximation}). One can further sharpen this restriction on $B$, especially constant factors, at the expense of inflating constants in \eqref{Eq: Uniform Power Condition}. See equation~\eqref{Eq: another bound for the quantile} in \Cref{Section: Proof of Theorem: Differential privacy of permutation tests}, \Cref{Lemma: Quantile Approximation} and \Cref{Remark: quantile approximation} for a further discussion.
\item It is worth mentioning that the first line of the condition~\eqref{Eq: Uniform Power Condition} relies on a polynomial dependence on $\alpha$ and $\beta$, which arise from the application of Chebyshev's and Markov's inequalities. If the considered test statistic has an exponential tail bound, these polynomial factors can be improved to logarithmic ones as we illustrate in \Cref{Section: Uniform separation and optimality}. 
\end{itemize}

Before moving on, let us briefly illustrate \Cref{Theorem: General uniform power condition} based on the plug-in IPM statistic considered in \Cref{Example: Sensitivity of Integral Probability Metric}. 

\begin{example}[Power Analysis against IPM alternatives] \normalfont \label{Example: Power Analysis against IPM}
	Continuing our discussion from \Cref{Example: Sensitivity of Integral Probability Metric}, denote the IPM between $P$ and $Q$ with a class of functions $\mathcal{F}$ as 
	\begin{align*}
		\mathrm{IPM}_{\mathcal{F}}(P,Q) = \sup_{f \in \mathcal{F}} \big| \mE_P[f(Y)] - \mE_Q[f(Z)] \big|.
	\end{align*}
	Without loss of generality, assume $n \leq m$ and write the pooled sample as $\mathcal{X}_{n+m} = \mathcal{Y}_n \cup \mathcal{Z}_m = \{X_1,\ldots,X_{n+m}\}$. Consider the maximum Rademacher complexity of $\mathcal{F}$ over all possible permuted samples given as
	\begin{align*}
		\mathcal{R}_n(\mathcal{F}) = \sup_{\bpi \in \boldsymbol{\Pi}_{n+m}}\mE\biggl[\sup_{f \in \mathcal{F}}\bigg| \frac{1}{n} \sum_{i=1}^n \omega_i f(X_{\pi_i})\bigg| \biggr],
	\end{align*}
	where $\{\omega_i\}_{i=1}^n$ are i.i.d.~Rademacher random variables independent of $\mathcal{X}_{n+m}$. Suppose that we implement \Cref{Algorithm: DP permutation test} using the plug-in IPM estimator $T(\mathcal{X}_{n+m})$ in \eqref{Eq: plug-in IPM}. Then \Cref{Theorem: General uniform power condition} yields that the resulting permutation test has the power greater than $1-\beta$ if 
	\begin{align*}
	   \mathrm{IPM}_{\mathcal{F}}(P,Q) \geq C_1\frac{\mathcal{R}_n(\mathcal{F})}{\sqrt{\alpha \beta}} + C_2 \frac{\sqrt{n}\Delta_T}{\sqrt{\alpha \beta}} + C_3 \frac{\Delta_T}{\xi_{\varepsilon,\delta}} \max\bigg\{ \! \log \biggl(\frac{1}{\alpha}\biggr),  \, \log \biggl(\frac{1}{\beta}\biggr) \bigg\},
    \end{align*}
	where $C_1,C_2,C_3$ are some positive constants. We defer a detailed analysis that leads to the above result to \Cref{Section: Details on Power against IPM}. We also refer to \citet{vanderVaart1996,bartlett2002rademacher,wainwright2019high} for additional information on the Rademacher complexity and illustrative examples.
\end{example}

So far we have examined the properties of the private permutation test in a general context. In the next section, we will apply our framework to the specific problem of kernel testing, and provide a detailed analysis. 

\section{Application: Differentially Private Kernel Tests} \label{Section: Applications to Kernel-Based Inference}
In recent years, there has been a growing trend in employing kernel-based methods for hypothesis testing problems, such as the MMD and the HSIC. This popularity is partly due to their ability to capture complex, non-linear relationships and to their straightforward implementation. Equipped with such benefits, the MMD~\citep{gretton2012kernel} is used to measure the difference between two probability distributions, while the HSIC~\citep{gretton2005measuring} is used to quantify the dependence between two random variables. In this and subsequent sections, we propose differentially private tests based on these two kernel-based measures, and provide an in-depth analysis of their theoretical properties.  

\paragraph{Terminology.} Before we begin, let us establish the terminology related to kernels. Consider a reproducing kernel $k: \mathbb{S} \times \mathbb{S} \mapsto \mathbb{R}$ defined on a separable topological space $\mathbb{S}$. Let $\mathcal{H}_k$ be a reproducing kernel Hilbert space (RKHS) endowed with kernel $k$. A kernel $k$ is said to be \emph{characteristic} if the kernel mean embedding 
\begin{align*}
	\mu_P  = \int_{\mathbb{S}} k(\cdot, x) \mathrm{d}P(x) \in \mathcal{H}_k
\end{align*}
is injective. In addition, a kernel $k : \mathbb{S} \times \mathbb{S} \mapsto \mathbb{R}$ is said to be \emph{translation invariant} if there exists a symmetric positive definite function $\kappa$ such that $k(x,y) = \kappa(x-y)$ for all $x,y \in \mathbb{S}$. Assuming that $0 \leq k(x,y) \leq K$ for all $x,y \in \mathbb{S}$, we say that the kernel $k$ has \emph{non-empty level sets} on $\mathbb{S}$ if, for any $\epsilon \in (0, K)$, there exist $x,y \in \mathbb{S}$, such that $k(x,y) \leq \epsilon$. Some popular examples of kernels include the Gaussian kernel $k(x,y) = e^{-\sigma \|x - y\|_2^2}$ and the Laplacian kernel $k(x,y) = e^{-\sigma \|x - y\|_1}$ for $\sigma >0$. These two kernels are translation invariant and known to be characteristic on $\mathbb{R}^d$ \citep[\emph{e.g.},][]{sriperumbudur2011universality}. They also have non-empty level sets on $\mathbb{R}^d$, which can be deduced from the continuity of the kernel function.

\subsection{Differentially Private MMD Test} \label{Section: Differentially Private MMD Test}
Starting with the MMD, suppose we are given mutually independent samples $\mathcal{Y}_n\coloneqq  \{Y_1,\ldots,Y_n\} \iid P$ and $\mathcal{Z}_m\coloneqq  \{Z_1,\ldots,Z_m\} \iid Q$ on a domain $\mathbb{S}$. Based on these samples, the two-sample problem aims to determine whether two probability distributions $P$ and $Q$ coincide. A majority of two-sample methodologies target a certain metric between $P$ and $Q$, and use their empirical counterpart as a test statistic. One such method is the non-private MMD test~\citep{gretton2012kernel} where the difference between $P$ and $Q$ is quantified in terms of MMD. To elaborate, consider the unit ball in a RKHS $\mathcal{H}_k$ denoted by $\mathcal{F}_k \coloneqq  \{f \in \mathcal{H}_k : \|f\|_{\mathcal{H}_k} \leq 1\}$. The maximum mean discrepancy between $P$ and $Q$ is defined as 
\begin{align*}
	\mathrm{MMD}_k(P,Q) \coloneqq  \sup_{f \in \mathcal{F}_k} \bigl\{ \mE_P[f(Y)] - \mE_Q[f(Z)] \bigr\}. 
\end{align*}
The empirical MMD is a plug-in estimator of MMD that replaces $P$ and $Q$ with the corresponding empirical probability measures. Formally, letting $\mathcal{X}_{n+m} = \mathcal{Y}_n \cup \mathcal{Z}_m$ be the pooled sample as before, the empirical MMD is given as
\begin{align} \label{Eq: empirical MMD}
	\widehat{\mathrm{MMD}}(\mathcal{X}_{n+m}) \coloneqq  \sup_{f \in \mathcal{F}_k} \biggl\{ \frac{1}{n} \sum_{i=1}^n f(Y_i) - \frac{1}{m} \sum_{j=1}^m \blue{f(Z_j)} \biggr\}.
\end{align}
Thanks to the reproducing kernel property, the empirical MMD can be computed straightforwardly. In particular, the squared empirical MMD can be calculated in quadratic time using the kernel-based expression
\begin{align} \label{Eq: closed form MMD}
	\widehat{\mathrm{MMD}}^2(\mathcal{X}_{n+m}) = \frac{1}{n^2} \sum_{i,j=1}^{n} k(Y_i,Y_j) +  \frac{1}{m^2} \sum_{i,j=1}^{m} k(Z_i,Z_j) - \frac{2}{nm} \sum_{i=1}^n\sum_{j=1}^m k(Y_i,Z_j).
\end{align}
In order to propose a private version of the MMD test, we begin with the global sensitivity of the empirical MMD.
\begin{lemma}[Sensitivity of Empirical MMD] \label{Lemma: Sensitivity of MMD}
	Assume that the kernel $k$ is bounded as $0 \leq k(x,y) \leq K$ for all $x,y \in \mathbb{S}$. Then the global sensitivity of the empirical MMD satisfies 
	\begin{align*}
		\sup_{\bpi \in \boldsymbol{\Pi}_{n+m}} \sup_{\substack{\mathcal{X}_{n+m},\tilde{\mathcal{X}}_{n+m}:\\d_{\mathrm{ham}}(\mathcal{X}_{n+m},\tilde{\mathcal{X}}_{n+m}) \leq 1}} \bigl| \widehat{\mathrm{MMD}}(\mathcal{X}_{n+m}^{\bpi}) - \widehat{\mathrm{MMD}}(\tilde{\mathcal{X}}_{n+m}^{\bpi}) \bigr| \leq \frac{\sqrt{2K}}{\min\{n,m\}}. 
	\end{align*}
    Moreover assume that $k$ is translation invariant, and has non-empty level sets in $\mathbb{S}$. Then the inequality becomes an equality.
\end{lemma}

The proof of \Cref{Lemma: Sensitivity of MMD} can be found in \Cref{Section: Proof of Lemma: Sensitivity of MMD}. Note that the sensitivity of the empirical MMD in \Cref{Lemma: Sensitivity of MMD} can be equivalently defined without the supremum over the permutations $\bpi \in \boldsymbol{\Pi}_{n+m}$ as $d_{\mathrm{ham}}(\mathcal{X}_{n+m},\tilde{\mathcal{X}}_{n+m})$ is the same as $d_{\mathrm{ham}}(\mathcal{X}_{n+m}^{\bpi},\tilde{\mathcal{X}}_{n+m}^{\bpi})$ for any $\bpi \in \boldsymbol{\Pi}_{n+m}$. As we will see in \Cref{Section: Differentially Private HSIC Test}, however, this property does not hold for independence testing. We also highlight our lower bound result, indicating that the upper bound $\sqrt{2K}/\min\{n,m\}$ cannot be improved for translation invariant kernels with non-empty level sets on their domain.

With \Cref{Lemma: Sensitivity of MMD} in place, we set the sensitivity parameter $\Delta_T = \sqrt{2 K}/\min\{n,m\}$ and run Algorithm~\ref{Algorithm: DP permutation test} with the empirical MMD in \eqref{Eq: empirical MMD} as the test statistic. We refer to the resulting private permutation test as the $\mathttt{dpMMD}$ test and denote it as $\phi_{\mathttt{dpMMD}}$. The $\mathttt{dpMMD}$ test has the following properties, which are proven in \Cref{Section: Proof of Theorem: Properties of dpMMD}. 
\begin{theorem}[Properties of $\mathttt{dpMMD}$ test] \label{Theorem: Properties of dpMMD}
		Let $\alpha \in (0,1)$ be a fixed constant, the kernel $k$ be bounded as $0 \leq k(x,y) \leq K$ for all $x,y \in \mathbb{S}$ and $n \leq m$. Then $\phi_{\mathttt{dpMMD}}$ satisfies the following properties:
	\begin{enumerate}
		\item[P1.] (Differential Privacy) For $\varepsilon >0$ and $\delta\in[0,1)$, $\phi_{\mathttt{dpMMD}}$ is $(\varepsilon,\delta)$-differentially private.
		\item[P2.] (Validity) The type I error of $\phi_{\mathttt{dpMMD}}$ is controlled at level $\alpha$ non-asymptotically.
		\item[P3.] (Consistency) Suppose that $\mathrm{MMD}_k(P,Q)$ is independent of the sample sizes and strictly positive for a fixed pair of $(P,Q)$. Moreover assume that $n^{-1} \xi_{\varepsilon,\delta}^{-1} \rightarrow 0$ as $n \rightarrow \infty$. Then for any sequence $B_n$ such that $\min_{n \geq 1} B_n > \alpha^{-1}-1$, we have $\lim_{n \rightarrow \infty} \mE_{P,Q} [\phi_{\mathttt{dpMMD}}] = 1$.
	\end{enumerate}
\end{theorem}
The first two properties on differential privacy and validity are clear in view of \Cref{Theorem: validity of private permutation tests} and \Cref{Theorem: Differential privacy of permutation tests}. It is well-known that the population MMD is strictly positive under the alternative if the kernel is characteristic~\citep{gretton2012kernel}. Hence the condition on the MMD metric in (\emph{P3}) is satisfied for characteristic kernels against any fixed alternative. Another highlight is that consistency holds irrespective of the relationship between $n$ and $m$. The power converges to one as long as the minimum sample size goes to infinity. We also remark that the condition $n^{-1} \xi_{\varepsilon,\delta}^{-1} \rightarrow 0$ is critical for obtaining consistency. If not, the empirical MMD is overwhelmed by the Laplace noise, which leads to a significant loss of power. 

We note in passing the recent work of \cite{yang2023differentially} that also utilizes the MMD in differentially private data analysis. Despite the fact that both \cite{yang2023differentially} and ours consider the differentially private MMD, their primary focus is on differentially private data generation, which is different from our focus on hypothesis testing.

\subsection{Differentially Private HSIC Test} \label{Section: Differentially Private HSIC Test}
Turning to the second application, suppose that we are given an i.i.d.~sample $\mathcal{X}_n=(X_1,\ldots,X_n)$, where $X_i=(Y_i,Z_i)$ follows a joint distribution $P_{YZ}$ on domain $\mathbb{Y} \times \mathbb{Z}$. Given $\mathcal{X}_n$, the aim of independence testing is to assess whether $Y$ and $Z$ are statistically independent. As a kernel dependence measure, the HSIC compares the joint probability measure $P_{YZ}$ to the product of marginals $P_YP_Z$. To formally define it, let $k$ and $\ell$ be kernels on $\mathbb{Y}$ and $\mathbb{Z}$, and let $k \otimes \ell$ be the product kernel given by
\[
	k \otimes \ell\bigl((y,z),(y',z')\bigr)=k(y,y')\ell(z,z'),
\]
for all $y,y' \in \mathbb{Y}$ and $z,z' \in \mathbb{Z}$. Further denoting the unit ball in the RKHS associated with $k \otimes \ell$ by $\mathcal{F}_{k \otimes \ell}$, HSIC is defined as\footnote{We remark that in the literature HSIC is often defined as the square of this quantity. However, for consistency with MMD, we define it without the square.}
\begin{align*}
	\mathrm{HSIC}_{k \otimes \ell} (P_{YZ}) \coloneqq  \sup_{f \in \mathcal{F}_{k \otimes \ell}} \big\{ \mE_{P_{YZ}}[f(Y,Z)] - \mE_{P_YP_Z}[f(Y,Z)] \big\}.
\end{align*}
In other words, the HSIC of $Y$ and $Z$ is simply the MMD between $P_{YZ}$ and $P_YP_Z$ with the product kernel $k \otimes \ell$. The empirical HSIC is a plug-in estimator given as
\begin{align} \label{Eq: empirical HSIC}
	\widehat{\mathrm{HSIC}}(\mathcal{X}_n) \coloneqq  \sup_{f \in \mathcal{F}_{k \otimes \ell}} \bigg\{ \frac{1}{n} \sum_{i=1}^n f(Y_i,Z_i) - \frac{1}{n^2}\sum_{i,j=1}^n f(Y_i,Z_j) \bigg\}.
\end{align}
Similarly to the empirical MMD, the squared empirical HSIC also has an explicit form in terms of the kernels $k$ and $\ell$ as
\begin{equation}
\begin{aligned} \label{Eq: closed form HSIC}
	\widehat{\mathrm{HSIC}}^2\!(\mathcal{X}_n) ~=~ & \frac{1}{n^2} \sum_{i,j=1}^n k(Y_i,Y_j) \ell(Z_i,Z_j) + \frac{1}{n^4} \sum_{i_1,i_2,j_1,j_2=1}^n k(Y_{i_1},Y_{j_1}) \ell(Z_{i_2},Z_{j_2})\\
	 -& \frac{2}{n^3} \sum_{i,j_1,j_2=1}^n k(Y_{i},Y_{j_1}) \ell(Z_{i},Z_{j_2}),
\end{aligned}
\end{equation}
which can be computed in quadratic time as described in \citet[][Theorem 1]{song2012feature}. For independence testing, the permutation test proceeds by randomly permuting either the $Y$ observations or the $Z$ observations. Here, we permute the $Z$ observations and denote $\mathcal{X}_n^{\bpi} = \{(Y_i,Z_{\pi_i})\}_{i=1}^n$. With this notation in place, the next lemma explores the global sensitivity of the empirical HSIC. 

\begin{lemma}[Sensitivity of Empirical HSIC] \label{Lemma: Sensitivity of HSIC}
	Assume that the kernels $k$ and $\ell$ are bounded as $0 \leq k(y,y') \leq K$ and $0 \leq \ell(z,z') \leq L$ for all $y,y' \in \mathbb{Y}$ and $z,z' \in \mathbb{Z}$. Then the global sensitivity of the empirical HSIC satisfies 
	\begin{align*}
		\sup_{\bpi \in \boldsymbol{\Pi}_n} \sup_{\substack{\mathcal{X}_{n},\tilde{\mathcal{X}}_{n}:\\d_{\mathrm{ham}}(\mathcal{X}_{n},\tilde{\mathcal{X}}_{n}) \leq 1}} \bigl| \widehat{\mathrm{HSIC}}(\mathcal{X}_{n}^{\bpi}) - \widehat{\mathrm{HSIC}}(\tilde{\mathcal{X}}_{n}^{\bpi}) \bigr| \leq \frac{4(n-1)}{n^2} \sqrt{KL}.
	\end{align*}
    Moreover assume that $k$ and $\ell$ are translation invariant, and have non-empty level sets on $\mathbb{Y}$ and $\mathbb{Z}$, respectively. Then the global sensitivity is lower bounded by $4(n-2.5)n^{-2} \sqrt{KL}$.
\end{lemma}

The proof of \Cref{Lemma: Sensitivity of MMD} can be found in \Cref{Section: Proof of Lemma: Sensitivity of HSIC}. Contrary to the MMD case, we observe that the two hamming distances, namely $d_{\mathrm{ham}}(\mathcal{X}^{\bpi}_{n},\tilde{\mathcal{X}}^{\bpi}_{n})$ and $d_{\mathrm{ham}}(\mathcal{X}_{n},\tilde{\mathcal{X}}_{n})$, can differ for independence testing. Consequently, the supremum over the permutations $\bpi \in \boldsymbol{\Pi}_{n}$ plays a non-trivial role in the sensitivity of the empirical HSIC. We mention the work of \cite{kusner2016private} that also examines the global sensitivity of the empirical HSIC. Our upper bound result improves theirs by replacing the constant factor $12$ to $4$ with a tighter analysis. In fact, as the lower bound result states, the proposed upper bound is asymptotically tight under mild conditions for $k$ and $\ell$. 

In view of the above lemma, we set $\Delta_T = 4(n-1) n^{-2} \sqrt{KL}$ and run Algorithm~\ref{Algorithm: DP permutation test} with the empirical HSIC in~\eqref{Eq: empirical HSIC} as the test statistic. We refer to the resulting permutation test as the $\mathttt{dpHSIC}$ test and denote it as $\phi_{\mathttt{dpHSIC}}$. Similar to \Cref{Theorem: Properties of dpMMD}, the $\mathttt{dpHSIC}$ test has the following properties, which are proven in \Cref{Section: Proof Theorem: Properties of dpHSIC}.

\begin{theorem}[Properties of $\mathttt{dpHSIC}$ test] \label{Theorem: Properties of dpHSIC}
	Let $\alpha \in (0,1)$ be a fixed constant and assume that the kernels $k$ and $\ell$ are bounded as $0 \leq k(y,y') \leq K$ and $0 \leq \ell(z,z') \leq L$ for all $y,y' \in \mathbb{Y}$ and $z,z' \in \mathbb{Z}$. Then $\phi_{\mathttt{dpHSIC}}$ satisfies the following properties:
	\begin{enumerate}
		\item[P1.] (Differential Privacy) For $\epsilon >0$ and $\delta\in[0,1)$, $\phi_{\mathttt{dpHSIC}}$ is $(\varepsilon,\delta)$-differentially private.
		\item[P2.] (Validity) The type I error of $\phi_{\mathttt{dpHSIC}}$ is controlled at level $\alpha$ non-asymptotically.
		\item[P3.] (Consistency) Suppose that $\mathrm{HSIC}_{k \otimes \ell}(P_{YZ})$ is independent of the sample sizes and strictly positive for a fixed distribution $P_{YZ}$. Moreover assume that $n^{-1} \xi_{\varepsilon,\delta}^{-1} \rightarrow 0$ as $n \rightarrow \infty$. Then for any sequence $B_n$ such that $\min_{n \geq 1} B_n > \alpha^{-1} - 1$, we have $\lim_{n \rightarrow \infty} \mE_{P_{YZ}} [\phi_{\mathttt{dpHSIC}}] = 1$.
	\end{enumerate}
\end{theorem}

As for the $\mathttt{dpMMD}$ test, the first two properties on differential privacy and validity are direct consequences of \Cref{Theorem: validity of private permutation tests} and \Cref{Theorem: Differential privacy of permutation tests}. The condition for consistency is ensured under any alternative when the kernels are characteristic~\citep{gretton2015simpler}. Therefore the $\mathttt{dpHSIC}$ test equipped with a characteristic kernel is pointwise consistent against any fixed alternative, provided that $n^{-1}\xi_{\varepsilon,\delta}^{-1} \rightarrow 0$ and $\min_{n \geq 1} B_n > \alpha^{-1}-1$. 

Before moving on and studying uniform power properties, let us briefly remark on the asymptotic null distributions of private kernel test statistics.

\begin{remark}[Asymptotic null distributions] \normalfont
	As mentioned earlier, when the null distribution is tractable, we can improve the power by eliminating the factor of 2 in the noise level. However, characterizing the limiting distribution is not a trivial task, even for non-private kernel statistics. In \Cref{Section: Limiting null distributions}, we show that a private kernel statistic converges in distribution to a mixture of Gaussian chaos and Laplace distributions, which is even more intricate than the limiting distribution of a non-private kernel statistic. A recent line of work~\citep{shekhar2022ind,shekhar2022two} propose cross MMD and cross HSIC that have a tractable limiting distribution with competitive power.  We leave the exploration of extending these variants to the private setting and comparing their power performance with our proposed methods as an avenue for future research.
\end{remark}

\section{Uniform Power and Optimality} \label{Section: Uniform separation and optimality}
In the previous section, we examined the fundamental properties of the private kernel tests, including their asymptotic power against fixed alternatives. This section delves into a more challenging setting where the alternative can shrink to the null as the sample size increases, and develops uniform power results. Moreover, we highlight an intrinsic trade-off between privacy and statistical power through the lens of minimax analysis, and explore optimality of the proposed private tests under the differential privacy constraint. In the main text, we focus on the analysis of the $\mathttt{dpMMD}$ test, and defer analogous results for the $\mathttt{dpHSIC}$ test to \Cref{Section: Separation in HSIC metric} and \Cref{Section: Separation in L2 for dpHSIC}.

\subsection{Separation in MMD Metric} \label{Section: Separation in MMD metric}
Consider the setting described in \Cref{Section: Differentially Private MMD Test}, and denote by $\mathcal{P}_{\mathbb{S}}$ the class of distributions defined on $\mathbb{S}$. Our first goal is to determine the minimum separation for $\phi_{\mathttt{dpMMD}}$ based on the MMD metric with kernel $k$. To this end, for $\rho > 0$, we define a class of paired distributions $(P,Q)$ such that
\begin{align*}
	\mathcal{P}_{\mathrm{MMD}_k}\!(\rho) \coloneqq  \big\{(P,Q) \in \mathcal{P}_{\mathbb{S}} \times \mathcal{P}_{\mathbb{S}} : \mathrm{MMD}_k(P,Q) \geq \rho \big\}.
\end{align*}
For a given target type II error $\beta \in (0,1 - \alpha)$, the minimum separation for the $\mathttt{dpMMD}$ test against $\mathcal{P}_{\mathrm{MMD}_k}\!(\rho)$ is given by
\begin{align}
	\label{Eq: separation MMD rate}
	\rho_{\phi_{\mathttt{dpMMD}}}(\alpha, \beta, \varepsilon, \delta, m,n) \coloneqq  \inf \biggl\{ \rho >0 : \sup_{(P,Q) \in \mathcal{P}_{\mathrm{MMD}_k}\!(\rho)} \mE_{P,Q}[1 - \phi_{\mathttt{dpMMD}}] \leq \beta  \biggr\}.
\end{align}
In simpler terms, the minimum separation $\rho_{\phi_{\mathttt{dpMMD}}}$ refers to the smallest MMD metric between $P$ and $Q$ that can be correctly detected by the $\mathttt{dpMMD}$ test with probability at least $1-\beta$. The next theorem provides an upper bound for $\rho_{\phi_{\mathttt{dpMMD}}}$ as a function of the parameters $\alpha$, $\beta$, $\varepsilon$, $\delta$, $m$, and $n$. The proof can be found in \Cref{Section: Proof of Theorem: Uniform separation for MMD}.

\begin{theorem}[Minimum Separation of $\mathttt{dpMMD}$ over $\mathcal{P}_{\mathrm{MMD}_k}$] \label{Theorem: Uniform separation for MMD}
	Assume that the kernel $k$ is bounded as $0 \leq k(x,y) \leq K$ for all $x,y \in \mathbb{S}$, and $n \leq m \leq \tau n$ for some fixed constant $\tau \geq 1$. Then for all values of $\alpha \in (0,1)$, $\beta \in (0,1-\alpha)$, $\varepsilon >0$, $\delta\in[0,1)$ and $B \geq 6\alpha^{-1} \log(2\beta^{-1})$, the minimum separation for $\phi_{\mathttt{dpMMD}}$ satisfies
	\begin{align*}
		\rho_{\phi_{\mathttt{dpMMD}}} \leq C_{K,\tau} \max \Biggl\{ \sqrt{\frac{\max\!\big\{\!\log(1/\alpha), \, 	
		\log(1/\beta)\big\}}{n}}, \, \frac{\max\!\big\{\!\log(1/\alpha), \, \log(1/\beta)\big\}}{n\xi_{\varepsilon,\delta}}  \Biggr\},
	\end{align*}
	where $C_{K,\tau}$ is a positive constant that depends only on $K$ and $\tau$, and $\xi_{\varepsilon,\delta}$ can be recalled in \eqref{Eq: definition of xi}.
\end{theorem}

We present several comments on the upper bound result.
\begin{itemize}
	\item \Cref{Theorem: Uniform separation for MMD} states that the separation rate for the $\mathttt{dpMMD}$ test becomes $n^{-1/2}$ in low privacy regimes (\emph{i.e.},~$\xi_{\varepsilon,\delta} \gtrsim n^{-1/2}$), whereas it becomes $n^{-1}\xi_{\varepsilon,\delta}^{-1}$ in high privacy regimes (\emph{i.e.},~$\xi_{\varepsilon,\delta} \lesssim n^{-1/2}$). Notably, this upper bound result allows the parameters $\alpha,\beta,\xi_{\varepsilon,\delta}$ to vary freely within the constraints in the theorem statement. We also mention that the minimum separation is meaningful only when $n^{-1}\xi_{\varepsilon,\delta}^{-1}\to 0$, which coincides with the condition for consistency established in \Cref{Theorem: Properties of dpMMD}. 
	\item We point out that $\phi_{\mathttt{dpMMD}}$ is equivalent to the non-DP MMD test~\citep{gretton2012kernel} when $\varepsilon \to \infty$ or $\delta \to 1$ (\emph{i.e.},~$\xi_{\varepsilon,\delta} \to \infty$). Thus, our result also yields the minimum separation rate for the non-DP MMD test as a byproduct.
	\item One can prove \Cref{Theorem: Uniform separation for MMD} by verifying the general conditions in \Cref{Theorem: General uniform power condition}. However this strategy results in polynomial factors of $\alpha$ and $\beta$ instead of logarithmic ones. To obtain logarithmic dependence in both $\alpha$ and $\beta$, we modify the proof of \Cref{Theorem: General uniform power condition} and utilize exponential concentration inequalities for the empirical MMD statistic (\Cref{Lemma: Concentration for MMD}) and permuted MMD statistic (\Cref{Lemma: Bobkovs inequality}). As we will see in \Cref{Theorem: Minimax separation in MMD}, these logarithmic factors cannot be improved further when $\alpha \asymp \beta$. 
	\item The constraint on the sample size ratio can be completely removed by using the Markov inequality for a permuted MMD statistic~(\Cref{Lemma: Markov for permuted MMD}). Nevertheless, this alternative approach yields a polynomial factor of $\alpha$ instead of a logarithmic one. See \Cref{Remark: Markov inequality instead of Exponential inequality}. It is currently unknown whether the constraint on $m$ and $n$ can be eliminated, while preserving the logarithmic factors.
\end{itemize}

We next investigate minimax optimality of $\phi_{\mathttt{dpMMD}}$ under certain regimes in $\mathbb{S} = \mathbb{R}^d$. To set the stage, let $\phi: \mathcal{Y}_n \cup \mathcal{Z}_m \mapsto \{0,1\}$ be a test function, and denote the set of $(\varepsilon,\delta)$-DP level $\alpha$ tests as 
\begin{align*}
	\Phi_{\alpha,\varepsilon,\delta} \coloneqq  \Big\{\phi : \sup_{P \in \mathcal{P}_{\mathbb{S}}} \mE_{P,P}[\phi] \leq \alpha \ \text{and} \ \text{$\phi$ is $(\varepsilon,\delta)$-DP} \Big\}.
\end{align*}
From a theoretical point of view, it is of interest to figure out an information-theoretic lower bound on the minimum separation for any test. This is often called the minimax separation or critical radius in the literature \citep{ingster1994minimax,ingster2003nonparametric,baraud2002non}. Formally, the minimax separation in terms of the MMD metric is defined as  
\begin{align*}
	\rho^\star_{\mathrm{MMD}}(\alpha,\beta,\varepsilon,\delta,m,n) \coloneqq  \inf\Bigl\{ \rho > 0: \inf_{\phi \in \Phi_{\alpha,\varepsilon,\delta}} \sup_{(P,Q) \in \mathcal{P}_{\mathrm{MMD}_k}\!(\rho)} \mE_{P,Q}[1 - \phi] \leq \beta \Big\}.
\end{align*}
In simpler terms, the minimax separation $\rho^\star_{\mathrm{MMD}}$ refers to the largest MMD metric between $P$ and $Q$ that cannot be correctly detected with probability at least $1-\beta$ by any level $\alpha$ test.
We say that a test $\phi$ is minimax rate optimal in terms of the MMD metric if the minimum separation of $\phi$ is equivalent to $\rho^\star_{\mathrm{MMD}}$ up to constant factors. The next theorem, proved in \Cref{Section: Proof of Theorem: Minimax separation in MMD}, establishes a lower bound for the minimax separation under the DP constraint, from which we demonstrate minimax optimality of the $\mathttt{dpMMD}$ test.

\begin{theorem}[Minimax Separation over $\mathcal{P}_{\mathrm{MMD}_k}$] \label{Theorem: Minimax separation in MMD}
	Let $\alpha$ and $\beta$ be real numbers in the interval $(0,1/5)$, $\varepsilon>0$, $\delta\in[0,1)$ and $n \leq m$. Assume that the kernel function $k$ is translation invariant on $\mathbb{R}^d$. In particular there exists some function $\kappa$ such that $k(x,y) = \kappa(x-y)$ for all $x,y \in \mathbb{R}^d$. Moreover, the kernel is non-constant in the sense that there exists a positive constant $\eta$ such that $\kappa(0) - \kappa(z) \geq \eta$ for some $z \in \mathbb{R}^d$. Then the minimax separation over $\mathcal{P}_{\mathrm{MMD}_k}$ is lower bounded as 
	\begin{align*} 
		\rho^\star_{\mathrm{MMD}} \geq C_{\eta} \max \Biggl\{ \min \! \Bigg( \! \sqrt{\frac{\log(1/(\alpha+\beta))}{n}}, \, 1 \Bigg), \, \min \Biggl( \frac{\log(1/\beta)}{n\xi_{\varepsilon,\delta}}, \, 1 \Biggr) \Biggr\},
	\end{align*}
	where $C_{\eta}$ is a positive constant that only depends on $\eta$, and $\xi_{\varepsilon,\delta}$ can be recalled in \eqref{Eq: definition of xi}.
\end{theorem}

Several remarks are in order. 

\begin{itemize}
	\item First of all, the restriction on $\alpha$ and $\beta$ is mild as we are typically interested in small values of $\alpha$ and $\beta$. In fact, the same result holds for any $\alpha,\beta$ such that $\alpha + \beta \leq C$ where $C$ is some fixed constant strictly smaller than $1/2$.
	We also note that for a bounded kernel ranging from $0$ and $K$, the MMD as well as the corresponding minimax separation cannot exceed $\sqrt{2K}$. Our lower bound result captures this restriction through the minimum operator.
	\item Second, our proof builds on the $(\varepsilon,\delta)$-DP Le~Cam's method outlined in~\citet{acharya2018differentially,acharya2021differentially}. This technique generalizes classical Le~Cam's two-point method~\citep{lecam1973convergence} to private settings via coupling argument. 
	As pointed out by \citet[][Lemma 5]{acharya2018differentially}, one can obtain a lower bound result for $(\varepsilon,\delta)$-DP by replacing $\varepsilon$ with $\varepsilon+\delta$ in the lower bound result for $\varepsilon$-DP. However, this method fails to yield a tight lower bound in terms of $\beta$. Our approach differs from \citet[][Lemma 5]{acharya2018differentially} and returns a sharp lower bound for all parameters of interest, namely $\beta,n,\varepsilon,\delta$.
	\item \Cref{Theorem: Minimax separation in MMD} holds for translation invariant kernels. Indeed, many kernels commonly used in practice are translation invariant including the Gaussian, Laplacian, inverse multiquadrics and Mat\'{e}rn kernels. Moreover, as discussed in \cite{tolstikhin2017minimax}, if we further assume that the kernel $k$ is characteristic, it guarantees the existence of $z \in \mathbb{R}^d$ and $\eta>0$ that satisfy the conditions of \Cref{Theorem: Minimax separation in MMD}. For instance, for the Gaussian kernel $k(x,y) = e^{-\sigma \|x - y\|_2^2}$, one can take $\eta = \frac{\sigma}{2} \|z\|_2^2$ for any non-zero $z$ such that $\|z\|_2^2 \leq \sigma^{-1}$. 
	\item The last point worth highlighting is that our lower bound permits varying values of $\alpha$ and $\beta$, which is in contrast to most existing research on minimax testing. A notable exception is the recent work by \cite{diakonikolas2021optimal}, which examines the sample complexity of testing for discrete distributions with high probability. 
\end{itemize}

We now compare the results of \Cref{Theorem: Uniform separation for MMD} and \Cref{Theorem: Minimax separation in MMD}, and observe that the lower bound for $\rho^\star_{\mathrm{MMD}}$ matches the upper bound for $\rho_{\phi_{\mathttt{dpMMD}}}$ in the regime where $m \asymp n$ and $\alpha \asymp \beta$ for $\mathbb{S} = \mathbb{R}^d$. This shows that the proposed $\mathttt{dpMMD}$ test is minimax rate optimal against the class of alternatives determined by the MMD metric in the considered regime. It is noteworthy that there is no restriction on the privacy parameters $\varepsilon>0$ and $\delta\in[0,1)$, and hence the $\mathttt{dpMMD}$ test achieves optimal separation rates in all privacy regimes.

\subsection{Separation in \texorpdfstring{$L_2$}{L2} Metric} \label{Section: Separation in L2 metric}
We next investigate the minimum separation of the $\mathttt{dpMMD}$ test in terms of the $L_2$ metric. Let $p$ and $q$ denote the Lebesgue density functions of $P$ and $Q$, respectively, defined on $\mathbb{R}^d$. As in \cite{schrab2021mmd} and \cite{li2019optimality}, we restrict our attention to a smooth class of density functions defined over a Sobolev ball. In particular, for a smoothness parameter $s>0$ and a radius $R>0$, the Sobolev ball $\mathcal{S}_d^s(R)$ is given as
\begin{align} \label{Eq: Sobolev ball}
	\mathcal{S}_d^s(R) := \Bigg\{ f \in L_1(\mathbb{R}^d) \cap L_2(\mathbb{R}^d) : \int_{\mathbb{R}^d} \|w\|_2^{2s} |\hat{f}(w) |^2 \mathrm{d}w \leq (2\pi)^{d} R^2 \Bigg\},
\end{align}
where $\hat{f}$ is the Fourier transform of $f$, \emph{i.e.}, $\hat{f}(w) = \int_{\mathbb{R}^d} f(x) e^{-ix^\top w} \mathrm{d}x$ for $w\in\mathbb{R}^d$. The condition $f \in L_1(\mathbb{R}^d) \cap L_2(\mathbb{R}^d)$ simply requires the function $f \colon \mathbb{R}^d\to\mathbb{R}$ to be both integrable and square-integrable with respect to the Lebesgue measure.
For $\rho>0$, let $\mathcal{P}_{L_2}\!(\rho)$ be the collection of paired distributions $(P,Q)$ on $\mathbb{R}^d \times \mathbb{R}^d$ where $P$ and $Q$ are equipped with the Lebesgue density functions $p$ and $q$, respectively, such that $\|p-q\|_{L_2} \geq \rho$. The target class of distributions is a subset of $\mathcal{P}_{L_2}\!(\rho)$ defined as
\begin{align*}
	\mathcal{P}_{L_2}^{s}\!(\rho) \coloneqq  \big\{(P,Q) \in\mathcal{P}_{L_2}\!(\rho): p-q \in \mathcal{S}_d^s(R), \ \max(\|p\|_{L_\infty},\|q\|_{L_\infty}) \leq M \big\}.
\end{align*}
The aim of this subsection is to characterize the minimum value of $\rho$ for which the $\mathttt{dpMMD}$ test has significant power uniformly over $\mathcal{P}_{L_2}^{s}(\rho)$. For simplicity, we focus on the $\mathttt{dpMMD}$ test with a Gaussian kernel. This choice is motivated by the observation that the population MMD with the Gaussian kernel approximates $\|p-q\|_{L_2}^2$ for small bandwidth values~\citep{li2019optimality}, and a similar result can be derived using other kernels in view of \cite{schrab2021mmd}. For $x = (x_1,\ldots,x_d)^\top \in \mathbb{R}^d$ and $y = (y_1,\ldots,y_d)^\top \in \mathbb{R}^d$, the Gaussian kernel with bandwidth $\boldsymbol{\lambda} = (\lambda_1,\ldots,\lambda_d)^\top \in (0, \infty)^d$ is given as
\begin{align*}
	k_{\boldsymbol{\lambda}}(x,y) = \prod_{i=1}^d \frac{1}{\sqrt{2\pi}\lambda_i} e^{-\frac{(x_i - y_i)^2}{2\lambda_i^2}}.
\end{align*}
Let us denote the minimum separation of the $\mathttt{dpMMD}$ test with the Gaussian kernel against $L_2$ alternatives as 
\begin{align} \label{Eq: L2 separation}
	\rho_{\phi_{\mathttt{dpMMD}},L_2}(\alpha, \beta, \varepsilon, \delta, m, n, d, s, R, M) \coloneqq  \inf \Biggl\{ \rho >0 : \sup_{(P,Q) \in \mathcal{P}_{L_2}^{s}\!(\rho)} \mE_{P,Q}[1 - \phi_{\mathttt{dpMMD}}] \leq \beta \Biggr\}.
\end{align}
\blue{For comparison with arbitrary private tests, define the minimax $L_2$ separation by
\begin{align} \label{Eq: minimax L2 separation}
	\rho^\star_{L_2}(\alpha,\beta,\varepsilon,\delta,m,n,d,s,R,M)
	\coloneqq
	\inf \Biggl\{\rho>0:
	\inf_{\phi\in\Phi_{\alpha,\varepsilon,\delta}}
	\sup_{(P,Q)\in\mathcal{P}_{L_2}^{s}\!(\rho)}
	\mE_{P,Q}[1-\phi]\leq\beta
	\Biggr\}.
\end{align}}
The next theorem, proved in \Cref{Section: Proof of Theorem: Minimax Separation over L2}, provides an upper bound for $\rho_{\phi_{\mathttt{dpMMD}},L_2}$ in terms of a set of parameters, including the bandwidth $\blambda$ and sample sizes.
\begin{theorem}[Minimum Separation of $\mathttt{dpMMD}$ over $\mathcal{P}_{L_2}^{s}$] \label{Theorem: Minimax Separation over L2}
	Assume that $n \leq m \leq \tau n$ for some fixed constant $\tau \geq 1$, and that $\alpha \in (0,e^{-1})$, $\beta \in (0,1-\alpha)$, $\varepsilon >0$, $\delta \in [0,1)$, $B \geq 6\alpha^{-1}  \log(2\beta^{-1})$ and $\prod_{i=1}^d \lambda_i \leq 1$. The minimum separation of the $\mathttt{dpMMD}$ test with the Gaussian kernel over $\mathcal{P}_{L_2}^{s}$ is upper bounded as 
	\begin{align*}
		\rho_{\phi_{\mathttt{dpMMD}},L_2}^2 \leq C_{\tau,\beta,s,R,M,d} \Biggl\{ \sum_{i=1}^d \lambda_i^{2s} & +  \frac{\log(1/\alpha)}{n\sqrt{\lambda_1\cdots \lambda_d}} +\frac{\log(1/\alpha)}{n^{3/2}\lambda_1\cdots \lambda_d \xi_{\varepsilon,\delta}}  \\
		& +  \frac{\log^2(1/\alpha)}{n^2\lambda_1\cdots \lambda_d \xi_{\varepsilon,\delta}^2}+ \frac{\log^{3/2}(1/\alpha)}{n^{3/2}(\lambda_1\cdots \lambda_d )^{3/4}\xi_{\varepsilon,\delta}} \Bigg\},
	\end{align*}
	where $\xi_{\varepsilon,\delta}$ is as in \eqref{Eq: definition of xi} and $C_{\tau,\beta,s,R,M,d}$ is a positive constant, depending only on $\tau,\beta,s,R,M,d$.
\end{theorem}

Several points merit emphasis. To facilitate our discussion, assume that $\alpha$ is a fixed number and write
\begin{align*}
	& (\mathrm{I}) = \frac{\log(1/\alpha)}{n\sqrt{\lambda_1\cdots \lambda_d}}, \ (\mathrm{II}) = \frac{\log(1/\alpha)}{n^{3/2}\lambda_1\cdots \lambda_d \xi_{\varepsilon,\delta}}, \ (\mathrm{III}) =  \frac{\log^2(1/\alpha)}{n^2\lambda_1\cdots \lambda_d \xi_{\varepsilon,\delta}^2}, \ (\mathrm{IV}) = \frac{\log^{3/2}(1/\alpha)}{n^{3/2}(\lambda_1\cdots \lambda_d )^{3/4}\xi_{\varepsilon,\delta}}.
\end{align*}
When $\alpha$ is fixed, we can absorb the term $(\mathrm{IV})$ into the term $(\mathrm{II})$ as $\lambda_1\cdots \lambda_d \leq 1$, and simplify the interpretation of the result as follows.
\begin{itemize}
	\item In the low privacy regime where the first term $(\mathrm{I})$ dominates the others, our result recovers \citet[][Theorem 6]{schrab2021mmd}, which studies the minimum separation of the non-private MMD test against $\mathcal{P}_{L_2}^{s}$. In this low privacy regime, by setting bandwidths $\lambda_i = n^{-2/(4s+d)}$ for $i \in [d]$, we can achieve the optimal separation rate over the Sobolev ball, that is $n^{-2s/(4s+d)}$.
	\item In the mid privacy regime where the term $(\mathrm{II})$ becomes a leading term, equating $(\mathrm{II})$ with $\sum_{i=1}^d \lambda_i^{2s}$ yields the optimal choice of bandwidths $\lambda_i = n^{-3/(4s+2d)} \xi_{\varepsilon,\delta}^{-1/(2s+d)}$ for $i \in [d]$. The resulting separation rate is  $n^{-3s/(4s+2d)} \xi_{\varepsilon,\delta}^{-s/(2s+d)}$. Similarly, in the high privacy regime where the term $(\mathrm{III})$ dominates the others, equating $(\mathrm{III})$ with $\sum_{i=1}^d \lambda_i^{2s}$ yields the optimal choice of bandwidths $\lambda_i = (n\xi_{\varepsilon,\delta})^{-2/(2s+d)}$ for $i \in [d]$. This yields the separation rate $(n\xi_{\varepsilon,\delta})^{-2s/(2s+d)}$. By tracking conditions for each term dominating the others, the minimum separation rate that one can achieve using different bandwidths is summarized as
	\begin{equation}
		\begin{aligned}
		\rho_{\phi_{\mathttt{dpMMD}},L_2} \lesssim \begin{cases}
			n^{-\frac{2s}{4s+d}}, & \textrm{if } n^{-\frac{2s-d/2}{4s+d}} \lesssim  \xi_{\varepsilon,\delta} \ \textrm{(low privacy)}, \\[2mm]
			( n^{\frac{3}{2}} \xi_{\varepsilon,\delta})^{-\frac{s}{2s+d}}, & \textrm{if } n^{-\frac{1}{2}} \lesssim\xi_{\varepsilon,\delta} \lesssim n^{-\frac{2s-d/2}{4s+d}}  \ \textrm{(mid privacy)}, \\[2mm]
			\left(n\xi_{\varepsilon,\delta}\right)^{-\frac{2s}{2s+d}}, & \textrm{if }  \xi_{\varepsilon,\delta} \lesssim n^{-\frac{1}{2}}  \ \textrm{(high privacy)}.
		\end{cases}
		\end{aligned} \label{Eq: L2 separation rate}
	\end{equation}
	In particular, when $ \xi_{\varepsilon,\delta} \asymp n^{-1/2}$, the separation rate becomes $n^{-s/(2s+d)}$, which is known to be the minimax optimal rate of density estimation under the $L_2$ loss. We refer to \Cref{Section: Interpretation} for a detailed discussion of the separation rate. 
	\item It is important to note that all these separation rates are achieved by using different bandwidths, which requires knowledge of the smoothness parameter $s$. Building on the idea of \cite{ingster2000adaptive,schrab2021mmd,biggs2023mmdfuse}, one can develop an aggregated $\mathttt{dpMMD}$ test that is adaptive to $s$ without losing much power. The main idea would be to consider a wide range of private MMD statistics with different bandwidths and aggregate them properly. A detailed analysis of this approach is left for future research.
\end{itemize} 

{
We next develop a lower bound for the minimax separation over $\mathcal{P}_{L_2}^{s}$ in the high-privacy regime.
\begin{theorem}[High-Privacy Minimax Separation over $\mathcal{P}_{L_2}^{s}$]
\label{Theorem: High-privacy minimax separation over L2}
	Assume that $d\geq1$, $s\geq d/2$, $R,M>0$, $n\leq m$, $\varepsilon>0$, and $\delta\in[0,1)$. Let $\alpha,\beta\in(0,1)$ satisfy $\alpha+\beta<0.9$. Then there exists a constant $c_{d,s,R,M,\alpha,\beta}>0$ such that
	\begin{align*}
		\rho^\star_{L_2}
		\geq
		c_{d,s,R,M,\alpha,\beta}
		\min\Bigl\{1,\,
		(n\xi_{\varepsilon,\delta})^{-\frac{2s}{2s+d}}
		\Bigr\}.
	\end{align*}
\end{theorem}

Combining \Cref{Theorem: High-privacy minimax separation over L2} with the high-privacy upper bound in \eqref{Eq: L2 separation rate} and the trivial bound $\|p-q\|_{L_2}\leq 2\sqrt{M}$ yields, for fixed $\alpha$ and $\beta$, $s\geq d/2$, $n\leq m\leq\tau n$, and $B$ as in \Cref{Theorem: Minimax Separation over L2},
\begin{align*}
	\rho^\star_{L_2}
	\asymp
	\rho_{\phi_{\mathttt{dpMMD}},L_2}
	\asymp
	\min\Bigl\{1,\,
	(n\xi_{\varepsilon,\delta})^{-\frac{2s}{2s+d}}
	\Bigr\}
	\qquad
	\text{whenever }
	\xi_{\varepsilon,\delta}\lesssim n^{-1/2}.
\end{align*}
Consequently, the $\mathttt{dpMMD}$ test is minimax rate-optimal over the Sobolev class throughout the high-privacy regime. Together with its low-privacy optimality established above, this proves minimax rate-optimality in both the low- and high-privacy regimes. We emphasize, however, that \Cref{Theorem: High-privacy minimax separation over L2} requires $s\geq d/2$. Establishing high-privacy optimality when $s<d/2$, as well as optimality in the mid-privacy regime, remains open and is left for future work.
}

In contrast to the prior work~\citep{li2019optimality,schrab2021mmd} that utilizes a U-statistic for minimax two-sample testing, our approach is based on a plug-in estimator, also known as a V-statistic, of the MMD. While plug-in estimators can often exhibit suboptimal performance in estimation problems due to their inherent bias, they can still achieve optimal results in testing problems. This can be explained by the interplay between the test statistic and the critical value in a testing procedure, where the bias terms in these components may offset each other. \Cref{Theorem: Minimax Separation over L2} demonstrates this phenomenon by showing that the test based on the plug-in estimator of the MMD attains the minimax separation rate over $\mathcal{P}_{L_2}^s$ in the low privacy regime. Perhaps more interestingly, the plug-in estimator can outperform the U-statistic by having lower sensitivity and thus leading to greater power in high privacy regimes. This aspect of plug-in estimators has not been noticed in the literature, and we provide a more detailed discussion in the next subsection.

\subsection{Private Test based on the MMD U-statistic} \label{Section: Private Test based on the MMD U-statistic}
It has been shown that kernel tests based on U-statistics often produce optimal separation rates in non-DP settings~\citep{li2019optimality,schrab2021mmd,albert2019adaptive,kim2020minimax}. Therefore, one can naturally expect that their private extensions perform similarly well across different privacy regimes. In this section, we prove that this is not necessarily the case. In particular, we illustrate that the private MMD permutation test based on a U-statistic is provably outperformed by our approach based on the plug-in MMD estimate in high privacy regimes. A similar result for HSIC can be found in \Cref{Section: HSIC U-statistic}.

We begin with the explicit form of the MMD U-statistic, which is an unbiased estimator of $\mathrm{MMD}^2_k$, given as
\begin{align*}
	U_{\mathrm{MMD}}(\mathcal{X}_{n+m}) \coloneqq  \frac{1}{n(n-1)} \sum_{1 \leq i \neq j \leq n} k(Y_i,Y_j)  + \frac{1}{m(m-1)} \sum_{1 \leq i \neq j \leq m} k(Z_i,Z_j)  - \frac{2}{nm}\sum_{i=1}^n \sum_{j=1}^m k(Y_i,Z_j).
\end{align*}
The following lemma calculates the global sensitivity of $U_{\mathrm{MMD}}$, which is proved in \Cref{Section: Proof of Lemma: Global sensitivity of U_MMD}.

\begin{lemma}[Global Sensitivity of $U_{\mathrm{MMD}}$] \label{Lemma: Global sensitivity of U_MMD}
	Assume that the kernel $k$ is bounded as $0 \leq k(x,y) \leq K$ for all $x,y \in \mathbb{S}$. In addition, assume that $k$ is translation invariant, and have non-empty level sets on $\mathbb{S}$. Then there exists a positive sequence $c_{m,n} \in [4,8]$ such that for all $2 \leq n \leq m$,
	\begin{align*}
		\sup_{\bpi \in \boldsymbol{\Pi}_{n+m}} \sup_{\substack{\mathcal{X}_{n+m},\tilde{\mathcal{X}}_{n+m}:\\d_{\mathrm{ham}}(\mathcal{X}_{n+m},\tilde{\mathcal{X}}_{n+m}) \leq 1}} \bigl| U_{\mathrm{MMD}}(\mathcal{X}_{n+m}^{\bpi}) - U_{\mathrm{MMD}}(\tilde{\mathcal{X}}_{n+m}^{\bpi}) \bigr| = \frac{c_{m,n} K}{n}.
	\end{align*}
\end{lemma}
The lemma above indicates that the global sensitivity of the U-statistic has the same dependence on $n$ as that of the plug-in MMD in \Cref{Lemma: Sensitivity of MMD}. However, it is important to mention that their target parameters are different. The U-statistic is an estimator of $\mathrm{MMD}^2_k$, whereas the plug-in estimator given in \eqref{Eq: empirical MMD} estimates $\mathrm{MMD}_k$ without squaring. This key difference can lead to a significant gap in their power performance in privacy regimes as explored below.  

Given the sensitivity of $U_{\mathrm{MMD}}$ in \Cref{Lemma: Global sensitivity of U_MMD}, we consider the private permutation test in Algorithm~\ref{Algorithm: DP permutation test} using the test statistic $U_{\mathrm{MMD}}$ and the global sensitivity $\Delta_T = c_{m,n}Kn^{-1}$. Let us denote the resulting private test by $\phi_{\mathttt{dpMMD}}^u$. We analyze the minimum separations of $\phi_{\mathttt{dpMMD}}^u$ over $\mathcal{P}_{\mathrm{MMD}_k}$ and $\mathcal{P}_{L_2}^s$ in \Cref{Theorem: Suboptimality of U-MMD} and \Cref{Theorem: Minimum separation of U-stat over L2}, respectively, and compare them with those of $\phi_{\mathttt{dpMMD}}$ based on the plug-in estimator. Starting with the MMD alternative, the following theorem demonstrates that $\phi_{\mathttt{dpMMD}}^u$ fails to achieve the minimax separation rate over $\mathcal{P}_{\mathrm{MMD}_k}$. 
\begin{theorem}[Suboptimality of $\phi_{\mathttt{dpMMD}}^u$ against MMD Alternatives] \label{Theorem: Suboptimality of U-MMD}
	Assume that the kernel $k$ fulfills the conditions specified in \Cref{Lemma: Global sensitivity of U_MMD}. Moreover, assume that if $P,Q \in \mathcal{P}_{\mathbb{S}}$, then $w P + (1-w) Q \in \mathcal{P}_{\mathbb{S}}$ for all $w \in [0,1]$, and there exist $P_0,Q_0 \in \mathcal{P}_{\mathbb{S}}$ such that $\mathrm{MMD}_k(P_0,Q_0) = \varrho_0$ for some fixed $\varrho_0 >0$. Let $\alpha \in \bigl((B+1)^{-1},1\bigr)$, $\beta \in (0,1-\alpha)$ be fixed values and $n \leq m$. 
	Consider the high privacy regime where $\xi_{\varepsilon,\delta} \asymp n^{-1/2 - r}$ with fixed $r \in (0, 1/2)$, for $\xi_{\varepsilon,\delta}$ as in \eqref{Eq: definition of xi}.
	Then the uniform power of $\phi_{\mathttt{dpMMD}}^u$ is asymptotically at most $\alpha$ over $ \mathcal{P}_{\mathrm{MMD}_k}\!(\rho)$ where
	\begin{align} \label{Eq: separation MMD rate for U-statistics}
		\rho  =  \log(n) \times \max \Biggl\{ \sqrt{\frac{\max\!\big\{\!\log(1/\alpha), \, \log(1/\beta)\big\}}{n}}, \, \frac{\max\!\big\{\!\log(1/\alpha), \, \log(1/\beta)\big\}}{n \xi_{\varepsilon,\delta}}  \Biggr\}.
	\end{align}
	In other words, it holds that 
	\begin{align*}
	\limsup_{n \rightarrow \infty} \inf_{(P,Q) \in \mathcal{P}_{\mathrm{MMD}_k}\!(\rho)} \mE_{P,Q}[\phi_{\mathttt{dpMMD}}^u]  \leq \alpha.
	\end{align*}
\end{theorem}

\Cref{Theorem: Suboptimality of U-MMD}, proven in \Cref{Section: Proof of Theorem: Suboptimality of U-MMD}, clearly shows that $\phi_{\mathttt{dpMMD}}^u$ is not minimax optimal in the MMD metric as $\rho/ \rho^\star_{\mathrm{MMD}} \rightarrow \infty$ as $n \rightarrow \infty$. We also mention that the factor $\log(n)$ in $\rho$ is chosen for convenience, and it can be replaced by any other positive sequence that increases slower than $n^r$ for $r \in (0,1/2)$. The suboptimal performance of $\phi_{\mathttt{dpMMD}}^u$ primarily stems from the relatively high noise level associated with the Laplace mechanism. Intuitively, we expect that $\phi_{\mathttt{dpMMD}}^u$ is powerful in the private regime when the target parameter $\mathrm{MMD}^2_k(P,Q)$ is larger than the Laplace noise level $(n\xi_{\varepsilon,\delta})^{-1}$, equivalently, $\mathrm{MMD}_k(P,Q)$ is larger than $(n\xi_{\varepsilon,\delta})^{-1/2}$. Otherwise, the test statistic will be dominated by the Laplace noise. Importantly, the minimax separation under high privacy regimes in \Cref{Theorem: Minimax separation in MMD} is associated with $\min\{(n\xi_{\varepsilon,\delta})^{-1},1\}$, which is smaller than $(n\xi_{\varepsilon,\delta})^{-1/2}$. 
This briefly explains the suboptimality of $\phi_{\mathttt{dpMMD}}^u$ against the MMD alternative. Nevertheless, our analysis is limited to the U-statistic with the Laplace mechanism, and it is unknown whether the U-statistic in conjunction with other DP mechanisms can lead to optimality.

Turning to the $L_2$ alternative, let us denote the minimum separation of $\phi_{\mathttt{dpMMD}}^u$ with the Gaussian kernel over $\mathcal{P}_{L_2}^s$ as $\rho_{\phi_{\mathttt{dpMMD}}^u,L_2}$, which is similarly defined as \eqref{Eq: L2 separation}. Our next concern is characterizing $\rho_{\phi_{\mathttt{dpMMD}}^u,L_2}$ and comparing it with the minimum separation of the $\mathttt{dpMMD}$ test established in \Cref{Theorem: Minimax Separation over L2}. 

\begin{theorem}[Minimum Separation of $\phi_{\mathttt{dpMMD}}^u$ over $\mathcal{P}_{L_2}^{s}$] \label{Theorem: Minimum separation of U-stat over L2} Assume that $n \leq m \leq \tau n$ for some fixed constant $\tau \geq 1$, and that $\alpha \in (0,e^{-1})$, $\beta \in (0,1)$, $\varepsilon >0$, $\delta \in [0,1)$, $B \geq 6 \alpha^{-1} \log(2\beta^{-1})$ and $\prod_{i=1}^d \lambda_i \leq 1$. The minimum separation of $\phi_{\mathttt{dpMMD}}^u$ with the Gaussian kernel over $\mathcal{P}_{L_2}^{s}$ is upper bounded as 
	\begin{align*}
		\rho_{\phi_{\mathttt{dpMMD}}^u,L_2}^2 \leq C_{\tau,\beta,s,R,M,d} \Biggl\{ \sum_{i=1}^d \lambda_i^{2s} & +  \frac{\log(1/\alpha)}{n\sqrt{\lambda_1\cdots \lambda_d}} +\frac{\log(1/\alpha)}{n\lambda_1\cdots \lambda_d \xi_{\varepsilon,\delta}} \Bigg\},
	\end{align*}
	where $C_{\tau,\beta,s,R,M,d}$ is a positive constant, depending only on $\tau,\beta,s,R,M,d$, and $\xi_{\varepsilon,\delta}$ is as in \eqref{Eq: definition of xi}.
\end{theorem}

The proof of \Cref{Theorem: Minimum separation of U-stat over L2} is given in \Cref{Section: Proof of Theorem: Minimum separation of U-stat over L2}. To simplify our discussion, assume that $\alpha$ is a fixed constant. In this case, by comparing \Cref{Theorem: Minimum separation of U-stat over L2} with \Cref{Theorem: Minimax Separation over L2}, the upper bound for $\rho_{\phi_{\mathttt{dpMMD}}^u,L_2}^2$ is smaller than that for $\rho_{\phi_{\mathttt{dpMMD}},L_2}^2$, up to a constant, only when $n\xi_{\varepsilon,\delta} \lesssim 1$. Since $\prod_{i=1}^d \lambda_i \leq 1$ and $\alpha$ is fixed, the condition $n\xi_{\varepsilon,\delta} \lesssim 1$ essentially means that $\|p - q\|_{L_2}$ needs to be sufficiently larger than a specific constant for significant power. However, this condition may be infeasible as we assume that $\|p\|_{L_\infty}$ and $\|q\|_{L_\infty}$ are bounded by $M$. In fact, our earlier result in \Cref{Theorem: Properties of dpMMD} suggests that the test is not even consistent in a pointwise sense when $n\xi_{\varepsilon,\delta} \lesssim 1$. Therefore, except for this boundary case, it is more beneficial to use $\phi_{\mathttt{dpMMD}}$ than $\phi_{\mathttt{dpMMD}}^u$ to achieve tighter separation rates over $\mathcal{P}_{L_2}^{s}$ in high privacy regimes. 

As we mentioned before, it remains an open question whether there exist alternative privacy mechanisms that could potentially yield an optimal test based on $U_{\mathrm{MMD}}$ in high privacy regimes. While we still advocate using the plug-in MMD estimate over $U_{\mathrm{MMD}}$ due to its smaller sensitivity, it would be interesting to explore this question in future work.

\section{Simulations} \label{Section: Simulations}

In this section, we compare the empirical power of dpMMD against other private two-sample tests, including 
the naive dpMMD introduced in \Cref{Eq: naive permutation p-value}
and the U-statistic dpMMD studied in \Cref{Section: Private Test based on the MMD U-statistic}.
We also implement two generic methods for privatizing the permutation MMD test, which we refer to as TOT MMD \citep{kazan2023test} and SARRM MMD \citep{pena2022differentially}. \blue{We additionally consider TOT--SARRM, which applies the SARRM hyperparameter-selection rule to TOT and therefore isolates the effect of the binomial test from that of hyperparameter selection.}
Finally, we also compare against the differentially private kernel two-sample test, TCS-ME, proposed by \citet{raj2020differentially}. A brief overview of these alternative methods is provided in \Cref{subsec:alternative_dp_tests}, with detailed information available in \Cref{Section: Alternative private tests}.

We empirically study the power attained by dpMMD on 
synthetic data sampled from perturbed uniform distributions in \Cref{subsec:mmd_uniform}, 
and on real-world high-dimensional CelebA image data in \Cref{subsec:mmd_celeba}. The latter scenario involves sensitive information on human faces, justifying the incorporation of differential privacy in the analysis. In both simulation settings, we observe consistent patterns in the power behavior, and a detailed discussion on the results is presented in \Cref{subsec:analysis_experiments}. For all our experiments, we use a Gaussian kernel, $B=2000$ permutations, and report power results averaged over $200$ repetitions.

Due to space constraints, we defer the dpHSIC simulations to \Cref{Section: Additional Simulations}, along with test-level analysis and additional low-privacy experiments. The code to run our tests and to reproduce the experiments is available at \url{https://github.com/antoninschrab/dpkernel}.

\subsection{Alternative differentially private tests}
\label{subsec:alternative_dp_tests}

We provide a brief introduction to the alternative differentially private tests, namely TCS-ME, TOT, and SARRM, with detailed implementation information available in \Cref{Section: Alternative private tests}.

\paragraph{TOT \citep{kazan2023test}.}
TOT is constructed based on the subsample-and-aggregate idea outlined in \citet{canonne2019structure}.
It is guaranteed to be differentially private and to correctly control the probability of type I error for any sample size, any number of partitioned subsets, and any sub-test significance level.
However, for non-parametric testing, there is no principled way to choose the last two parameters and one has to rely on heuristics in practice. This heuristic aspect presents a notable disadvantage of TOT, being highly sensitive to the choice of these parameters. 

\paragraph{SARRM \citep{pena2022differentially}.}
As another method based on the subsample-and-aggregate idea, SARRM also depends on the number of partitioned subsets and on the sub-test significance level. \citet{pena2022differentially} propose a method to select these parameters, which can be implemented for MMD/HSIC tests.  
However, the differential privacy constraint and type I error control for SARRM are only guaranteed for sufficiently large sample sizes (determined by a minimum number of partitioned subsets to use), which means that SARRM simply cannot be run in some settings depending on the values of $\varepsilon$, $\alpha$, and $n$.

\paragraph{TCS-ME \citep{raj2020differentially}.}
The TCS-ME test is a privatized version of the ME test~\citep{jitkrittum2016interpretable} that utilizes kernel mean embeddings. This method builds on a Hotelling-type test statistic privatized by the Gaussian mechanism, and requires a careful choice of test locations. The resulting test is $(\varepsilon, \delta)$-DP for $\delta>0$ (run with $\delta=10^{-5}$).
A major limitation of TCS-ME is its potential for significant miscalibration, particularly in high-dimensional settings. See our discussion on type I error control in \Cref{subsec:analysis_additional_experiments}.

\FloatBarrier

\subsection{Perturbed Uniform Distributions}
\label{subsec:mmd_uniform}

As recalled in \Cref{Section: Separation in L2 metric}, the minimax $L_2$ separation rate over the Sobolev ball $\mathcal{S}_d^s(R)$ is $n^{-2s/(4s+d)}$ in the non-privacy regime.
A lower bound for this minimax separation rate is derived by constructing two densities whose difference lies in the Sobolev ball with a small $L_2$ norm. As explained by \citet[Appendix D]{schrab2021mmd}, the uniform and perturbed uniform densities meet the requirements in the lower bound construction, and we adopt this setting in our two-sample experiments. In more detail, we compare the uniform distribution with its perturbed counterpart, varying the amplitudes. 

\begin{figure}[t!]
	\centering
	\includegraphics[width=\textwidth]{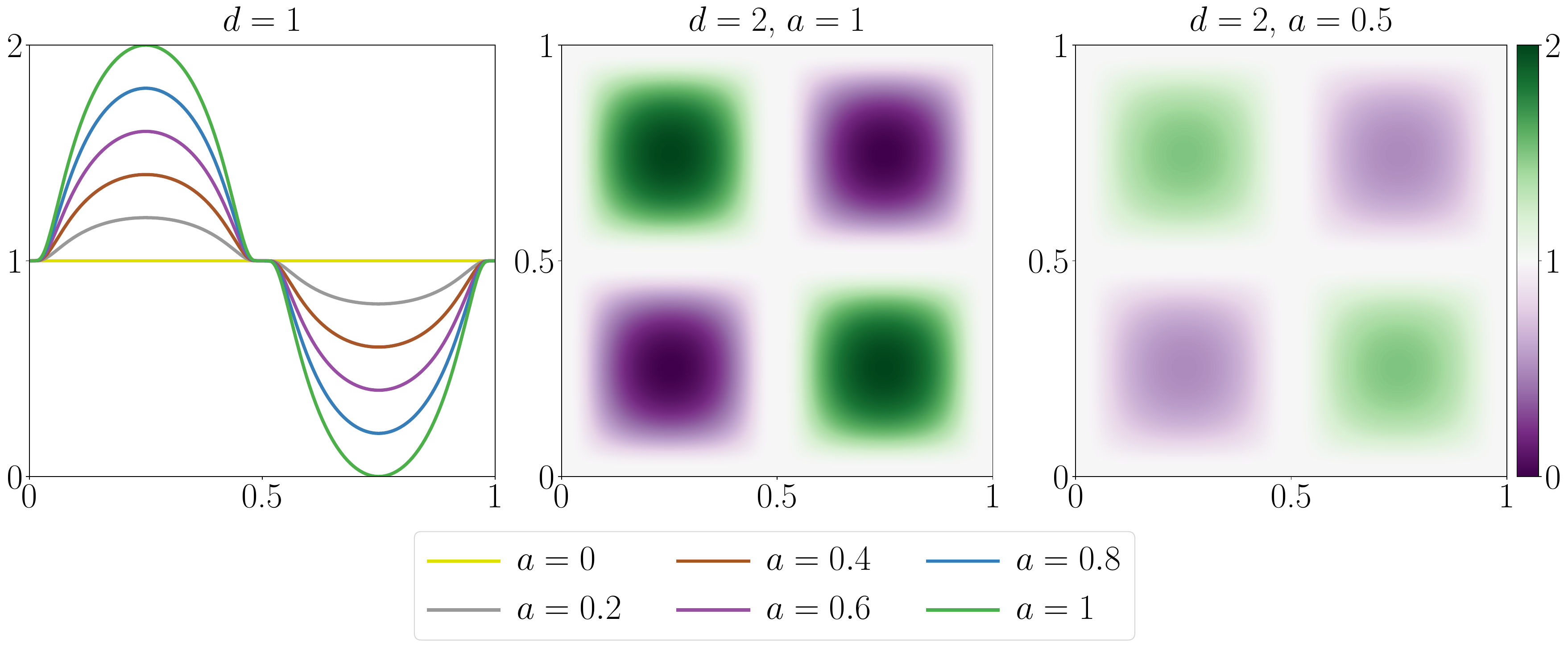}
	\captionsetup{format=hang}
	\caption{
		Perturbed uniform $d$-dimensional densities on $[0,1]^d$ with varying perturbation amplitude $a$.
	}
	\label{fig:perturbed_densities}
\end{figure}

\begin{figure}[t!]
	\centering
	\includegraphics[width=\textwidth]{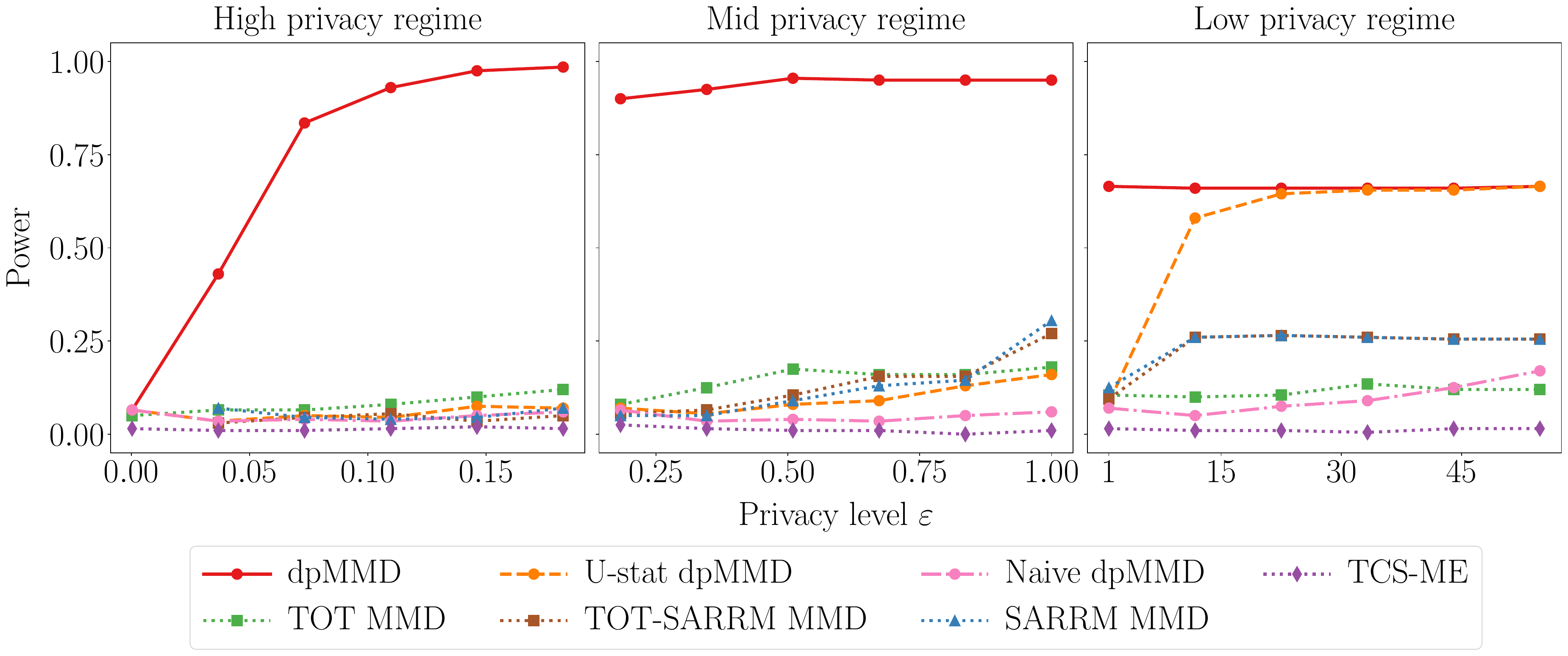}
	\captionsetup{format=hang}
	\caption{
		Comparing uniform vs.~perturbed uniform while varying the privacy level $\varepsilon$. We set the sample sizes $m = n = 3000$ and dimension $d=1$, and change the privacy level $\varepsilon$ and perturbation amplitude $a$ as follows: 
		\emph{(Left)} Privacy level $\varepsilon$ from $1/n$ to $10/\sqrt{n}$, perturbation amplitude $a=0.2$. 
		\emph{(Middle)} Privacy level $\varepsilon$ from $10/\sqrt{n}$ to $1$, perturbation amplitude $a=0.15$. 
		\emph{(Right)} Privacy level $\varepsilon$ from $1$ to $\sqrt{n}$, perturbation amplitude $a=0.1$. 
	}
	\label{fig:uniform_mmd_privacy}
\end{figure}

\begin{figure}[t!]
	\centering
	\includegraphics[width=\textwidth]{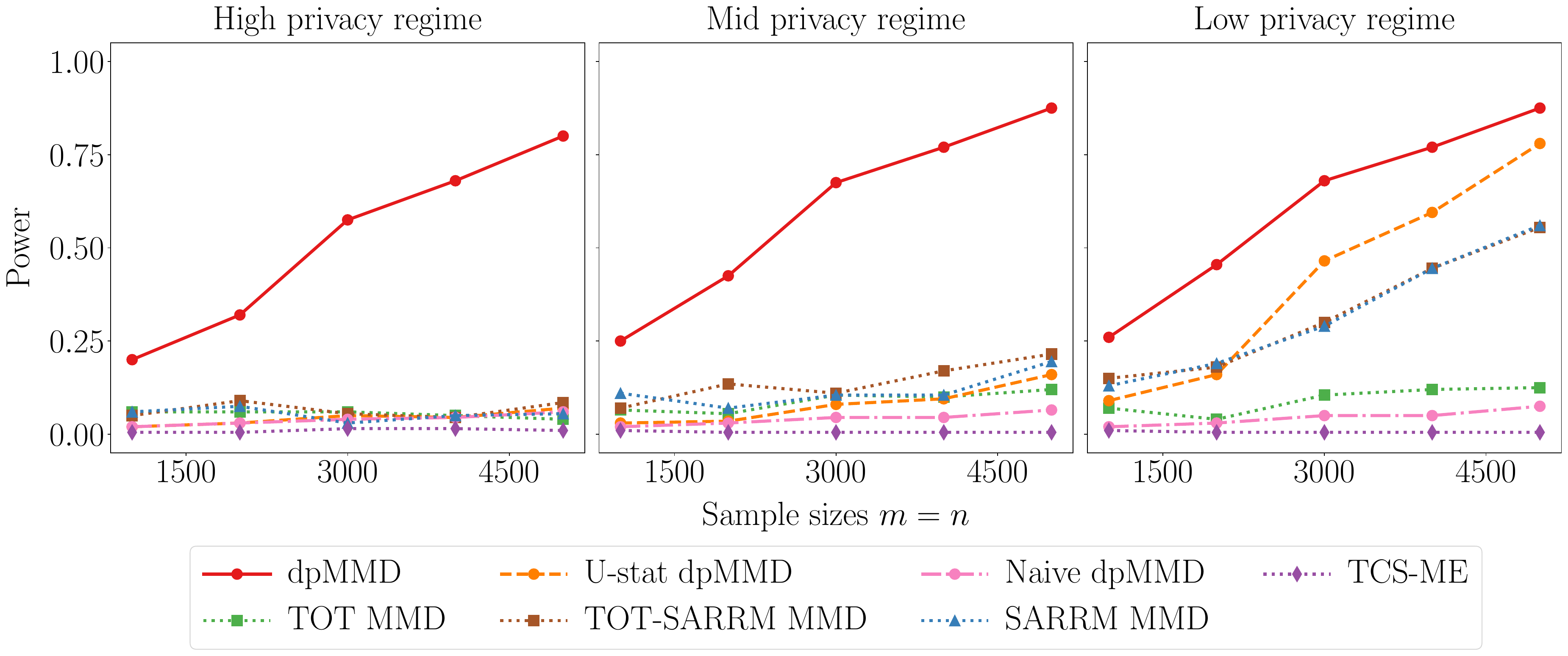}
	\captionsetup{format=hang}
	\caption{
		Comparing uniform vs.~perturbed uniform while varying the sample sizes $m=n$. We set the dimension $d=1$ and perturbation amplitude $a=0.1$. We change the privacy level as follows:
		\emph{(Left)} Privacy level $\varepsilon=10/\sqrt{n}$. 
		\emph{(Middle)} Privacy level $\varepsilon=1$. 
		\emph{(Right)} Privacy level $\varepsilon=\sqrt{n}/10$. 
	}
	\label{fig:uniform_mmd_sample_size}
\end{figure}

\begin{figure}[t!]
	\centering
	\includegraphics[width=\textwidth]{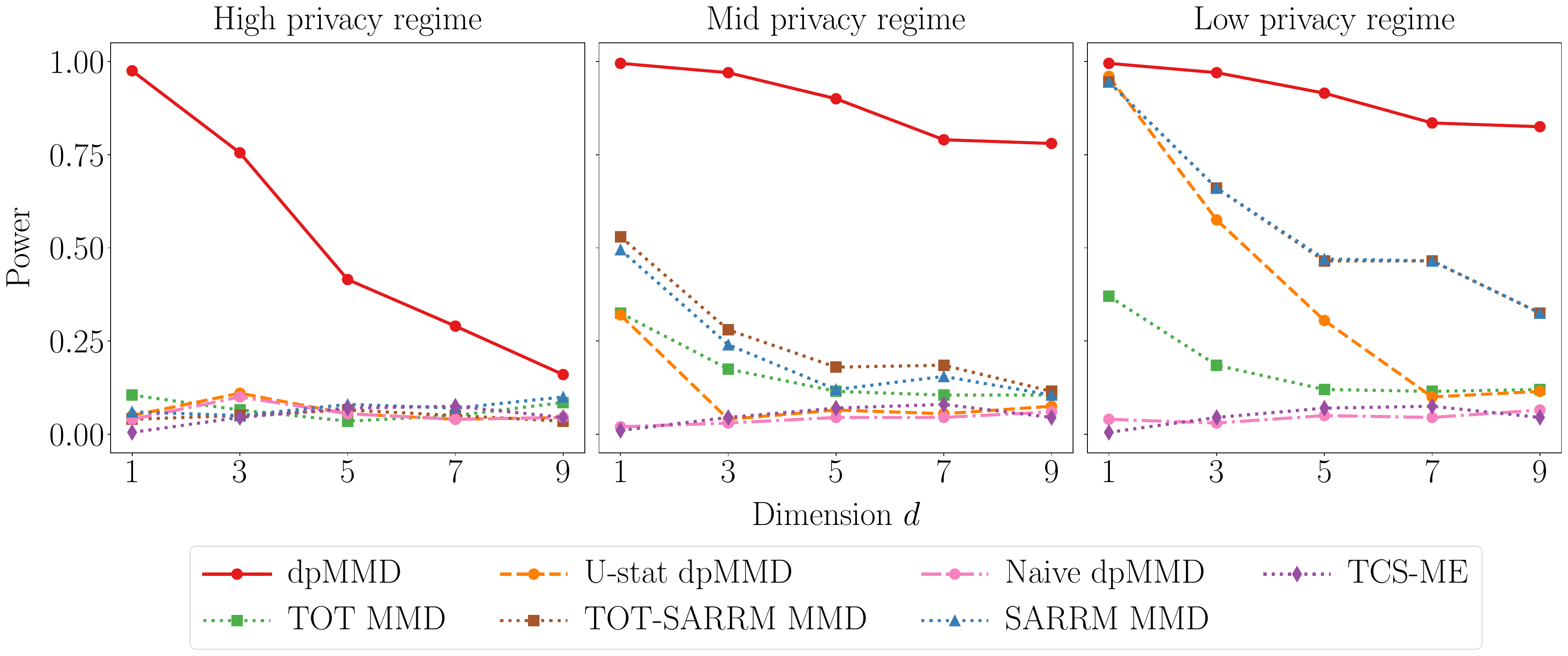}
	\captionsetup{format=hang}
	\caption{
		Comparing uniform vs.~perturbed uniform while varying the dimension $d$. We set the sample sizes $m = n = 3000$ and perturbation amplitude $a=0.2$. We change the privacy level as follows:
		\emph{(Left)} Privacy level $\varepsilon=10/\sqrt{n}$. 
		\emph{(Middle)} Privacy level $\varepsilon=1$. 
		\emph{(Right)} Privacy level $\varepsilon=\sqrt{n}/10$. 
	}
	\label{fig:uniform_mmd_dimension}
\end{figure}

Specifically, the considered uniform distribution on $[0,1]^d$ has density $\mathds{1}(x \in [0,1]^d)$ for $x\in\mathbb{R}^d$, while
the perturbed uniform density on $[0,1]^d$ with a perturbation amplitude $a\in[0,1]$ is
\begin{equation}
	\label{eq:perturbed_uniform_density}
	\mathds{1}\bigl(x \in [0,1]^d\bigr) + a\, \prod_{i=1}^dP(x_i), \quad \text{for $x = (x_1,\ldots,x_d)^\top \in\mathbb{R}^d$,}
\end{equation}
where the one-dimensional perturbation is defined as
\begin{align*}
	P(x_i) \coloneqq \exp\!\left(1-\frac{1}{1-(4x_i-1)^2} \right) \!\mathds{1}\bigl(x_i \in (0,1/2)\bigr)
	- \exp\!\left(1-\frac{1}{1-(4x_i-3)^2} \right) \!\mathds{1}\bigl(x_i \in (1/2,1)\bigr).
\end{align*}
This definition matches that of \citet[Equation 17]{schrab2021mmd} with only one (scaled) perturbation per dimension.
The one-dimensional and two-dimensional perturbed densities with various perturbation amplitudes are visualized in \Cref{fig:perturbed_densities}.

We run our perturbed uniform experiments under three different settings where we vary
the privacy level $\varepsilon$ (\Cref{fig:uniform_mmd_privacy}),
the sample sizes $m=n$ (\Cref{fig:uniform_mmd_sample_size}),
and the dimension $d$ (\Cref{fig:uniform_mmd_dimension}). Additionally, we provide experiments with `strong-signal' alternatives in the low privacy regime in \Cref{fig:uniform_mmd_appendix} in \Cref{subsec:low_privacy}, as well as a level analysis in \Cref{fig:uniform_mmd_level} in \Cref{subsec:level}.
We discuss all experimental results of \Cref{fig:uniform_mmd_privacy,fig:uniform_mmd_sample_size,fig:uniform_mmd_dimension} in \Cref{subsec:analysis_experiments}.

\begin{figure}[t!]
	\centering
	\includegraphics[width=0.6\textwidth]{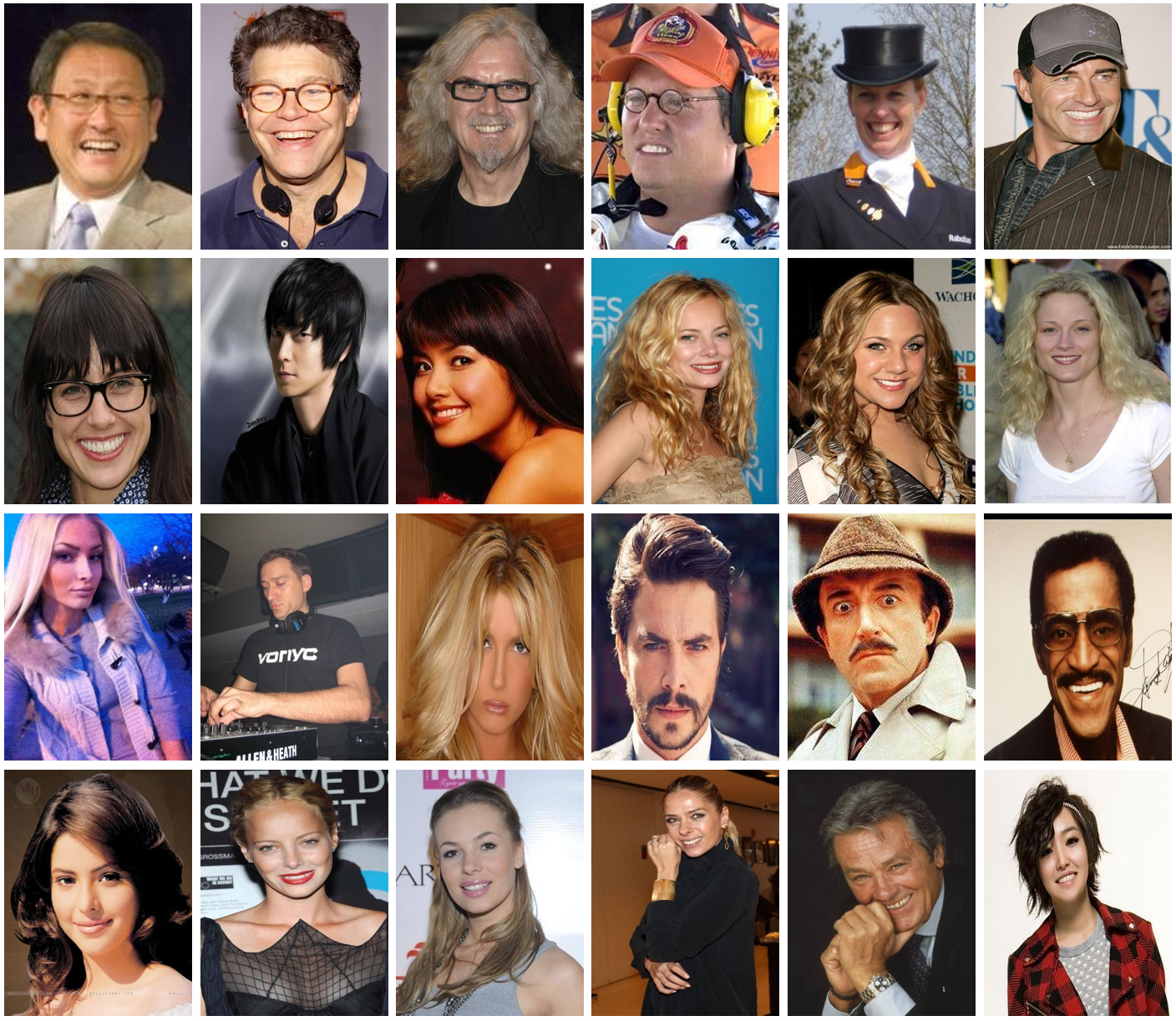}
	\captionsetup{format=hang}
	\caption{Selected 
		\href{https://mmlab.ie.cuhk.edu.hk/projects/CelebA.html}{CelebA} images in dimension $3\times 178\times 218$. 
	}
	\label{fig:celeba}
\end{figure}

\subsection{CelebA}
\label{subsec:mmd_celeba}

As some potential real-world applications of differentially private two-sample tests, we consider \href{https://mmlab.ie.cuhk.edu.hk/projects/CelebA.html}{CelebA} face images which in practice would be highly confidential, and hence the use of DP tests is thoroughly justified.
The CelebA dataset \citep{liu2015faceattributes} consists of $202,599$ face images of $10,177$ identities with a large diversity of face attributes, poses and backgrounds.
For illustration purposes, we display a selection of CelebA images in \Cref{fig:celeba}. It is worth highlighting that we run our tests on the original full-resolution images $(3\times 178\times 218)$ without any modifications.

In our experiments, one sample consists of uniformly-sampled face images of women, while the other is `corrupted' with corruption parameter $c\in[0,1]$ in the following sense:
we uniformly sample face images of women with probability $1-c$, and of men with probability $c$.

\begin{figure}[t!]
	\centering
	\includegraphics[width=\textwidth]{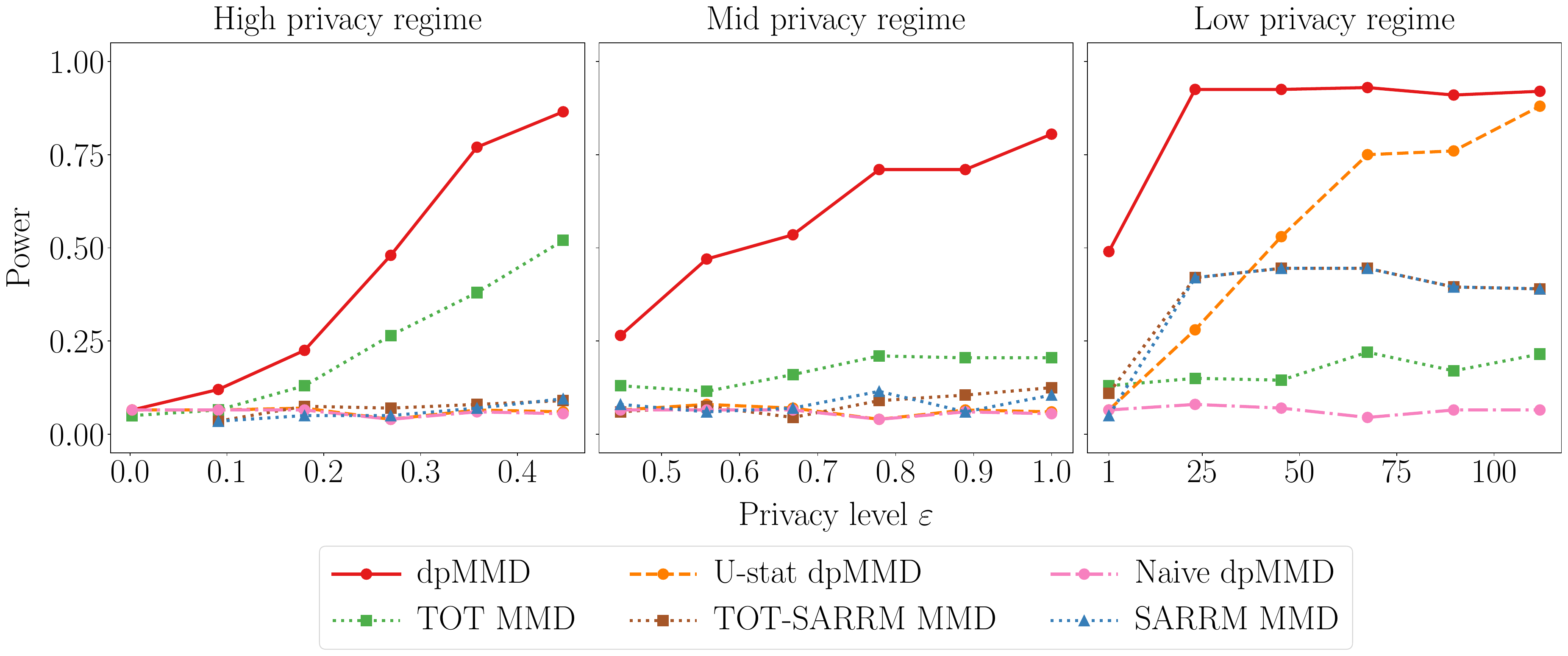}
	\captionsetup{format=hang}
	\caption{
		Comparing CelebA women/men images while varying the privacy level $\varepsilon$. We set the sample sizes $m = n = 500$, and change the parameters as follows:
		\emph{(Left)} Privacy level $\varepsilon$ from $1/n$ to $10/\sqrt{n}$, corruption $c=1$. 
		\emph{(Middle)} Privacy level $\varepsilon$ from $10/\sqrt{n}$ to $1$, corruption $c=0.6$. 
		\emph{(Right)} Privacy level $\varepsilon$ from $1$ to $\sqrt{n}$, corruption $c=0.5$.
	}
	\label{fig:celeba_mmd_privacy}
\end{figure}

\begin{figure}[t!]
	\centering
	\includegraphics[width=\textwidth]{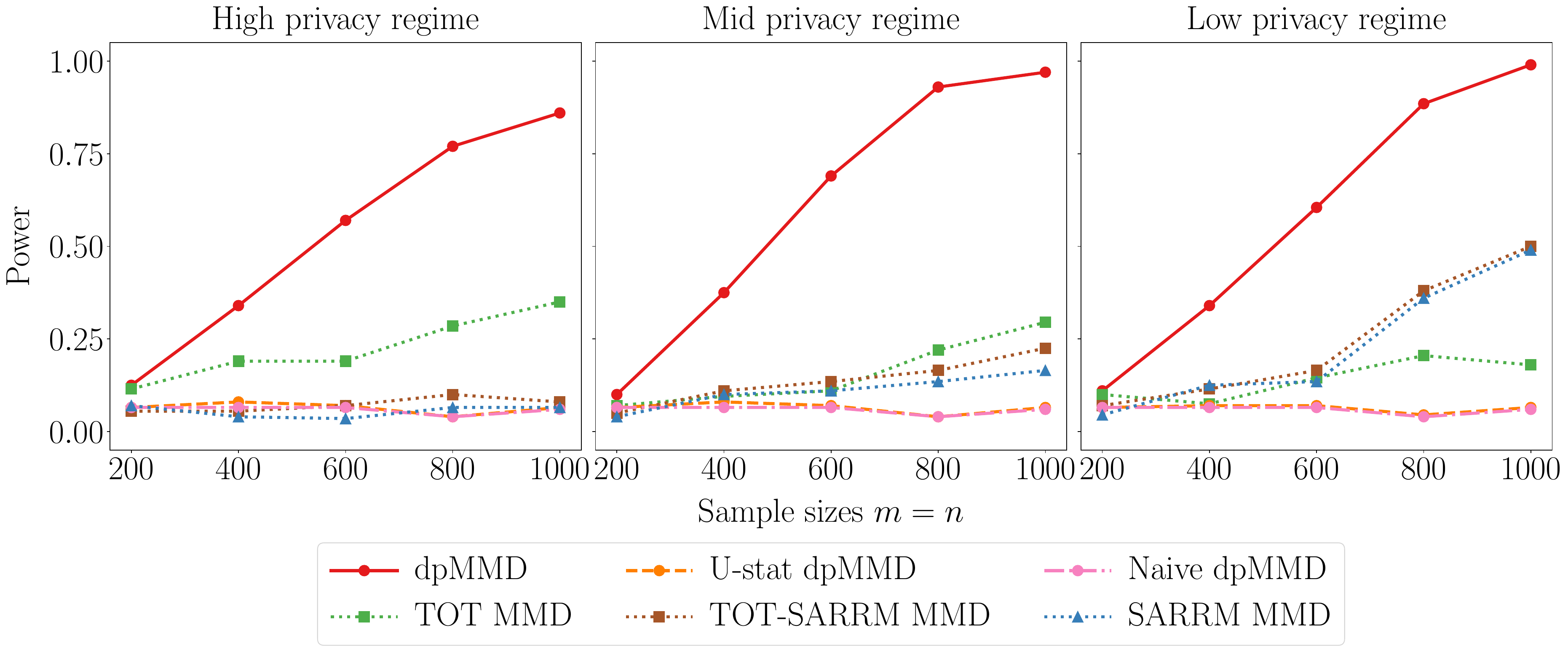}
	\captionsetup{format=hang}
	\caption{
		Comparing CelebA women/men images while varying the sample sizes $m=n$. We set the other parameters as follows:
		\emph{(Left)} Privacy level $\varepsilon=10/\sqrt{n}$, corruption $c=0.7$. 
		\emph{(Middle)} Privacy level $\varepsilon=1$, corruption $c=0.5$. 
		\emph{(Right)} Privacy level $\varepsilon=\sqrt{n}/10$, corruption $c=0.4$. 
	}
	\label{fig:celeba_mmd_sample_size}
\end{figure}

\begin{figure}[t!]
	\centering
	\includegraphics[width=\textwidth]{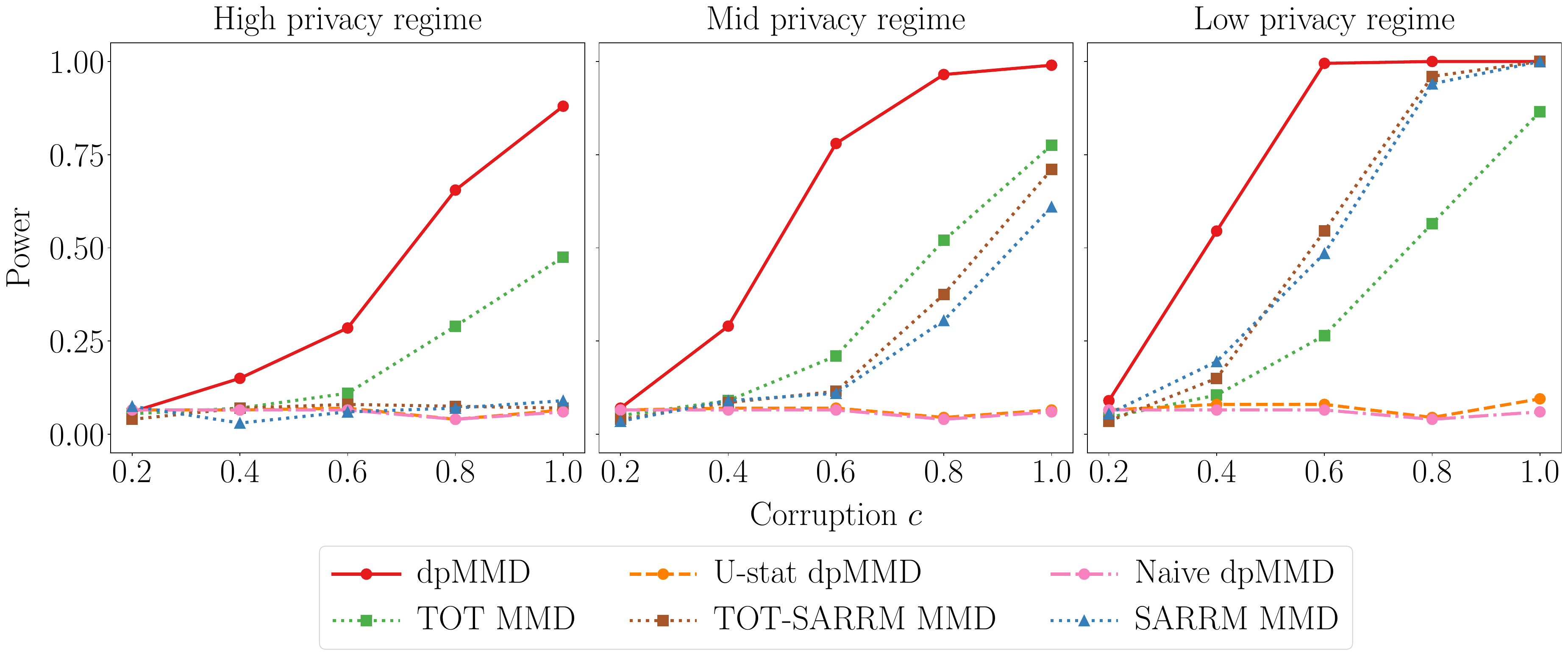}
	\captionsetup{format=hang}
	\caption{
		Comparing CelebA women/men images while varying the corruption parameter $c$. We set the sample sizes $m = n = 500$, and privacy parameter $\varepsilon$ as follows:
		\emph{(Left)} Privacy level $\varepsilon=10/\sqrt{n}$. 
		\emph{(Middle)} Privacy level $\varepsilon=1$. 
		\emph{(Right)} Privacy level $\varepsilon=\sqrt{n}/10$. 
	}
	\label{fig:celeba_mmd_corruption}
\end{figure}

We run several CelebA experiments while varying
the privacy level $\varepsilon$ (\Cref{fig:celeba_mmd_privacy}),
the sample sizes $m=n$ (\Cref{fig:celeba_mmd_sample_size}),
and the corruption $c$ (\Cref{fig:celeba_mmd_corruption}).
As in \Cref{subsec:mmd_uniform}, we consider the high/mid/low privacy regimes for each of these.
We also verify that all tests are well-calibrated in \Cref{fig:uniform_mmd_level} in \Cref{subsec:level}.
TCS-ME is excluded from our power analysis on the CelebA data as we empirically observed that this method is not well-calibrated for this dataset. 

Our results consistently demonstrate that kernel tests using a simple Gaussian kernel are able to capture complex image distribution shifts (see \Cref{fig:celeba}) even in this extremely high dimensional setting with $d = 3\times 178\times 218 = 116,412$. This result is surprising given that the Gaussian kernel simply compares the distance between pixels at the same location without using information about the image structure. We discuss this aspect more in \Cref{subsec:analysis_experiments}.

\subsection{Analysis of Main Experimental Results}
\label{subsec:analysis_experiments}

We now analyze the results of the perturbed uniform and CelebA experiments presented in  \Cref{subsec:mmd_uniform,subsec:mmd_celeba}, respectively. 

\paragraph{Overview.}
First and foremost, we observe that dpMMD achieves significantly higher power than all other tests across all privacy regimes. This trend remains consistent when varying the other parameters $\{$privacy, sample sizes, dimension, corruption$\}$. 
In the high and mid privacy regimes illustrated in \Cref{fig:uniform_mmd_privacy,fig:uniform_mmd_sample_size,fig:uniform_mmd_dimension}, only dpMMD is able to detect the perturbation on the uniform distribution.
In the high and mid privacy regimes presented in  \Cref{fig:celeba_mmd_privacy,fig:celeba_mmd_sample_size,fig:celeba_mmd_corruption}, dpMMD clearly outperforms all other tests but TOT MMD and SARRM MMD eventually manage to detect the CelebA image distributional shift.
In the low privacy regime, dpMMD also achieves the highest power, which is eventually matched by U-stat dpMMD as the privacy parameter increases. 

\paragraph{Varying privacy and free privacy.}
As can be seen in \Cref{fig:uniform_mmd_privacy,fig:celeba_mmd_privacy}, increasing $\varepsilon$ (\emph{i.e.}, lowering privacy) first leads to an increase in power. However, we observe that, after some threshold, further increasing $\varepsilon$ does not increase the power of dpMMD. Essentially, in this low privacy regime, dpMMD already attains the power of the non-private MMD test and the (low) privacy guarantee then comes for free. 
This empirical observation is also supported by our theoretical findings that the non-private MMD and $L_2$ optimal uniform separation rates are attained by dpMMD in the low privacy regime (see \Cref{Theorem: Uniform separation for MMD,Theorem: Minimax Separation over L2}).

\paragraph{Varying privacy for U-stat dpMMD.}
Furthermore, \Cref{fig:uniform_mmd_privacy,fig:celeba_mmd_privacy} shows that the dpMMD test using the U-statistic is powerless in the high and mid privacy regimes, but it eventually reaches the power of dpMMD (using the V-statistic) in the low privacy regime, which is theoretically justified by our suboptimality result~(\Cref{Theorem: Suboptimality of U-MMD}) for dpMMD based on the U-statistic and by the fact that the non-private MMD U-statistic test is minimax rate optimal \citep{schrab2021mmd}.

\paragraph{Varying the sample size.}
When varying the sample sizes $m=n$ in \Cref{fig:uniform_mmd_sample_size,fig:celeba_mmd_sample_size} with fixed high/mid/low privacy level $\varepsilon\in\{10/\sqrt{n},1,\sqrt{n}/10\}$, the power of all tests naturally increases, while the power of dpMMD increases faster than the others. In the low privacy regime of \Cref{fig:uniform_mmd_sample_size}, we see that the power of U-stat dpMMD approaches that of dpMMD as the sample size increases. This can be explained by the aforementioned reasoning along with the observation that the privacy level $\varepsilon=\sqrt{n}/10$ also increases (\emph{i.e.}, lower privacy) in this setting.

\paragraph{Varying the problem difficulty.}
In \Cref{fig:uniform_mmd_dimension,fig:celeba_mmd_corruption}, the sample sizes and privacy levels are fixed while the difficulty of the problem is varied.
For the perturbed uniform experiment, as the dimension of the problem increases, the perturbation becomes more difficult to detect and hence the power decreases for each test. 
We observe nonetheless that the power of dpMMD deteriorates at a much slower rate than that of the other tests.
For the CelebA experiment, the test power increases with the corruption parameter of the image sampler, with dpMMD always achieving the highest power, followed by either TOT MMD or SARRM MMD.

\paragraph{The power of kernel methods.}
We end this discussion with some remarks regarding the CelebA experiments of \Cref{subsec:mmd_celeba}.
First, we emphasize again that the quadratic-time kernel tests run swiftly on the full-resolution CelebA image data, which has $116,412$ pixels per image.
Second, given the large diversity of faces, poses and backgrounds (see \Cref{fig:celeba}), it is remarkable that dpMMD is able to detect such complicated differences in high-dimensional image distributions while using an off-the-shelf Gaussian kernel which entirely ignores the image structure and simply averages distances between pixel values at the same locations.
Third, the fact that such complex testing problems can now be solved while guaranteeing differential privacy is extremely important, especially when dealing with data as personal and sensitive as facial data.
The DP constraint essentially guarantees privacy in the sense that confidential information about a single face image cannot be recovered.
Fourth, we stress that the power results reported are meaningful and that dpMMD truly detects the difference between CelebA images of women and men, which is justified as type I error control is correctly retained (\Cref{subsec:level}).
We believe these CelebA experiments strongly advocate the use of tests leveraging kernel methods for differential privacy.

\subsection{Analysis of Additional Experimental Results}
\label{subsec:analysis_additional_experiments}

Before concluding the paper, we briefly summarize the results of the additional experiments in \Cref{Section: Additional Simulations}.

\paragraph{Independence testing.}
In the HSIC testing experiments of \Cref{subsec:hsic_uniform}, we consider the problem of detecting the dependence of variables with the perturbed uniform joint density introduced in \Cref{subsec:mmd_uniform}, and hence with uniform marginal densities.
This setting is of particular interest as this corresponds to the joint density used in deriving the non-private $L_2$ minimax independence lower bound over Sobolev balls by \citet{albert2019adaptive}.
The exact same power dynamics and aforementioned observations, which hold for the MMD-based tests, also apply to the HSIC variants of all tests, and indeed dpHSIC achieves substantially higher power than all the other tests.

\paragraph{High-signal \& low-privacy alternatives.}
We remark that Naive dpMMD and TCS-ME have almost no power against the alternatives considered in the previous subsections. This is not because these tests are faulty but because the signal is too weak to be detected by these tests. In fact, the Naive dpMMD test with privacy $\varepsilon$ is exactly equivalent to the dpMMD test with privacy $2/(B\varepsilon)$ (recall $B=2000$ permutations), which justifies the poor performance of Naive dpMMD. For a sanity check, we consider `high-signal' \& `low-privacy' alternatives  in \Cref{subsec:low_privacy} and show that dpMMD and TCS-ME are indeed able to detect the difference when the signal is large enough. The results can be found in \Cref{fig:uniform_mmd_appendix}.

\paragraph{Type I error control.}
In \Cref{subsec:level}, we run experiments under the null hypothesis.
This corresponds to no perturbation for the perturbed uniform two-sample and independence settings (amplitude $a=0$), and to no corruption in the CelebA sampler (corruption $c=0$).
All tests are well-calibrated with empirical level around $\alpha$ under all settings considered, except TCS-ME.
Indeed, we observe in \Cref{fig:uniform_mmd_level} that when testing two samples from a 100-dimensional uniform distribution, TCS-ME fails to control the type I error rate, which is estimated to be around $0.5=10\alpha$ instead of around $\alpha=0.05$.
This depicts a major limitation of TCS-ME especially in high-dimensional settings. However, we point out that this test is well-calibrated in the settings of \Cref{subsec:mmd_uniform}, ensuring a fair comparison of power therein.

\section{Discussion} \label{Section: Discussion}
In this work, we have proposed differentially private permutation tests and examined their theoretical and empirical performance. The prior work on differentially private testing has often been limited in its practical applicability, being restricted to discrete data or relying on asymptotic theory that does not offer confidence in finite-sample settings. Our permutation framework addresses these challenges by introducing practical tools, which are applicable to diverse settings with finite-sample guarantees for both type I error control and differential privacy. In addition to general power properties, we have provided a detailed power analysis of the proposed method in the context of kernel testing, and shown that the proposed private kernel tests achieve minimax optimal power in terms of kernel metrics in all privacy regimes. We have also analyzed the testing power against nonparametric $L_2$ alternatives, and established minimum separation rates in all privacy regimes. Finally, we have conducted an extensive simulation study to validate our theoretical findings as well as to highlight the practical value of our approach. 

Our work raises several intriguing open questions that deserve further investigation, as outlined below.
\begin{itemize}
	\item \blue{\textbf{Beyond Global DP.} Our main results focus on global $(\varepsilon,\delta)$-differential privacy with the Laplace mechanism, while \Cref{Section: Renyi Differentially Private Permutation Tests} extends the framework to R\'{e}nyi differential privacy. Developing private permutation tests under further privacy notions and studying their minimax properties remain promising directions.}
	\item \textbf{Other Applications.} We illustrated the proposed method in the context of two-sample and independence testing, with a focus on kernel-based tests. The permutation method has been employed successfully in other statistical problems such as testing for regression coefficients~\citep{diciccio2017robust} and conditional independence testing~\citep{kim2022local}. It is therefore compelling to broaden the application of our framework by tackling other statistical problems in privacy settings. Future work can also focus on conducting a detailed analysis of the private permutation test using non-kernel test statistics.
	\item \textbf{Variants of dpMMD and dpHSIC.} In our analysis, we utilized the plug-in estimators of MMD and HSIC based on single kernels. In recent years, significant progress has been made to reduce the computational complexity~\citep{schrab2022efficient,domingo2023compress} as well as to avoid the bandwidth selection issue~\citep{schrab2021mmd,biggs2023mmdfuse} of this standard approach. Considering these developments, a promising avenue for future research would be to extend these recent advances to privacy-preserving settings.
	\item \blue{\textbf{Computationally Efficient DP Tests.} The plug-in MMD and HSIC estimators have quadratic computational cost. Adapting recent sub-quadratic kernel tests~\citep{schrab2022efficient,domingo2023compress,choi2024computational} to the private setting could clarify the trade-off among privacy, computation, and statistical power.}
	\item \blue{\textbf{Kernel Choice.} A large Gaussian bandwidth is often effective for broad distributional differences, whereas a small bandwidth emphasizes localized discrepancies. Kernel selection and aggregation are already challenging without privacy~\citep{gretton2012optimal,sutherland2016generative,schrab2021mmd,biggs2023mmdfuse}; under differential privacy, adapting over more kernels additionally requires privacy accounting and can increase the injected noise. Developing principled private kernel-selection methods is therefore an important direction.}
	\item \blue{\textbf{Minimax Separation under DP.} The $L_2$ upper bounds now match lower bounds in the low-privacy regime and, for $s\geq d/2$, in the high-privacy regime. Whether the intermediate-privacy rate is optimal remains open. More broadly, establishing sharp minimax rates under differential privacy for other metrics would deepen our understanding of privacy-induced statistical difficulty.}
\end{itemize}

\paragraph{Acknowledgements.} The authors thank Aaditya Ramdas for his valuable feedback and encouragement, and the anonymous reviewers for their thoughtful comments and constructive suggestions. \blue{Ilmun Kim acknowledges support from the Korean government (RS-2023-00211073) and KAIST startup funding (KAIST-G04250059). Antonin Schrab acknowledges support from the U.K. Research and Innovation under the UCL CDT in Foundational AI (grant EP/S021566/1) and under the AI World-Leading Researcher Fellowship on Advancing Modern Data-Driven Robust AI (grant G111021).}

\bibliographystyle{apalike}
\bibliography{reference}

\clearpage

\appendix

\allowdisplaybreaks

\section{Overview of Appendices}

This appendix includes additional results, technical details, and proofs omitted from the main text. The remaining material is organized as follows. 
\begin{itemize}
	\item \blue{In \Cref{Section: Renyi Differentially Private Permutation Tests}, we extend the proposed framework to R\'{e}nyi differential privacy.}
	\item In \Cref{Section: Additional Results}, we present additional results including 
	\begin{enumerate}[(i)]
		\item Limiting null distributions of privatized kernel statistics~(\Cref{Section: Limiting null distributions}),
		\item General consistency results~(\Cref{Section: General Pointwise Consistency Result}), 		
		\item Detailed explanation of \Cref{Example: Sensitivity of Integral Probability Metric} (\Cref{Section: Details on IPM}), 
		\item Detailed explanation of \Cref{Example: Power Analysis against IPM} (\Cref{Section: Details on Power against IPM}), 
		\item Minimax separation rate in the HSIC metric (\Cref{Section: Separation in HSIC metric}), 
		\item Minimum separation of the $\mathttt{dpHSIC}$ test in the $L_2$ distance (\Cref{Section: Separation in L2 for dpHSIC}), 
		\item Negative results of HSIC U-statistic (\Cref{Section: HSIC U-statistic}) and 
		\item Analyses of minimum separation rates (\Cref{Section: Interpretation}).
	\end{enumerate}
	\item In \Cref{Section: Additional Simulations}, we present additional simulations including 
	\begin{enumerate}[(i)]
		\item \blue{two-sample testing on Parkinson's telemonitoring data~(\Cref{subsec:parkinsons}),}
		\item \blue{experiments with unequal sample sizes~(\Cref{subsec:different_sample_sizes}),}
		\item independence testing with dpHSIC~(\Cref{subsec:hsic_uniform}),
		\item additional experiments on power in the low-privacy regime~(\Cref{subsec:low_privacy}), and
		\item experiments on type I error rates~(\Cref{subsec:level}). 
	\end{enumerate}
	\item In \Cref{Section: Alternative private tests}, we explain in detail other private tests considered in our simulation studies. 
	\item In \Cref{Section: Proofs for the Main Text}, we provide the proofs of the results presented in the main text.
	\item In \Cref{Section: Proofs of the additional results}, we collect the proofs of the additional results provided in \Cref{Section: Additional Results}. 
	\item In \Cref{Section: Technical Lemmas}, we collect technical lemmas used in the main proofs. 
\end{itemize}
Throughout these appendices, we use an additional set of notation described below. 

\paragraph{Additional Notation.} For a sequence of random variables $X_n$, we say that $X_n = O_P(1)$ if for any $\epsilon >0$, there exists $M_\epsilon, N_\epsilon > 0$ such that $\mP(|X_n|  > M_\epsilon) <\epsilon$ for all $n \geq N_\epsilon$. Similarly, we say that $X_n = o_P(1)$ if for any $\epsilon > 0$, there exists $N_\epsilon>0$ such that $\mP(|X_n| > \epsilon) < \epsilon$ for all $n \geq N_\epsilon$. For a sequence of positive numbers $a_n$, $X_n = O_P(a_n)$ (resp.~$o_P(a_n)$) means $a_n^{-1} X_n = O_P(1)$ (resp.~$a_n^{-1} X_n = o_P(1)$). For $x \in \mathbb{R}$, $\ceil{x}$ denotes the least integer greater than or equal to $x$. We use the notation $X_n \convD X$ to denote the convergence in distribution and $X_n \convP X$ to denote the convergence in probability. Consider two probability distributions $P$ and $Q$ on a measurable space $(\Omega, \mathcal{F})$. The total variation (TV) distance between $P$ and $Q$ is defined as
\begin{align*}
	d_{\mathrm{TV}}(P,Q) = \sup_{A \in \mathcal{F}}\, \bigl| P(A) - Q(A) \bigr| = \frac{1}{2} \| P- Q\|_1 = \frac{1}{2} \int \bigg| \frac{\dd P}{\dd \nu} - \frac{\dd Q}{\dd \nu} \bigg| \dd\nu,
\end{align*}
where $\nu$ is a common dominating measure of $P$ and $Q$, and $\frac{\dd P}{\dd \nu}$, $\frac{\dd Q}{\dd \nu}$ are their density functions with respect to $\nu$. The Kullback--Leibler (KL) divergence of $P$ from $Q$ is given as 
\begin{align*}
	d_{\mathrm{KL}}(P,Q) = \int  \frac{\dd P}{\dd \nu} \log \biggl( \frac{\dd P/\dd\nu}{\dd Q/\dd\nu} \biggr) \dd\nu.
\end{align*}
We let $\mathbf{i}_p^q$ denote the set of all $p$-tuples drawn without replacement from $[q]$. The $n$-fold product distribution of a distribution $P$ is denoted as $P^{\otimes n}$. 

{
\section{R\'{e}nyi Differentially Private Permutation Tests}
\label{Section: Renyi Differentially Private Permutation Tests}

\citet{mironov2017renyi} introduced R\'{e}nyi differential privacy (RDP), which measures privacy loss using R\'{e}nyi divergence. RDP retains key properties of pure differential privacy, including post-processing and composition, while often simplifying the privacy analysis of iterative algorithms. This section shows that the proposed permutation framework extends directly to RDP.

For probability distributions $P$ and $Q$ with densities $p$ and $q$, the R\'{e}nyi divergence of order $\gamma>1$ is
\[
	D_\gamma(P\|Q)
	=
	\frac{1}{\gamma-1}
	\log
	\mE_{X\sim Q}
	\left[
	\left\{\frac{p(X)}{q(X)}\right\}^{\gamma}
	\right].
\]
We use $\gamma$ for the order because $\alpha$ already denotes the significance level.

\begin{definition}[R\'{e}nyi Differential Privacy]
\label{Definition: Renyi DP}
	A randomized algorithm $\mathcal{A}$ is $(\gamma,\varepsilon)$-RDP if, for every auxiliary input $w$ and all neighboring datasets $\mathcal{X}_n,\widetilde{\mathcal{X}}_n$,
	\[
		D_\gamma\!\left(
		\mathcal{A}(\mathcal{X}_n;w)
		\middle\|
		\mathcal{A}(\widetilde{\mathcal{X}}_n;w)
		\right)
		\leq\varepsilon.
	\]
\end{definition}

\begin{algorithm}[H]\raggedright
\caption{R\'{e}nyi Differentially Private Permutation Test}
\label{Algorithm: Renyi DP permutation test}
	\textbf{Input}: Data $\mathcal{X}_n$, significance level $\alpha\in(0,1)$, privacy parameters $\gamma>1$ and $\varepsilon>0$, test statistic $T$, global sensitivity (or an upper bound) $\Delta_T$, and number of permutations $B\in\mathbb{N}$.\\
	\vskip .3em
	\textbf{For} $i\in[B]$ \textbf{do}
	\begin{algorithmic}
		\State Generate a random permutation $\bpi_i$ of $[n]$ and $\zeta_i\sim N(0,1)$.
		\State Set $M_i\leftarrow T(\mathcal{X}_n^{\bpi_i})+\Delta_T\sqrt{2\gamma/\varepsilon}\,\zeta_i$.
	\end{algorithmic}
	\textbf{End For}
	\begin{algorithmic}
		\State Generate $\zeta_0\sim N(0,1)$ and set $M_0\leftarrow T(\mathcal{X}_n)+\Delta_T\sqrt{2\gamma/\varepsilon}\,\zeta_0$.
		\State Compute $\widehat p_{\mathrm{rdp}}$ from $\{M_i\}_{i=0}^B$ as in \eqref{Eq: permutation p-value}.
	\end{algorithmic}
	\textbf{Output}: Reject $H_0$ if $\widehat p_{\mathrm{rdp}}\leq\alpha$.
\end{algorithm}

The procedure uses the Gaussian mechanism, whose RDP calibration follows from \citet[][Corollary 3]{mironov2017renyi}. In our setting, the statistic is scalar, so its $\ell_1$ and $\ell_2$ sensitivities coincide.

\begin{lemma}[Gaussian Mechanism under RDP]
\label{Lemma: Gaussian Mechanism under RDP}
	If a scalar statistic has global sensitivity $S$, then adding a centered Gaussian variable with variance $\gamma S^2/(2\varepsilon)$ yields a $(\gamma,\varepsilon)$-RDP release. In particular, taking $S=2\Delta_T$ gives variance $2\gamma\Delta_T^2/\varepsilon$.
\end{lemma}

\begin{theorem}[Properties of the RDP Permutation Test]
\label{Theorem: Properties of RDP permutation test}
	Given $\alpha\in(0,1)$, the test $\mathds{1}(\widehat p_{\mathrm{rdp}}\leq\alpha)$ in \Cref{Algorithm: Renyi DP permutation test} has the following properties.
	\begin{enumerate}
		\item If $\mathcal{X}_n$ is exchangeable under $H_0$, then, for every $B,n\geq1$,
		\[
			\sup_{P\in\mathcal{P}_0}
			\mP_P(\widehat p_{\mathrm{rdp}}\leq\alpha)
			=
			\frac{\floor{(B+1)\alpha}}{B+1}
			\leq\alpha.
		\]
		\item The test is $(\gamma,\varepsilon)$-RDP.
		\item For a fixed alternative $P$ and fixed $\alpha$, if $\mP_P(M_0\leq M_1)\to0$, then the test is pointwise consistent for every sequence $B_n$ satisfying $\inf_n B_n>\alpha^{-1}-1$.
		\item Let $\beta\in(0,1-\alpha)$ and $B\geq6\alpha^{-1}\log(2\beta^{-1})$. There exist universal constants $C_1,C_2>0$ such that the uniform power is at least $1-\beta$ whenever, for every $P\in\mathcal{P}_1$,
		\begin{align*}
			\mE_P[T(\mathcal{X}_n)]
			-\mE_{P,\boldsymbol{\pi}}[T(\mathcal{X}_n^{\boldsymbol{\pi}})]
			\geq{}&
			C_1
			\sqrt{
			\frac{
			\mV_P[T(\mathcal{X}_n)]
			+\mV_{P,\boldsymbol{\pi}}[T(\mathcal{X}_n^{\boldsymbol{\pi}})]
			}{\alpha\beta}
			}\\
			&+
			C_2\Delta_T\sqrt{\frac{\gamma}{\varepsilon}}
			\max\left\{
			\log\left(\frac1\alpha\right),
			\log\left(\frac1\beta\right)
			\right\}.
		\end{align*}
	\end{enumerate}
\end{theorem}

\begin{proof}
	The validity and consistency arguments use only exchangeability and therefore follow the proofs of \Cref{Theorem: validity of private permutation tests,Theorem: Pointwise Consistency}. For privacy, \Cref{Lemma: Quantile representation} expresses the decision as a post-processing of a scalar comparison with sensitivity at most $2\Delta_T$; \Cref{Lemma: Gaussian Mechanism under RDP} then gives the stated guarantee. Finally, the proof of \Cref{Theorem: General uniform power condition} applies after replacing the Laplace tail bound by the Gaussian tail bound, yielding the displayed noise term.
\end{proof}

For finite $\gamma$, the privacy-noise contribution scales as $\varepsilon^{-1/2}$ rather than the $\varepsilon^{-1}$ dependence under pure DP. This improvement comes with a weaker privacy notion: by monotonicity of R\'{e}nyi divergence, the limit $\gamma=\infty$ corresponds to pure DP.

\begin{remark}[Gaussian Differential Privacy]\normalfont
	The construction can also be adapted to other privacy notions. In particular, $f$-differential privacy provides a unified framework for privacy-loss accounting, and Gaussian differential privacy (GDP) is naturally attained by the Gaussian mechanism~\citep{dong2022gaussian}. Applying \citet[][Theorem 2.7]{dong2022gaussian} to \Cref{Algorithm: Renyi DP permutation test} shows that the test is $\sqrt{4\varepsilon/\gamma}$-GDP.
\end{remark}
}

\section{Additional Results} \label{Section: Additional Results}
In this section, we provide additional technical results omitted in the main text. The proofs for these additional results are relegated to \Cref{Section: Proofs of the additional results}.

\subsection{Limiting Null Distributions} \label{Section: Limiting null distributions}
The following proposition derives the limiting distributions of privatized kernel test statistics. In particular, denote 
\begin{align*}
	M_{\mathrm{MMD}} := \widehat{\mathrm{MMD}}(\mathcal{X}_{n+m}) + \frac{2\sqrt{2K}}{n\xi_{\varepsilon,\delta}} \zeta \quad \text{and} \quad  M_{\mathrm{HSIC}} := \widehat{\mathrm{HSIC}}(\mathcal{X}_{n}) + \frac{8(n-1)\sqrt{KL}}{n^2\xi_{\varepsilon,\delta}} \zeta,
\end{align*}
where $\zeta \sim \mathsf{Laplace}(0,1)$ independent of everything else. The results developed in \Cref{Section: Applications to Kernel-Based Inference} guarantee that $M_{\mathrm{MMD}}$ and $M_{\mathrm{HSIC}}$ are $(\varepsilon,\delta)$-DP for bounded kernels. These two private statistics have the following limiting behavior under the null. The proof of \Cref{Proposition: Asymptotic null distributions} can be found in \Cref{Section: Proof of Proposition: Asymptotic null distributions}.
\begin{proposition}[Asymptotic Null Distributions] \label{Proposition: Asymptotic null distributions}
	\leavevmode 
	\begin{itemize}
		\item (MMD) Assume that the kernel $k$ is bounded as $0 \leq k(x,y) \leq K$ for all $x,y \in \mathbb{S}$, and $\frac{m}{n+m} \rightarrow \omega \in (0,1)$. Write $\sigma = 2\sqrt{2K} n^{-1}\xi_{\varepsilon,\delta}^{-1}$. Then there exists a deterministic sequence $\{\lambda_i\}_{i=1}^\infty$ such that 
		\begin{align*}
			\begin{cases}
				(n+m)^{1/2} M_{\mathrm{MMD}}  \convD  \sqrt{\sum_{i=1}^\infty \lambda_i Z_i^2} + \eta \zeta \quad & \text{if $\sqrt{n+m} \sigma  \rightarrow \eta \in [0, \infty)$,} \\[.5em]
				\sigma^{-1} M_{\mathrm{MMD}}  \convD \zeta \quad & \text{if $\sqrt{n+m} \sigma \rightarrow \infty$,}
			\end{cases}
		\end{align*}
		where $\{Z_i\}_{i=1}^\infty \iid N(0,1)$ and $\zeta \sim \mathsf{Laplace}(0,1)$ are independent.
		\item (HSIC) Assume that the kernels $k$ and $\ell$ are bounded as $0 \leq k(y,y') \leq K$ and $0 \leq \ell(z,z') \leq L$ for all $y,y' \in \mathbb{Y}$ and $z,z' \in \mathbb{Z}$. Write $\sigma = 8(n-1)n^{-2} \sqrt{KL} \xi_{\varepsilon,\delta}^{-1}$. Then there exist deterministic sequences $\{\lambda_i\}_{i=1}^\infty$ and $\{\eta_i\}_{i=1}^\infty$ such that 
		\begin{align*}
			\begin{cases}
				n^{1/2} M_{\mathrm{HSIC}}   \convD  \sqrt{\sum_{i=1}^\infty\sum_{j=1}^\infty \lambda_i \eta_j Z_{i,j}^2} + \eta \zeta \quad & \text{if $\sqrt{n} \sigma  \rightarrow \eta \in [0, \infty)$,} \\[.5em]
				\sigma^{-1} M_{\mathrm{HSIC}}  \convD \zeta \quad & \text{if $\sqrt{n} \sigma \rightarrow \infty$,}
			\end{cases}
		\end{align*}
		where $\{Z_{i,j}\}_{i,j=1}^\infty \iid N(0,1)$ and $\zeta \sim \mathsf{Laplace}(0,1)$ are independent.
	\end{itemize}
	\vskip 1em 
	\emph{\textbf{Remark.} Note that the sequences $\{\lambda_i\}_{i=1}^\infty$ and $\{\eta_i\}_{i=1}^\infty$ are associated with the eigenvalues of integral kernel operators. See \cite{gretton2012kernel} and \cite{zhang2018large} for details.}
\end{proposition}

\subsection{General Pointwise Consistency} \label{Section: General Pointwise Consistency Result}
One of the desiderata of nonparametric tests is their pointwise consistency: the power converges to one as the sample size increases against any fixed alternative of interest. The following lemma develops a general pointwise consistency result that can be applied broadly to resampling-based tests~(\emph{e.g.},~bootstrap and permutation tests). We then leverage this general result to derive conditions for pointwise consistency of the differentially private permutation test in \Cref{Theorem: Pointwise Consistency}. 

\begin{lemma}[General Conditions for Consistency] \label{Lemma: General conditions for consistency}
	Let $\alpha \in (0,1)$ be a fixed constant. Suppose that $\{W_{1,n},\ldots,W_{B_n,n}\}$ are i.i.d.~random variables conditional on a sigma field $\mathcal{G}$, and $W_{0,n}$ is constant conditional on the same sigma field $\mathcal{G}$. Suppose further that $\lim_{n \rightarrow \infty} \mP(W_{0,n} \leq W_{1,n}) = 0$. Then for any positive sequence of $B_n$ such that $\min_{n \geq 1} B_n > \alpha^{-1} - 1$, we have 
	\begin{align*}
		\lim_{n \rightarrow \infty} \mP \biggl( \frac{1}{B_n+1} \biggl\{ \sum_{i=1}^{B_n} \mathds{1}(W_{0,n} \leq W_{i,n}) + 1 \biggr\} \leq \alpha \biggr) = 1.
	\end{align*}
\end{lemma}

The proof of this result can be found in \Cref{Section: Proof of Lemma: General conditions for consistency}. When $B_n$ is a fixed quantity, \Cref{Lemma: General conditions for consistency} can be proved using a union bound. As mentioned in the main text, proving this consistency result for a general sequence of $B_n$ is non-trivial and thereby we highlight it as our contribution. 

\vskip 1em

To illustrate \Cref{Lemma: General conditions for consistency} in a simple setting, let us consider a wild bootstrap test for normal mean testing described as follows.

\begin{example}[Wild Bootstrap] \normalfont
	Suppose that we are given $\mathcal{X}_n = \{X_1,\ldots,X_n\} \iid N(\mu, \sigma^2)$, and our interest is in testing whether $H_0 : \mu = 0$ or $H_1 : \mu > 0$. Consider a test statistic $W_{0,n} = \frac{1}{n} \sum_{i=1}^n X_i$, and let $\mathcal{G}$ be the sigma field generated by $\mathcal{X}_n$. Then conditional on $\mathcal{G}$, the test statistic $W_{0,n}$ is constant. Now generate i.i.d.~Rademacher random variables $\epsilon_1,\ldots,\epsilon_n$, and define $W_{1,n} = \frac{1}{n} \blue{\sum_{i=1}^n} \epsilon_i X_i$. The other statistics $W_{2,n},\ldots,W_{B_n,n}$ are defined similarly using independent sets of i.i.d.~Rademacher random variables. Then conditional on $\mathcal{G}$, the sequence $\{W_{1,n},\ldots,W_{B_n,n}\}$ consists of i.i.d.~random variables. We reject the null when 
	\begin{align*}
		\frac{1}{B_n+1} \biggl\{ \sum_{i=1}^{B_n} \mathds{1}(W_{0,n} \leq W_{i,n}) + 1 \biggr\} \leq \alpha
	\end{align*}
	and the resulting test is called a (wild) bootstrap test. By the law of large numbers, it can be seen that $W_{0,n} \convP \mu$ and $W_{1,n}  \convP 0$ for fixed mean $\mu >0$ and finite variance $\sigma^2$, implying that $\lim_{n \rightarrow \infty} \mP(W_{0,n} \leq W_{1,n}) = 0$. Hence the bootstrap test is consistent according to \Cref{Lemma: General conditions for consistency}. The type I error is also controlled since $\{W_{0,n},W_{1,n},\ldots,W_{B_{n},n}\}$ are exchangeable under the null.
\end{example}

\subsection{Details on \Cref{Example: Sensitivity of Integral Probability Metric}} \label{Section: Details on IPM}
Let us denote the pooled sample as $\mathcal{X}_{n+m} = \mathcal{Y}_n \cup \mathcal{Z}_m = \{X_1,\ldots,X_{n+m}\}$, and the permutation counterpart permuted according to $\bpi$ as $\mathcal{X}_{n+m}^{\bpi}$. Then the plug-in estimator using $\mathcal{X}_{n+m}^{\bpi}$ is given as
\begin{align*}
	T_{\mathrm{IPM}}(\mathcal{X}_{n+m}^{\bpi}) = \sup_{f \in \mathcal{F}} \bigg| \frac{1}{n} \sum_{i=1}^n f(X_{\pi_i}) - \frac{1}{m}\sum_{i=1}^m f(X_{\pi_{n+i}}) \bigg|.
\end{align*}
By the reverse triangle inequality, it holds that
\begin{align*}
	\sup_{\substack{\mathcal{X}_{n+m}, \tilde{\mathcal{X}}_{n+m}:\\d_{\mathrm{ham}}(\mathcal{X}_{n+m}, \tilde{\mathcal{X}}_{n+m}) \leq 1}} \big|T_{\mathrm{IPM}}(\mathcal{X}_{n+m}^{\bpi}) - T_{\mathrm{IPM}}(\tilde{\mathcal{X}}_{n+m}^{\bpi}) \big| \leq \widetilde{\Delta}_T \coloneqq\frac{1}{\min\{n,m\}} \sup_{X,X' \in \mathbb{S}} \sup_{f \in \mathcal{F}}|f(X) - f(X')|.
\end{align*}
The above bound is independent of $\bpi$ and thus the global sensitivity is at most $\widetilde{\Delta}_T$. Indeed, this upper bound is tight. For $X,X' \in \mathbb{S}$, set $\mathcal{X}_{n+m} = \{X,\ldots,X\}$ and $\tilde{\mathcal{X}}_{n+m} = \{X',X,\ldots,X\}$. Then the global sensitivity $\Delta_{T}$ is lower bounded as
\begin{align*}
	\frac{1}{\min\{n,m\}} \sup_{f \in \mathcal{F}}|f(X) - f(X')|\leq \Delta_{T}, 
\end{align*}
which holds for all $X,X' \in \mathbb{S}$. Hence it holds $\widetilde{\Delta}_T \leq \Delta_{T}$, and the global sensitivity becomes $\Delta_{T} = \widetilde{\Delta}_T$.

\subsection{Details on \Cref{Example: Power Analysis against IPM}} \label{Section: Details on Power against IPM}

In view of condition~\eqref{Eq: Uniform Power Condition}, there are four terms that we need to investigate, namely $\mE[T(\mathcal{X}_{n+m})]$, $\mE[T(\mathcal{X}_{n+m}^{\boldsymbol{\pi}})]$, $\mV[T(\mathcal{X}_{n+m})]$ and $\mV[T(\mathcal{X}_{{n+m}}^{\boldsymbol{\pi}})]$. 

Starting with the expected value $\mE[T(\mathcal{X}_{n+m})]$, the triangle inequality yields
\begin{align*}
	\sup_{f \in \mathcal{F}} \bigg| \frac{1}{n} \sum_{i=1}^n f(Y_i) - \frac{1}{m}\sum_{i=1}^m f(Z_i) \bigg| ~\geq~  \mathrm{IPM}_{\mathcal{F}}(P,Q)  & -  \sup_{f \in \mathcal{F}} \bigg| \frac{1}{n} \sum_{i=1}^n f(Y_i) - \mE_P[f(Y)] \bigg| \\[.5em]
	& - \sup_{f \in \mathcal{F}} \bigg|\mE_Q[f(Z)] - \frac{1}{m}\sum_{i=1}^m f(Z_i) \bigg| 
\end{align*}
and thus 
\begin{align*}
	\mE[T(\mathcal{X}_{n+m})] ~\geq~ & \mathrm{IPM}_{\mathcal{F}}(P,Q) - \mE \biggl[\sup_{f \in \mathcal{F}} \bigg| \frac{1}{n} \sum_{i=1}^n f(Y_i) - \mE_P[f(Y)] \bigg|\biggr] \\[.5em]
	& ~~~~~~~~~~~~~~~~ - \mE\biggl[\sup_{f \in \mathcal{F}} \bigg|\mE_Q[f(Z)] - \frac{1}{m}\sum_{i=1}^m f(Z_i) \bigg| \biggr] \\[.5em]
	\geq ~ & \mathrm{IPM}_{\mathcal{F}}(P,Q) - 4\mathcal{R}_n(\mathcal{F}),
\end{align*}
where the last inequality uses the standard symmetrization trick~\citep[\emph{e.g.},][Chapter 2.3]{vanderVaart1996}. In particular, introducing i.i.d.~copies $\tilde{Y}_i$ of $Y_i$, Jensen's inequality along with the triangle inequality yields
\begin{align*}
	\mE \biggl[\sup_{f \in \mathcal{F}} \bigg| \frac{1}{n} \sum_{i=1}^n f(Y_i) - \mE_P[f(Y)] \bigg|\biggr] \leq \mE \biggl[\sup_{f \in \mathcal{F}} \bigg| \frac{1}{n} \sum_{i=1}^n \omega_i \{f(Y_i) - f(\tilde{Y}_i)\} \bigg|\biggr] \leq 2 \mathcal{R}_n(\mathcal{F}).
\end{align*} 
The other expectation can be handled exactly the same way, which allows us to obtain the lower bound for $\mE[T(\mathcal{X}_{n+m})] $. 

We next look at $\mE[T(\mathcal{X}_{n+m}^{\boldsymbol{\pi}})]$, which is equal to
\begin{align*}
	\mE[T(\mathcal{X}_{n+m}^{\boldsymbol{\pi}})] = \mE\biggl[ \sup_{f \in \mathcal{F}} \bigg| \frac{1}{n} \sum_{i=1}^n f(X_{\pi_i}) - \frac{1}{m}\sum_{i=1}^m f(X_{\pi_{n+i}}) \bigg| \biggr].
\end{align*}
Let $I_1,\ldots,I_{\binom{m}{n}}$ denote all subsets $(i_1,\ldots,i_n)$ of $[m]$ satisfying $1 \leq i_1 < \ldots < i_n \leq m$. Observe that the sample mean can be expressed as the U-statistic with a kernel of order $n$ as below:
\begin{align*}
	\frac{1}{m}\sum_{i=1}^m f(X_{\pi_{n+i}}) = \frac{1}{\binom{m}{n}} \sum_{1 \leq j \leq \binom{m}{n}} \bigg\{\frac{1}{n} \sum_{i \in I_j} f(X_{\pi_{n+i}}) \bigg\}.
\end{align*}
This observation along with Jensen's inequality yields 
\begin{align*}
	\mE[T(\mathcal{X}_{n+m}^{\boldsymbol{\pi}})]  \leq \frac{1}{\binom{m}{n}} \sum_{1 \leq j \leq \binom{m}{n}} \mE\biggl[ \sup_{f \in \mathcal{F}} \bigg| \frac{1}{n} \sum_{i=1}^n f(X_{\pi_i}) - \frac{1}{n}\sum_{i \in I_j} f(X_{\pi_{n+i}}) \bigg| \biggr].
\end{align*}
Without loss of generality, set $I_j = (1,\ldots,n)$ and then we have
\begin{align*}
	\mE\biggl[ \sup_{f \in \mathcal{F}} \bigg| \frac{1}{n} \sum_{i=1}^n f(X_{\pi_i}) - \frac{1}{n}\sum_{i \in I_j} f(X_{\pi_{n+i}}) \bigg| \biggr] = \mE\biggl[ \sup_{f \in \mathcal{F}} \bigg| \frac{1}{n} \sum_{i=1}^n \omega_i \{f(X_{\pi_i}) - f(X_{\pi_{n+i}})\} \bigg| \biggr],
\end{align*}
as randomly switching the order between $\pi_i$ and $\pi_{n+i}$ for $i \in [n]$ does not change the distribution of 
\begin{align*}
	\sup_{f \in \mathcal{F}} \bigg| \frac{1}{n} \sum_{i=1}^n f(X_{\pi_i}) - f(X_{\pi_{n+i}}) \bigg|.
\end{align*}
Therefore we may upper bound $\mE[T(\mathcal{X}_{n+m}^{\boldsymbol{\pi}})]$ as
\begin{align*}
	\mE[T(\mathcal{X}_{n+m}^{\boldsymbol{\pi}})] \leq 2 \mathcal{R}_n(\mathcal{F}).
\end{align*}

Moving to the variance terms, since $T(\mathcal{X}_{n+m})$ is a function of independent random variables, we can apply the Efron--Stein inequality for the variance \citep[][Corollary 3.2]{boucheron2013concentration} as
\begin{align*}
	\mV[T(\mathcal{X}_{n+m})] ~\leq~&  \frac{1}{2} \sum_{i=1}^{n+m} \mE\bigl[ \bigl(T(\mathcal{X}_{n+m}) - T(\mathcal{X}_{n+m}^{(i)})\bigr)^2 \bigr] \\[.5em]
	\leq ~ & \frac{1}{2n}\sup_{X,X' \in \mathbb{S}} \sup_{f \in \mathcal{F}}|f(X) - f(X')|^2 + \frac{1}{2m} \sup_{X,X' \in \mathbb{S}} \sup_{f \in \mathcal{F}}|f(X) - f(X')|^2 \\[.5em]
	\leq ~ & n\Delta_T^2,
\end{align*}
where $\mathcal{X}_{n+m}^{(i)} = (X_1,\ldots,X_i',\ldots,X_{n+m})$, and $X_i'$ is an i.i.d.~copy of $X_i$ independent of everything else. This proves that $\mV[T(\mathcal{X}_{n+m})] \leq n\Delta_T^2$. 

For the last term, the law of total variance shows
\begin{align*}
	\mV[T(\mathcal{X}_{{n+m}}^{\bpi})] = \mV[\mE\{ T(\mathcal{X}_{{n+m}}^{\bpi}) \given \bpi\}] +  \mE[\mV\{ T(\mathcal{X}_{{n+m}}^{\bpi}) \given \bpi \}].
\end{align*}
Conditional on $\bpi$, the same analysis as above yields $\mV\{ T(\mathcal{X}_{{n+m}}^{\bpi}) \given \bpi \} \leq  n\Delta_T^2$ based on the Efron--Stein inequality, which gives $\mE[\mV\{ T(\mathcal{X}_{{n+m}}^{\bpi}) \given \bpi \}] \leq  n\Delta_T^2$. For the first term, using the same trick used in the analysis of $\mE[T(\mathcal{X}_{n+m}^{\bpi})]$, we have
\begin{align*}
	 \mV[\mE\{ T(\mathcal{X}_{{n+m}}^{\bpi}) \given \bpi\}] ~\leq~ & \mE[(\mE\{ T(\mathcal{X}_{{n+m}}^{\bpi}) \given \bpi\})^2] \\[.5em]
	 \leq ~ &  \mE\biggl[ \bigg( \mE \biggl\{ \sup_{f \in \mathcal{F}} \bigg| \frac{1}{n} \sum_{i=1}^n \omega_i \{f(X_{\pi_i}) - f(X_{\pi_{n+i}})\} \bigg| \,\bigg|\, \bpi\biggr\} \bigg)^2 \biggr] \\[.5em]
	 \leq ~ & 4 \mathcal{R}_{n}^2(\mathcal{F}).
 \end{align*}

Having these ingredients, one can directly check the condition~\eqref{Eq: Uniform Power Condition} is fulfilled if 
\begin{align*}
	\mathrm{IPM}_{\mathcal{F}}(P,Q) \geq C_1\frac{\mathcal{R}_n(\mathcal{F})}{\sqrt{\alpha \beta}} + C_2 \frac{\sqrt{n}\Delta_T}{\sqrt{\alpha \beta}} + C_3 \frac{\Delta_T}{\xi_{\varepsilon,\delta}} \max\bigg\{ \! \log \biggl(\frac{1}{\alpha}\biggr),  \, \log \biggl(\frac{1}{\beta}\biggr) \bigg\},
\end{align*}
where $C_1,C_2,C_3$ are some positive constants.

\subsection{Separation in HSIC metric} \label{Section: Separation in HSIC metric}
In this subsection, we develop results similar to those in \Cref{Section: Separation in MMD metric} in terms of the HSIC metric. We start by discussing the minimum separation for the $\mathttt{dpHSIC}$ test in \Cref{Theorem: Uniform separation for HSIC} and then establish the matching lower bound in \Cref{Theorem: Minimax separation in HSIC}. Letting $\mathcal{P}_{\mathbb{Y} \times \mathbb{Z}}$ denote the class of distributions on $\mathbb{Y} \times \mathbb{Z}$ and $\rho > 0$, we define the set of alternative distributions as
\begin{align*}
	\mathcal{P}_{\mathrm{HSIC}_{k \otimes \ell}}\!(\rho) \coloneqq  \big\{ P_{YZ} \in \mathcal{P}_{\mathbb{Y} \times \mathbb{Z}} : \mathrm{HSIC}_{k \otimes \ell}(P_{YZ}) \geq \rho \big\}.
\end{align*}
For a given target type II error $\beta \in (0,1-\alpha)$, the minimum separation for the $\mathttt{dpHSIC}$ test against $\mathcal{P}_{\mathrm{HSIC}_{k \otimes \ell}}\!(\rho)$ is given as 
\begin{align}
	\label{Eq: separation HSIC rate}
	\rho_{\phi_{\mathttt{dpHSIC}}}(\alpha,\beta,\varepsilon,\delta,n) \coloneqq  \inf \biggl\{ \rho > 0 : \sup_{P_{YZ} \in \mathcal{P}_{\mathrm{HSIC}_{k \otimes \ell}}\!(\rho)} \mE_{P_{YZ}} [1 - \phi_{\mathttt{dpHSIC}}] \leq \beta \biggr\}.
\end{align}
The next theorem, which is analogous to \Cref{Theorem: Uniform separation for MMD} for the $\mathttt{dpMMD}$ test, provides an upper bound for the minimum separation $\rho_{\phi_{\mathttt{dpHSIC}}}(\alpha,\beta,\varepsilon,\delta,n)$ as a function of $\alpha$, $\beta$, $\varepsilon$, $\delta$ and $n$. 

\begin{theorem}[Minimum Separation of $\phi_{\mathttt{dpHSIC}}$] \label{Theorem: Uniform separation for HSIC}
	Assume that the kernels $k$ and $\ell$ are bounded as $0 \leq k(y,y') \leq K$ and $0 \leq \ell(z,z') \leq L$ for all $y, y' \in \mathbb{Y}$ and $z, z' \in \mathbb{Z}$. For all values of $\alpha \in (0,1)$, $\beta \in (0, 1-\alpha)$, $\varepsilon > 0$, $\delta\in[0,1)$ and $B \geq 6\alpha^{-1} \log(2\beta^{-1})$, the minimum separation for $\phi_{\mathttt{dpHSIC}}$ satisfies 
	\begin{align*}
		\rho_{\phi_{\mathttt{dpHSIC}}}(\alpha,\beta,\varepsilon,\delta,n) \leq C_{K,L} \max \Biggl\{ \sqrt{\frac{\max\!\big\{\!\log(1/\alpha), \, \log(1/\beta)\big\}}{n}}, \, \frac{\max\!\big\{\!\log(1/\alpha), \, \log(1/\beta)\big\}}{n\xi_{\varepsilon,\delta}}  \Biggr\},
	\end{align*}
	where $\xi_{\varepsilon,\delta}$ is given in \eqref{Eq: definition of xi} and $C_{K,L}$ is a positive constant that only depends on $K$ and $L$.
\end{theorem}
The proof of \Cref{Theorem: Uniform separation for HSIC}, given in \Cref{Section: Proof of Theorem: Uniform separation for HSIC}, follows a similar approach to that of \blue{\Cref{Theorem: Uniform separation for MMD}}. A notable difference, however, is that we need to use exponential tail bounds for the HSIC statistic instead of the MMD statistic. To this end, we develop exponential concentration results for the empirical HSIC (\Cref{Lemma: Exponential inequality for the empirical HSIC}) and the permuted HSIC~(\Cref{Lemma: Concentration inequality for permuted HSIC}), by leveraging the recent result of \cite{kim2020minimax} and McDiarmid's inequality. As clearly demonstrated in the proof, the use of these tools is crucial in obtaining the logarithmic factors in $\alpha$ and $\beta$.

We now examine minimax optimality of $\phi_{\mathttt{dpHSIC}}$ with a focus on the cases where $\mathbb{Y} = \mathbb{R}^{d_Y}$ and $\mathbb{Z} = \mathbb{R}^{d_Z}$. For this direction, we establish a lower bound for the minimax separation $\rho^\star_{\mathrm{HSIC}}$ defined below and compare it with $\rho_{\phi_{\mathttt{dpHSIC}}}$. Let $\phi : \mathcal{X}_n \mapsto \{0,1\}$ be a test function and $\mathcal{P}_0 \coloneqq  \{ P_{YZ} \in \mathcal{P}_{\mathbb{Y} \times \mathbb{Z}}: P_{YZ} = P_YP_Z \}$ be the set of null distributions on $\mathbb{R}^{d_Y} \times \mathbb{R}^{d_Z}$. We denote the set of $(\varepsilon,\delta)$-DP level $\alpha$ tests as
\begin{align*}
	\Phi_{\alpha,\varepsilon,\delta} \coloneqq   \Big\{\phi : \sup_{P_{YZ} \in \mathcal{P}_0} \mE_{P_{YZ}}[\phi] \leq \alpha \ \text{and} \ \text{$\phi$ is $(\varepsilon,\delta)$-DP} \Big\},
\end{align*}
and define the minimax separation in terms of the HSIC metric as
\begin{align*}
	\rho^\star_{\mathrm{HSIC}}(\alpha,\beta,\varepsilon,\delta,n)  \coloneqq  \inf\biggl\{ \rho > 0: \inf_{\phi \in \Phi_{\alpha,\varepsilon,\delta}}  \sup_{P_{YZ} \in \mathcal{P}_{\mathrm{HSIC}_{k \otimes \ell}}\!(\rho)} \mE_{P_{YZ}} [1 - \phi] \leq \beta \bigg\}.
\end{align*}
Similar to \Cref{Theorem: Minimax separation in MMD}, the next result establishes a lower bound for the minimax separation $\rho^\star_{\mathrm{HSIC}}$ under the DP constraint. 
\begin{theorem}[Minimax Separation in HSIC] \label{Theorem: Minimax separation in HSIC}
	Let $\alpha$ and $\beta$ be real numbers in the interval $(0,1/5)$, $\varepsilon>0$ and $\delta\in[0,1)$. Assume that the kernel functions $k$ and $\ell$ are translation invariant on $\mathbb{R}^{d_Y}$ and $\mathbb{R}^{d_Z}$, respectively. In particular, there exist some functions $\kappa_Y, \kappa_Z$ such that $k(y,y') = \kappa_Y(y-y')$ for all $y,y' \in \mathbb{R}^{d_Y}$ and $\ell(z,z') = \kappa_Z(z-z')$ for all $z,z' \in \mathbb{R}^{d_Z}$. Moreover, assume that the kernels are non-constant in the sense that there exist positive constants $\eta_Y,\eta_Z$ such that $\kappa_Y(0) - \kappa_Y(y_0) \geq \eta_Y$ and $\kappa_Z(0) - \kappa_Z(z_0) \geq \eta_Z$ for some $y_0 \in \mathbb{R}^{d_Y}$ and $z_0 \in \mathbb{R}^{d_Z}$. Then the minimax separation over $\mathcal{P}_{\mathrm{HSIC}_{k \otimes \ell}}\!(\rho)$ is lower bounded as 
	\begin{align*} 
		\rho^\star_{\mathrm{HSIC}}(\alpha,\beta,\varepsilon,\delta,n)  \geq C_{\eta_Y,\eta_Z} \max \Biggl\{ \min \bigg\{  \sqrt{\frac{\log(1/(\alpha+\beta))}{n}}, \, 1 \bigg\}, \, \min \biggl\{ \frac{\log(1/\beta)}{n\xi_{\varepsilon,\delta}}, \, 1 \biggr\} \Biggr\},
	\end{align*}
	where $\xi_{\varepsilon,\delta}$ is given in \eqref{Eq: definition of xi} and $C_{\eta_Y,\eta_Z}$ is a positive constant that only depends on $\eta_Y$ and $\eta_Z$.
\end{theorem}
The proof of this result can be found in \Cref{Section: Proof of Theorem: Minimax separation in HSIC}. By comparing the result of \Cref{Theorem: Uniform separation for HSIC} with the above result, it becomes evident that $\rho_{\phi_{\mathttt{dpHSIC}}} \asymp	\rho^\star_{\mathrm{HSIC}}$ under the conditions $\alpha \asymp \beta$ and $\mathbb{Y} \times \mathbb{Z} = \mathbb{R}^{d_Y} \times \mathbb{R}^{d_Z}$. Consequently, $\phi_{\mathttt{dpHSIC}}$ achieves minimax rate optimality across all parameters, namely $(\alpha,\beta,\varepsilon,n)$, provided that the type I error and the type II error are comparable. Other comments on \Cref{Theorem: Minimax separation in MMD} can similarly be applied to \Cref{Theorem: Minimax separation in HSIC}. For instance, numerous kernels frequently used in practice are translation invariant and the existence of $\eta_Y$ and $\eta_Z$ in the theorem statement is guaranteed as long as the kernels $k$ and $\ell$ are characteristic~\citep{tolstikhin2017minimax}.

The minimax results for both $\mathttt{dpMMD}$ and $\mathttt{dpHSIC}$ tests indicate the existence of an inherent trade-off between privacy and statistical power. More concretely, one cannot expect to achieve the parametric $\sqrt{n}$-separation rate in both MMD and HSIC metrics when the privacy parameter $\xi_{\varepsilon,\delta}$ is much smaller than $n^{-1/2}$. This intrinsic trade-off has been observed in a variety of statistical problems~\citep{diakonikolas2015differentially,acharya2018differentially,kamath2019privately,kamath2020private,acharya2021differentially} and our work serves to reaffirm this trade-off in the context of kernel-based testing problems.

\subsection{Separation in \texorpdfstring{$L_2$}{L2} Metric} \label{Section: Separation in L2 for dpHSIC}
We next investigate the minimum separation of the $\mathttt{dpHSIC}$ test against a class of alternatives measured in terms of the $L_2$ metric. For $\rho > 0$, let $\tilde{\mathcal{P}}_{L_2}(\rho)$ be the collection of distributions $P_{YZ}$ on $\mathbb{R}^{d_Y} \times \mathbb{R}^{d_Z}$, where $P_{YZ}$ is equipped with the Lebesgue density function $p_{YZ}$ and the product of the marginals $p_Yp_Z$ such that $\|p_{YZ} - p_Yp_Z\|_{L_2} \geq \rho$. The class of alternative distributions of interest is a subset of $\tilde{\mathcal{P}}_{L_2}(\rho)$ given as
\begin{align*}
	\tilde{\mathcal{P}}_{L_2}^s(\rho) := \bigg\{ P_{YZ} \in  \tilde{\mathcal{P}}_{L_2}(\rho) : p_{YZ} - p_Y p_Z \in \mathcal{S}_{d_Y + d_Z}^s(R), \ \max(\|p_{YZ}\|_{L_\infty}, \|p_{Y} p_Z\|_{L_\infty}) \leq M \bigg\},
\end{align*}
where $\mathcal{S}_{d_Y + d_Z}^s(R)$ is the Sobolev ball defined in~\eqref{Eq: Sobolev ball}. As for the $\mathttt{dpMMD}$ test in \Cref{Section: Separation in L2 metric}, we focus on the use of the Gaussian kernels for simplicity.  In particular, for $y,y' \in \mathbb{R}^{d_Y}$ and $z,z' \in \mathbb{R}^{d_Z}$, consider the Gaussian kernels with bandwidths $\blambda = (\lambda_1,\ldots,\lambda_{d_Y}) \in (0,\infty)^{d_Y}$ and $\bmu = (\mu_1,\ldots,\mu_{d_Z}) \in (0, \infty)^{d_Z}$ given as
\begin{align*}
	k_{\blambda}(y,y') = \prod_{i=1}^{d_Y} \frac{1}{\sqrt{2\pi}\lambda_i} e^{-\frac{(y_i-y_i')^2}{2\lambda_i^2}} \quad \text{and} \quad 	\ell_{\bmu}(z,z') = \prod_{i=1}^{d_Z} \frac{1}{\sqrt{2\pi}\mu_i} e^{-\frac{(z_i-z_i')^2}{2\mu_i^2}}.
\end{align*}
Let us denote the minimum separation of the $\mathttt{dpHSIC}$ test with the above Gaussian kernels against $L_2$ alternatives as 
\begin{align*} 
	\rho_{\phi_{\mathttt{dpHSIC}},L_2}(\alpha, \beta, \varepsilon, \delta, n, d_Y, d_Z, s, R, M) \coloneqq  \inf \Biggl\{ \rho >0 : \sup_{P_{YZ} \in \tilde{\mathcal{P}}_{L_2}^{s}\!(\rho)} \mE_{P_{YZ}}[1 - \phi_{\mathttt{dpHSIC}}] \leq \beta \Biggr\}.
\end{align*}
The next theorem presents an upper bound for $\rho_{\phi_{\mathttt{dpHSIC}},L_2}$ in terms of a set of parameters including $\blambda$, $\bmu$ and $n$. The proof of \Cref{Theorem: Minimax Separation over L2 for HSIC} is given in \Cref{Section: Proof of Theorem: Minimax Separation over L2 for HSIC}.

\begin{theorem}[Minimum Separation of $\mathttt{dpHSIC}$ over $\tilde{\mathcal{P}}_{L_2}^{s}$] \label{Theorem: Minimax Separation over L2 for HSIC}
	Assume that $\alpha \in (0,1)$, $\beta \in (0,1-\alpha)$, $\varepsilon >0$, $\delta \in [0,1)$, $B \geq 6\alpha^{-1} \log(2\beta^{-1})$, $\prod_{i=1}^{d_Y} \lambda_i \leq 1$ and $\prod_{i=1}^{d_Z} \mu_i \leq 1$. The minimum separation of the $\mathttt{dpHSIC}$ test with the Gaussian kernel over $\tilde{\mathcal{P}}_{L_2}^{s}$ is upper bounded as
	\begin{align*}
		\rho_{\phi_{\mathttt{dpHSIC}},L_2}^2 \leq C_{\alpha,\beta,s,R,M,d_Y,d_Z} \Bigg\{  & \sum_{i=1}^{d_Y} \lambda_i^{2s} + \sum_{i=1}^{d_Z} \mu_i^{2s}  + \frac{1}{n\sqrt{\lambda_1 \cdots \lambda_{d_Y} \mu_1 \cdots \mu_{d_Z}}} \\[.5em]
		& + \frac{1}{n^2 \lambda_1 \cdots \lambda_{d_Y} \mu_1 \cdots \mu_{d_Z} \xi_{\varepsilon,\delta}^2}  + \frac{1}{n^{3/2} \lambda_1 \cdots \lambda_{d_Y} \mu_1 \cdots \mu_{d_Z} \xi_{\varepsilon,\delta}} \Bigg\},
	\end{align*}
	where $C_{\alpha,\beta,s,R,M,d_Y,d_Z} $ is a positive constant, depending only on $\alpha,\beta,s,R,M,d_Y,d_Z$, and $\xi_{\varepsilon,\delta}$ is given in \eqref{Eq: definition of xi}.  
\end{theorem}

We make a few comments on \Cref{Theorem: Minimax Separation over L2 for HSIC}.

\begin{itemize}
	\item In contrast to the prior work~\citep{berrett2021optimal,albert2019adaptive,kim2020minimax,schrab2022efficient} that considers U-statistic for independence testing against $L_2$ alternatives, our result is based on the V-statistic of the HSIC, which requires additional effort in dealing with bias terms. As emphasized before, the V-statistic is more favorable than the U-statistic in high privacy regimes as the former has smaller global sensitivity than the latter.
	\item The main idea behind the proof is similar to that of \Cref{Theorem: Minimax Separation over L2} for the $\mathttt{dpMMD}$ test where we first analyze the difference between the U- and V-statistics of the MMD, and then leverage the existing results of the U-statistic. Extending this strategy to the HSIC requires substantial effort as the U- and V-statistics of the HSIC are based on the fourth-order kernel (whereas, the corresponding MMD statistics are based on the second-order kernel).
	\item In the proof, we do not keep track of the dependence on the type I error rate $\alpha$ in deriving the minimum separation. We believe that the dependence on $\alpha$ can be improved in view of \citet[][Theorem 3]{schrab2022efficient} where they establish a logarithmic dependence on $\alpha$ using U-statistics. Nevertheless, extending this result to the V-statistic poses a non-trivial challenge, and we defer this investigation to future research.
	\item We remark that if $\xi_{\varepsilon,\delta}$ is sufficiently large so that the following inequality holds
	\begin{align*}
		\rho_{\phi_{\mathttt{dpHSIC}},L_2}^2 \leq C_{\alpha,\beta,s,R,M,d_Y,d_Z}  \Biggl\{\sum_{i=1}^{d_Y} \lambda_i^{2s} + \sum_{i=1}^{d_Z} \mu_i^{2s}  + \frac{1}{n\sqrt{\lambda_1 \cdots \lambda_{d_Y} \mu_1 \cdots \mu_{d_Z}}}\Biggr\},
	\end{align*}
	then we obtain the same separation rate as the non-private HSIC test studied in \cite{albert2019adaptive}. In particular, letting $\lambda_1 = \cdots \lambda_{d_Y} = \mu_{1} = \cdots = \mu_{d_Z} = n^{-2/(4s+d_Y+d_Z)}$, we see that the $\mathttt{dpHSIC}$ test achieves minimax optimal rate $n^{-2s/(4s+d_Y+d_Z)}$ against the Sobolev ball.
	\item As for the result~\eqref{Eq: L2 separation rate} of the $\mathttt{dpMMD}$ test, one can achieve different minimum separation rates by varying the bandwidth parameters of the Gaussian kernel. In particular, a similar analysis given in  \Cref{Section: Interpretation} verifies that 
		\begin{equation}
		\begin{aligned}
			\rho_{\phi_{\mathttt{dpHSIC}},L_2} \lesssim \begin{cases}
				n^{-\frac{2s}{4s+d}}, & \textrm{if } n^{-\frac{2s-d/2}{4s+d}} \lesssim  \xi_{\varepsilon,\delta} \ \textrm{(low privacy)}, \\[2mm]
				( n^{\frac{3}{2}} \xi_{\varepsilon,\delta})^{-\frac{s}{2s+d}}, & \textrm{if } n^{-\frac{1}{2}} \lesssim\xi_{\varepsilon,\delta} \lesssim n^{-\frac{2s-d/2}{4s+d}}  \ \textrm{(mid privacy)}, \\[2mm]
				\left(n\xi_{\varepsilon,\delta}\right)^{-\frac{2s}{2s+d}}, & \textrm{if }  \xi_{\varepsilon,\delta} \lesssim n^{-\frac{1}{2}}  \ \textrm{(high privacy)}.
			\end{cases}
		\end{aligned} \label{Eq: L2 separation rate for HSIC}
	\end{equation}
	\item As in the $\mathttt{dpMMD}$ test, the optimal choice of the bandwidths assumes that the smoothness parameter $s$ is known. In this regard, it would be interesting to develop an adaptive test that does not require the knowledge of $s$, while retaining (nearly) optimal power. We leave this direction for future work. We also note that establishing lower bounds in high privacy regimes is another interesting avenue for future work.
\end{itemize}

We next turn to the minimum separation of the private independence test based on a U-statistic, and compare it with the one established in \Cref{Theorem: Minimax Separation over L2 for HSIC}.

\subsection{Private Test based on the HSIC U-statistic}  \label{Section: HSIC U-statistic}
Mirroring \Cref{Section: Private Test based on the MMD U-statistic}, this subsection develops a similar negative result for the $\mathttt{dpHSIC}$ test based on a U-statistic. As an unbiased estimator of the square of $\mathrm{HSIC}_{k \otimes \ell}(P_{YZ})$, the U-statistic of HSIC is given as
\begin{equation}
\begin{aligned} \label{Eq: U-HSIC}
	U_{\mathrm{HSIC}} ~=~ & \frac{1}{n(n-1)} \sum_{(i,j) \in \mathbf{i}_2^n} k(Y_i,Y_j) \ell(Z_i,Z_j) + \frac{(n-4)!}{n!} \sum_{(i_1,i_2,j_1,j_2) \in \mathbf{i}_4^n} k(Y_{i_1},Y_{j_1}) \ell(Z_{i_2},Z_{j_2}) \\
	-~ & \frac{2}{n(n-1)(n-2)} \sum_{(i,j_1,j_2) \in \mathbf{i}_3^n} k(Y_{i},Y_{j_1}) \ell(Z_{i},Z_{j_2}),
\end{aligned}
\end{equation}
where we recall that $\mathbf{i}_p^q$ denotes the set of all $p$-tuples drawn without replacement from $[q]$. To privatize $U_{\mathrm{HSIC}}$ as well as the test function through the Laplace mechanism, the following lemma computes the global sensitivity of $U_{\mathrm{HSIC}}$.

\begin{lemma} \label{Lemma: Global sensitivity of U_HSIC}
	Assume that the kernels $k$ and $\ell$ are bounded as $0 \leq k(y,y') \leq K$ and $0 \leq \ell(z,z') \leq L$ for all $y,y' \in \mathbb{Y}$ and $z,z' \in \mathbb{Z}$. In addition, assume that $k$ and $\ell$ are translation invariant, and have non-empty level sets on $\mathbb{Y}$ and $\mathbb{Z}$, respectively. Then there exists a positive sequence $c_n \in [2,24]$ such that for all $n \geq 4$,
	\begin{align*}
		\sup_{\bpi \in \boldsymbol{\Pi}_n} \sup_{\substack{\mathcal{X}_{n},\tilde{\mathcal{X}}_{n}:\\d_{\mathrm{ham}}(\mathcal{X}_{n},\tilde{\mathcal{X}}_{n}) \leq 1}} \bigl| U_{\mathrm{HSIC}}(\mathcal{X}_{n}^{\bpi}) - U_{\mathrm{HSIC}}(\tilde{\mathcal{X}}_{n}^{\bpi}) \bigr| = \frac{c_n KL}{n}.
	\end{align*}
\end{lemma}

The proof of \Cref{Lemma: Global sensitivity of U_HSIC} is given in \Cref{Section: Proof of Lemma: Global sensitivity of U_HSIC}. We would like to highlight that determining the precise global sensitivity of $U_{\mathrm{HSIC}}$ is more challenging than that of $U_{\mathrm{MMD}}$ as the former is a fourth-order U-statistic, whereas the latter is a second order one. Having established the global sensitivity of $U_{\mathrm{HSIC}}$, consider the $\mathttt{dpHSIC}$ test in Algorithm~\ref{Algorithm: DP permutation test} using the test statistic $U_{\mathrm{HSIC}}$ and the global sensitivity $\Delta_T = \frac{c_n KL}{n}$ where $c_n$ can be an arbitrary sequence of constants between $2$ and $24$. We denote the resulting test as $\phi_{\mathttt{dpHSIC}}^u$ and show its suboptimal power property. The proof of the result below can be found in \Cref{Section: Proof of Theorem: Suboptimality of U_HSIC}

\begin{theorem}[Suboptimality of $\phi_{\mathttt{dpHSIC}}^u$] \label{Theorem: Suboptimality of U_HSIC}
	Assume that the kernels $k$ and $\ell$ fulfill the conditions specified in \Cref{Lemma: Global sensitivity of U_HSIC}. Moreover, assume that if $P_{YZ} \in \mathcal{P}_{\mathbb{Y} \times \mathbb{Z}}$, then $w P_{YZ} + (1-w) \blue{P_YP_Z} \in  \mathcal{P}_{\mathbb{Y} \times \mathbb{Z}}$ for all $w \in [0,1]$, and there exists $P_{YZ} \in \mathcal{P}_{\mathbb{Y} \times \mathbb{Z}}$ such that $\mathrm{HSIC}_{k \otimes \ell}(P_{YZ}) = \varrho_0$ for some fixed $\varrho_0>0$. Let $\alpha > \frac{1}{B+1}$ and $\alpha \in (0,1)$, $\beta \in (0,1-\alpha)$ be fixed values. Consider the high privacy regime where $\xi_{\varepsilon,\delta} \asymp  n^{-1/2 - r}$ with fixed $r \in (0,1/2)$, for $\xi_{\varepsilon,\delta}$ as in \eqref{Eq: definition of xi}. Then the uniform power of $\phi_{\mathttt{dpHSIC}}^u$ is asymptotically at most $\alpha$ over $ \mathcal{P}_{\mathrm{HSIC}_{k \otimes \ell}}\!(\rho)$ where $\rho$ is given in \eqref{Eq: separation MMD rate for U-statistics}. In other words, it holds that 
	\begin{align*}
		\limsup_{n \rightarrow \infty} \inf_{P_{YZ} \in \mathcal{P}_{\mathrm{HSIC}_{k \otimes \ell}}\!(\rho)} \mE_{P_{YZ}}[\phi_{\mathttt{dpHSIC}}^u]  \leq \alpha.
	\end{align*}
\end{theorem}

The comments made for \Cref{Theorem: Suboptimality of U-MMD} also apply to \Cref{Theorem: Suboptimality of U_HSIC}. Specifically, as in the MMD case, the U-statistic of the HSIC requires a higher noise level than the V-statistic to ensure differential privacy through the Laplace mechanism. To be more precise, $\phi_{\mathttt{dpHSIC}}^u$ becomes powerful in the privacy regime when the target parameter $\mathrm{HSIC}_{k \otimes \ell}^2$ is greater than the Laplace noise level $(n\xi_{\varepsilon,\delta})^{-1}$ or equivalently $\mathrm{HSIC}_{k \otimes \ell}$ is greater than $(n\xi_{\varepsilon,\delta})^{-1/2}$. However, the second term $\min\{(n\xi_{\varepsilon,\delta})^{-1},1\}$ in the minimax separation rate (\Cref{Theorem: Minimax separation in HSIC}) is smaller than $(n\xi_{\varepsilon,\delta})^{-1/2}$, and thus we may see that $\phi_{\mathttt{dpHSIC}}^u$ becomes powerless when $\min\{(n\xi_{\varepsilon,\delta})^{-1},1\} \lesssim \mathrm{HSIC}_{k \otimes \ell} \lesssim (n\xi_{\varepsilon,\delta})^{-1/2}$. 
We make this intuition more precise in the proof of \Cref{Theorem: Suboptimality of U_HSIC} and establish that the worse-case power is essentially bounded by significance level $\alpha$.

Our next concern is the minimum separation of $\phi_{\mathttt{dpHSIC}}^u$ test over $\tilde{\mathcal{P}}_{L_2}^{s}$, which mirrors \Cref{Theorem: Minimum separation of U-stat over L2} for the U-statistic of the MMD. As in \Cref{Theorem: Minimum separation of U-stat over L2}, we focus on the Gaussian kernels for the HSIC, and establish the following result. The proof can be found in \Cref{Section: Proof of Theorem: Minimum separation of HSIC U-stat over L2}. 

\begin{theorem}[Minimum Separation of $\phi_{\mathttt{dpHSIC}}^u$ over $\tilde{\mathcal{P}}_{L_2}^{s}$] \label{Theorem: Minimum separation of HSIC U-stat over L2} 
	Assume that $\alpha \in (0, 1)$, $\beta \in (0,1-\alpha)$, $\varepsilon > 0$, $\delta \in [0,1)$, $B \geq 6\alpha^{-1}\log(2\beta^{-1})$, $\prod_{i=1}^{d_Y} \lambda_i \leq 1$ and $\prod_{i=1}^{d_Z} \mu_i \leq 1$. Then the minimum separation of $\phi_{\mathttt{dpHSIC}}^u$ with Gaussian kernels over $\tilde{\mathcal{P}}_{L_2}^{s}$ is upper bounded as 
	\begin{align*}
		\rho_{\phi_{\mathttt{dpHSIC}}^u,L_2}^2 \leq C_{\alpha,\beta,s,R,M,d_Y,d_Z} \Bigg\{  \sum_{i=1}^{d_Y} \lambda_i^{2s} + \sum_{i=1}^{d_Z} \mu_i^{2s} & + \frac{1}{n\sqrt{\lambda_1 \cdots \lambda_{d_Y} \mu_1 \cdots \mu_{d_Z}}} \\[.5em]
		& + \frac{1}{n \lambda_1 \cdots \lambda_{d_Y} \mu_1 \cdots \mu_{d_Z} \xi_{\varepsilon,\delta}} \Bigg\},
	\end{align*}
	where $C_{\alpha,\beta,s,R,M,d_Y,d_Z}$ is a positive constant, depending only on $\alpha,\beta,R,M,d_Y,d_Z$, and $\xi_{\varepsilon,\delta}$ is given in \eqref{Eq: definition of xi}.
\end{theorem}

By comparing \Cref{Theorem: Minimum separation of HSIC U-stat over L2} with \Cref{Theorem: Minimax Separation over L2 for HSIC}, we see that the upper bound for $\rho_{\phi_{\mathttt{dpHSIC}}^u,L_2}^2$ is smaller than that of $\rho_{\phi_{\mathttt{dpHSIC}},L_2}^2$ if $n\xi_{\varepsilon,\delta} \lesssim 1$. Since we assume $\prod_{i=1}^{d_Y} \lambda_i \leq 1$ and $\prod_{i=1}^{d_Z} \mu_i \leq 1$, the condition~$n\xi_{\varepsilon,\delta} \lesssim 1$ implies that $\|p_{YZ} - p_Yp_Z\|_{L_2}$ needs to be sufficiently large to ensure significant power. However, the $L_2$ distance cannot be made sufficiently large as $\|p_{YZ}\|_{L_\infty}$ and $\|p_Yp_Z\|_{L_\infty}$ are assumed to be bounded. Therefore, the proposed test $\phi_{\mathttt{dpHSIC}}$ is more favorable than the one based on the U-statistic in terms of obtaining a tight separation rate in the $L_2$ distance. Nevertheless, when $\xi_{\varepsilon,\delta}$ is sufficiently large, the minimum separation satisfies 
	\begin{align*}
	\rho_{\phi_{\mathttt{dpHSIC}}^u,L_2}^2 \leq C_{\alpha,\beta,s,R,M,d_Y,d_Z} \Bigg\{  \sum_{i=1}^{d_Y} \lambda_i^{2s} + \sum_{i=1}^{d_Z} \mu_i^{2s}  + \frac{1}{n\sqrt{\lambda_1 \cdots \lambda_{d_Y} \mu_1 \cdots \mu_{d_Z}}} \Bigg\}.
\end{align*}
Therefore, in low privacy regimes, both $\phi_{\mathttt{dpHSIC}}$ and $\phi_{\mathttt{dpHSIC}}^u$ tests achieve the same separation rate in terms of the $L_2$ distance, which can be optimal when $\lambda_1 = \cdots = \lambda_{d_Y} = \mu_1 =\cdots = \mu_{d_Z} \asymp n^{-2/(4s+d_Y+d_Z)}$.

Like the negative results of $U_{\mathrm{MMD}}$, our negative results of $U_{\mathrm{HSIC}}$ are based solely on the Laplace mechanism. It is currently unknown whether there is an alternative privacy mechanism, enabling a test based on $U_{\mathrm{HSIC}}$ to be optimal in high privacy regimes. We leave this important direction for future work.

\subsection{Analyses of Minimum Separation Rates} \label{Section: Interpretation}
In this section, we provide details on the separation rate stated in \eqref{Eq: L2 separation rate}. Throughout our discussion, we treat the level $\alpha$ as a fixed constant. 
Let us recall the three terms in the separation rate: 
\begin{align*}
	(\mathrm{I}) &= \frac{\log(1/\alpha)}{n\sqrt{\lambda_1\cdots \lambda_d}}, \\ 
	(\mathrm{II}) &= \frac{\log(1/\alpha)}{n^{3/2}\lambda_1\cdots \lambda_d \xi_{\varepsilon,\delta}}, \\ 
	(\mathrm{III}) &=  \frac{\log^2(1/\alpha)}{n^2\lambda_1\cdots \lambda_d \xi_{\varepsilon,\delta}^2}.
\end{align*}
Among these three terms, the dominant term becomes 
\begin{align*}
	\begin{cases}
		\mathrm{(I)},  \quad & \text{if $n^{-1/2}(\lambda_1 \dots \lambda_d)^{-1/2} \lesssim \xi_{\varepsilon,\delta}$ (low privacy regime)}, \\[.5em]
		\mathrm{(II)}, \quad & \text{if $n^{-1/2} \lesssim \xi_{\varepsilon,\delta} \lesssim n^{-1/2}(\lambda_1 \dots \lambda_d)^{-1/2}$ (mid privacy regime),} \\[.5em]
		\mathrm{(III)}, \quad & \text{if $\xi_{\varepsilon,\delta} \lesssim n^{-1/2}$ (high privacy regime).}
	\end{cases}
\end{align*}
We choose the bandwidth parameters so that each dominating term matches the order of $\sum_{i=1}^d \lambda_i^{2s}$ and compute the resulting rate explained below.
\begin{itemize}
	\item  \textbf{Low privacy regime.} First, in the low privacy regime, we have the separation rate satisfying 
	\begin{align*}
		\rho_{\phi_{\mathttt{dpMMD}},L_2}^2 \leq C_{\tau,\beta,s,R,M,d} \Biggl\{ \sum_{i=1}^d \lambda_i^{2s} + \frac{\log(1/\alpha)}{n\sqrt{\lambda_1\cdots \lambda_d}}  \Bigg\}.
	\end{align*}
	Equating $\mathrm{(I)}$ with $\sum_{i=1}^{d} \lambda_i^{2s}$ yields optimal bandwidths $\lambda_i = n^{-\frac{2}{4s+d}}$ for $i=1,\ldots,d$, and the resulting rate is $n^{-\frac{4s}{4s+d}}$, \emph{i.e.}, $\rho_{\phi_{\mathttt{dpMMD}},L_2}^2 \lesssim n^{-\frac{4s}{4s+d}}$. This rate is known to be minimax optimal over the Sobolev ball under the non-DP setting~\citep{li2019optimality,schrab2021mmd}. Therefore, when $\lambda_i = n^{-\frac{2}{4s+d}}$ for $i \in [d]$, we can achieve an optimal separation rate whenever $\xi_{\varepsilon,\delta} \gtrsim n^{- \frac{2s-d/2}{4s+d}}$. 
	\item  \textbf{Mid privacy regime.}  Second, in the mid privacy regime, we have the separation rate satisfying 
	\begin{align*}
		\rho_{\phi_{\mathttt{dpMMD}},L_2}^2 \leq C_{\tau,\beta,s,R,M,d} \Biggl\{ \sum_{i=1}^d \lambda_i^{2s} + \frac{\log(1/\alpha)}{n^{3/2}\lambda_1\cdots \lambda_d \xi_{\varepsilon,\delta}} \Bigg\}.
	\end{align*}
	Equating $\mathrm{(II)}$ with $\sum_{i=1}^{d} \lambda_i^{2s}$ yields optimal bandwidths $\lambda_i =(n^{3/2} \xi_{\varepsilon,\delta})^{- \frac{1}{2s+d}}$ for $i=1,\ldots,d$, and the resulting rate becomes $(n^{3/2} \xi_{\varepsilon,\delta})^{- \frac{2s}{2s+d}}$. Therefore, we can achieve the separation rate $(n^{3/2} \xi_{\varepsilon,\delta})^{- \frac{2s}{2s+d}}$ using bandwidths $\lambda_i =(n^{3/2} \xi_{\varepsilon,\delta})^{- \frac{1}{2s+d}}$ for $i \in [d]$ whenever $n^{-1/2} \lesssim \xi_{\varepsilon,\delta} \lesssim n^{-1/2}(\lambda_1 \dots \lambda_d)^{-1/2}$, which is equivalent to 
	\begin{align*}
		n^{-1/2} \lesssim \xi_{\varepsilon,\delta} \lesssim n^{-\frac{2s-d/2}{4s+d}}.
	\end{align*}
	\item \textbf{High privacy regime.} Lastly, in the high privacy regime, the separation rate satisfies 
	\begin{align*}
	\rho_{\phi_{\mathttt{dpMMD}},L_2}^2 \leq C_{\tau,\beta,s,R,M,d} \Biggl\{ \sum_{i=1}^d \lambda_i^{2s} + \frac{\log^2(1/\alpha)}{n^2\lambda_1\cdots \lambda_d \xi_{\varepsilon,\delta}^2}  \Bigg\}.
	\end{align*}
	Equating $\mathrm{(III)}$ with $\sum_{i=1}^{d} \lambda_i^{2s}$ yields optimal bandwidths $\lambda_i =(n \xi_{\varepsilon,\delta})^{- \frac{2}{2s+d}}$ for $i \in [d]$, and the resulting rate becomes $(n\xi_{\varepsilon,\delta})^{-\frac{4s}{2s+d}}$. Consequently, we have the separation rate $(n\xi_{\varepsilon,\delta})^{-\frac{4s}{2s+d}}$ whenever $\xi_{\varepsilon,\delta} \lesssim n^{-1/2}$. 
\end{itemize}
To summarize, we can achieve the separation rate using different bandwidths in different privacy regimes as
\begin{align*}
	\rho_{\phi_{\mathttt{dpMMD}},L_2}^2 \lesssim \begin{cases}
		n^{-\frac{4s}{4s+d}},  \quad & \text{in the low privacy regime with $n^{- \frac{2s-d/2}{4s+d}} \lesssim \xi_{\varepsilon,\delta}$}, \\[.5em]
		(n^{3/2} \xi_{\varepsilon,\delta})^{- \frac{2s}{2s+d}},  \quad &  \text{in the mid privacy regime with $n^{-1/2} \lesssim \xi_{\varepsilon,\delta} \lesssim n^{-\frac{2s-d/2}{4s+d}}$},\\[.5em]
		(n\xi_{\varepsilon,\delta})^{-\frac{4s}{2s+d}},  \quad &  \text{in the high privacy regime with $\xi_{\varepsilon,\delta} \lesssim n^{-1/2}$}.
	\end{cases}
\end{align*}
Each separation rate can be attained by using the bandwidths $\lambda_1 = \dots = \lambda_d = n^{-\frac{2}{4s+d}}$ for the low privacy regime, $\lambda_1 = \dots = \lambda_d=(n^{3/2}  \xi_{\varepsilon,\delta})^{- \frac{1}{2s+d}}$ for the mid privacy regime and $\lambda_1 = \dots = \lambda_d = (n \xi_{\varepsilon,\delta})^{- \frac{2}{2s+d}}$ for the high privacy regime.

\section{Additional Simulations} 
\label{Section: Additional Simulations}

In this section, we present further experiments on:
\begin{enumerate}[(i)]
	\item \blue{two-sample testing on Parkinson's telemonitoring data (\Cref{subsec:parkinsons}),}
	\item \blue{robustness to unequal sample sizes (\Cref{subsec:different_sample_sizes}),}
	\item independence testing (\Cref{subsec:hsic_uniform}),
	\item power in low privacy regimes (\Cref{subsec:low_privacy}) and
	\item test significance levels (\Cref{subsec:level}).
\end{enumerate}
These results are all analyzed in \Cref{subsec:analysis_experiments}.

{
\subsection{Parkinson's Telemonitoring Data}
\label{subsec:parkinsons}

To further demonstrate the practical utility of our approach, we consider the \href{https://archive.ics.uci.edu/dataset/189/parkinsons+telemonitoring}{Parkinson's telemonitoring dataset}~\citep{tsanas2009accurate}, which consists of $19$ features extracted from $5{,}875$ biomedical voice measurements of patients with early-stage Parkinson's disease. Because speech data are sensitive and personal, tests with differential privacy guarantees are particularly relevant in this setting.

We partition the measurements into `mild' and `severe' groups according to whether the total Unified Parkinson's Disease Rating Scale score is below or above $30$.\footnote{The UPDRS scores are used only to form the groups and are not among the $19$ tested features.} As in the CelebA experiments, one sample is drawn from the mild group, while observations in the other sample are drawn from the severe and mild groups with probabilities $c$ and $1-c$, respectively. We vary the privacy level, the common sample size, and the corruption parameter in \Cref{fig:parkinsons_mmd_privacy,fig:parkinsons_mmd_sample_size,fig:parkinsons_mmd_corruption}. The dpMMD test outperforms the competing tests across all three privacy regimes and experimental settings. The qualitative trends agree with those observed for the perturbed-uniform and CelebA experiments.

\begin{figure}[!htbp]
	\centering
	\includegraphics[width=0.9\textwidth]{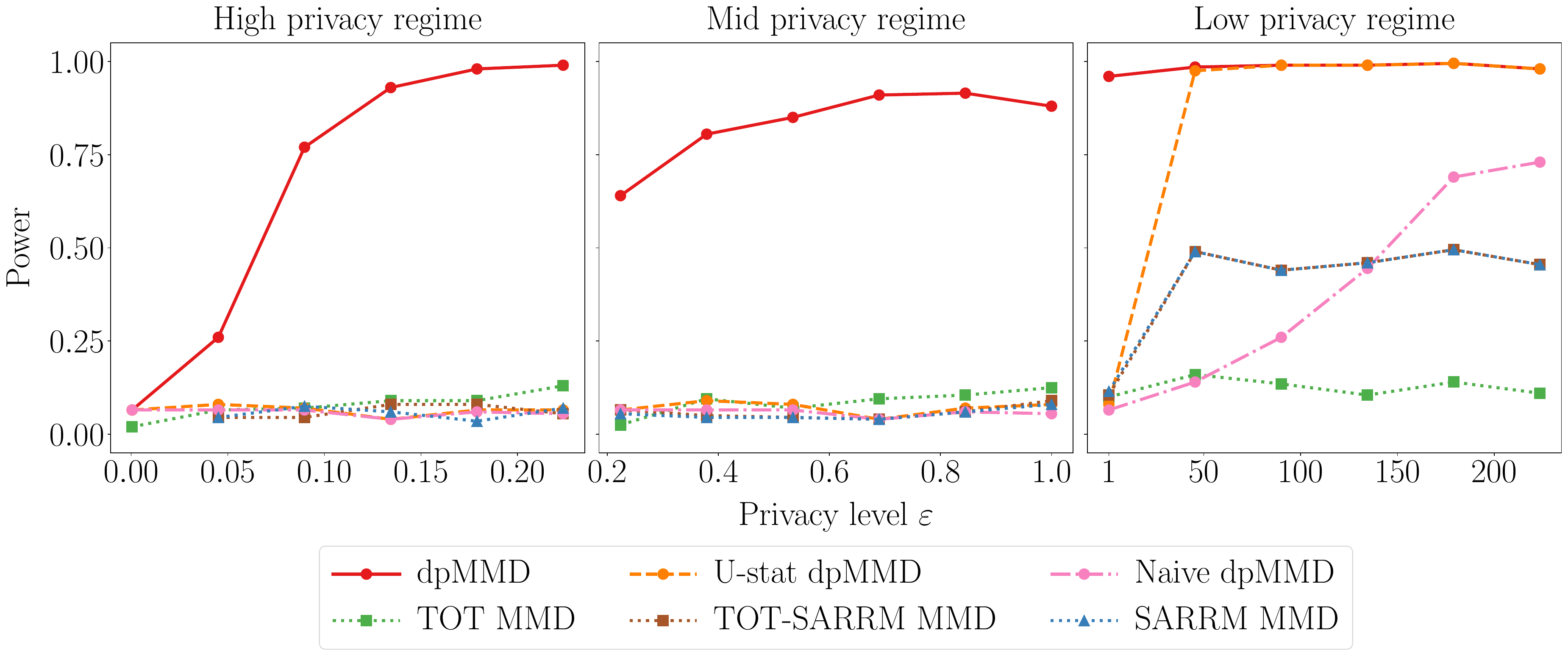}
	\captionsetup{format=hang}
	\caption{Comparing speech data of patients with mild/severe Parkinson's disease while varying $\varepsilon$. We set $m=n=2000$ and use: \emph{(Left)} $\varepsilon$ from $1/n$ to $10/\sqrt n$, $c=1$; \emph{(Middle)} $\varepsilon$ from $10/\sqrt n$ to $1$, $c=0.6$; \emph{(Right)} $\varepsilon$ from $1$ to $\sqrt n$, $c=0.65$.}
	\label{fig:parkinsons_mmd_privacy}
\end{figure}

\begin{figure}[!htbp]
	\centering
	\includegraphics[width=0.9\textwidth]{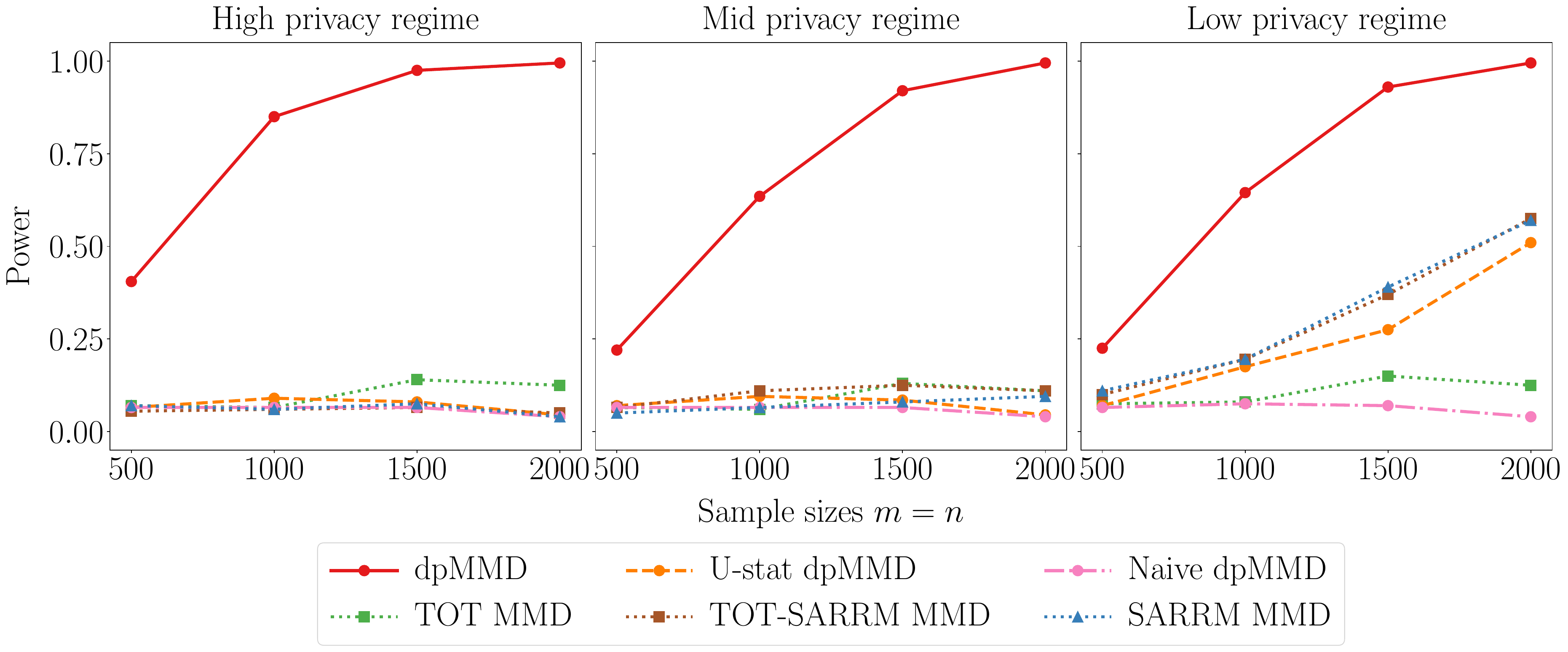}
	\captionsetup{format=hang}
	\caption{Comparing speech data of patients with mild/severe Parkinson's disease while varying $m=n$. We use: \emph{(Left)} $\varepsilon=10/\sqrt n$, $c=1$; \emph{(Middle)} $\varepsilon=1$, $c=0.7$; \emph{(Right)} $\varepsilon=\sqrt n/10$, $c=0.7$.}
	\label{fig:parkinsons_mmd_sample_size}
\end{figure}

\begin{figure}[!htbp]
	\centering
	\includegraphics[width=0.9\textwidth]{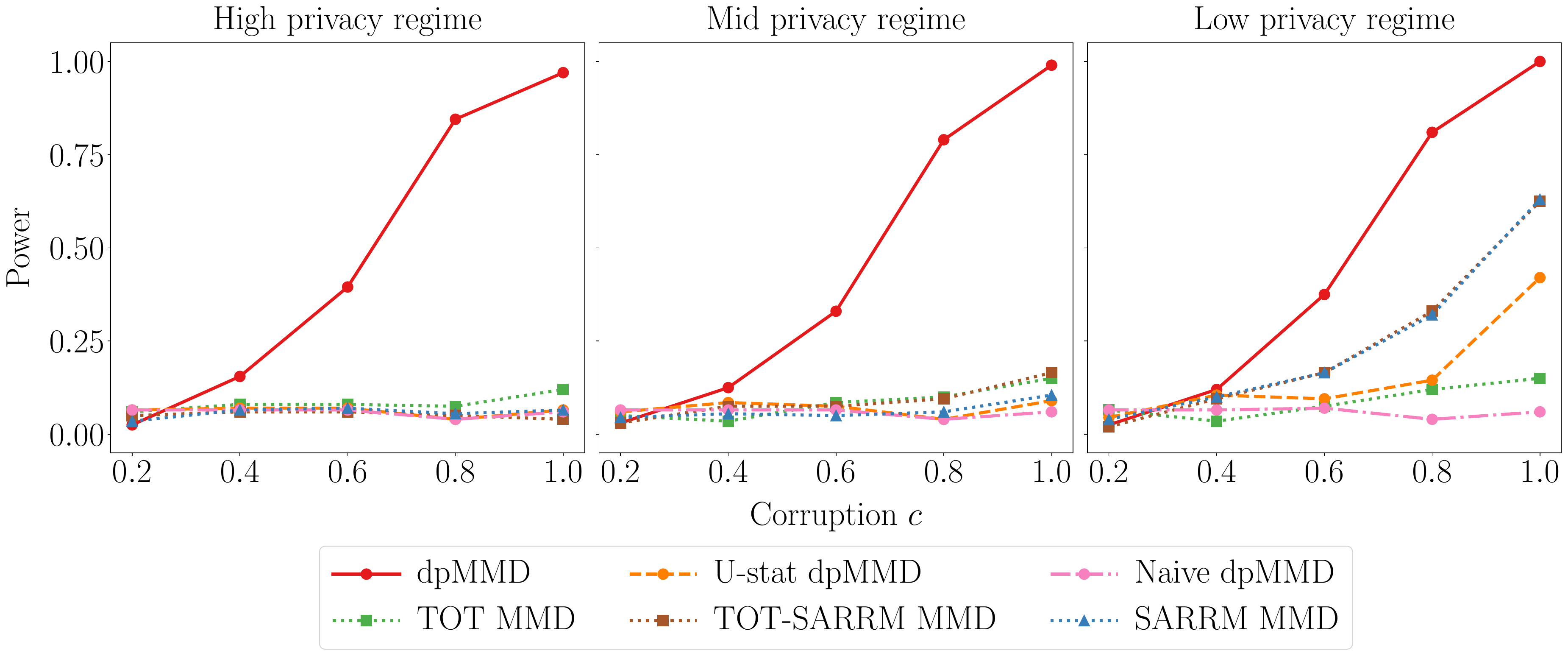}
	\captionsetup{format=hang}
	\caption{Comparing speech data of patients with mild/severe Parkinson's disease while varying $c$. We use: \emph{(Left)} $\varepsilon=10/\sqrt n$, $n=1500$; \emph{(Middle)} $\varepsilon=1$, $n=1000$; \emph{(Right)} $\varepsilon=\sqrt n/10$, $n=1000$.}
	\label{fig:parkinsons_mmd_corruption}
\end{figure}

\FloatBarrier

\subsection{Two-Sample Testing with Unequal Sample Sizes}
\label{subsec:different_sample_sizes}

We next fix one sample size $n$ and vary the other sample size $m$. In \Cref{fig:uniform_mmd_sample_size_different,fig:celeba_mmd_sample_size_different,fig:parkinsons_mmd_sample_size_different}, the power of dpMMD increases substantially as $m$ approaches $n$ from below. This behavior reflects the $1/\min(m,n)$ scaling of the privatization noise. Once $m>n$, the noise scale remains fixed at $1/n$, so further increases in $m$ yield much smaller gains, especially in the high-privacy regime.

\begin{figure}[!htbp]
	\centering
	\includegraphics[width=\textwidth]{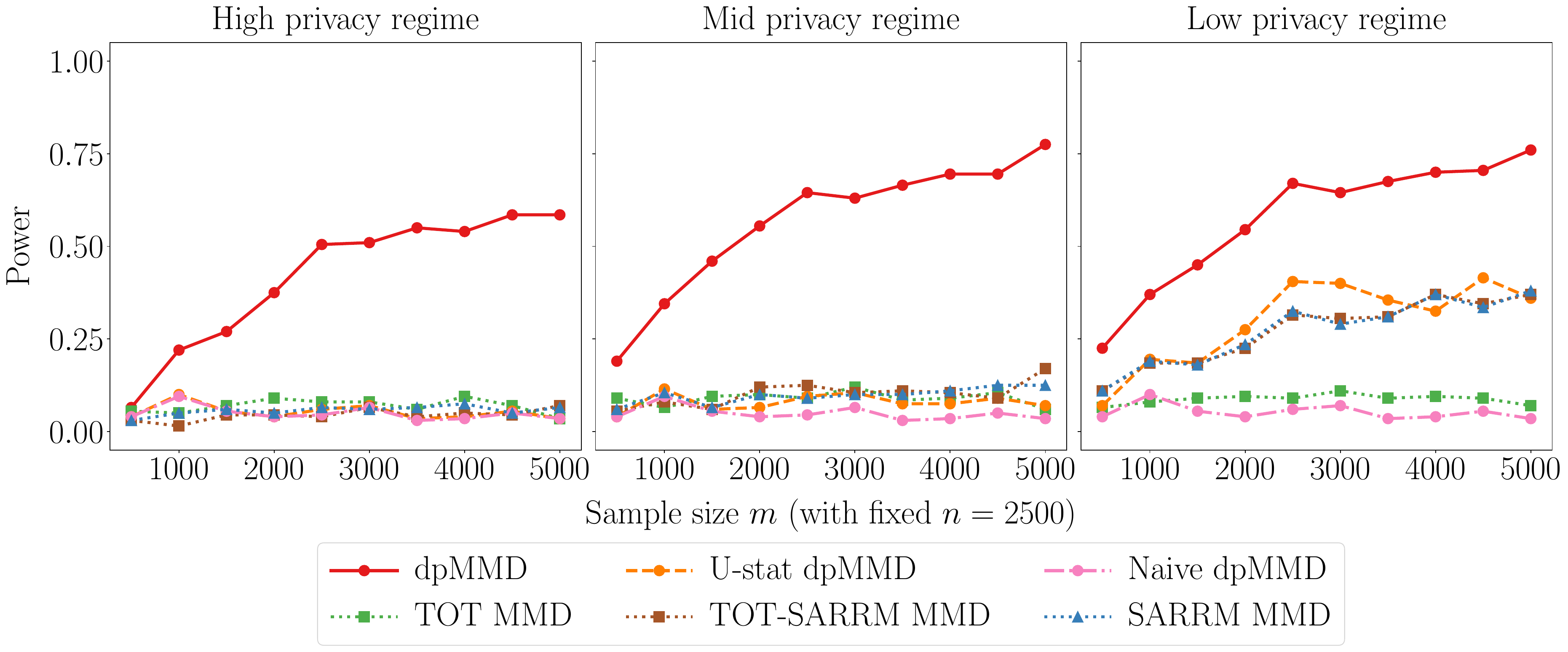}
	\captionsetup{format=hang}
	\caption{Comparing uniform vs.~perturbed uniform while varying $m$ and fixing $n=2500$, $d=1$, and $a=0.1$. We use: \emph{(Left)} $\varepsilon=10/\sqrt n$; \emph{(Middle)} $\varepsilon=1$; \emph{(Right)} $\varepsilon=\sqrt n/10$.}
	\label{fig:uniform_mmd_sample_size_different}
\end{figure}

\begin{figure}[!htbp]
	\centering
	\includegraphics[width=\textwidth]{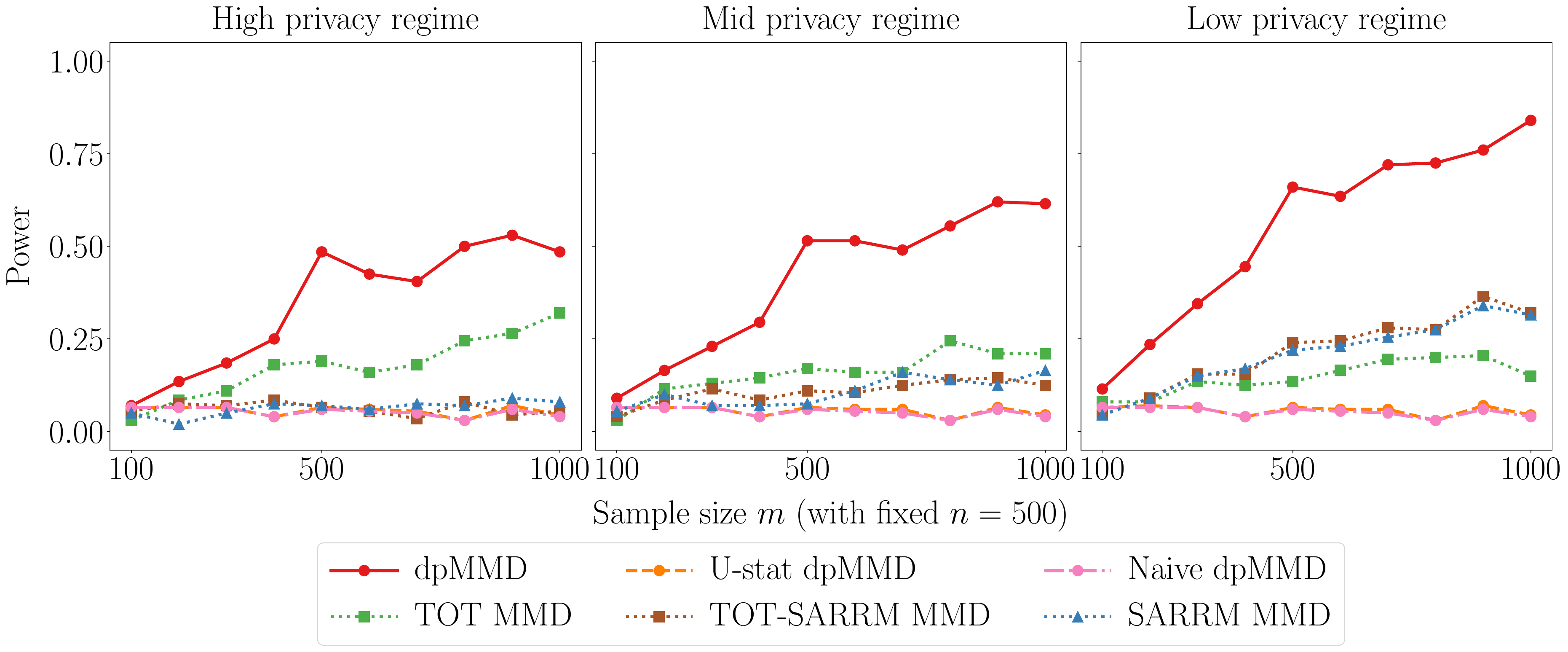}
	\captionsetup{format=hang}
	\caption{Comparing CelebA images while varying $m$ and fixing $n=500$. We use: \emph{(Left)} $\varepsilon=10/\sqrt n$, $c=0.7$; \emph{(Middle)} $\varepsilon=1$, $c=0.5$; \emph{(Right)} $\varepsilon=\sqrt n/10$, $c=0.45$.}
	\label{fig:celeba_mmd_sample_size_different}
\end{figure}

\begin{figure}[!htbp]
	\centering
	\includegraphics[width=\textwidth]{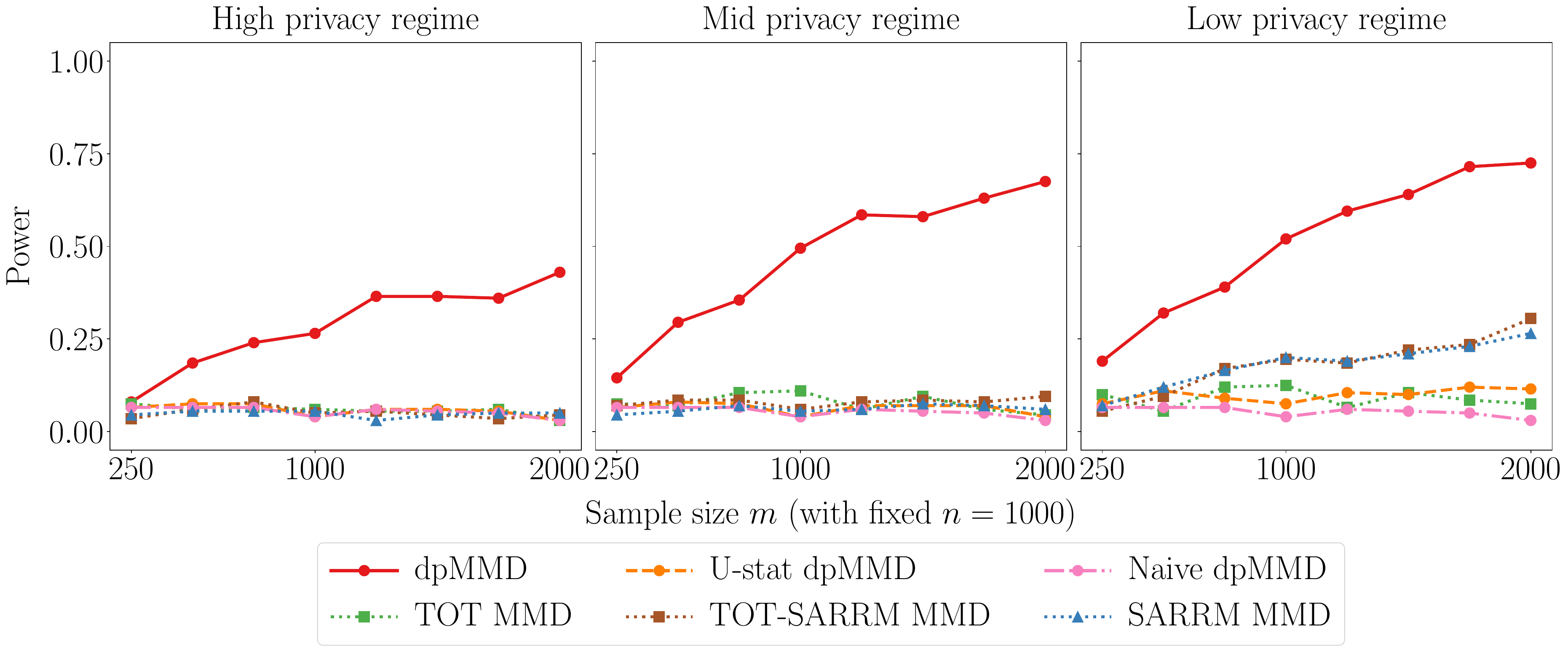}
	\captionsetup{format=hang}
	\caption{Comparing speech data of patients with mild/severe Parkinson's disease while varying $m$ and fixing $n=1000$ and $c=0.65$. We use: \emph{(Left)} $\varepsilon=10/\sqrt n$; \emph{(Middle)} $\varepsilon=1$; \emph{(Right)} $\varepsilon=\sqrt n/10$.}
	\label{fig:parkinsons_mmd_sample_size_different}
\end{figure}

\FloatBarrier
}

\subsection{Perturbed Uniform Distributions for Independence Testing}
\label{subsec:hsic_uniform}

\begin{figure}[!htbp]
	\centering
	\includegraphics[width=0.72\textwidth]{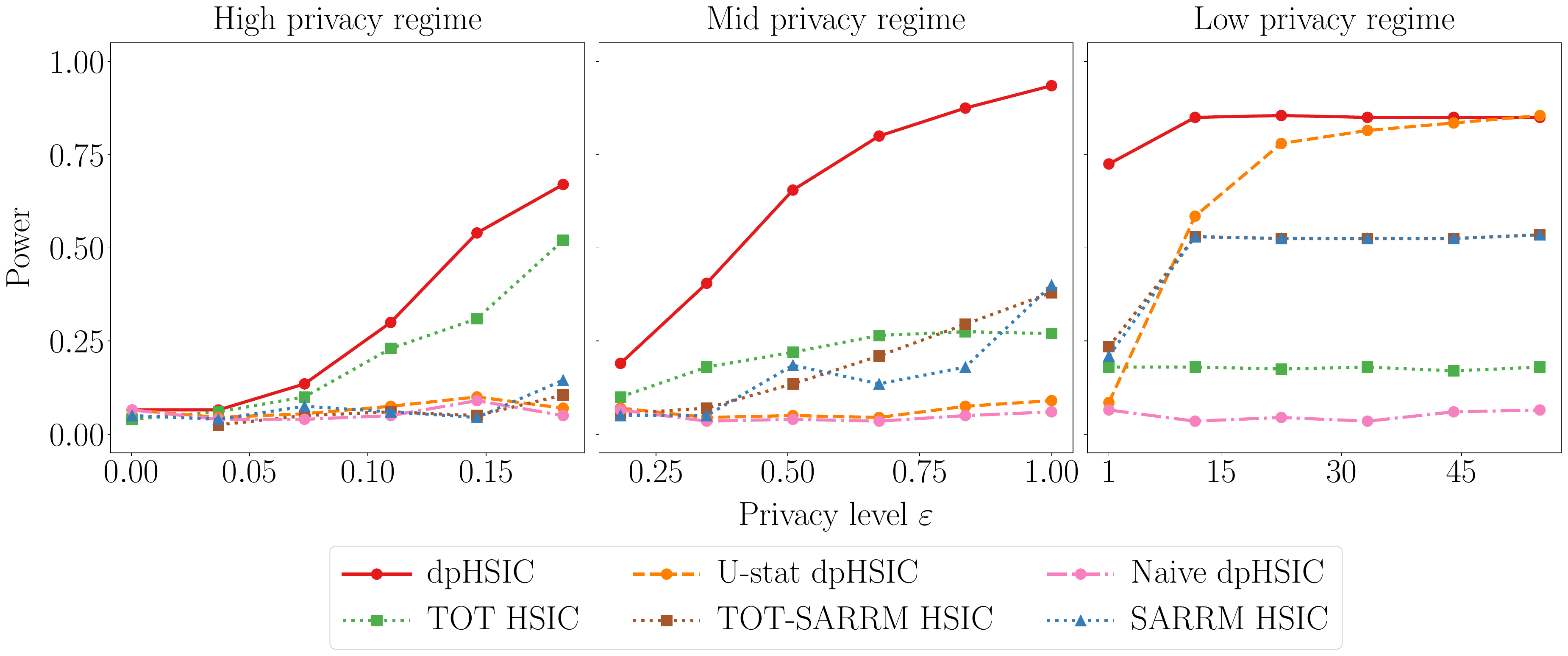}
	\captionsetup{format=hang}
	\caption{
		Measuring the dependence on joint perturbed uniform distributions while varying the privacy level $\varepsilon$. We set the sample size $n = 3000$ and dimensions $d_X=d_Y=1$. We change the privacy level and perturbation amplitude as follows:
		\emph{(Left)} Privacy level $\varepsilon$ from $1/n$ to $10/\sqrt{n}$, perturbation amplitude $a=0.4$. 
		\emph{(Middle)} Privacy level $\varepsilon$ from $10/\sqrt{n}$ to $1$, perturbation amplitude $a=0.2$. 
		\emph{(Right)} Privacy level $\varepsilon$ from $1$ to $\sqrt{n}$, perturbation amplitude $a=0.15$. 
	}
	\label{fig:uniform_hsic_privacy}
\end{figure}

\begin{figure}[!htbp]
	\centering
	\includegraphics[width=0.72\textwidth]{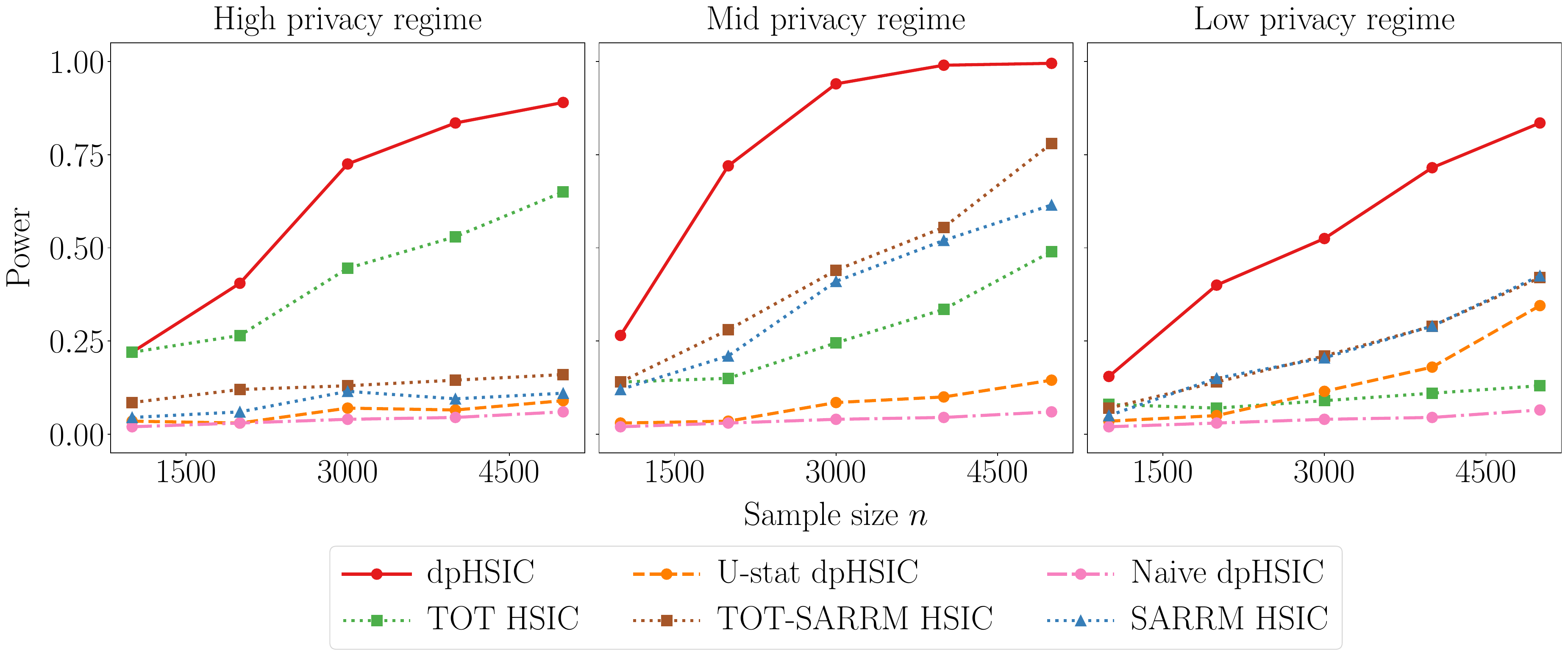}
	\captionsetup{format=hang}
	\caption{
		Measuring the dependence on joint perturbed uniform distributions while varying the sample size $n$. We set the dimensions $d_X=d_Y=1$, and change other parameters as follows: 
		\emph{(Left)} Privacy level $\varepsilon=10/\sqrt{n}$, perturbation amplitude $a=0.4$. 
		\emph{(Middle)} Privacy level $\varepsilon=1$, perturbation amplitude $a=0.2$. 
		\emph{(Right)} Privacy level $\varepsilon=\sqrt{n}/10$, perturbation amplitude $a=0.1$. 
	}
	\label{fig:uniform_hsic_sample_size}
\end{figure}

\begin{figure}[!htbp]
	\centering
	\includegraphics[width=0.8\textwidth]{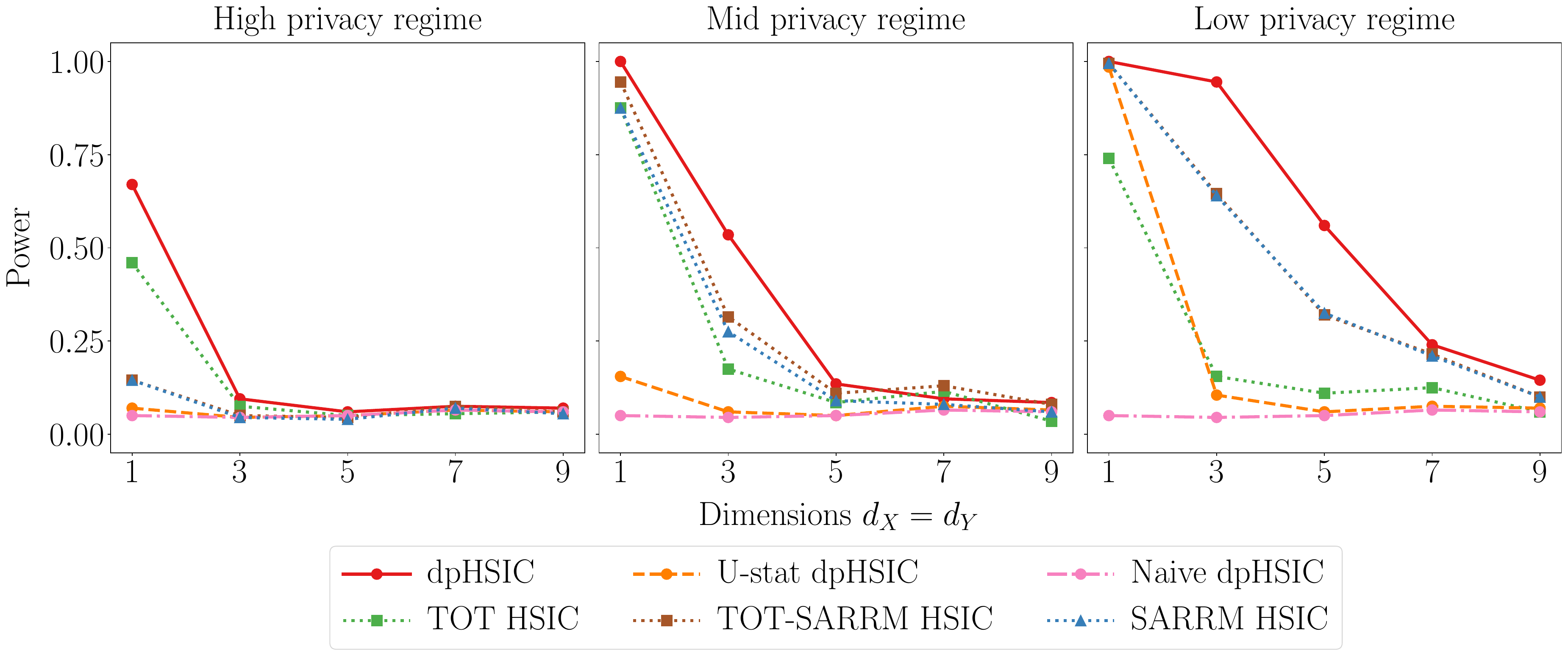}
	\captionsetup{format=hang}
	\caption{
		Measuring the dependence on joint perturbed uniform distributions while varying the dimensions $d_X=d_Y$. We set the sample size $n = 3000$, and change other parameters as follows:
		\emph{(Left)} Privacy level $\varepsilon=10/\sqrt{n}$, perturbation amplitude $a=0.4$. 
		\emph{(Middle)} Privacy level $\varepsilon=1$, perturbation amplitude $a=0.35$. 
		\emph{(Right)} Privacy level $\varepsilon=\sqrt{n}/10$, perturbation amplitude $a=0.3$. 
	}
	\label{fig:uniform_hsic_dimension}
\end{figure}

We run independence testing simulations comparing the power of dpHSIC, U-Stat dpHSIC, Naive dpHSIC, TOT HSIC and SARRM HSIC, \blue{as well as TOT--SARRM HSIC,} where the samples are drawn from the perturbed uniform joint density defined in \Cref{eq:perturbed_uniform_density}.
The privacy level, sample size, and dimension are varied in \Cref{fig:uniform_hsic_privacy,fig:uniform_hsic_sample_size,fig:uniform_hsic_dimension}, respectively.
As analyzed in \Cref{subsec:analysis_experiments}, the same power trends are observed as in the two-sample MMD case of \Cref{subsec:mmd_uniform}.
Overall, dpHSIC achieves the highest power of all independence DP tests in all experimental settings and across all privacy regimes.

\FloatBarrier

\begin{figure}[!htbp]
	\centering
	\includegraphics[width=\textwidth]{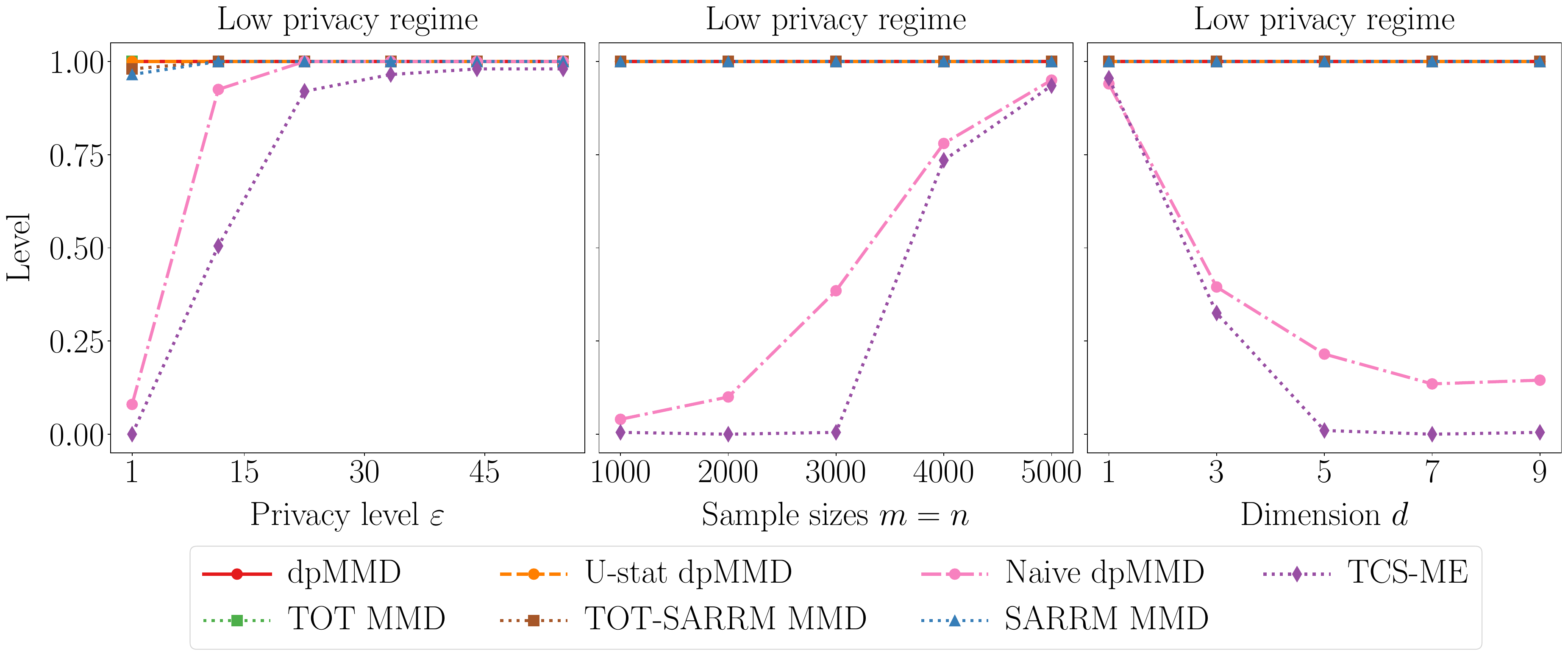}
	\captionsetup{format=hang}
	\caption{
		Comparing uniform vs.~perturbed uniform in low privacy regimes. We set the perturbation amplitude $a=1$ and change parameters as follows:
		\emph{(Left)} Privacy level $\varepsilon$ from $1$ to $\sqrt{n}$, sample sizes $m = n = 3000$, dimension $d=1$.
		\emph{(Middle)} Privacy level $\varepsilon=\sqrt{n}/10$, dimension $d=1$.
		\emph{(Right)} Privacy level $\varepsilon=\sqrt{n}/10$, sample sizes $m = n = 5000$.
	}
	\label{fig:uniform_mmd_appendix}
\end{figure}

\subsection{High-Signal \& Low-Privacy}
\label{subsec:low_privacy}

\begin{figure}[!htbp]
	\centering
	\includegraphics[width=\textwidth]{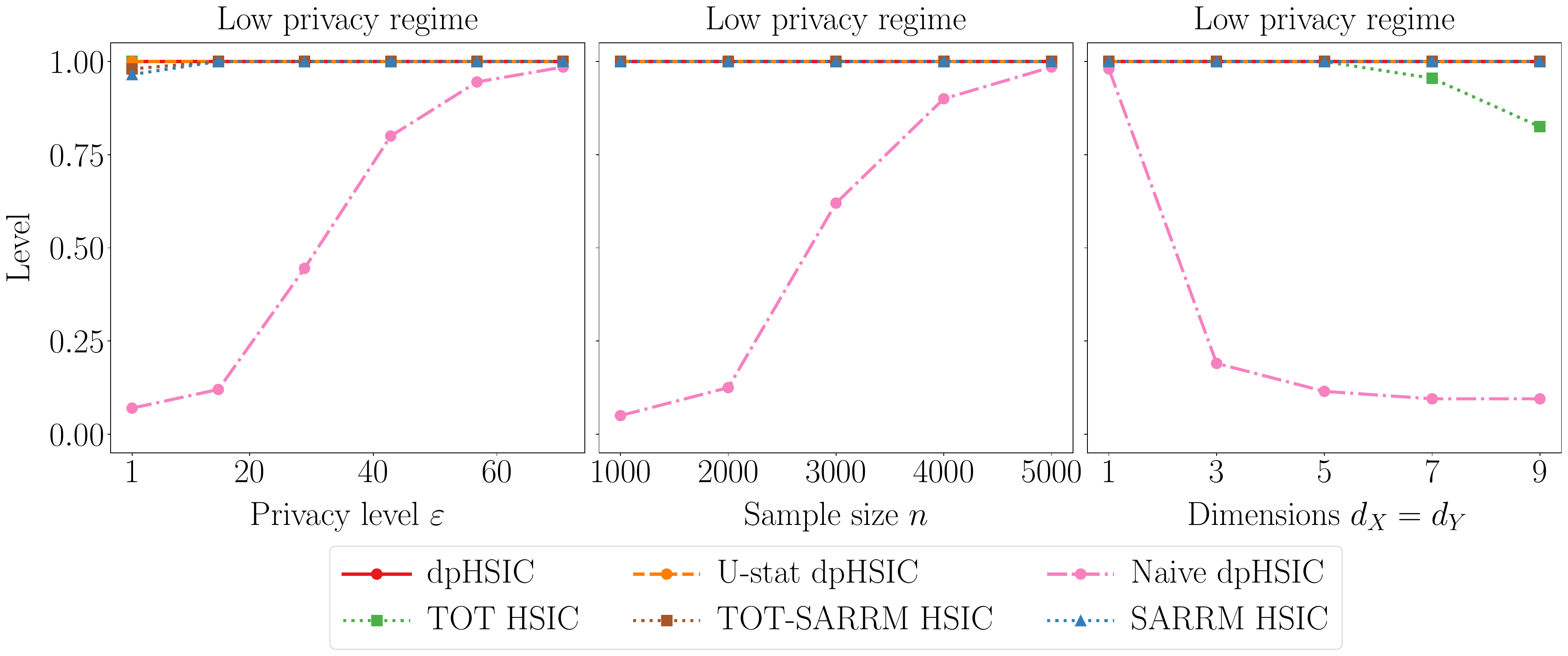}
	\captionsetup{format=hang}
	\caption{
		Measuring the dependence on joint perturbed uniform distributions in low privacy regimes. We set the perturbation amplitude $a=1$ and change parameters as follows:
		\emph{(Left)} Privacy level $\varepsilon$ from $1$ to $\sqrt{n}$, sample size $n = 5000$, dimensions $d_X=d_Y=1$.
		\emph{(Middle)} Privacy level $\varepsilon=\sqrt{n}$, dimensions $d_X=d_Y=1$.
		\emph{(Right)} Privacy level $\varepsilon=\sqrt{n}$, sample size $n = 5000$.
	}
	\label{fig:uniform_hsic_appendix}
\end{figure}

We present power results on perturbed uniform two-sample and independence simulations in the low privacy regime with a high signal $(a=1)$ in \Cref{fig:uniform_mmd_appendix,fig:uniform_hsic_appendix}.
This illustrates that TCS-ME, Naive dpMMD and Naive dpHSIC can indeed detect easier alternatives (while being almost powerless against the more challenging alternatives) and eventually attain power one.
We observe that TCS-ME and Naive dpMMD attain similar power, with Naive dpMMD actually being slightly more powerful.

\subsection{Level}
\label{subsec:level}

\begin{figure}[!htbp]
	\centering
	\includegraphics[width=\textwidth]{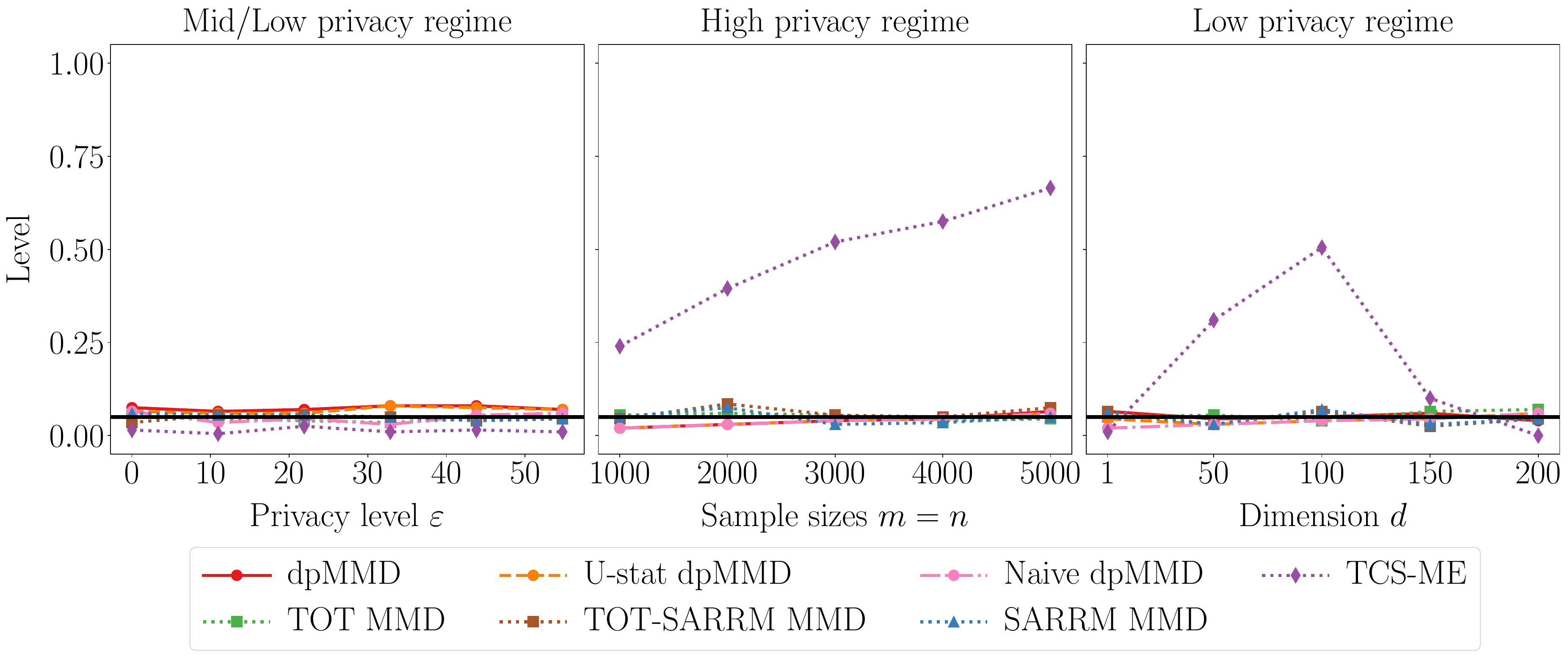}
	\captionsetup{format=hang}
	\caption{
		Type I error rates for two-sample testing under the uniform null distribution, \emph{i.e.}, with perturbation amplitude $a=0$. We vary parameters as follows:
		\emph{(Left)} Privacy level $\varepsilon$ from $1/\sqrt{n}$ to $\sqrt{n}$, sample sizes $m = n = 3000$, dimension $d=1$.
		\emph{(Middle)} Privacy level $\varepsilon=10/\sqrt{n}$, dimension $d=100$.
		\emph{(Right)} Privacy level $\varepsilon=\sqrt{n}/10$, sample sizes $m = n = 3000$.
	}
	\label{fig:uniform_mmd_level}
\end{figure}

\begin{figure}[!htbp]
	\centering
	\includegraphics[width=\textwidth]{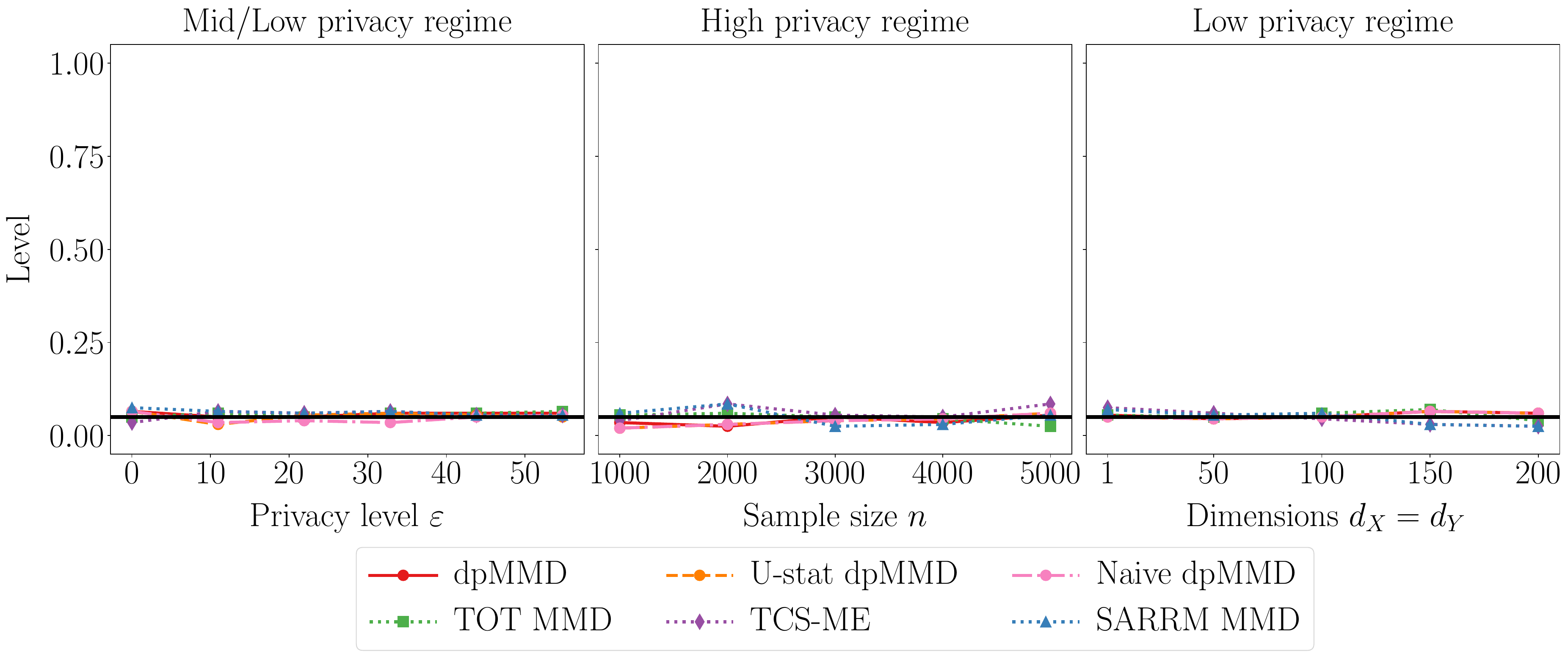}
	\captionsetup{format=hang}
	\caption{
		Type I error rates for independence testing under the uniform joint null distribution, \emph{i.e.}, with perturbation amplitude $a=0$. We vary parameters as follows:
		\emph{(Left)} Privacy level $\varepsilon$ from $1/\sqrt{n}$ to $\sqrt{n}$, sample size $n = 3000$, dimensions $d_X=d_Y=1$.
		\emph{(Middle)} Privacy level $\varepsilon=10/\sqrt{n}$, dimensions $d_X=d_Y=1$.
		\emph{(Right)} Privacy level $\varepsilon=\sqrt{n}/10$, sample size $n = 3000$.
	}
	\label{fig:uniform_hsic_level}
\end{figure}

\begin{figure}[!htbp]
	\centering
	\includegraphics[width=0.78\textwidth,trim=0 0 550bp 0,clip]{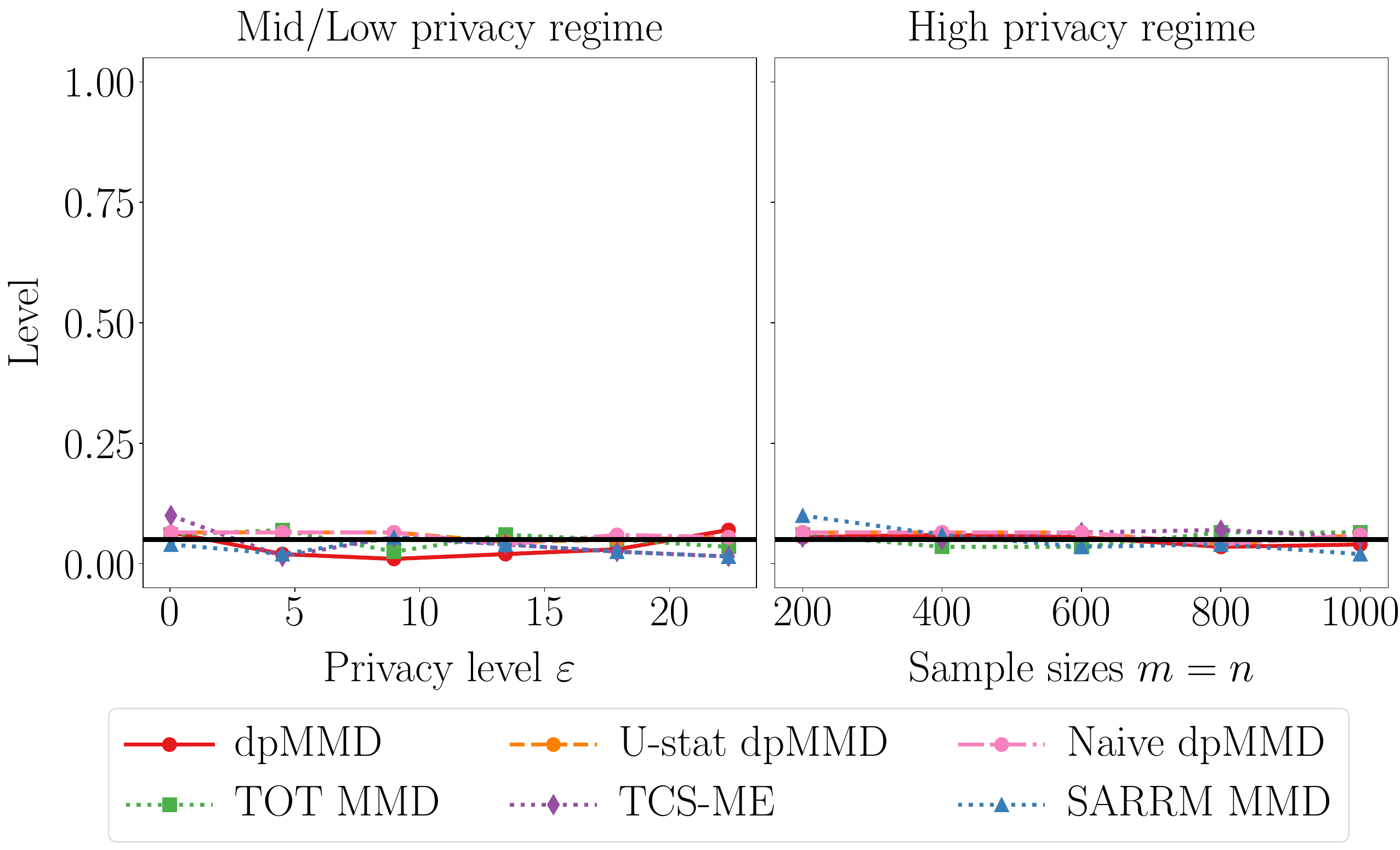}
	\captionsetup{format=hang}
	\caption{
		Type I error rates for two-sample testing for CelebA women images with zero-corruption $(c=0)$. We vary parameters as follows: \emph{(Left)} Privacy level $\varepsilon$ from $1/\sqrt{n}$ to $\sqrt{n}$, sample sizes $m = n = 500$.
		\emph{(Right)} Privacy level $\varepsilon=10/\sqrt{n}$.
	}
	\label{fig:celeba_mmd_level}
\end{figure}

We verify in \Cref{fig:uniform_mmd_level,fig:uniform_hsic_level,fig:celeba_mmd_level,fig:parkinsons_mmd_level} whether the tests are well-calibrated, \emph{i.e.}, their type I error rates are well-controlled at the significance level $\alpha=0.05$.
To this end, we run experiments under the null where the underlying distribution is uniform without any perturbation (\emph{i.e.}, the amplitude parameter $a=0$ in the settings described in \Cref{subsec:mmd_uniform} and \Cref{subsec:hsic_uniform}) for both two-sample and independence testing. We also evaluate the type I error rates of private MMD tests by drawing CelebA face images of women (\Cref{subsec:mmd_celeba}) for both samples (\emph{i.e.}, the corruption parameter $c=0$). \blue{For the Parkinson's experiment, both samples are drawn from the mild group.} \Cref{fig:uniform_mmd_level,fig:uniform_hsic_level,fig:celeba_mmd_level,fig:parkinsons_mmd_level} indicate that all tests are well-calibrated in all experimental settings, except TCS-ME.
Indeed, TCS-ME appears to be extremely miscalibrated, for example in dimension $d=100$ where its level is around $10\alpha=0.5$ instead of being around $\alpha=0.05$. We note, nonetheless, that TCS-ME controls the type I error correctly at level $\alpha$ in the one-dimensional case considered in \Cref{subsec:mmd_uniform}.

\FloatBarrier

{
\begin{figure}[!htbp]
	\centering
	\includegraphics[width=0.78\textwidth,trim=0 0 550bp 0,clip]{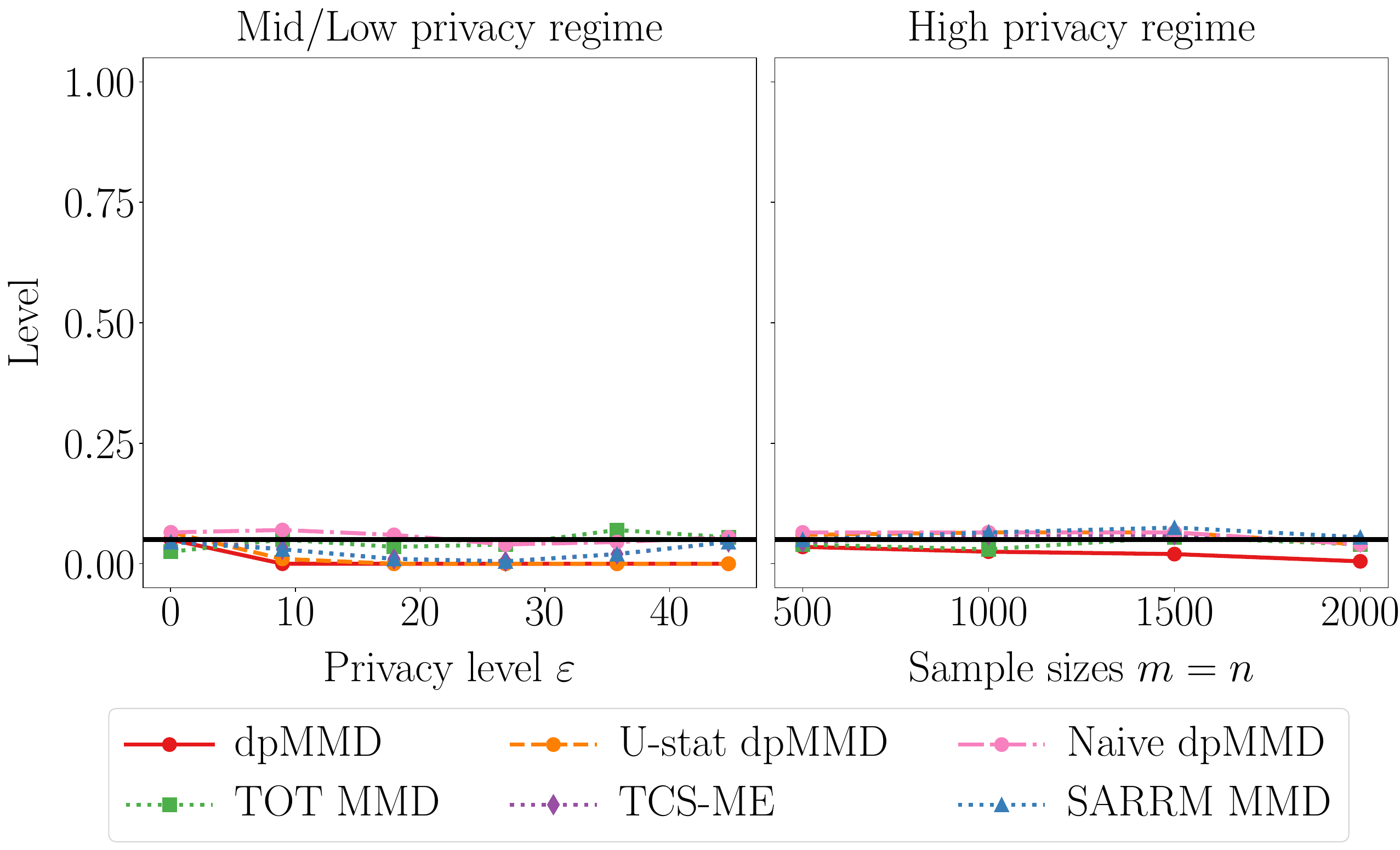}
	\captionsetup{format=hang}
	\caption{Type I error rates for two-sample testing on speech data from patients with mild Parkinson's disease, corresponding to $c=0$. \emph{(Left)} We vary $\varepsilon$ from $1/\sqrt n$ to $\sqrt n$ with $m=n=2000$. \emph{(Right)} We fix $\varepsilon=10/\sqrt n$. The nominal level $\alpha=0.05$ is shown as a black line.}
	\label{fig:parkinsons_mmd_level}
\end{figure}
}

\FloatBarrier

\section{Alternative private tests} 
\label{Section: Alternative private tests}

We first explain how standard MMD and HSIC permutation tests can be privatized using the methods of \citet{kazan2023test} in \Cref{subsec: tot} (TOT) and of \citet{pena2022differentially} in \Cref{subsec: sarrm} (SARRM),
both relying on the subsample-and-aggregate procedure of \citet{canonne2019structure}.
\blue{We also describe TOT--SARRM in \Cref{subsec: TOT-SARRM}, which combines the TOT test with the SARRM hyperparameter-selection strategy.}
We then present the approach of \citet{raj2020differentially} to constructing a private MMD test in \Cref{subsec: Raj} (TCS-ME).
We provide our own \href{https://github.com/antoninschrab/dpkernel-paper}{implementation} of TOT and SARRM in JAX, and use the \href{https://github.com/hcllaw/private_tst}{implementation} of \citet{raj2020differentially} directly for TCS-ME.

\subsection{TOT: Test of Tests}
\label{subsec: tot}

\citet{kazan2023test} propose a procedure called `Test of Tests' which allows us to privatize any existing test. 
The resulting test, which we refer to as TOT, is $(\varepsilon, 0)$-DP and has well-calibrated level \citep[Theorems 3.1 and 3.2]{kazan2023test}.
The test relies on the Truncated-Uniform-Laplace distribution $\textrm{Tulap}(b)$ \citep[Definition 4.1]{awan2018differentially} 
for positive real scale parameter $b$.
It is defined as the distribution with the following cumulative distribution function 
\begin{equation}
\label{eq: Tulap CDF}
F^{\textrm{Tulap}}_{b}(x) = \begin{cases}
(1+b)^{-1}\, b^{-[x]_{\textrm{near}}} \Big(b+\big(x-[x]_{\textrm{near}}+\frac 12\big)\big(1-b\big)\Big)&\textrm{if }x\leq 0,\\[2mm]
1-(1+b)^{-1}\, b^{[x]_{\textrm{near}}}\Big(b+ \big([x]_{\textrm{near}}-x+\frac 12 \big) \big(1-b\big)\Big) &\textrm{if }x>0.
\end{cases}
\end{equation}

We summarize the procedure of the Test of Tests \citep[Algorithm 1]{kazan2023test} tailored to our permutation setting in \Cref{alg:tot}, and remark that their private test does not depend on the global sensitivity of the test statistic.

\begin{algorithm}[!t]
\caption{Permutation Test of Tests TOT \citep{kazan2023test}}
\label{alg:tot}
\begin{algorithmic}
    \State \textbf{Input:} 
    Data $\mathcal{X}_n$, 
    significance level $\alpha\in(0,1)$, 
    privacy level $\varepsilon>0$, 
    test statistic function $T(\cdot)$, 
    number of permutations $B\in\mathbb N$
    .\\
    \State \textbf{Procedure:} 
    \State Choose $\alpha_0$ (sub-test level) and $S$ (number of subsets) according to \Cref{condition tot}.
    \State  Partition $\mathcal{X}_n$ into $S$ disjoint subsets $\mathcal{X}_{n,1},\dots,\mathcal{X}_{n,S}$ (of equal sizes if possible)\;
    \For{$s \in [S]$}
    	\For{$b \in [B]$}
			\State Generate a random permutation $\bpi_{b,s}$ of $[|\mathcal{X}_{n,s}|]$. 
			\State Compute $T(\mathcal{X}_{n,s}^{\bpi_{b,s}})$. 
    	\EndFor
    	\State Compute $T(\mathcal{X}_{n,s})$.
		\State Compute $\hat{p}_s$ the permutation $p$-value for the subset $\mathcal{X}_{n,s}$ as in \eqref{Eq: permutation p-value representation}.
    \EndFor
    \State Compute the number of rejects $a = |\{s : \hat{p}_s < \alpha_0 \}|$.
    \State Privatize $a$ using Truncated-Uniform-Laplace noise \citep[Algorithm 2]{awan2018differentially}:
    \State \indent Generate random quantities $g_1,g_2\sim \textrm{Geom}(1-e^{-\varepsilon})$ and $u\sim \textrm{Unif}(-0.5,0.5)$.
    \State \indent Compute $z = a + g_1 - g_2 + u$.
    \State Compute the $p$-value \citep[Algorithm 1]{awan2018differentially}:
    \State \indent Compute $\hat{p} = \sum_{s=0}^S \binom{S}{s} \alpha_0^s (1-\alpha_0)^{S-s} F^{\textrm{Tulap}}_{e^{-\varepsilon}}(s-z)$ with $F^{\textrm{Tulap}}_{e^{-\varepsilon}}$ as in \eqref{eq: Tulap CDF}. \\
    \State \textbf{Output:} Reject $H_0$ if $\hat{p} \leq \alpha$.
\end{algorithmic}
\end{algorithm}

\begin{procedure}[{\citealt[Section 3.4]{kazan2023test}}]
\label{condition tot}
To the best of our understanding, in the nonparametric MMD and HSIC settings considered, the power of the permutation tests does not admit a closed-form expression. Hence, it is not possible to perform the optimization procedure proposed by \citet[Section 3.4]{kazan2023test} to select the optimal combination of the number of subsets $S$ and the sub-test significance threshold $\alpha_0$. 
The authors point out that \emph{`the optimization tends to favor high values of $S$, with very small subsamples and high significance thresholds $\alpha_0$ on the subtests.'}
We verified these observations empirically and have found on separate data that setting the values to $S=\sqrt{n}$ and $\alpha_0=5\alpha$ consistently provide high power, we use these heuristics in our implementation.
\end{procedure}

Note that the privatized statistic $z$ in \Cref{alg:tot} is the sum of two terms: ``$a$'' which under the null follows a $\textrm{Binomial}(S,\alpha_0)$ distribution, and a noise quantity which follows a $\textrm{Tulap}(e^{-\varepsilon})$ distribution. Hence, the $p$-value is $\mathbb{P}(B+N \geq z)$ where $B\sim\textrm{Binomial}(S,\alpha_0)$ and $N\sim\textrm{Tulap}(e^{-\varepsilon})$, and it can be estimated as in the last step of \Cref{alg:tot} \citep[Algorithm 1]{awan2018differentially}.

While TOT is guaranteed to be $(\varepsilon,0)$-DP and able to maintain a significance level $\alpha$ for any number of subsets $S$ for the data partitioning, any sub-test significance level $\alpha_0\in(0,1)$, and any sample size $n$ (which is not the case of SARRM introduced below), this comes at the cost of having to rely on heuristics for $\alpha_0$ and $S$ in the non-parametric testing setting (as explained in \Cref{condition tot}).

\subsection{SARRM: Subsampled and Aggregated Randomized Response Mechanism}
\label{subsec: sarrm}

The privatization procedure of \citet{pena2022differentially} is similar to that of \citet{kazan2023test}: both methods split the data into subsets on which the non-private test is run. \citet{kazan2023test} then runs the optimal binomial test of \citet{awan2018differentially}, while \citet{pena2022differentially} relies on a binomial test based on a randomized response mechanism.

For binary input $x\in\{0,1\}$ and $p\in[0,1]$, the randomized response mechanism is defined as
\begin{equation}
\label{eq: randomized response mechanism}
r_p(x) \coloneqq
\begin{cases}
	x &\textrm{with probability } p, \\
	1-x &\textrm{with probability } 1-p.
\end{cases}
\end{equation}

We present the Subsampled and Aggregated Randomized Response Mechanism (SARRM) Test of \citet{pena2022differentially} in \Cref{alg:sarrm}.
Again, unlike our procedure, we point out that this test does not rely on the global sensitivity of the test statistic.

\begin{algorithm}[H]
\caption{Permutation SARRM Test \citep{pena2022differentially}}
\label{alg:sarrm}
\begin{algorithmic}
    \State \textbf{Input:} 
    Data $\mathcal{X}_n$, 
    significance level $\alpha\in(0,1)$, 
    privacy level $\varepsilon>0$, 
    test statistic function $T(\cdot)$, 
    number of permutations $B\in\mathbb N$
    .\\
    \State \textbf{Procedure:} 
    \State Compute $p = e^\varepsilon / (1 + e^\varepsilon)$.
    \State Choose $\alpha_0$ (sub-test level) and $k$ (determining the number of subsets) according to \Cref{condition sarrm}.
    \State Partition $\mathcal{X}_n$ into $S=2k+1$ disjoint subsets $\mathcal{X}_{n,1},\dots,\mathcal{X}_{n,S}$ (of equal sizes if possible)\;
    \For{$s \in [S]$}
    	\For{$b \in [B]$}
			\State Generate a random permutation $\bpi_{b,s}$ of $[|\mathcal{X}_{n,s}|]$. 
			\State Compute $T(\mathcal{X}_{n,s}^{\bpi_{b,s}})$. 
    	\EndFor
    	\State Compute $T(\mathcal{X}_{n,s})$.
		\State Compute $\hat{p}_s$ the permutation $p$-value for the subset $\mathcal{X}_{n,s}$ as in \eqref{Eq: permutation p-value representation}.
		\State Compute sub-test output $t_s = \mathds{1}(\hat{p}_s\leq \alpha_0)$.

    \EndFor
    \State Compute statistic $T = \sum_{s=1}^S r_p(t_s)$ with $r_p$ as in \eqref{eq: randomized response mechanism}. \\

    \State \textbf{Output:} Reject $H_0$ if $T > k$.
\end{algorithmic}
\end{algorithm}

Recall the meaning of the test parameters of SARRM: $k$ determines the number of subsets $2k+1$, $p$ is the randomized response mechanism probability, $\alpha_0$ is the sub-test significance level, $\alpha$ is the significance level, and $\varepsilon$ is the privacy level.
\citet[Proposition 2]{pena2022differentially} show that SARRM is $(\varepsilon,0)$-DP with
\begin{equation}
\label{eq: condition dp sarrm}
\varepsilon = \log\!\big(\mathbb{P}(B_1>k)\big/\mathbb{P}(B_0>k)\big)
\end{equation}
where $B_0\sim \mathrm{Binomial}(2k+1, 1-p)$ and $B_1\sim \mathrm{Bernoulli}(p) + \mathrm{Binomial}(2k, 1-p)$.
\citet[Proposition 7]{pena2022differentially} use the fact that the test statistic $T = \sum_{s=1}^{2k+1} r_p(t_s)$ defined in \Cref{alg:sarrm} follows a $\mathrm{Binomial}(2k+1,q_{p,\alpha_0})$ distribution under the null, where $q_{p,\alpha_0}\coloneqq p\alpha_0+(1-p)(1-\alpha_0)$, from which it can be deduced that SARRM has significance level $\alpha$ provided that 
\begin{equation}
\label{eq: condition level sarrm}
\mathrm{level}(k, p, \alpha_0) \coloneqq \sum_{\ell=k+1}^{2k+1} \binom{2k+1}{\ell} q_{p,\alpha_0}^\ell (1-q_{p,\alpha_0})^{2k+1-\ell} \leq \alpha.
\end{equation}

\begin{procedure}[{\citealt[Section 3.3]{pena2022differentially}}]
\label{condition sarrm}
The aim is to determine the parameters $k\in\mathbb{N}\setminus\{0\}$, $p\in(0.5, 1)$ and $\alpha_0\geq\alpha_{0,\mathrm{min}}$ for $\alpha_{0,\mathrm{min}}=0.0025$ (user-specified), such that SARRM achieves the given privacy level $\varepsilon$ and significance level $\alpha$.
To start, set the value of $k$ to 1.
The parameter $p$ can then be chosen according to \eqref{eq: condition dp sarrm} to guarantee the $(\varepsilon,0)$-DP property.
If $\mathrm{level}(k, p, \alpha_{0,\mathrm{min}}) \leq \alpha$, then there exists some $\alpha_0\geq\alpha_{0,\mathrm{min}}$ such that $\mathrm{level}(k, p, \alpha_0) = \alpha$ which guarantees test level $\alpha$, all required conditions are then satisfied, and the values $k$, $p$ and $\alpha_0$ are used for SARRM. 
If $\mathrm{level}(k, p, \alpha_{0,\mathrm{min}}) > \alpha$, then level $\alpha$ cannot be attained, the value $k+1$ is then considered instead of $k$, and the procedure is repeated.\footnote{Following this procedure, we are able to replicate the results of \citet[Table 1]{pena2022differentially}.} The procedure is guaranteed to terminate. If an upper bound on $k$ is known (as in our experiments since the sample size is fixed), then the search of the parameter $k$ can be performed using a bisection method which runs considerably faster.
\end{procedure}

Note that, in order to achieve a given privacy level $\varepsilon$ and significance level $\alpha$, a certain number of subsets $2k+1$ (determined by \Cref{condition sarrm}) must be used. This is only possible if the sample size is greater than the number of subsets, that is, if $n \geq 2k+1$. If this is not the case, then SARRM simply cannot be run with that privacy level, significance level and sample size, which is a considerable limitation of SARRM.

{
\subsection{TOT--SARRM}
\label{subsec: TOT-SARRM}

TOT uses the sub-test level $\alpha_0$ and number of subsets $S$ as hyperparameters. The guidance of \citet{kazan2023test} does not directly apply to our nonparametric setting, so our original TOT implementation used the empirically chosen values $\alpha_0=5\alpha$ and $S=\sqrt n$. SARRM also uses $\alpha_0$ and $S=2k+1$, together with the randomized-response probability $p$, and \citet{pena2022differentially} provide the selection rule in \Cref{condition sarrm}.

To distinguish the effect of the binomial test from that of hyperparameter selection, we apply the SARRM rule to choose $\alpha_0$ and $S$ and then run TOT with those values. We call the resulting procedure TOT--SARRM. Although the selection rule was derived specifically for SARRM, TOT remains a valid level-$\alpha$ test for every $\alpha_0$ and $S$. Moreover, TOT uses the uniformly most powerful binomial test of \citet{awan2018differentially}; hence, with common $\alpha_0$ and $S$, TOT--SARRM has power at least as high as SARRM~\citep[Theorem 4.1]{kazan2023test}. This is what we observe empirically: the two procedures have identical power in most settings, indicating that their hyperparameter choices often matter more than the distinction between their binomial tests.
}

\subsection{TCS-ME: Trusted-Curator, perturbed Statistic, Mean Embedding}
\label{subsec: Raj}

\citet{raj2020differentially} propose several variants of a DP kernel two-sample test: 
\begin{itemize}
	\item considering both the trusted-curator and the no-trusted-entity privacy settings,
	\item relying on the analytic Gaussian mechanism \citep{balle2018improving} to inject privacy noise either in the statistic or in the means and covariances of feature vectors,
	\item using the two kernel statistics of \citet{jitkrittum2016interpretable} based on mean embeddings and on smooth characteristic functions,
	\item either optimizing the test locations and bandwidth, or sampling locations and using the median heuristic bandwidth,
	\item using the asymptotic $\chi^2$ null distribution, or an approximation to the null distribution.
\end{itemize}

We compare against the most powerful of these tests (as shown in \citealt[Figure 2]{raj2020differentially}), which the authors call TCS-ME. This test uses the Trusted-Curator (TC) privacy setting, consistent with the framework considered for dpMMD. It directly perturbs a test statistic (S) based on kernel mean embeddings (ME), uses 20\% of the data to optimize the test locations and kernel bandwidth for power, and approximates the null distribution.

TCS-ME is guaranteed to be $(\varepsilon, \delta)$-DP for $\delta>0$ as it relies on the (analytic) Gaussian mechanism.
In our simulations, we set $\delta=10^{-5}$ as in the experiments of \citet{raj2020differentially}, and compare against $(\varepsilon, 0)$-DP tests (a slightly more restrictive constraint).

As can be seen in \Cref{subsec:level}, while TCS-ME controls the type I error at level $\alpha$ in the one-dimensional case (which is the setting considered \Cref{subsec:mmd_uniform} where we compared against TCS-ME), we observe that it can be extremely poorly calibrated (type I error up to $10\alpha$) when working in higher dimensions.

\section{Proofs for the Main Text} \label{Section: Proofs for the Main Text}
This section collects the proofs of the results provided in the main text. 

\subsection{Proof of \Cref{Theorem: validity of private permutation tests}}
When $\mathcal{X}_n = (X_1,\ldots,X_n)$ are exchangeable, $M_0,\ldots,M_B$ are exchangeable by construction. Moreover, since $M_0,M_1,\ldots,M_B$ are distinct with probability one due to i.i.d.~Laplace noises, the second result of \Cref{Lemma: permutation p-value} proves the equality.

\subsection{Proof of \Cref{Theorem: Differential privacy of permutation tests}} \label{Section: Proof of Theorem: Differential privacy of permutation tests}
Let us denote the $1-\alpha$ quantile of $\{M_i\}_{i=0}^B$ and the $1 - \alpha_\star$ quantile of $\{M_i\}_{i=1}^B$ by $q_{1-\alpha}$ and $r_{1-\alpha_\star}$, respectively, where
\begin{align*}
	\alpha_\star = \max\bigg\{ \biggl(\frac{B+1}{B} \alpha - \frac{1}{B}\biggr), \, 0 \bigg\} \in [0,1).
\end{align*}
We then claim that 
\begin{align*}
	\mathds{1}(\hat{p}_{\mathrm{dp}} \leq \alpha) \overset{\mathrm{(i)}}{=} \mathds{1}(M_0 > q_{1-\alpha}) \overset{\mathrm{(ii)}}{=}  \mathds{1}(M_0 > r_{1-\alpha_\star}) \mathds{1} \biggl( \alpha \geq \frac{1}{B+1} \biggr),
\end{align*}
where identity~(i) holds by the quantile representation of the permutation test in \Cref{Lemma: Quantile representation},
and identity~(ii) holds by \Cref{Lemma: Alternative expression}. Moreover, \Cref{Lemma: Sensitivity of quantiles} proves that $r_{1-\alpha_\star}$ has the global sensitivity at most $\Delta_T$, and the triangle inequality yields that $T_0 -  r_{1-\alpha_\star}$ has the global sensitivity at most $2\Delta_T$ under condition~\eqref{Eq: global sensitivity}. 
Therefore the Laplace mechanism ensures that $M_0 - r_{1-\alpha_\star} = T_0 - r_{1-\alpha_\star} + 2\Delta_T \xi_{\varepsilon,\delta}^{-1} \zeta_0$ is $(\varepsilon,\delta)$-DP.
Finally, noting that the test function $\mathds{1}(\hat{p}_{\mathrm{dp}} \leq \alpha) = \mathds{1}(M_0 - r_{1-\alpha_\star}  > 0) \,\mathds{1} ( \alpha \geq 1/(B+1))$ is a function of $(\varepsilon,\delta)$-DP statistic, the post processing property of differential privacy~(\Cref{Lemma: Post-Processing}) proves that the given permutation test is differentially private at privacy levels $\varepsilon$ and $\delta$. This completes the proof of \Cref{Theorem: Differential privacy of permutation tests}.

\subsection{Proof of \Cref{Theorem: Pointwise Consistency}}
The result of \Cref{Theorem: Pointwise Consistency} follows as a corollary of \Cref{Lemma: General conditions for consistency}. In more detail, let $\mathcal{G}$ be the sigma field generated by $\{\mathcal{X}_n,\zeta_0\}$. Conditional on $\mathcal{G}$, observe that $\{M_1,\ldots,M_B\}$ are i.i.d.~where the randomness arises from $\{\bpi_i, \zeta_i\}_{i=1}^B$, and $M_0$ is constant. Therefore, the proposed private test is pointwise consistent as the conditions in \Cref{Lemma: General conditions for consistency} are satisfied.

\subsection{Proof of \Cref{Theorem: General uniform power condition}}  \label{Section: Proof Theorem: General uniform power condition}
By the quantile representation of the permutation test~(\Cref{Lemma: Quantile representation}), the type II error can be written as
\begin{align*}
	\mP_{P} (\hat{p}_{\mathrm{dp}} > \alpha) = \mP_{P} \bigl(M_0 \leq q_{1-\alpha,B} \bigr),
\end{align*}
where $q_{1-\alpha,B}$ is the $1-\alpha$ quantile of $M_0,M_1,\ldots,M_B$. Let us split the proof into two steps:
\begin{itemize}
	\item In the first step, we present an upper bound for $q_{1-\alpha,B}$ holding with high probability. In particular, the following holds
	\begin{align*}
		q_{1-\alpha,B} \leq \mE_P[T(\mathcal{X}_n^{\bpi})]  + 10 \sqrt{\frac{\mV_P[T(\mathcal{X}_n^{\bpi})]}{\alpha\beta}} + 2 \Delta_T \xi_{\varepsilon,\delta}^{-1}  F_{\zeta}^{-1}(1-\alpha/4)
	\end{align*}
	with high probability (say $1- 3\beta/4$) under the condition that $B \geq 2\alpha^{-1} \max \{3 \log(2\beta^{-1}), 1- \alpha\}$. In fact, we can improve the condition for $B$ at the cost of having a worse dependence on $\alpha$ and $\beta$. Specifically, we prove at the end of this subsection that:  
	\begin{align} \label{Eq: another bound for the quantile}
		q_{1-\alpha,B} \leq  \mE_{P,\boldsymbol{\pi}}[T(\mathcal{X}_n^{\boldsymbol{\pi}})] +8\sqrt{\frac{\mV_{P,\boldsymbol{\pi}}[T(\mathcal{X}_n^{\boldsymbol{\pi}})]}{\gamma \alpha \beta}} + \frac{16}{\sqrt{\gamma \alpha \beta}} \Delta_T \xi_{\varepsilon,\delta}^{-1} 
	\end{align}
	holds with probability at least $1-\beta/4$, provided that $B \geq \frac{1-\alpha}{\alpha(1-\gamma)}$ for any $\gamma \in (0,1)$.
	\item In the second step, we use the above quantile bound to show that
 control of the type II error is guaranteed under the condition in \eqref{Eq: Uniform
Power Condition}.
\end{itemize} 
For simplicity, we omit the subscript $P$ in the expectation, variance and probability operators throughout the proof.

\vskip 1em 

\noindent \textbf{Step 1 (Bounding the quantile).} We first focus on the quantile  $q_{1-\alpha,B}$ and present a high probability upper bound. To this end, we apply \Cref{Lemma: Quantile Approximation} by setting $\mathcal{G}$ therein to be the sigma field generated by $\mathcal{X}_n$. In particular, under the condition on $B$ and by letting $q_{1-\alpha/6,\infty}$ be the $1-\alpha/6$ quantile of the conditional distribution of $T(\mathcal{X}_n^{\bpi})+ 2 \Delta_T \xi_{\varepsilon,\delta}^{-1} \zeta$ given $\mathcal{X}_n$ (or equivalently the limiting value of $q_{1-\alpha/6,B}$ with $B = \infty$), \Cref{Lemma: Quantile Approximation} along with \Cref{Remark: quantile approximation} guarantees that the following event
\begin{align*}
	E_1 := \big\{q_{1-\alpha,B} \leq q_{1-\alpha/6,\infty} \big\}
\end{align*}
holds with probability at least $1-\beta/2$. We further define $q_{1-\alpha/12,\infty}^a$ and $q_{1-\alpha/12,\infty}^b$ as
\begin{align*}
	& q_{1-\alpha/12,\infty}^a \coloneqq  \inf \Biggl\{ x \in \mathbb{R} : \frac{1}{|\boldsymbol{\Pi}_{n}|} \sum_{\bpi \in \boldsymbol{\Pi}_{n}} \mathds{1}\bigl( T(\mathcal{X}_n^{\bpi}) \leq x \bigr) \geq 1 - \alpha/12 \Biggr\} \quad \text{and} \\[.5em]
	& q_{1-\alpha/12,\infty}^b \coloneqq  2 \Delta_T \xi_{\varepsilon,\delta}^{-1} F_{\zeta}^{-1}(1 - \alpha/12),
\end{align*}
where $|\boldsymbol{\Pi}_{n}|$ denotes the cardinality of the set $\boldsymbol{\Pi}_{n}$. Then \Cref{Lemma: Quantile inequality} gives the inequality that $q_{1-\alpha/6,\infty} \leq q_{1-\alpha/12,\infty}^a + q_{1-\alpha/12,\infty}^b$.

Next we further upper bound $q_{1-\alpha/12,\infty}^a$ and $q_{1-\alpha/12,\infty}^b$ by more manageable terms. To begin with $q_{1-\alpha/12,\infty}^a$, Chebyshev's inequality conditional on $\mathcal{X}_n$ yields
\begin{align*}
	\mP_{\bpi} \bigl\{ T(\mathcal{X}_n^{\bpi})   \geq  \mE_{\bpi}[T(\mathcal{X}_n^{\bpi}) \given \mathcal{X}_n] + t \given \mathcal{X}_n \bigr\} \leq \frac{\mV_{\bpi}[T(\mathcal{X}_n^{\bpi}) \given \mathcal{X}_n]}{t^2} \quad \text{for any $t > 0$.}
\end{align*}
Consequently, $q_{1-\alpha/12,\infty}^a$ is bounded as
\begin{align*}
	q_{1-\alpha/12,\infty}^a \leq \sqrt{\frac{12\mV_{\bpi}[T(\mathcal{X}_n^{\bpi}) \given \mathcal{X}_n]}{\alpha}} + \mE_{\bpi}[T(\mathcal{X}_n^{\bpi}) \given \mathcal{X}_n].
\end{align*}
Now define events $E_2$ and $E_3$,
\begin{align*}
	& E_2 := \Big\{ \mV_{\bpi}[T(\mathcal{X}_n^{\bpi}) \given \mathcal{X}_n] < 8\beta^{-1} \mE\bigl( \mV_{\bpi}[T(\mathcal{X}_n^{\bpi}) \given \mathcal{X}_n] \bigr) \Big\}, \\[.5em]
	& E_3 := \bigg\{ \big| \mE_{\bpi}[T(\mathcal{X}_n^{\bpi}) \given \mathcal{X}_n] - \mE[T(\mathcal{X}_n^{\bpi})] \big|  < \sqrt{8\beta^{-1} \mV\bigl(\mE_{\bpi}[T(\mathcal{X}_n^{\bpi}) \given \mathcal{X}_n]\bigr) } \bigg\},
\end{align*}
where each of the events holds with probability at least $1-\beta/8$ by Markov's inequality and Chebyshev's inequality, respectively. 
We emphasize that the expectation $\mE[T(\mathcal{X}_n^{\bpi})]$ is taken with respect to both $\mathcal{X}_n$ and $\bpi$.
Under these events $E_2$ and $E_3$, it can be seen that 
\begin{align*}
	q_{1-\alpha/12,\infty}^a ~\leq~ & \mE[T(\mathcal{X}_n^{\bpi})]  + \sqrt{\frac{96\mE\bigl(\mV_{\bpi}[T(\mathcal{X}_n^{\bpi}) \given \mathcal{X}_n]\bigr)}{\alpha\beta}} + \sqrt{\frac{ 8\mV\bigl(\mE_{\bpi}[T(\mathcal{X}_n^{\bpi}) \given \mathcal{X}_n]\bigr)}{\beta}} \\[.5em]
	\leq ~ & \mE[T(\mathcal{X}_n^{\bpi})]  + 10 \sqrt{\frac{\mV[T(\mathcal{X}_n^{\bpi})]}{\alpha\beta}}
\end{align*}
where the second inequality uses the law of total variance. In summary, we have shown that with probability at least $1- 3\beta/4$, 
\begin{align*}
	q_{1-\alpha,B} \leq \mE[T(\mathcal{X}_n^{\bpi})]  + 10 \sqrt{\frac{\mV[T(\mathcal{X}_n^{\bpi})]}{\alpha\beta}} + 2 \Delta_T \xi_{\varepsilon,\delta}^{-1}  F_{\zeta}^{-1}(1-\alpha/12).
\end{align*}

\vskip 1em

\noindent \textbf{Step 2 (Bounding the type II error).} Given the upper bound for $q_{1-\alpha,B}$ from the previous step, we now examine the type II error of the differentially private permutation test and show that it is bounded by $\beta$ under the given condition. In particular, writing 
\begin{align*}
	& R_1 := 10 \sqrt{\frac{\mV[T(\mathcal{X}_n^{\bpi})]}{\alpha\beta}} + 2 \Delta_T \xi_{\varepsilon,\delta}^{-1}  F_{\zeta}^{-1}(1-\alpha/12), \\[.5em]
	& R_2 := 2 \Delta_T \xi_{\varepsilon,\delta}^{-1}  F_{\zeta}^{-1}(1-\beta/8),
\end{align*}
the result from Step 1 yields
\begin{align} \nonumber
	\mP(M_0 \leq q_{1-\alpha,B})  ~\leq~ & \mP\bigl(M_0 \leq \mE[T(\mathcal{X}_n^{\bpi})] + R_1 \bigr)  + \frac{3\beta}{4} \\[.5em] \nonumber
	= ~ &  \mP\bigl(T_0 + 2 \Delta_T \xi_{\varepsilon,\delta}^{-1} \zeta_0 \leq \mE[T(\mathcal{X}_n^{\bpi})] + R_1 \bigr) + \frac{3\beta}{4} \\[.5em] \label{Eq: last step of power bound}
	\leq  ~ &  \mP\bigl(T_0 \leq \mE[T(\mathcal{X}_n^{\bpi})] + R_1 + R_2 \bigr)  + \frac{7\beta}{8},
\end{align}
where the last inequality holds since
\begin{align*}
	\mP\bigl( -\xi_{\varepsilon,\delta}^{-1} \zeta_0 > F_{\zeta}^{-1}(1-\beta/8) \bigr) \leq \frac{\beta}{8}.
\end{align*}
We also remark that for any $\alpha \in (0,1)$ and $\beta \in (0,1-\alpha)$ (implying that $\min\{\alpha,\beta\} < 1/2$), 
\begin{align} \label{Eq: quantile bound for Laplace}
	F_{\zeta}^{-1}(1-\alpha/12) +  F_{\zeta}^{-1}(1-\beta/8) \leq 7 \max \{\log(1/\alpha), \log(1/\beta)\}.
\end{align}
In more detail, the quantile function of $\zeta$ is given as 
\begin{align*}
	F_{\zeta}^{-1}(t) = 
	\begin{cases}
		\log(2t) < 0, \quad & \text{if $t < 0.5$}, \\
		-\log\bigl(2(1-t)\bigr) \geq 0, \quad & \text{if $t \geq 0.5$.}
	\end{cases}
\end{align*}
Based on this expression along with the fact that $\min\{\alpha,\beta\} < 1/2$, we have
\begin{align*}
		F_{\zeta}^{-1}(1-\alpha/12) +  F_{\zeta}^{-1}(1-\beta/8) ~ = ~ & \log(6/\alpha) + \log(4/\beta) = \log(24) + \log(1/\alpha) + \log(1/\beta) \\[.5em]
		\leq ~ &  7 \max \{\log(1/\alpha), \log(1/\beta)\}.
\end{align*}
Hence, under the given condition of \eqref{Eq: Uniform Power Condition}, we can ensure that 
\begin{align*}
	\mE[T(\mathcal{X}_n^{\bpi})] + R_1 + R_2  \leq \mE[T(\mathcal{X}_n)]  - \sqrt{8\beta^{-1}\mV[T(\mathcal{X}_n)]}.
\end{align*}
For clarity, we emphasize the notational difference between $\mE[T(\mathcal{X}_n)] = \mE_{P}[T(\mathcal{X}_n)]$ and  $\mE[T(\mathcal{X}_n^{\bpi})] = \mE_{P,\bpi}[T(\mathcal{X}_n^{\bpi})]$ where the former expectation is taken over $\mathcal{X}_n$, whereas the latter expectation is taken over both $\mathcal{X}_n$ and $\bpi$. This bound together with Chebyshev's inequality yields
\begin{align*}
	\mP\bigl(T_0 \leq \mE[T(\mathcal{X}_n^{\bpi})] + R_1 + R_2 \bigr) ~\leq~ & \mP\bigl(T_0 \leq \mE[T(\mathcal{X}_n)]  - \sqrt{8\beta^{-1}\mV[T(\mathcal{X}_n)]} \bigr) \\[.5em]
	= ~ &  \mP\bigl( \sqrt{8\beta^{-1}\mV[T(\mathcal{X}_n)]}  \leq  \mE[T(\mathcal{X}_n)] - T_0 \bigr) \\[.5em]
	\leq ~ & \frac{\beta}{8}.
\end{align*}
Now combining the inequality~\eqref{Eq: last step of power bound} with the above yields that 
\begin{align*}
	\mP(M_0 \leq q_{1-\alpha,B}) \leq  \beta.
\end{align*}
This bound holds for any $P \in \mathcal{P}_1$ and the upper bound is independent of $P$. Hence the uniform guarantee stated in \Cref{Theorem: General uniform power condition} holds.

\vskip 1em

\noindent \textbf{Proof of \eqref{Eq: another bound for the quantile}.} We now prove the inequality~\eqref{Eq: another bound for the quantile}, which holds under a more relaxed condition for $B$. Observe that $q_{1-\alpha,B}$ is the $1-\alpha$ quantile of $\{M_0,M_1,\ldots,M_B\}$ satisfying
\begin{align*}
	q_{1-\alpha,B} ~=~& \inf\bigg\{ x \in \mathbb{R} : \frac{1}{B+1} \sum_{i=0}^B \mathds{1}(M_i \leq x) \geq 1 - \alpha  \bigg\} \\[.5em]
	\leq ~& \inf\bigg\{ x \in \mathbb{R} : \frac{1}{B} \sum_{i=1}^B \mathds{1}(M_i \leq x) \geq \frac{B+1}{B} (1-\alpha) \bigg\} \\[.5em]
	\leq ~ & \inf\bigg\{ x \in \mathbb{R} : \frac{1}{B} \sum_{i=1}^B \mathds{1}(M_i \leq x) \geq 1 - \gamma \alpha \bigg\}:= \tilde{q}_{1- \gamma \alpha, B},
\end{align*}
where the last inequality holds under the condition that $B \geq \frac{1-\alpha}{\alpha(1-\gamma)}$ for $\gamma \in (0,1)$. Letting 
\begin{align*}
	\hat{\mu}_B = \frac{1}{B} \sum_{i=1}^B M_i \quad \text{and} \quad \hat{\sigma}^2_B = \frac{1}{B} \sum_{i=1}^B (M_i - \hat{\mu}_B)^2,
\end{align*} 
Chebyshev's inequality yields
\begin{align*}
	\frac{1}{B} \sum_{i=1}^B \mathds{1}(M_i - \hat{\mu}_B \geq t) \leq \frac{1}{t^2} \hat{\sigma}^2_B \overset{\mathrm{set}}{=} \gamma \alpha,
\end{align*}
which implies that the quantile $\tilde{q}_{1- \gamma \alpha, B}$ is bounded above by
\begin{align} \label{Eq: preliminary inequallity}
	\tilde{q}_{1- \gamma \alpha, B} \leq \hat{\mu}_B + \sqrt{\frac{\hat{\sigma}^2_B}{\gamma \alpha}}.
\end{align}
Our strategy is to identify a high-probability upper bound for $\hat{\mu}_B + \sqrt{\hat{\sigma}^2_B / (\gamma \alpha)}$, which is again an upper bound for $q_{1-\alpha,B}$ based on the previous observation. To this end, note that $M_1,\ldots,M_B$ are i.i.d.~random variables conditional on $\mathcal{X}_n$. Hence, it holds by the law of total expectation that  
\begin{align*}
	\mE[(\hat{\mu}_B -  \mE[M_1])^2] ~=~&  \mE[(\hat{\mu}_B - \mE[M_1 \given \mathcal{X}_n] + \mE[M_1 \given \mathcal{X}_n] - \mE[M_1])^2] \\[.5em]
	=~&  \mE[(\hat{\mu}_B - \mE[M_1 \given \mathcal{X}_n])^2] + \mE[(\mE[M_1 \given \mathcal{X}_n] - \mE[M_1])^2] \\[.5em]
	=~& \frac{1}{B} \mE [\mV(M_1 \given \mathcal{X}_n)]  + \mV[ \mE[M_1 \given \mathcal{X}_n)] \\[.5em]
	\leq~ & \mV[M_1]. 
\end{align*}
This together with Chebyshev's inequality yields that 
\begin{align*}
	\mP\bigl( |\hat{\mu}_B -  \mE[M_1] | \geq t  \bigr) \leq \frac{1}{t^2} \mV[M_1] \overset{\mathrm{set}}{=} \frac{\beta}{8}.
\end{align*}
Moreover, by Markov's inequality, we have
\begin{align*}
	\mP(\hat{\sigma}_B \geq t) \leq \frac{1}{t^2} \mE[\hat{\sigma}^2_B] \leq \frac{1}{t^2} \mE[\mV(M_1 \given \mathcal{X}_n)] \overset{\mathrm{set}}{=} \frac{\beta}{8}.
\end{align*} 
Combining the above two inequalities along with \eqref{Eq: preliminary inequallity} yields that 
\begin{align*}
	q_{1-\alpha,B} \leq \mE[M_1] + \sqrt{\frac{8}{\beta}\mV[M_1]} + \sqrt{\frac{8}{\gamma \alpha \beta} \mE[\mV(M_1 \given \mathcal{X}_n)]}
\end{align*}
holding with probability at least $1-\beta/4$. Since $\mE[\mV(M_1 \given \mathcal{X}_n)] \leq \mV[M_1]$ and $M_1$ has the same distribution as $T(\mathcal{X}_n^{\bpi}) + 2\Delta_T \xi_{\varepsilon,\delta}^{-1} \zeta$, we further have 
\begin{align*}
	q_{1-\alpha,B} \leq  \mE_{P,\boldsymbol{\pi}}[T(\mathcal{X}_n^{\boldsymbol{\pi}})] +8\sqrt{\frac{\mV_{P,\boldsymbol{\pi}}[T(\mathcal{X}_n^{\boldsymbol{\pi}})]}{\gamma \alpha \beta}} + \frac{16}{\sqrt{\gamma \alpha \beta}} \frac{\Delta_T}{\xi_{\varepsilon,\delta}}
\end{align*}
with probability at least $1-\beta/4$, provided that $B \geq \frac{1-\alpha}{\alpha(1-\gamma)}$ and $\gamma \in (0,1)$. This proves the inequality~\eqref{Eq: another bound for the quantile}.
 
\subsection{Proof of \Cref{Lemma: Sensitivity of MMD}} \label{Section: Proof of Lemma: Sensitivity of MMD}

We start with the proof of the upper bound, followed by the proof of the lower bound.

\subsubsection{Upper Bound}
To simplify the proof, we may exploit the permutation invariance property of $\widehat{\mathrm{MMD}}(\mathcal{X}_{n+m}^{\bpi})$ within $\{X_{\pi_1} ,\ldots, X_{\pi_n}\}$ and $\{Y_{\pi_{n+1}},\ldots,Y_{\pi_{n+m}}\}$, and focus on two specific cases of neighboring dataset: 
\begin{align*}
	\tilde{\mathcal{X}}_{n+m}^{\bpi,a}=\{X_{\pi_1}',\ldots, X_{\pi_n},Y_{\pi_{n+1}},\ldots,Y_{\pi_{n+m}}\} \quad  \text{and} \quad \tilde{\mathcal{X}}_{n+m}^{\bpi,b}=\{X_{\pi_1},\ldots, X_{\pi_n},Y_{\pi_{n+1}}',\ldots,Y_{\pi_{n+m}}\},
\end{align*}
where $X_{\pi_1}'$ is an i.i.d.~copy of $X_{\pi_1}$ and $Y_{\pi_{n+1}}'$ is an i.i.d.~copy of $Y_{\pi_{n+1}}$. The bound for the other neighboring datasets can be proven by following the same lines of the proof. 

Starting with the neighboring dataset $\tilde{\mathcal{X}}_{n+m}^{\bpi,a}$, note that 
\begin{align*}
	& \widehat{\mathrm{MMD}}(\mathcal{X}_{n+m}^{\bpi})  \\[.5em]
	=~  & \sup_{f \in \mathcal{F}_k} \biggl\{ \frac{1}{n} \sum_{i=1}^n f(X_{\pi_i}) - \frac{1}{m} \sum_{j=1}^m f(Y_{\pi_{n+j}}) \biggr\} \\[.5em]
	=~ & \sup_{f \in \mathcal{F}_k} \biggl\{ \frac{1}{n} \sum_{i=1}^n f(X_{\pi_i}) - \frac{1}{m} \sum_{j=1}^m f(Y_{\pi_{n+j}}) + \frac{1}{n}f(X_{\pi_1}') - \frac{1}{n}f(X_{\pi_1}') \biggr\} \\[.5em]
	=~ &  \sup_{f \in \mathcal{F}_k} \biggl\{ \frac{1}{n} \Bigl(f(X_{\pi_1}') + \sum_{i=2}^n f(X_{\pi_i}) \Bigr) - \frac{1}{m} \sum_{j=1}^m f(Y_{\pi_{n+j}}) + \frac{1}{n}f(X_{\pi_1}) - \frac{1}{n}f(X_{\pi_1}') \biggr\} \\[.5em]
	\leq ~ &  \sup_{f \in \mathcal{F}_k} \biggl\{ \frac{1}{n} \Bigl(f(X_{\pi_1}') + \sum_{i=2}^n f(X_{\pi_i})\Bigr) - \frac{1}{m} \sum_{j=1}^m f(Y_{\pi_{n+j}})\biggr\} + \sup_{f \in \mathcal{F}_k} \biggl\{ \frac{1}{n}f(X_{\pi_1}) - \frac{1}{n}f(X_{\pi_1}') \biggr\} \\[.5em]
	= ~ &  \widehat{\mathrm{MMD}}(\mathcal{X}_{n+m}^{\bpi,a}) + \sup_{f \in \mathcal{F}_k} \biggl\{ \frac{1}{n}f(X_{\pi_1}) - \frac{1}{n}f(X_{\pi_1}') \biggr\},
\end{align*}
where the inequality uses the triangle inequality. Since $|f(x)-f(y)| = |\langle f, k(x,\cdot) - k(y,\cdot)\rangle_{\mathcal{H}_k}| \leq \|f\|_{\mathcal{H}_k} \sqrt{k(x,x) + k(y,y) - 2 k(x,y)}$ by the Cauchy--Schwarz inequality as well as the reproducing kernel property, we obtain
\begin{align*}
	\sup_{f \in \mathcal{F}_k} \biggl\{ \frac{1}{n}f(X_{\pi_1}) - \frac{1}{n}f(X_{\pi_1}') \biggr\} \leq \frac{1}{n} \sup_{f \in \mathcal{F}_k} |f(X_{\pi_1}) - f(X_{\pi_1}')| \leq \frac{\sqrt{2K}}{n}, 
\end{align*}
where the second inequality uses the fact that $\|f\|_{\mathcal{H}_k} \leq 1$ and $0 \leq k(x,y) \leq K$ for all $x,y \in \mathbb{S}$. We therefore observe $\widehat{\mathrm{MMD}}(\mathcal{X}_{n+m}^{\bpi})- \widehat{\mathrm{MMD}}(\mathcal{X}_{n+m}^{\bpi,a}) \leq \sqrt{2K}/n$. Since the same argument holds with the roles of $\mathcal{X}_{n+m}^{\bpi}$ and $\mathcal{X}_{n+m}^{\bpi,a}$ reversed, we conclude that
\begin{align*}
	\bigl| \widehat{\mathrm{MMD}}(\mathcal{X}_{n+m}^{\bpi}) - \widehat{\mathrm{MMD}}(\mathcal{X}_{n+m}^{\bpi,a}) \bigr| \leq \frac{\sqrt{2K}}{n}.
\end{align*} 
Next, for the dataset $\tilde{\mathcal{X}}_{n+m}^{\bpi,b}$, we can use essentially the same argument to show that 
\begin{align*}
	\bigl| \widehat{\mathrm{MMD}}(\mathcal{X}_{n+m}^{\bpi}) - \widehat{\mathrm{MMD}}(\mathcal{X}_{n+m}^{\bpi,b}) \bigr| \leq \frac{\sqrt{2K}}{m}.
\end{align*}
Since $n \leq m$, the sensitivity is bounded by $\sqrt{2K}/n$ for any neighboring datasets and the upper bound is independent of $\bpi$. Therefore it holds that 
\begin{align*}
	\sup_{\bpi \in \boldsymbol{\Pi}_{n+m}} \sup_{\substack{\mathcal{X}_{n+m},\tilde{\mathcal{X}}_{n+m}:\\d_{\mathrm{ham}}(\mathcal{X}_{n+m},\tilde{\mathcal{X}}_{n+m}) \leq 1}} \bigl| \widehat{\mathrm{MMD}}(\mathcal{X}_{n+m}^{\bpi}) - \widehat{\mathrm{MMD}}(\tilde{\mathcal{X}}_{n+m}^{\bpi}) \bigr| \leq \frac{\sqrt{2K}}{n}, 
\end{align*}
which in turn proves the first claim in \Cref{Lemma: Sensitivity of MMD}.

\subsubsection{Lower Bound}
We now further assume that the kernel $k$ is translation invariant, and has non-empty level sets in $\mathbb{S}$. In this case, we show that the inequality becomes an equality by considering appropriate sets $\mathcal{X}_{n+m}$ and $\tilde{\mathcal{X}}_{n+m}$. More specifically, for a translation invariant kernel $k$ and a constant $\epsilon >0$, consider values $x_a$ and $x_b$ in $\mathbb{S}$ such that $k(x_a,x_b) = \epsilon_\star \leq \epsilon$, which is guaranteed by our assumption that $k$ has non-empty level sets. Let $Y_1=\ldots = Y_n = Z_1 = \ldots = Z_m = x_a$. In this scenario, the empirical MMD becomes zero for any permutation $\bpi \in \boldsymbol{\Pi}_{n+m}$. Let $Y_1' = x_b$ so that $k(Y_1,Y_1') = \epsilon_\star$. Then based on the closed form expression of the empirical MMD~\eqref{Eq: closed form MMD}, we have 
\begin{align*}
    \widehat{\mathrm{MMD}}^2(\tilde{\mathcal{X}}_{n+m}^{\bpi}) ~=~&  \frac{(n-1)^2 + 1}{n^2}K + \frac{2(n-1)}{n^2} \epsilon_\star + K  - \frac{2m(n-1)}{nm} K - \frac{2m}{nm}\epsilon_\star \\[.5em]
    =~ & \frac{2K}{n^2} - \frac{2\epsilon_\star}{n^2}.
\end{align*}
Therefore 
\begin{align*}
    \sup_{\bpi \in \boldsymbol{\Pi}_{n+m}} \sup_{\substack{\mathcal{X}_{n+m},\tilde{\mathcal{X}}_{n+m}:\\d_{\mathrm{ham}}(\mathcal{X}_{n+m},\tilde{\mathcal{X}}_{n+m}) \leq 1}} \bigl| \widehat{\mathrm{MMD}}(\mathcal{X}_{n+m}^{\bpi}) - \widehat{\mathrm{MMD}}(\tilde{\mathcal{X}}_{n+m}^{\bpi}) \bigr| \geq \frac{\sqrt{2(K-\epsilon_\star)}}{n}.
\end{align*}
Since $\epsilon$ is an arbitrary positive number and $\epsilon_\star \leq \epsilon$, we may take $\epsilon \rightarrow 0$ and thus the lower bound becomes $\sqrt{2K}/n$, indicating that the upper bound $\sqrt{2K}/n$ is tight under the given conditions. This proves the second claim in \Cref{Lemma: Sensitivity of MMD}.

\subsection{Proof of \Cref{Theorem: Properties of dpMMD}} \label{Section: Proof of Theorem: Properties of dpMMD}
Let us prove the three claims made in \Cref{Theorem: Properties of dpMMD} below.
\begin{itemize}
	\item \emph{Proof of Differential Privacy.} This result follows by \Cref{Theorem: Differential privacy of permutation tests} along with \Cref{Lemma: Sensitivity of MMD}. 
	\item \emph{Proof of Validity.} This result follows by \Cref{Theorem: validity of private permutation tests}.
	\item\emph{Proof of Consistency.} Having \Cref{Lemma: General conditions for consistency} in place,  the only condition we need to verify for consistency is $\lim_{n \rightarrow \infty} \mP(W_{0,n} \leq W_{1,n}) = 0$. Let $\mathcal{G}$ be the sigma field generated by $\{\zeta_0,\mathcal{X}_{n+m}\}$, and let 
	\begin{align*}
		W_{i,n} = M_i =  \widehat{\mathrm{MMD}}(\mathcal{X}_{n+m}^{\bpi_i}) + \frac{2\sqrt{2K}}{n\xi_{\varepsilon,\delta}}\zeta_i, \quad \text{for each $i \in \{0,1\}$,}
	\end{align*}
	where we recall $\xi_{\varepsilon,\delta} = \varepsilon+\log(1/(1-\delta))$.
	Our goal is to show that $\lim_{n \rightarrow \infty} \mP(M_0 \leq M_1) = 0$ under the conditions. Since the pair of distributions $(P,Q)$ and kernel $k$ do not vary with $n$, we consider the kernel bound $K$ to be a fixed constant throughout. Since $K$ is a constant and $(n\xi_{\varepsilon,\delta})^{-1} \rightarrow 0$, the scale of the Laplace noise $\frac{2\sqrt{2K}}{n\xi_{\varepsilon,\delta}}$ goes to zero as well. Thus Slutsky's theorem yields $\frac{2\sqrt{2K}}{n\xi_{\varepsilon,\delta}} \times \zeta_i = o_P(1)$ for each $i \in \{0,1\}$. Next, for $i=0$, \citet[][Theorem 7]{gretton2012kernel} indicate that the unpermuted MMD statistic satisfies $ \widehat{\mathrm{MMD}}(\mathcal{X}_{n+m}^{\bpi_0}) \convP \mathrm{MMD}_k(P,Q)$ where the symbol $\convP$ denotes convergence in probability. Hence another application of Slutsky's theorem yields $M_0 \convP \mathrm{MMD}_k(P,Q)$. On the other hand, we have $\widehat{\mathrm{MMD}}(\mathcal{X}_{n+m}^{\bpi_1}) \convP 0$ by \Cref{Lemma: Markov for permuted MMD}, and again Slutsky's theorem gives $M_1  \convP 0$. Combining these two results, we further have $M_0 - M_1 \convP \mathrm{MMD}_k(P,Q)$. To complete the proof, note that $\mathrm{MMD}_k(P,Q)\eqqcolon  \epsilon >0$ is assumed to be a fixed positive constant and therefore
	\begin{align*}
		\lim_{n \rightarrow \infty} \mP(M_0 \leq M_1) \leq \lim_{n \rightarrow \infty} \mP(|M_0 - M_1| \geq \epsilon) = 0,
	\end{align*}
	by the definition of convergence in probability. This completes the proof of \Cref{Theorem: Properties of dpMMD}.
\end{itemize}

\subsection{Proof of \Cref{Lemma: Sensitivity of HSIC}} \label{Section: Proof of Lemma: Sensitivity of HSIC}
We prove the upper bound result and the lower bound result in order. 
\subsubsection{Upper Bound}
We begin by considering the upper bound result of \Cref{Lemma: Sensitivity of HSIC}. Let us consider the dataset $\mathcal{X}_n = \{(Y_1,Z_1),\ldots,(Y_n,Z_n)\}$ and its neighboring dataset denoted by
\begin{align*}
	\tilde{\mathcal{X}}_n = \{(Y_1',Z_1'),(Y_2,Z_2),\ldots,(Y_n,Z_n)\} \coloneqq  \{(\tilde{Y}_1,\tilde{Z}_1),(\tilde{Y}_2,\tilde{Z}_2),\ldots,(\tilde{Y}_n,\tilde{Z}_n)\}
\end{align*}
where $(Y_1',Z_1')$ is an i.i.d.~copy of $(X_1,Y_1)$ independent of everything else. That is, $\mathcal{X}_n$ and $\tilde{\mathcal{X}}_n$ differ in their first component. We only prove the result of \Cref{Lemma: Sensitivity of HSIC} focusing on $\tilde{\mathcal{X}}_n$, and the proof for the other neighboring datasets can be derived similarly due to the symmetric structure of the empirical HSIC.

For a given permutation $\bpi$, denote by $\tilde{\mathcal{X}}_{n}^{\bpi}$ the neighboring dataset $\tilde{\mathcal{X}}_{n}$ whose $Z$ values are permuted based on $\bpi$. Let us divide the cases into two: (i) $\pi_1 = 1$ and (ii) $\pi_1 \neq 1$, and provide the proofs separately. For the first case where $\pi_1 = 1$, we write $k(y,\cdot) = \psi_Y(y)$ and $\ell(z,\cdot) = \psi_Z(z)$. We also write the sample mean of $\psi_Y(Y_1),\ldots,\psi_Y(Y_n)$ as $\overline{\psi}_Y$ and the sample mean of $\psi_Z(Z_1),\ldots,\psi_Z(Z_n)$ as $\overline{\psi}_Z$. Similarly write the sample mean of $\psi_Y(\tilde{Y}_1),\ldots,\psi_Y(\tilde{Y}_n)$ as $\tilde{\psi}_Y$ and the sample mean of $\psi_Z(\tilde{Z}_1),\ldots,\psi_Z(\tilde{Z}_{n})$ as $\tilde{\psi}_Z$. Then by adding and subtracting the same terms, we can connect $ \widehat{\mathrm{HSIC}}(\mathcal{X}_{n}^{\bpi})$ with $\widehat{\mathrm{HSIC}}(\tilde{\mathcal{X}}_{n}^{\bpi}) $ as follows:
\begin{align*}
	& \widehat{\mathrm{HSIC}}(\mathcal{X}_{n}^{\bpi}) \\[.5em] 
	= ~ &  \sup_{f \in \mathcal{F}_{k \otimes \ell}} \Biggl\{ \frac{1}{n} \sum_{i=1}^n f(Y_i,Z_{\pi_i}) - \frac{1}{n^2} \sum_{i=1}^n \sum_{j=1}^n f(Y_i,Z_{\pi_j}) \Biggr\} \\[.5em]
	= ~ & \sup_{f \in \mathcal{F}_{k \otimes \ell}}  \Bigg\langle f, \, \frac{1}{n} \sum_{i=1}^n \big\{ \psi_Y(Y_i) - \overline{\psi}_Y \big\}\big\{ \psi_Z(Z_{\pi_i}) - \overline{\psi}_Z \big\} \Bigg\rangle_{\mathcal{H}_{k \otimes \ell}} \\[.5em]
	= ~ & \sup_{f \in \mathcal{F}_{k \otimes \ell}}  \Bigg\langle f, \, \frac{1}{n} \sum_{i=1}^n \big\{ \psi_Y(Y_i) - \psi_Y(\tilde{Y}_i) + \psi_Y(\tilde{Y}_i) - \overline{\psi}_Y + \tilde{\psi}_Y - \tilde{\psi}_Y \big\}\big\{ \psi_Z(Z_{\pi_i}) - \overline{\psi}_Z \big\} \Bigg\rangle_{\mathcal{H}_{k \otimes \ell}} \\[.5em]
	\overset{(\mathrm{i})_a}{=}  ~ & \sup_{f \in \mathcal{F}_{k \otimes \ell}}  \Bigg\langle f, \, \frac{1}{n} \sum_{i=1}^n \big\{ \psi_Y(Y_i) - \psi_Y(\tilde{Y}_i) + \psi_Y(\tilde{Y}_i) - \tilde{\psi}_Y \big\}\big\{ \psi_Z(Z_{\pi_i}) - \overline{\psi}_Z \big\} \Bigg\rangle_{\mathcal{H}_{k \otimes \ell}} \\[.5em] 
	\leq~ & \sup_{f \in \mathcal{F}_{k \otimes \ell}}  \Bigg\langle f, \, \frac{1}{n} \big\{ \psi_Y(Y_1) - \psi_Y(\tilde{Y}_1) \big\}\big\{ \psi_Z(Z_{\pi_1}) - \overline{\psi}_Z \big\} \Bigg\rangle_{\mathcal{H}_{k \otimes \ell}}  \\[.5em]
	& \hskip 5em + \sup_{f \in \mathcal{F}_{k \otimes \ell}}  \Bigg\langle f, \, \frac{1}{n} \sum_{i=1}^n \big\{ \psi_Y(\tilde{Y}_i) - \tilde{\psi}_Y \big\}\big\{ \psi_Z(Z_{\pi_i}) - \overline{\psi}_Z \big\} \Bigg\rangle_{\mathcal{H}_{k \otimes \ell}} \\[.5em] 
	= ~ & \sup_{f \in \mathcal{F}_{k \otimes \ell}}  \Bigg\langle f, \, \frac{1}{n} \big\{ \psi_Y(Y_1) - \psi_Y(\tilde{Y}_1) \big\}\big\{ \psi_Z(Z_{\pi_1}) - \overline{\psi}_Z \big\} \Bigg\rangle_{\mathcal{H}_{k \otimes \ell}}  \\[.5em]
	& + \sup_{f \in \mathcal{F}_{k \otimes \ell}}  \Bigg\langle f, \, \frac{1}{n} \sum_{i=1}^n \big\{ \psi_Y(\tilde{Y}_i) - \tilde{\psi}_Y \big\}\big\{ \psi_Z(Z_{\pi_i})- \psi_Z(\tilde{Z}_{\pi_i}) + \psi_Z(\tilde{Z}_{\pi_i}) - \tilde{\psi}_Z + \tilde{\psi}_Z- \overline{\psi}_Z \big\} \Bigg\rangle_{\mathcal{H}_{k \otimes \ell}} \\[.5em] 
	\overset{(\mathrm{ii})_a}{\leq} ~ &  \widehat{\mathrm{HSIC}}(\tilde{\mathcal{X}}_{n}^{\bpi}) + \sup_{f \in \mathcal{F}_{k \otimes \ell}}  \Bigg\langle f, \, \frac{1}{n} \big\{ \psi_Y(Y_1) - \psi_Y(\tilde{Y}_1) \big\}\big\{ \psi_Z(Z_{\pi_1}) - \overline{\psi}_Z \big\} \Bigg\rangle_{\mathcal{H}_{k \otimes \ell}}  \\[.5em]
	& \hskip 5em +  \sup_{f \in \mathcal{F}_{k \otimes \ell}}  \Bigg\langle f, \, \frac{1}{n} \big\{ \psi_Y(\tilde{Y}_1) - \tilde{\psi}_Y \big\}\big\{ \psi_Z(Z_{\pi_1}) - \psi_Z(\tilde{Z}_{\pi_1}) \big\} \Bigg\rangle_{\mathcal{H}_{k \otimes \ell}},
\end{align*}
where step~$(\mathrm{i})_a$ uses $\frac{1}{n} \sum_{i=1}^n (\overline{\psi}_Y - \tilde{\psi}_Y) \{ \psi_Z(Z_{\pi_i}) - \overline{\psi}_Z\} = 0$, and step~$(\mathrm{ii})_a$ follows similarly. Therefore 
\begin{align*}
	 \bigl| \widehat{\mathrm{HSIC}}(\mathcal{X}_{n}^{\bpi}) - \widehat{\mathrm{HSIC}}(\tilde{\mathcal{X}}_{n}^{\bpi}) \bigr| ~\leq~ & \sup_{f \in \mathcal{F}_{k \otimes \ell}}  \Bigg| \Bigg\langle f, \, \frac{1}{n} \big\{ \psi_Y(Y_1) - \psi_Y(\tilde{Y}_1) \big\}\big\{ \psi_Z(Z_{\pi_1}) - \overline{\psi}_Z \big\} \Bigg\rangle_{\mathcal{H}_{k \otimes \ell}} \Bigg| \\[.5em]
	 + & \sup_{f \in \mathcal{F}_{k \otimes \ell}}  \Bigg| \Bigg\langle f, \, \frac{1}{n} \big\{ \psi_Y(\tilde{Y}_1) - \tilde{\psi}_Y \big\}\big\{ \psi_Z(Z_{\pi_1}) - \psi_Z(\tilde{Z}_{\pi_1}) \big\} \Bigg\rangle_{\mathcal{H}_{k \otimes \ell}} \Bigg| \\[.5em]
	 = ~ & \mathrm{(I)} + \mathrm{(II)}.
\end{align*}
For the term (I), the Cauchy--Schwarz inequality along with the fact $\|f\|_{\mathcal{H}_{k \otimes \ell}} \leq 1$ yields
\begin{align*}
	\mathrm{(I)} ~\leq ~ & \bigg\|   \frac{1}{n} \big\{ \psi_Y(Y_1) - \psi_Y(\tilde{Y}_1) \big\}\big\{ \psi_Z(Z_{\pi_1}) - \overline{\psi}_Z \big\} \bigg\|_{\mathcal{H}_{k \otimes \ell}}  \\[.5em]
	= ~ &   \frac{1}{n} \big\|   \psi_Y(Y_1) - \psi_Y(\tilde{Y}_1) \big\|_{\mathcal{H}_{k}}  \big\|  \psi_Z(Z_{\pi_1}) - \overline{\psi}_Z  \big\|_{\mathcal{H}_{\ell}} \\[.5em]
	\overset{\mathrm{(i)}_b}{\leq} ~ & \frac{1}{n^2} \sum_{i=1}^n  \big\|   \psi_Y(Y_1) - \psi_Y(\tilde{Y}_1) \big\|_{\mathcal{H}_{k}} \big\| \psi_Z(Z_{\pi_1}) - \psi_Z(Z_{\pi_i}) \big\|_{\mathcal{H}_{\ell}} \\[.5em]
	\overset{\mathrm{(ii)}_b}{\leq} ~ & \frac{n-1}{n^2} \sqrt{2K} \sqrt{2L}.
\end{align*}
In the above, step~$(\mathrm{i})_b$ uses Jensen's inequality, and step~$(\mathrm{ii})_b$ follows since 
\begin{align*}
	& \| \psi_Y(y) - \psi_Y(y') \|_{\mathcal{H}_k} = \sqrt{k(y,y) + k(y',y') - 2k(y,y')} \leq \sqrt{2K} \quad \text{and} \\[.5em]
	& \| \psi_Z(z) - \psi_Z(z') \|_{\mathcal{H}_\ell} = \sqrt{\ell(z,z) + \ell(z',z') - 2\ell(z,z')} \leq \sqrt{2L},	
\end{align*}
for any $y,y' \in \mathbb{Y}$ and $z,z' \in \mathbb{Z}$. The second term (II) can be similarly handled and shown to be
\begin{align*}
	\mathrm{(II)} ~\leq ~ \frac{n-1}{n^2} \sqrt{2K} \sqrt{2L}.
\end{align*}
Thus it holds that 
\begin{align*}
	\bigl| \widehat{\mathrm{HSIC}}(\mathcal{X}_{n}^{\bpi}) - \widehat{\mathrm{HSIC}}(\tilde{\mathcal{X}}_{n}^{\bpi}) \bigr| \leq \frac{4(n-1)\sqrt{KL}}{n^2}.
\end{align*}

The other case where $\pi_1 \neq 1$ can be proven similarly. Without loss of generality, we may assume that $\pi_1 = 2$. Then by following the same calculations as before, we have
\begin{align*}
	\bigl| \widehat{\mathrm{HSIC}}(\mathcal{X}_{n}^{\bpi}) - \widehat{\mathrm{HSIC}}(\tilde{\mathcal{X}}_{n}^{\bpi}) \bigr| ~\leq~ & \sup_{f \in \mathcal{F}_{k \otimes \ell}}  \Bigg| \Bigg\langle f, \, \frac{1}{n} \big\{ \psi_Y(Y_1) - \psi_Y(\tilde{Y}_1) \big\}\big\{ \psi_Z(Z_{\pi_1}) - \overline{\psi}_Z \big\} \Bigg\rangle_{\mathcal{H}_{k \otimes \ell}} \Bigg| \\[.5em]
	+ & \sup_{f \in \mathcal{F}_{k \otimes \ell}}  \Bigg| \Bigg\langle f, \, \frac{1}{n} \big\{ \psi_Y(\tilde{Y}_2) - \tilde{\psi}_Y \big\}\big\{ \psi_Z(Z_{\pi_2}) - \psi_Z(\tilde{Z}_{\pi_2}) \big\} \Bigg\rangle_{\mathcal{H}_{k \otimes \ell}} \Bigg| \\[.5em]
	\leq ~ & \frac{4(n-1)\sqrt{KL}}{n^2},
\end{align*}
which completes the proof of the first claim of \Cref{Lemma: Sensitivity of HSIC}.

\subsubsection{Lower Bound} 
We now construct an example where the sensitivity becomes $4(n-2.5)n^{-2}\sqrt{KL}$. For a given $\epsilon \in (0, \min\{K,L\})$, we may assume that there exist $y_a,y_b \in \mathbb{Y}$ and $z_a,z_b \in \mathbb{Z}$ such that $k(y_a,y_b) = \epsilon_{\star} \leq \epsilon$ and $\ell(z_a,z_b) = \epsilon_{\star}' \leq \epsilon$, since $k$ and $\ell$ are assumed to have non-empty level sets on $\mathbb{Y}$ and $\mathbb{Z}$. Consider the two datasets
\begin{align*}
	\mathcal{X}_n = \begin{bmatrix}
		y_a & z_a \\
		y_a & z_a \\
		y_b & z_b \\
		\vdots & \vdots \\
		y_b & z_b
	\end{bmatrix} \quad \text{and} \quad 
	\tilde{\mathcal{X}}_n = \begin{bmatrix}
	y_a & z_a \\
	y_b & z_b \\
	y_b & z_b \\
	\vdots & \vdots \\
	y_b & z_b
	\end{bmatrix},
\end{align*}
where $d_{\mathrm{ham}}(\mathcal{X}_n, \tilde{\mathcal{X}}_n) = 1$. We then consider a permutation $\bpi = (2,1,3,4,\ldots,n)$ and denote the corresponding permuted datasets as
\begin{align*}
	\mathcal{X}_n^{\bpi} = \begin{bmatrix}
		y_a & z_a \\
		y_a & z_a \\
		y_b & z_b \\
		\vdots & \vdots \\
		y_b & z_b
	\end{bmatrix} \quad \text{and} \quad 
	\tilde{\mathcal{X}}_n^{\bpi} = \begin{bmatrix}
		y_a & z_b \\
		y_b & z_a \\
		y_b & z_b \\
		\vdots & \vdots \\
		y_b & z_b
	\end{bmatrix}.
\end{align*}
Using the closed form expression of the squared HSIC in \eqref{Eq: closed form HSIC} and the property of translation invariant kernels, in particular $k(y,y) = K$ and $\ell(z,z) = L$ for all $y,z$, it can be seen that 
\begin{align*}
	\widehat{\mathrm{HSIC}}^2(\mathcal{X}_{n}^{\bpi}) ~=~ & \frac{4+(n-2)^2}{n^2} KL + \frac{4(n-2)\epsilon_\star \epsilon_\star'}{n^2} \\[.5em]
	+ ~ &   \bigg\{ \frac{4 + (n-2)^2}{n^2}K + \frac{4(n-2)}{n^2}\epsilon_\star \bigg\}\bigg\{ \frac{4 + (n-2)^2}{n^2}L + \frac{4(n-2)}{n^2}\epsilon_\star' \bigg\} \\[.5em]
	- ~ 4& \bigg\{ \frac{4}{n^3}KL  + \frac{2K \epsilon_\star'(n-2)}{n^3}  + \frac{2L \epsilon_\star(n-2)}{n^3} + \frac{(n-2)^2\epsilon_\star \epsilon_\star'}{n^3} \bigg\} \\[.5em]
	- ~ 2&(n-2) \bigg\{ \frac{4\epsilon_\star \epsilon_\star'}{n^3} + \frac{2\epsilon_\star L(n-2)}{n^3} + \frac{2 K\epsilon_\star'(n-2)}{n^3} + \frac{KL(n-2)^2}{n^3}\bigg\} \\[.5em]
	= ~ & \frac{16(n-2)^2}{n^4}KL +  C_1 \epsilon_\star \epsilon_\star' + C_2 \epsilon_\star + C_3 \epsilon_\star',
\end{align*}
where $C_1,C_2,C_3$ are constants that only depend on $K,L,n$. A similar calculation shows 
\begin{align*}
		\widehat{\mathrm{HSIC}}^2(\tilde{\mathcal{X}}_{n}^{\bpi}) ~=~ & \frac{4}{n^4} KL + C_1' \epsilon_\star \epsilon_\star' + C_2' \epsilon_\star + C_3' \epsilon_\star',
\end{align*}
where $C_1',C_2',C_3'$ are constants that only depend on $K,L,n$. These results together with the reverse triangle inequality yields
\begin{align*}
	\bigl| \widehat{\mathrm{HSIC}}(\mathcal{X}_{n}^{\bpi}) - \widehat{\mathrm{HSIC}}(\tilde{\mathcal{X}}_{n}^{\bpi}) \bigr| ~\geq~&  \bigg|  \sqrt{\frac{16(n-2)^2}{n^4}KL} - \sqrt{\frac{4}{n^4} KL}  \bigg| - h(\epsilon_\star, \epsilon_\star',K,L,n)\\[.5em]
	= ~ & \frac{4(n-2.5)}{n^2} \sqrt{KL} - h(\epsilon_\star, \epsilon_\star',K,L,n),
\end{align*}
where $h(\epsilon_\star, \epsilon_\star',K,L,n)$ is some function of $\epsilon_\star, \epsilon_\star',K,L,n$, which goes to zero as $\epsilon_\star \rightarrow 0$ and $\epsilon_\star' \rightarrow 0$ for each fixed $K,L,n$. Therefore it holds that 
\begin{align*}
	\sup_{\bpi \in \boldsymbol{\Pi}_n} \sup_{\substack{\mathcal{X}_{n},\tilde{\mathcal{X}}_{n}:\\d_{\mathrm{ham}}(\mathcal{X}_{n},\tilde{\mathcal{X}}_{n}) \leq 1}} \bigl| \widehat{\mathrm{HSIC}}(\mathcal{X}_{n}^{\bpi}) - \widehat{\mathrm{HSIC}}(\tilde{\mathcal{X}}_{n}^{\bpi}) \bigr| \geq \frac{4(n-2.5)}{n^2} \sqrt{KL} - h(\epsilon_\star, \epsilon_\star',K,L,n),
\end{align*}
and the lower bound result follows by letting $\epsilon \rightarrow 0$ given that $\max\{\epsilon_\star, \epsilon_\star'\} \leq \epsilon$. This completes the proof of \Cref{Lemma: Sensitivity of HSIC}.

\subsection{Proof \Cref{Theorem: Properties of dpHSIC}} \label{Section: Proof Theorem: Properties of dpHSIC}
Let us prove the three claims made in \Cref{Theorem: Properties of dpHSIC} below.
\begin{itemize}
	\item \emph{Proof of Differential Privacy.} This result follows by \Cref{Theorem: Differential privacy of permutation tests} along with \Cref{Lemma: Sensitivity of HSIC}. 
	\item \emph{Proof of Validity.} This result follows by \Cref{Theorem: validity of private permutation tests}.
	\item\emph{Proof of Consistency.} The proof of consistency for $\phi_{\mathttt{dpHSIC}}$ is essentially the same as that for $\phi_{\mathttt{dpMMD}}$ in \Cref{Theorem: Properties of dpMMD}. We only need to verify that $\widehat{\mathrm{HSIC}}(\mathcal{X}_{n}^{\bpi_0}) \convP \mathrm{HSIC}_{k \otimes \ell}(P_{YZ})$ and $\widehat{\mathrm{HSIC}}(\mathcal{X}_{n}^{\bpi_1}) \convP 0$. Indeed, the first claim follows by \Cref{Lemma: Exponential inequality for the empirical HSIC} and the second claim follows by \Cref{Lemma: Concentration inequality for permuted HSIC}, given that $\alpha, K, L$ are fixed quantities. The remaining steps are the same as those for the proof of \Cref{Theorem: Properties of dpMMD} and thus we omit the details.
\end{itemize}

\subsection{Proof of \Cref{Theorem: Uniform separation for MMD}} \label{Section: Proof of Theorem: Uniform separation for MMD}
The proof of \Cref{Theorem: Uniform separation for MMD} follows the same structure as that of \Cref{Theorem: General uniform power condition}. The main difference is that we make use of exponential concentration inequalities of the MMD statistic in order to obtain the logarithmic dependence on both $\alpha$ and $\beta$. Throughout the proof, we denote positive constants that only depend on $K$ and $\tau$ by $C_1,C_2,\ldots$ whose values may vary in different places. As in the proof of \Cref{Theorem: General uniform power condition}, we build on the quantile representation of the permutation test~(\Cref{Lemma: Quantile representation}) from which we see that the type II error can be written as
\begin{align*}
	\mE[1 - \phi_{\mathttt{dpMMD}}] =  \mP(M_0 \leq q_{1-\alpha,B}), 
\end{align*}
where $q_{1-\alpha,B}$ denotes the $1-\alpha$ quantile of $\{M_0,M_1,\ldots,M_B\}$, and $M_i$ is given as $\widehat{\mathrm{MMD}}(\mathcal{X}_{n+m}^{\bpi_i}) + 2\Delta_T \xi_{\varepsilon,\delta}^{-1} \zeta_i$ for $i \in \{0\} \cup [B]$.

\vskip 1em 

\noindent \textbf{Step 1 (Bounding the quantile).}  As in the proof of \Cref{Theorem: General uniform power condition}, denote by $q_{1-\alpha/12,\infty}^a$ and $q_{1-\alpha/12,\infty}^b$, the $1-\alpha/12$ quantile of the conditional distribution of $\widehat{\mathrm{MMD}}(\mathcal{X}_{n+m}^{\bpi})$ given $\mathcal{X}_{n+m}$ and that of $2\Delta_T \xi_{\varepsilon,\delta}^{-1} \zeta$, respectively. Then following the argument given in \Cref{Section: Proof Theorem: General uniform power condition} and under the condition for $B$, we see that the inequality $q_{1-\alpha,B} \leq q_{1-\alpha/12,\infty}^a + q_{1-\alpha/12,\infty}^b$ holds with probability at least $1-\beta/2$. 

Next we further upper bound $q_{1-\alpha/12,\infty}^a$ and $q_{1-\alpha/12,\infty}^b$ by more manageable terms. Applying \Cref{Lemma: Bobkovs inequality} with the condition $n \leq m \leq \tau n$ yields the following inequality:
\begin{align*}
	q_{1-\alpha/12,\infty}^a \leq C_1 \sqrt{\frac{1}{n} \log \biggl(\frac{12}{\alpha}\biggr)}.
\end{align*}
On the other hand, as noted in \Cref{Section: Proof Theorem: General uniform power condition}, $q_{1-\alpha/12,\infty}^b$ is simply the $1-\alpha/12$ quantile of $2 \Delta_T\xi_{\varepsilon,\delta}^{-1} \zeta$, namely 
\begin{align*}
	q_{1-\alpha/12,\infty}^b = 2\Delta_T \xi_{\varepsilon,\delta}^{-1} F_{\zeta}^{-1}(1-\alpha/12).
\end{align*}
Therefore, the following inequality holds with probability at least $1-\beta/2$:
\begin{align} \label{Eq: upper bound for q}
	q_{1-\alpha,B} \leq C_1 \sqrt{\frac{1}{n} \log \biggl(\frac{12}{\alpha}\biggr)} + 2\Delta_T \xi_{\varepsilon,\delta}^{-1} F_{\zeta}^{-1}(1-\alpha/12). 
\end{align} 

\vskip .5em

\begin{remark}[Condition for the sample-size ratio] \normalfont  \label{Remark: Markov inequality instead of Exponential inequality}
	The condition $n \leq m \leq \tau n$ is required to apply \Cref{Lemma: Bobkovs inequality} under which we can obtain the logarithmic factor of $\alpha$ in the upper bound for $q_{1-\alpha/12,\infty}^a$. This condition can be eliminated by using \Cref{Lemma: Markov for permuted MMD}, which shows that  
	\begin{align*}
		q_{1-\alpha/12,\infty}^a \leq C_2 \sqrt{\frac{1}{n\alpha}},
	\end{align*}
	without the condition for the sample-size ratio constraint. It remains an open question whether we can achieve the logarithmic factor of $\alpha$ in the upper bound without the constraint $n \leq m \leq \tau n$, which we leave for future research.
\end{remark}

\vskip 1em

\noindent \textbf{Step 2 (Bounding the type II error).} Let $(P,Q)$ be a pair of distributions in $\mathcal{P}_{\mathrm{MMD}_k}\!(\rho)$. Given $(P,Q)$, \Cref{Lemma: Concentration for MMD} (concentration inequality for $\widehat{\mathrm{MMD}}$) ensures that the following event $E_1$:
\begin{align*}
	E_1 \coloneqq   \biggl\{ \bigg|  \mathrm{MMD}_k(P,Q) - \widehat{\mathrm{MMD}}(\mathcal{X}_{n+m}) \bigg| \leq  C_3 \sqrt{\frac{\log(8/\beta)}{n}} \biggr\}
\end{align*}
holds with probability at least $1 - \beta/4$. Notice that the inverse cumulative distribution function of $\zeta \sim \mathsf{Laplace}(0,1)$ is given by 
\begin{align*}
	F^{-1}_{\zeta}(p) \coloneqq  - \mathrm{sign}(p - 0.5) \times \log(1 - 2|p - 0.5|) \quad \text{for $p \in (0,1)$.}
\end{align*}
This yields that the probability of the event 
\begin{align*}
	E_2\coloneqq  \big\{ 2\Delta_T \xi_{\varepsilon,\delta}^{-1} \zeta_0 > 2\Delta_T \xi_{\varepsilon,\delta}^{-1}F^{-1}_{\zeta}(\beta/4)\big\}
\end{align*}
is equal to $1 - \beta/4$. 
Using these preliminary results, it can be seen that the type II error of $\phi_{\mathttt{dpMMD}}$ satisfies
\begin{align*} 
	& \mE[1 - \phi_{\mathttt{dpMMD}}]  ~ \overset{(\mathrm{i})}{=} ~ \mP(M_0 \leq q_{1-\alpha,B}) ~ \overset{(\mathrm{ii})}{=} ~  \mP\bigl(\widehat{\mathrm{MMD}}(\mathcal{X}_{n+m}) + 2\Delta_T \xi_{\varepsilon,\delta}^{-1} \zeta_0 \leq q_{1-\alpha,B}\bigr) \\[.5em] \nonumber
	\overset{(\mathrm{iii})}{\leq} ~ &  \mP\biggl(\mathrm{MMD}_k(P,Q) + 2\Delta_T \xi_{\varepsilon,\delta}^{-1}F^{-1}_{\zeta}(\beta/4) \leq q_{1-\alpha,B} + C_3 \sqrt{\frac{\log(8/\beta)}{n}} \biggr) + \mP(E_1^c) + \mP(E_2^c) \\[.5em] 
	\overset{\mathrm{(iv)}}{\leq} ~ &  \mP\biggl(\mathrm{MMD}_k(P,Q) \leq q_{1-\alpha,B} + C_3 \sqrt{\frac{\log(8/\beta)}{n}} - 2\Delta_T \xi_{\varepsilon,\delta}^{-1}F^{-1}_{\zeta}(\beta/4) \biggr)  + \frac{\beta}{2} \\[.5em]
	\overset{(\mathrm{v})}{\leq} ~ & \mP\biggl(\mathrm{MMD}_k(P,Q) \leq C_1 \sqrt{\frac{\log(12/\alpha)}{n}} + 2\Delta_T \xi_{\varepsilon,\delta}^{-1} F_{\zeta}^{-1}(1-\alpha/12)  \\[.5em]
	& ~~~~~~~~~~~~~~~~~~~~~~~~ + C_3 \sqrt{\frac{\log(8/\beta)}{n}} - 2\Delta_T \xi_{\varepsilon,\delta}^{-1}F^{-1}_{\zeta}(\beta/4) \biggr) + \beta,
\end{align*}
where step~(i) uses \Cref{Lemma: Quantile representation}, step~(ii) uses the definition of $M_0$, step~(iii) follows by the union bound, step~(iv) uses the properties of the events $E_1$ and $E_2$, and step~(v) uses the inequality given in \eqref{Eq: upper bound for q} and the union bound. 

Since $\beta/4 < 1/2$ and $1-\alpha/12 > 1/2$, it follows that 
\begin{align*}
	F^{-1}_{\zeta}(\beta/4) = \log \biggl( \frac{\beta}{2} \biggr) \quad \text{and} \quad F^{-1}_{\zeta}(1-\alpha/12) = \log\biggl( \frac{6}{\alpha} \biggr),
\end{align*}
which gives 
\begin{align*}
	2 \Delta_T\xi_{\varepsilon,\delta}^{-1}  F_{\zeta}^{-1}(1-\alpha/12) - 2\Delta_T \xi_{\varepsilon,\delta}^{-1}F^{-1}_{\zeta}(\beta/4) ~=~&  2\Delta_T \xi_{\varepsilon,\delta}^{-1} \biggl\{  \log \biggl( \frac{6}{\alpha} \biggr) +  \log \biggl( \frac{2}{\beta} \biggr) \biggr\} \\[.5em]
	\leq ~ & C_4 \Delta_T \xi_{\varepsilon,\delta}^{-1} \max \biggl\{  \log \biggl( \frac{1}{\alpha} \biggr), \,  \log \biggl( \frac{1}{\beta} \biggr) \biggr\}.
\end{align*}
Consequently, the type II error of the $\mathttt{dpMMD}$ test is bounded by 
\begin{align*}
	& \mE[1 - \phi_{\mathttt{dpMMD}}] \\[.5em] 
	\leq ~ & \mP\biggl(\mathrm{MMD}_k(P,Q) \leq C_5 \sqrt{\frac{\max\bigl\{\log(1/\alpha),\log(1/\beta)\}}{n}}   + C_6 \frac{\max\bigl\{\log(1/\alpha),\log(1/\beta)\}}{n\xi_{\varepsilon,\delta}} \biggr) +  \beta \\[.5em]
	\leq ~ & \beta,
\end{align*}
where the last inequality holds by taking $C_{K,\tau}$ to be larger than, for instance, $2 \max \{C_5, C_6\} + 1$ in the theorem statement. Finally, we note that the above upper bound is independent of $(P,Q)$, and thus complete the proof by taking the supremum on both sides over $\mathcal{P}_{\mathrm{MMD}_k}\!(\rho)$.

\subsection{Proof of \Cref{Theorem: Minimax separation in MMD}} \label{Section: Proof of Theorem: Minimax separation in MMD}
Under the conditions of \Cref{Theorem: Minimax separation in MMD}, we analyze that the minimax separation $\rho^\star_{\mathrm{MMD}} (\alpha,\beta,\varepsilon,\delta,m,n) =\rho^\star_{\mathrm{MMD}}$ in the non-privacy regime and in the privacy regime, separately. In both regimes, we use a common trick that reduces the two-sample problem to the one-sample problem (also known as goodness-of-fit testing) stated below. 

\paragraph{Reducing the two-sample problem to the one-sample problem.} We can think of the one-sample problem as a special case of the two-sample problem by assuming that we have access to an infinite number of observations from one of the two distributions (say $Q$). That is, we know the entire information of the distribution $Q$ and $m = \infty$. \citet[][Lemma 1]{arias2018remember} builds on this simple observation and formalizes that the minimax risk of the two-sample problem is greater than or equal to that of the one-sample problem. We follow this strategy and work with the one-sample problem. 

More formally, pick one distribution $Q_0$ (specified later in the proof) from $\mathcal{P}$ and fix it throughout. Let $\phi: \mathcal{X}_n \mapsto \{0,1\}$ be a test function and define the set of $(\varepsilon,\delta)$-DP level $\alpha$ one-sample tests as 
\begin{align*}
	\Phi_{\alpha,\varepsilon,\delta,Q_0} \coloneqq   \Big\{\phi : \mE_{Q_0}[\phi] \leq \alpha \ \text{and} \ \text{$\phi$ is $(\varepsilon,\delta)$-DP} \Big\}.
\end{align*}
Then, letting $\mathcal{P}_{1,\mathrm{MMD}_k}(\rho; Q_0) \coloneqq  \{ P \in \mathcal{P} : \mathrm{MMD}_k(P,Q_0) \geq \rho \}$, \citet[][Lemma 1]{arias2018remember} yields
\begin{align*}
	\rho^\star_{\mathrm{MMD}} \geq \inf\Bigl\{ \rho > 0:  \inf_{\phi \in \Phi_{\alpha,\varepsilon,\delta,Q_0}} \sup_{P  \in \mathcal{P}_{1,\mathrm{MMD}_k}(\rho; Q_0)} \mE_{P}[1 - \phi] \leq \beta \Big\}.
\end{align*}
Therefore, once we prove 
\begin{subequations}
\begin{align}
	& \inf\Bigl\{ \delta > 0 :  \inf_{\phi \in \Phi_{\alpha,\varepsilon,Q_0}} \sup_{P  \in \mathcal{P}_{1,\mathrm{MMD}_k}(\rho;Q_0)} \mE_{P}[1 - \phi] \leq \beta \Big\} \nonumber\\[.5em]
	\geq ~ & C_{\eta} \max \Biggl\{ \min \bigg\{  \sqrt{\frac{\log(1/(\alpha+\beta))}{n}}, \, 1 \bigg\}, \, \min \biggl\{ \frac{\log(1/\beta)}{n(\varepsilon+\delta)}, \, 1 \biggr\} \Biggr\} \label{Equation: MMD rate (i)}\\[.5em]
	\geq ~ & C_{\eta} \max \Biggl\{ \min \bigg\{  \sqrt{\frac{\log(1/(\alpha+\beta))}{n}}, \, 1 \bigg\}, \, \min \biggl\{ \frac{\log(1/\beta)}{n\big(\varepsilon+\log(1/(1-\delta))\big)}, \, 1 \biggr\} \Biggr\}, \label{Equation: MMD rate (ii)}
\end{align}
\end{subequations}
for some $Q_0 \in \mathcal{P}$, then the desired result of \Cref{Theorem: Minimax separation in MMD} follows.
The two lower bounds in \eqref{Equation: MMD rate (ii)} are proved in \Cref{Section: Non-privacy regime MMD,Section: Privacy regime MMD}.
The inequality \eqref{Equation: MMD rate (i)} $\geq$ \eqref{Equation: MMD rate (ii)} holds as $\log(1/(1-\delta))\geq\delta$ for all $\delta\in[0,1)$, and is sufficient to prove \Cref{Theorem: Minimax separation in MMD}. 
Nonetheless, to provide intuition on the tightness of the lower bound, it is interesting to understand that when $\alpha\asymp\beta$ the two rates in \Cref{Equation: MMD rate (i),Equation: MMD rate (ii)} are the same, which we prove in \Cref{Section: Equivalence MMD}.

\subsubsection{Non-privacy regime}  \label{Section: Non-privacy regime MMD}

We first show that $\rho^\star_{\mathrm{MMD}} \geq C_\eta \min\big\{\sqrt{\log(1/(\alpha+\beta))/n}, 1\big\}$ in the non-privacy regime. Writing the set of non-private level $\alpha$ tests (\emph{i.e.}, $\varepsilon \to \infty$) as
\begin{align*}
	\Phi_{\alpha, \infty, Q_0} \coloneqq  \Big\{\phi : \mE_{Q_0}[\phi] \leq \alpha \Big\},
\end{align*}
we clearly have $\Phi_{\alpha, \varepsilon, \delta, Q_0} \subset \Phi_{\alpha, \infty, Q_0}$. Hence, in this non-private regime, it suffices to show 
\begin{align*}
	\inf\Bigl\{ \rho > 0 :  \inf_{\phi \in \Phi_{\alpha,\infty,Q_0}} \sup_{P \in \mathcal{P}_{1,\mathrm{MMD}_k}(\rho;Q_0)} \mE_{P}[1 - \phi] \leq \beta \Big\} \geq C_\eta \min \Bigg\{ \sqrt{\frac{\log(1/(\alpha+\beta))}{n}}, \, 1 \Bigg\},
\end{align*}
as the infimum term above is smaller than or equal to $\rho^\star_{\mathrm{MMD}}$.
To this end, we use the standard Le Cam's two point method~\citep{lecam1973convergence,le2012asymptotic}. In particular, for any $P_0 \in \mathcal{P}_{1,\mathrm{MMD}_k}(\widetilde{\rho};Q_0)$ with
\begin{align*}
	\widetilde{\rho} = C_\eta \min \Bigg\{ \sqrt{\frac{\log(1/(\alpha+\beta))}{n}}, \, 1 \Bigg\},
\end{align*}
we have
\begin{align*}
	\inf_{\phi \in \Phi_{\alpha,\infty,Q_0}} \sup_{P  \in \mathcal{P}_{1,\mathrm{MMD}_k}(\widetilde{\rho};Q_0)} \mE_{P}[1 - \phi]  ~\geq~ & \inf_{\phi \in \Phi_{\alpha,\infty,Q_0}} \mE_{P_0}[1 - \phi]  = 1 - \sup_{\phi \in \Phi_{\alpha,\infty,Q_0}} \mE_{P_0}[\phi] \\[.5em]
	\geq ~ & 1 - \alpha - d_{\mathrm{TV}}(P_0^{\otimes n}, Q_0^{\otimes n})  \\[.5em]
	\overset{(\mathrm{i})}{\geq} ~ & 1 - \alpha - 1 + \frac{1}{2} e^{- d_{\mathrm{KL}}(P_0^{\otimes n} \| Q_0^{\otimes n})}  ~ \overset{(\mathrm{ii})}{=} ~ \frac{1}{2} e^{-n \times d_{\mathrm{KL}}(P_0 \| Q_0)} - \alpha,
\end{align*} 
where step~(i) holds by Bretagnolle--Huber inequality~\citep[][Lemma B.4]{canonne2022topics} and step~(ii) uses the chain rule of the KL divergence. Therefore the minimax type II error is bounded by $\beta$, that is
\begin{align*}
	\inf_{\phi \in \Phi_{\alpha,\infty,Q_0}} \sup_{P  \in \mathcal{P}_{1,\mathrm{MMD}_k}(\widetilde{\rho};Q_0)} \mE_{P}[1 - \phi] \geq \beta,
\end{align*}
provided that $\alpha + \beta < 0.4$ and 
\begin{align*}
	d_{\mathrm{KL}}(P_0 \| Q_0)  \leq \frac{1}{n} \log\biggl( \frac{1}{2(\alpha + \beta)} \biggr).
\end{align*}
Hence it is enough to find an instance of $(P_0,Q_0)$ in $\mathbb{R}^d$ such that 
\begin{subequations}
	\begin{align} \label{Eq: sufficient condition in non-private regime a}
		& \mathrm{MMD}_k(P_0, Q_0) \geq C_\eta \min \Bigg\{ \sqrt{\frac{\log(1/(\alpha+\beta))}{n}}, \, 1 \Bigg\} \quad \text{and} \\[.5em]   \label{Eq: sufficient condition in non-private regime b}
		& d_{\mathrm{KL}}(P_0 \| Q_0) \leq \frac{1}{n} \log\biggl( \frac{1}{2(\alpha + \beta)} \biggr). 
	\end{align}
\end{subequations}
Motivated by the proof of \citet[][Theorem 1]{tolstikhin2017minimax}, we pick two discrete distributions $P_0 = p_0 \delta_x + (1 - p_0) \delta_v$ and $Q_0 = q_0 \delta_x + (1-q_0) \delta_v$, where $x,v \in \mathbb{R}^d$, $0 < p_0, q_0 < 1$ and $\delta_x$ is a Dirac measure at $x$. In this setting, the calculations in \cite{tolstikhin2017minimax} show that 
\begin{align} \label{Eq: Explicit expression}
	\mathrm{MMD}_k(P_0, Q_0) = \sqrt{2 (p_0 - q_0)^2 \bigl( \kappa(0) - \kappa(x-v) \bigr)}.
\end{align} 
In addition, under the same setting, the KL divergence of $P_0$ from $Q_0$ is bounded by
\begin{align*}
	d_{\mathrm{KL}}(P_0 \| Q_0) \leq \frac{(p_0 - q_0)^2}{q_0(1-q_0)}.
\end{align*}
Now set
\begin{align*}
	p_0 = \frac{1}{2} + \min \Bigg\{  \sqrt{\frac{1}{4n} \log\biggl( \frac{1}{2(\alpha + \beta)} \biggr)}, \frac{1}{2} \Bigg\} \quad \text{and} \quad q_0 = \frac{1}{2},
\end{align*}
so that both $P_0$ and $Q_0$ are valid probability measures, which leads to 
\begin{align*}
	d_{\mathrm{KL}}(P_0 \| Q_0) \leq \min \Biggl\{\frac{1}{n} \log\biggl( \frac{1}{2(\alpha + \beta)} \biggr), \, 1\Biggr\} \leq \frac{1}{n} \log\biggl( \frac{1}{2(\alpha + \beta)} \biggr).
\end{align*}
We also choose $x$ and $v$ such that $\kappa(0) - \kappa(z) \geq \eta$ for $z = x - v$ where $\eta$ is given in the theorem assumption. In this setting, the explicit expression of the MMD metric~\eqref{Eq: Explicit expression} can be used to show 
\begin{align*}
	\mathrm{MMD}_k(P_0, Q_0) ~\geq~ & \sqrt{2\eta}  \min \Bigg\{  \sqrt{\frac{1}{4n} \log\biggl( \frac{1}{2(\alpha + \beta)} \biggr)}, \, \frac{1}{2} \Bigg\} \\[.5em]
	\geq~ & C_\eta \min \Bigg\{  \sqrt{\frac{1}{n} \log\biggl( \frac{1}{\alpha + \beta} \biggr)}, \, 1 \Bigg\}
\end{align*}
where the second inequality holds by taking $C_\eta$ sufficiently small under the condition $\alpha + \beta < 0.4$. Hence we have shown that the sufficient conditions in \eqref{Eq: sufficient condition in non-private regime a} and \eqref{Eq: sufficient condition in non-private regime b} are fulfilled, and thus the minimax separation satisfies
\begin{align*}
	\rho^\star_{\mathrm{MMD}} \geq C_\eta \min \Bigg\{  \sqrt{\frac{\log(1/(\alpha+\beta))}{n}}, \, 1 \Bigg\}. 
\end{align*}
Next we focus on the privacy regime and prove the second term in the lower bound.

\subsubsection{Privacy regime}
\label{Section: Privacy regime MMD}
The second term in the lower bound can be proved based on the coupling idea described in \citet[][Theorem 11]{acharya2018differentially} and \citet[][Theorem 1]{acharya2021differentially}. In detail, we pick $P_0$ from $\mathcal{P}_{1,\mathrm{MMD}_k}(\delta;Q_0)$, and let $(\mathcal{X}_n, \mathcal{Y}_n)$ be a coupling between $P_0^{\otimes n}$ and $Q_0^{\otimes n}$ with $D = \mE[d_{\mathrm{ham}}(\mathcal{X}_n, \mathcal{Y}_n)]$. Then for any $\phi \in \Phi_{\alpha,\varepsilon,Q_0}$, the proof of \citet[][Theorem 1]{acharya2021differentially} shows 
\begin{align*}
	1 - \alpha ~\leq~ \mE_{Q_0}[1 -  \phi] &~\leq~ \mE_{P_0}[1 - \phi] \cdot e^{10D \varepsilon} + 0.1 + 10D\delta e^{10D \varepsilon}\\
	\mE_{P_0}[1 - \phi]  &~\geq~ e^{-10D \varepsilon } (0.9 - \alpha) - 10D\delta.
\end{align*}
For completeness, we give details here. By Markov's inequality, the event $\mathcal{E}^c =\{ d_{\mathrm{ham}}(\mathcal{X}_n, \mathcal{Y}_n) > 10D\}$ holds with probability at most $1/10$ (here one can replace the constant $10$ with some other positive number). Moreover,
\begin{align*}
	\mE_{Q_0}[1 -  \phi] = \mP_{Q_0}(\phi = 0) ~\leq~& \mP_{Q_0}(\phi = 0 \given \mathcal{E}) \mP(\mathcal{E}) + \mP(\mathcal{E}^c)  \\[.5em]
	\leq ~ & \mP_{Q_0}(\phi = 0 \given \mathcal{E})\mP(\mathcal{E}) + 0.1 \\[.5em] 
	\overset{(\dagger)}{\leq} ~ & \mP_{P_0}(\phi = 0 \given \mathcal{E}) \mP(\mathcal{E}) \cdot e^{10D \varepsilon}  + 0.1 + 10D\delta e^{(10D-1)\varepsilon}\\[.5em]
	\leq ~ & \mE_{P_0}[1 - \phi] \cdot e^{10D \varepsilon} + 0.1 + 10D\delta e^{10D\varepsilon},
\end{align*}
where step~$(\dagger)$ holds due to the fact that $\phi$ is $\varepsilon$-differentially private and inequality~\eqref{Eq: implication of DP definition}. Letting $E\coloneqq 10D$ for notation purposes, and using the fact that $\alpha \in (0,1/5)$, we obtain
\begin{align*}
	\inf_{\phi \in \Phi_{\alpha,\varepsilon,Q_0}} \sup_{P  \in \mathcal{P}_{1,\mathrm{MMD}_k}(\rho;Q_0)} \mE_{P}[1 - \phi] ~\geq~ &  \inf_{\phi \in \Phi_{\alpha,\varepsilon,Q_0}} \mE_{P_0}[1 - \phi] \\[.5em]
	\geq ~ & e^{-E \varepsilon} (0.9 - \alpha) - E \delta \\[.5em]
	\geq ~ & 0.5 e^{-E \varepsilon} - E \delta  \\[.5em]
	\geq ~ & 0.25 e^{-E \varepsilon},
\end{align*}
provided that the condition $E\delta \leq 0.25e^{-E\varepsilon}$ holds.
Hence, we require 
\begin{align*}
	\inf_{\phi \in \Phi_{\alpha,\varepsilon,Q_0}} \sup_{P  \in \mathcal{P}_{1,\mathrm{MMD}_k}(\rho;Q_0)} \mE_{P}[1 - \phi]  ~\geq~ 0.25 e^{-E\varepsilon} ~\geq~ \beta,
\end{align*}
which is fulfilled when 
\begin{equation}
	\label{cond1}
	E ~\leq~ \frac{1}{\varepsilon} \log \biggl( \frac{1}{4\beta} \biggr)
	\qquad \text{and} \qquad
	E\delta ~\leq~ \frac{1}{4}e^{-E\varepsilon}.
\end{equation} 
The second condition, equivalently $E\delta e^{E\varepsilon} \leq 1/4$, is trivial when $\delta=0$. Hence by assuming $\delta >0$, we verify that
\begin{align*}
    E \leq \frac{1}{4(\varepsilon+\delta)} \quad \text{implies} \quad E\delta e^{E\varepsilon} \leq \frac{1}{4}.
\end{align*}
To see this, under the condition $E \leq \frac{1}{4(\varepsilon+\delta)}$ and letting $x = \varepsilon/\delta$,
\begin{align*}
    E\delta e^{E\varepsilon} \leq \frac{\delta}{4(\varepsilon+\delta)} e^{\frac{\varepsilon}{4(\varepsilon+\delta)}} = \frac{1}{4(x+1)}e^{\frac{x}{4(x+1)}}
\end{align*}
Hence it suffices to show that for all $x>0$ we have
\begin{align*}
    \frac{1}{4(x+1)}e^{\frac{x}{4(x+1)}} \leq \frac{1}{4} \qquad \Longleftrightarrow \qquad f(x) \coloneqq  4(x+1) \log(x+1) - x \geq 0.
\end{align*}
The first derivative of $f$ is given as $f'(x) = 4\log(x+1)+3 > 0$, which implies that $f$ is monotone increasing and $f(x) \geq 0$ for all $x>0$. 
Hence, for all $\varepsilon>0$ and $\delta\in[0,1)$, the conditions in \Cref{cond1} are satisfied when 
\begin{equation*}
	E ~\leq~ \frac{1}{\varepsilon} \log \biggl( \frac{1}{4\beta} \biggr)
	\qquad \text{and} \qquad
	E ~\leq~ \frac{1}{4(\varepsilon+\delta)},
\end{equation*} 
which are in particular satisfied when
\begin{equation}
	\label{cond3}
	E ~\leq~ \frac{1}{4(\varepsilon+\delta)} \min \Biggl\{ \log \biggl( \frac{1}{4\beta} \biggr) , 1 \Biggl\}.
\end{equation} 
Furthermore, given $P_0$ and $Q_0$, \citet[][Lemma 20]{acharya2021differentially} ensures the existence of a coupling between $P_0^{\otimes n}$ and $Q_0^{\otimes n}$ such that 
\begin{align*}
	D = \mE[d_{\mathrm{ham}}(\mathcal{X}_1^n, \mathcal{Y}_1^n)] = \frac{n}{2} \| P_0 - Q_0 \|_1,
\end{align*}
so $E=10D=5n\| P_0 - Q_0 \|_1$. The condition in Equation~(\ref{cond3}) becomes
\begin{equation}
	\label{cond4}
	\| P_0 - Q_0 \|_1 ~\leq~ \frac{1}{20n(\varepsilon+\delta)} \min \Biggl\{ \log \biggl( \frac{1}{4\beta} \biggr) , 1 \Biggl\}.
\end{equation}
As in the non-private case, we can pick two discrete distributions $P_0 = p_0 \delta_x + (1 - p_0) \delta_v$ and $Q_0 = q_0 \delta_x + (1-q_0) \delta_v$, where $q_0 = 1/2$ and 
\begin{align*}
	p_0 = \frac{1}{2} + \min \Biggl\{ \frac{1}{40n(\varepsilon+\delta)} \log \biggl( \frac{1}{4\beta} \biggr), \, \frac{1}{2} \Biggr\}.
\end{align*}
This choice ensures that 
\begin{align*}
	\| P_0 - Q_0 \|_1  = 2|p_0 - q_0|  = 2 \min \Biggl\{ \frac{1}{40n(\varepsilon+\delta)} \log \biggl( \frac{1}{4\beta} \biggr), \, \frac{1}{2} \Biggr\}. 
\end{align*}
We also choose $x$ and $v$ such that $\kappa(0) - \kappa(z) \geq \eta$ for $z = x - v$ under which 
\begin{align*}
	\mathrm{MMD}_k(P_0, Q_0) 
	&\geq \sqrt{2\eta} \min \Biggl\{ \frac{1}{40n(\varepsilon+\delta)} \log \biggl( \frac{1}{4\beta} \biggr), \, \frac{1}{2} \Biggr\} \\[.5em]	
	&\geq  C_\eta \min \Biggl\{ \frac{1}{n(\varepsilon+\delta)} \log \biggl( \frac{1}{\beta} \biggr), \, 1 \Biggr\},
\end{align*}
where the second inequality holds under the condition $\beta \in (0,1/5)$. Therefore, under privacy regime, for all $\varepsilon>0$ and all $\delta\in[0,1)$, it holds that
\begin{align*}
	\rho^\star_{\mathrm{MMD}} \geq  C_\eta \min \Biggl\{ \frac{\log(1/\beta)}{n(\varepsilon+\delta)}, \, 1 \Biggr\},
\end{align*}
which proves the inequality in \eqref{Equation: MMD rate (i)}.

\subsubsection{Equivalence of rates}  \label{Section: Equivalence MMD}

We now prove the equivalence in \eqref{Equation: MMD rate (ii)} when $\alpha\asymp\beta$.
As previously noted, we have \eqref{Equation: MMD rate (i)} $\geq$ \eqref{Equation: MMD rate (ii)} as 
$$
	\log(1/(1-\delta))\geq\delta \quad\textrm{ for all $\delta\in[0,1)$.}
$$
Hence, it remains to show that \eqref{Equation: MMD rate (i)} $\leq$ \eqref{Equation: MMD rate (ii)}, so that, up to constants, we have
\begin{align*}
	&\max \Biggl\{ \min \bigg\{  \sqrt{\frac{\log(1/\beta)}{n}}, \, 1/2 \bigg\}, \, \min \biggl\{ \frac{\log(1/\beta)}{n(\varepsilon+\delta)}, \, 1/2 \biggr\} \Biggr\} \\[.5em]
	\overset{(\star)}{\asymp} ~ &\max \Biggl\{ \min \bigg\{  \sqrt{\frac{\log(1/\beta)}{n}}, \, 1/2 \bigg\}, \, \min \biggl\{ \frac{\log(1/\beta)}{n\big(\varepsilon+\log(1/(1-\delta))\big)}, \, 1/2 \biggr\} \Biggr\}.
\end{align*}
Here, we have used the fact that $\alpha\asymp\beta$ and have absorbed some terms in the constants to replace the $1$'s with $1/2$'s in the minima.
For the result to be non-trivial we need to assume
\begin{equation}
	\label{Equation: Equivalence condition 1/2}
	\sqrt{\frac{\log(1/\beta)}{n}} ~<~ 1/2
	\qquad\textrm{and}\qquad 
	\frac{\log(1/\beta)}{n\big(\varepsilon+\log(1/(1-\delta))\big)} ~<~ 1/2.
\end{equation}

\begin{itemize}
	\item If $\delta\in[0,1/2)$, then $2\log(2) \delta \geq \log(1/(1-\delta)) \geq \delta$, and so we get
	$$
	\frac{\log(1/\beta)}{n\big(\varepsilon+\log(1/(1-\delta))\big)}
	~\leq~\frac{\log(1/\beta)}{n(\varepsilon+\delta)}
	~\leq~ 2\log(2)\frac{\log(1/\beta)}{n\big(\varepsilon+\log(1/(1-\delta))\big)}
	$$
	which proves $(\star)$ for $\delta\in[0,1/2)$.
	\item If $\delta\in[1/2,1)$, then $\delta\geq 1/2 > \sqrt{\log(1/\beta)/n}$, hence we get
	$$
	\frac{\log(1/\beta)}{n \big(\varepsilon + \log(1/(1-\delta))\big)}
	~\leq~\frac{\log(1/\beta)}{n\delta}
	~<~ \sqrt{\frac{\log(1/\beta)}{n}}
	$$ 
	and 
	$$
	\frac{\log(1/\beta)}{n (\varepsilon+ \delta)}
	~<~ \sqrt{\frac{\log(1/\beta)}{n}},
	$$
	so both sides of $(\star)$ are equal to $\sqrt{{\log(1/\beta)}/{n}}$ when $\delta\in[1/2,1)$. So, when $\delta \geq \sqrt{\log(1/\beta) / n}$, the non-DP rate dominates.
\end{itemize}

We have shown that $(\star)$ holds for all $\delta\in[0,1)$, which completes the proof.

\subsection{Proof of \Cref{Theorem: Minimax Separation over L2}} \label{Section: Proof of Theorem: Minimax Separation over L2}
To simplify the notation, we denote the square of the empirical MMD based on $\mathcal{X}_{n+m}$ as $V_{\mathrm{MMD}}$ (or equivalently the V-statistic of MMD), and that based on permuted data $\mathcal{X}_{n+m}^{\bpi}$ as $V_{\bpi,\mathrm{MMD}}$. We also denote the U-statistic of the MMD as $U_{\mathrm{MMD}}$, \emph{i.e.},
\begin{align*}
	U_{\mathrm{MMD}} = \frac{1}{n(n-1)} \sum_{1 \leq i \neq j \leq n} k_{\blambda}(Y_i,Y_j)  + \frac{1}{m(m-1)} \sum_{1 \leq i \neq j \leq m} k_{\blambda}(Z_i,Z_j)  - \frac{2}{nm}\sum_{i=1}^n \sum_{j=1}^m k_{\blambda}(Y_i,Z_j),
\end{align*}
and $U_{\bpi,\mathrm{MMD}}$ is similarly defined based on $\mathcal{X}_{n+m}^{\bpi}$. Remark that the difference between the V-statistic and the U-statistic of MMD is
\begin{equation}
	\begin{aligned} \label{Eq: difference between V and U}
		V_{\mathrm{MMD}} - U_{\mathrm{MMD}} ~=~ &  \frac{1}{n} \prod_{i=1}^d \frac{1}{\sqrt{2 \pi} \lambda_i} + \frac{1}{m} \prod_{i=1}^d \frac{1}{\sqrt{2 \pi} \lambda_i} \\[.5em]
		& - \frac{1}{n^2(n-1)} \sum_{1 \leq i \neq j \leq n} k_{\blambda}(Y_i,Y_j) - \frac{1}{m^2(m-1)} \sum_{1 \leq i \neq j \leq m} k_{\blambda}(Z_i,Z_j). 
	\end{aligned}
\end{equation}
Given this connection, our proof strategy is to leverage known results of the U-statistic along with a careful analysis of the difference between $V_\mathrm{MMD}$ and $U_\mathrm{MMD}$. Writing the sensitivity of $\sqrt{V_{\mathrm{MMD}}}$ as $\Delta_{V^{1/2}}$, \Cref{Lemma: Quantile representation} ensures that the $\mathttt{dpMMD}$ test rejects the null if and only if 
\begin{align*}
	\sqrt{V_\mathrm{MMD}} + \frac{2\Delta_{V^{1/2}}}{\xi_{\varepsilon,\delta}} \zeta_0 > q_{1-\alpha,B},
\end{align*}
where $q_{1-\alpha,B}$ is the $1-\alpha$ quantile of $\{ \sqrt{V_{\bpi_i,\mathrm{MMD}}} + 2 \Delta_{V^{1/2}}\xi_{\varepsilon,\delta}^{-1} \zeta_i \}_{i=0}^B$. We then closely follow the proof steps given in \Cref{Section: Proof Theorem: General uniform power condition} to prove the claim.

\vskip 1em

\noindent \textbf{Bounding the type II error.} \sloppy As in the proof of \Cref{Theorem: General uniform power condition} in \Cref{Section: Proof Theorem: General uniform power condition}, denote by $q_{1-\alpha/12,\infty}^a$ and $q_{1-\alpha/12,\infty}^b$, the $1-\alpha/12$ quantile of the conditional distribution of $V_{\bpi,\mathrm{MMD}}^{1/2}$ given $\mathcal{X}_{n+m}$ and that of $2\Delta_{V^{1/2}} \xi_{\varepsilon,\delta}^{-1} \zeta$, respectively. Then for $B \geq 6\alpha^{-1} \log(2/\beta)$, the analysis given in the proof of \Cref{Theorem: General uniform power condition} shows that the type II error of the $\mathttt{dpMMD}$ test can be bounded as
\begin{align*}
	\mP\bigl(\sqrt{V_\mathrm{MMD}} + 2 \Delta_{V^{1/2}} \xi_{\varepsilon,\delta}^{-1} \zeta_0 \leq q_{1-\alpha,B}\bigr) ~ \leq ~ & \mP\bigl( \sqrt{V_\mathrm{MMD}} + 2 \Delta_{V^{1/2}} \xi_{\varepsilon,\delta}^{-1} \zeta_0 \leq q_{1-\alpha/12,\infty}^a + q_{1-\alpha/12,\infty}^b \bigr) + \beta/2 \\[.5em]
	\leq ~ & \mP\bigl( \sqrt{V_\mathrm{MMD}}  \leq q_{1-\alpha/12,\infty}^a + q_{1-\alpha/12,\infty}^b + R_2 \bigr) + 5\beta/8,
\end{align*}
where we recall $R_2 = 2 \Delta_{V^{1/2}} \xi_{\varepsilon,\delta}^{-1} F_{\zeta}^{-1}(1-\beta/8)$. We also note that the quantile $ q_{1-\alpha/12,\infty}^b$ can be explicitly written as
\begin{align*}
	q_{1-\alpha/12,\infty}^b = 2 \Delta_{V^{1/2}} \xi_{\varepsilon,\delta}^{-1} F_{\zeta}^{-1}(1-\alpha/12),
\end{align*}
and using the inequality~\eqref{Eq: quantile bound for Laplace}, the type II error is upper bounded as 
\begin{align} \nonumber
	& \mP\bigl(\sqrt{V_\mathrm{MMD}} + 2 \Delta_{V^{1/2}} \xi_{\varepsilon,\delta}^{-1} \zeta_0 \leq q_{1-\alpha,B}\bigr) \\[.5em]
	\leq ~ & \mP\bigl( \sqrt{V_\mathrm{MMD}} \leq q_{1-\alpha/12,\infty}^a + 14 \Delta_{V^{1/2}} \xi_{\varepsilon,\delta}^{-1} \max\{\log(1/\alpha),\log(1/\beta)\} \bigr) + 5\beta/ 8. \label{Eq: type II error bound}
\end{align}
Given this bound for the type II error, our next effort lies in bounding $q_{1-\alpha/12,\infty}^a$. 

\vskip 1em

\noindent \textbf{Bounding $q_{1-\alpha/12,\infty}^a$.} To obtain an upper bound for $q_{1-\alpha/12,\infty}^a$, we leverage the exponential inequality for permuted U-statistics studied in \cite{kim2020minimax} and \cite{schrab2021mmd}. To this end, denote the probability function, which is taken over $\bpi$ conditional on $\mathcal{X}_{m+n}$, as $\mP_{\bpi}$ and note that for any $t>0$,
\begin{align*}
	& \mP_{\bpi} \Biggl( \sqrt{V_{\bpi,\mathrm{MMD}}} \geq \sqrt{\left(\frac{1}{n} + \frac{1}{m}\right) \prod_{i=1}^d \frac{1}{\sqrt{2 \pi} \lambda_i} + t} \Biggr) \\[.5em]
	=  ~ & \mP_{\bpi} \Biggl( V_{\bpi,\mathrm{MMD}} \geq \left(\frac{1}{n} + \frac{1}{m}\right) \prod_{i=1}^d \frac{1}{\sqrt{2 \pi} \lambda_i} + t  \Biggr) \\[.5em]
	= ~ &    \mP_{\bpi} \Biggl( V_{\bpi,\mathrm{MMD}} - U_{\bpi,\mathrm{MMD}} + U_{\bpi,\mathrm{MMD}} \geq \left(\frac{1}{n} + \frac{1}{m}\right) \prod_{i=1}^d \frac{1}{\sqrt{2 \pi} \lambda_i} + t  \Biggr) \\[.5em]
	=~ &  \mP_{\bpi} \Biggl( \biggl\{V_{\bpi,\mathrm{MMD}} - U_{\bpi,\mathrm{MMD}} - \left(\frac{1}{n} + \frac{1}{m}\right) \prod_{i=1}^d \frac{1}{\sqrt{2 \pi} \lambda_i}\biggr\}+ U_{\bpi,\mathrm{MMD}} \geq t  \Biggr) \\[.5em]
	\leq~ &  \mP_{\bpi} \bigl( U_{\bpi,\mathrm{MMD}} \geq t \bigr),
\end{align*}
where the last equality holds since for any permutation $\bpi$
\begin{align*}
	V_{\bpi,\mathrm{MMD}} - U_{\bpi,\mathrm{MMD}}  - \left(\frac{1}{n} + \frac{1}{m}\right) \prod_{i=1}^d \frac{1}{\sqrt{2 \pi} \lambda_i} \leq 0,
\end{align*}
which can be seen from the previous expression of $V_\mathrm{MMD}-U_\mathrm{MMD}$ in \eqref{Eq: difference between V and U} and non-negativity of $k_{\blambda}$. Hence $q_{1-\alpha/12,\infty}^a$ is bounded by 
\begin{align*}
	q_{1-\alpha/12,\infty}^a \leq \sqrt{\left(\frac{1}{n} + \frac{1}{m}\right) \prod_{i=1}^d \frac{1}{\sqrt{2 \pi} \lambda_i} + q_{1-\alpha/12,\infty}^U},
\end{align*}
where $q_{1-\alpha/12,\infty}^U$ denotes the $1-\alpha/12$ quantile of the permutation distribution of $U_{\bpi,\mathrm{MMD}}$. Now using the exponential bound for $U_{\bpi,\mathrm{MMD}}$ in \citet[][Equation 59]{kim2020minimax}, and following the proof of \citet[][Proposition 4]{schrab2021mmd}, we see that there exists a constant $C_{\tau,\beta,M,d}$ such that 
\begin{align*}
	\mP\biggl( q_{1-\alpha/12,\infty}^U \leq C_{\tau,\beta,M,d} \frac{\log(1/\alpha)}{n \sqrt{\lambda_1 \cdots \lambda_d}} \biggr) \geq 1 - \beta/8,
\end{align*}
under the conditions of \Cref{Theorem: Minimax Separation over L2}. For simplicity, let us write 
\begin{align} \label{Eq: definition of b_{n,m,blamda}}
	b_{n,m,\blambda} =  \left(\frac{1}{n} + \frac{1}{m}\right) \prod_{i=1}^d \frac{1}{\sqrt{2 \pi} \lambda_i} \leq  \frac{C_{\tau,d}}{n\lambda_1 \dots \lambda_d}.
\end{align}
Observe that we have $\max\{\log(1/\alpha),\log(1/\beta)\} \leq C_\beta \log(1/\alpha)$ for $\alpha \in (0,e^{-1})$. Putting these pieces together and continuing from \eqref{Eq: type II error bound}, the type II error is further bounded by 
\begin{align}  \nonumber
	& \mP\bigl(\sqrt{V_\mathrm{MMD}} + 2 \Delta_{V^{1/2}} \xi_{\varepsilon,\delta}^{-1} \zeta_0 \leq q_{1-\alpha,B}\bigr) \\[.5em]  \nonumber
	\leq ~ & \mP\Biggl( \sqrt{V_\mathrm{MMD}} \leq \sqrt{b_{n,m,\blambda} + C_{\tau,\beta,M,d} \frac{\log(1/\alpha)}{n\sqrt{\lambda_1 \cdots \lambda_d}}} + C_{\beta} \Delta_{V^{1/2}} \xi_{\varepsilon,\delta}^{-1} \log(1/\alpha) \Biggr) + 6\beta/ 8 \\[.5em]  \nonumber
	= ~ & \mP\Biggl( V_\mathrm{MMD} \leq b_{n,m,\blambda} + C_{\tau,\beta,M,d} \frac{\log(1/\alpha)}{n\sqrt{\lambda_1 \cdots \lambda_d}} + C_{\beta}^2 \Delta_{V^{1/2}}^2 \xi_{\varepsilon,\delta}^{-2} \log^2(1/\alpha) \\[.5em] \label{Eq: upper bound for type II error}
	& ~~~~~~~~ + 2 \sqrt{b_{n,m,\blambda} + C_{\tau,\beta,M,d} \frac{\log(1/\alpha)}{n\sqrt{\lambda_1 \cdots \lambda_d}}} C_{\beta} \Delta_{V^{1/2}} \xi_{\varepsilon,\delta}^{-1} \log(1/\alpha) \Biggr) + 6\beta/ 8. 
\end{align}

\vskip 1em

\noindent \textbf{Connecting $V_\mathrm{MMD}$ with $U_\mathrm{MMD}$.} In order to express the condition for type II error control in terms of the $L_2$ distance, we make use of the existing results for the U-statistic in \cite{schrab2021mmd}. For simplicity, we use the notation $C_1, C_2, \ldots$ to represent constants that may depend on $\tau, \beta, s, R, M, d$. The specific values of these constants may vary in different places. To proceed, let us make an observation that 
\begin{align*}
	\mE\biggl[ \frac{1}{n^2(n-1)} \sum_{(i,j) \in \mathbf{i}_2^n} k_{\blambda}(Y_i,Y_j) \biggr] ~=~& \frac{1}{n} \mE \bigl[ k_{\blambda}(Y_1,Y_2) \bigr] = \frac{1}{n} \int \int p(y_1)p(y_2) k_{\blambda}(y_1,y_2) dy_1 dy_2 \\[.5em]
	\leq ~ & \frac{\|p\|_{L_\infty}}{n}
\end{align*}
and 
\begin{align} \label{Eq: upper bound for the variance of U-statistic}
	\mV\biggl[ \frac{1}{n^2(n-1)} \sum_{(i,j) \in \mathbf{i}_2^n} k_{\blambda}(Y_i,Y_j) \biggr]  \leq   \frac{C_1}{n^3} \mE[k_{\blambda}^2(Y_1,Y_2)] \leq  \frac{C_2}{n^3\lambda_1 \dots \lambda_d}.
\end{align}
To explain the last display, note that the U-statistic 
\begin{align*}
	\frac{1}{n(n-1)} \sum_{(i,j) \in \mathbf{i}_2^n} k_{\blambda}(Y_i,Y_j) 
\end{align*}
achieves the minimum variance among all unbiased estimators of $k_{\blambda}(Y_1,Y_2)$, including a linear estimator given as
\begin{align*}
	\frac{1}{\floor{n/2}} \sum_{i=1}^{\floor{n/2}} k_{\blambda}(Y_{2i-1},Y_{2i}).
\end{align*}
The variance of the above linear estimator is bounded by $n^{-1} \mE[k_{\blambda}^2(Y_1,Y_2)]$, up to a constant factor, thereby the inequality~\eqref{Eq: upper bound for the variance of U-statistic} holds. A similar calculation shows that 
\begin{align*}
	& \mE\biggl[ \frac{1}{m^2(m-1)} \sum_{(i,j) \in \mathbf{i}_2^m} k_{\blambda}(Z_i,Z_j) \biggr] \leq \frac{\|p\|_{L_\infty}}{n}\quad \text{and} \\[.5em]
	& \mV\biggl[ \frac{1}{m^2(m-1)} \sum_{(i,j) \in \mathbf{i}_2^m} k_{\blambda}(Z_i,Z_j) \biggr] \leq \frac{C_3}{n^3 \lambda_1 \dots \lambda_d}.
\end{align*}
Based on these results combined with Chebyshev's inequality, we have the following statement, which holds with probability at least $1-\beta/8$,
\begin{align*}
	& V_\mathrm{MMD} - U_\mathrm{MMD} + U_\mathrm{MMD} \\[.5em]
	=~ & b_{n,m,\blambda}  - \frac{1}{n^2(n-1)} \sum_{(i,j) \in \mathbf{i}_2^n} k_{\blambda}(Y_i,Y_j) - \frac{1}{m^2(m-1)} \sum_{(i,j) \in \mathbf{i}_2^m} k_{\blambda}(Z_i,Z_j) + U_\mathrm{MMD} \\[.5em]
	\geq ~ &  b_{n,m,\blambda}  - \frac{C_4}{n} - \frac{C_5}{n^{3/2}\sqrt{\lambda_1 \dots \lambda_d}} + U_\mathrm{MMD}.
\end{align*}
Hence, continuing from \eqref{Eq: upper bound for type II error}, the type II error is upper bounded by 
\begin{equation}
\begin{aligned} \label{Eq: upper bound for type II error 2}
	&\mP\Biggl( U_\mathrm{MMD}  \leq \frac{C_1}{n} + \frac{C_2}{n^{3/2}\sqrt{\lambda_1 \cdots \lambda_d}} + \frac{C_3 \log(1/\alpha)}{n\sqrt{\lambda_1 \cdots \lambda_d}} + \frac{C_4 \Delta_{V^{1/2}}^2 \log^2(1/\alpha)}{\xi_{\varepsilon,\delta}^{-2}} \\[.5em]
	& ~~~~~~~~~~~~~~~ + 2 C_5 \sqrt{b_{n,m,\blambda} + \frac{C_6 \log(1/\alpha)}{n\sqrt{\lambda_1 \cdots \lambda_d}}} \frac{\Delta_{V^{1/2}} \log(1/\alpha)}{\xi_{\varepsilon,\delta}} \Biggr) + 7\beta/8. 
\end{aligned}
\end{equation}
Recall that we assume $\alpha \in (0,e^{-1})$ and $\lambda_1\cdots \lambda_d \leq 1$. Moreover, for the Gaussian kernel, \Cref{Lemma: Sensitivity of MMD} shows that the sensitivity of $\sqrt{V_\mathrm{MMD}}$ satisfies
\begin{align*}
	\Delta_{V^{1/2}} \leq \frac{C_6}{n\sqrt{\lambda_1 \dots \lambda_d}}.
\end{align*}
Using these conditions along with the inequality~\eqref{Eq: definition of b_{n,m,blamda}} for $b_{n,m\blambda}$, the previous bound~\eqref{Eq: upper bound for type II error 2} can be upper bounded as
\begin{align*}
	&\mP\Biggl( U_\mathrm{MMD} \leq  \frac{C_1 \log(1/\alpha)}{n\sqrt{\lambda_1 \cdots \lambda_d}} + \frac{C_2 \log^2(1/\alpha)}{n^2 \lambda_1 \cdots \lambda_d \xi_{\varepsilon,\delta}^{2}} \\[.5em]
	& ~~~~~~~~~~~~~~~ + \frac{C_3 \log(1/\alpha)}{n^{3/2} \lambda_1 \cdots \lambda_d \xi_{\varepsilon,\delta}} +  \frac{C_4\log^{3/2}(1/\alpha)}{n^{3/2}(\lambda_1\cdots\lambda_d)^{3/4}\xi_{\varepsilon,\delta}} \Biggr) + 7\beta/8. 
\end{align*}

\vskip 1em

\noindent \textbf{Condition in terms of $L_2$ distance.} As shown in \citet[][Lemma 2]{schrab2021mmd}, a sufficient condition for the first probability term in the above display to be less than $\beta/8$ is
\begin{align*}
	\mathrm{MMD}_{k_{\blambda}}^2 ~\geq~ & C_5 \sqrt{\mV[U_\mathrm{MMD}]} +  \frac{C_1 \log(1/\alpha)}{n\sqrt{\lambda_1 \cdots \lambda_d}} + \frac{C_2 \log^2(1/\alpha)}{n^2 \lambda_1 \cdots \lambda_d \xi_{\varepsilon,\delta}^{2}} \\[.5em]
	& +  \frac{C_3 \log(1/\alpha)}{n^{3/2} \lambda_1 \cdots \lambda_d \xi_{\varepsilon,\delta}} +  \frac{C_4\log^{3/2}(1/\alpha)}{n^{3/2}(\lambda_1\cdots\lambda_d)^{3/4}\xi_{\varepsilon,\delta}}.
\end{align*}
Moreover, writing the difference of two densities as $\psi = p-q$ and the convolution of $\psi$ and $k_{\blambda}$ as $\psi \ast k_{\blambda}$, the proof of \citet[][Theorem 5]{schrab2021mmd} yields that the previous condition is implied by
\begin{align*}
	\|\psi\|_{L_2}^2 ~\geq ~  \|\psi - \psi \ast \varphi_\lambda\|_{L_2}^2 + C_6 & \Biggl\{ \frac{\log(1/\alpha)}{n\sqrt{\lambda_1 \cdots \lambda_d}} + \frac{\log^2(1/\alpha)}{n^2 \lambda_1 \cdots \lambda_d \xi_{\varepsilon,\delta}^{2}} \\[.5em]
	& +  \frac{ \log(1/\alpha)}{n^{3/2} \lambda_1 \cdots \lambda_d \xi_{\varepsilon,\delta}} +  \frac{\log^{3/2}(1/\alpha)}{n^{3/2}(\lambda_1\cdots\lambda_d)^{3/4}\xi_{\varepsilon,\delta}} \Biggr\}.
\end{align*}
Lastly, the proof of \citet[][Theorem 6]{schrab2021mmd} yields that over the Sobolev ball, a sufficient condition for the previous inequality is
\begin{align*}
	\|\psi\|_{L_2}^2 ~ \geq ~  C_7 \Biggl\{ \sum_{i=1}^d \lambda_i^{2s} & +  \frac{\log(1/\alpha)}{n\sqrt{\lambda_1 \cdots \lambda_d}} + \frac{\log^2(1/\alpha)}{n^2 \lambda_1 \cdots \lambda_d \xi_{\varepsilon,\delta}^{2}} \\[.5em]
	& +  \frac{ \log(1/\alpha)}{n^{3/2} \lambda_1 \cdots \lambda_d \xi_{\varepsilon,\delta}} +  \frac{\log^{3/2}(1/\alpha)}{n^{3/2}(\lambda_1\cdots\lambda_d)^{3/4}\xi_{\varepsilon,\delta}} \Biggr\}
\end{align*}
as claimed.

{
\subsection{Proof of \Cref{Theorem: High-privacy minimax separation over L2}}
\label{Section: Proof of high-privacy minimax separation over L2}

We construct a localized pair of alternatives and apply the same coupling argument used for the privacy-regime lower bound in the proof of \Cref{Theorem: Minimax separation in MMD}.

\paragraph{Construction.}
Let
\[
	\ell=\left(\frac{2}{M}\right)^{1/d},
	\qquad
	Q_0=[0,\ell]^d,
	\qquad
	q_0(x)=\frac{M}{2}\mathds{1}_{Q_0}(x).
\]
Choose a nonzero function $\psi\in C_c^\infty((-1,1)^d)$ such that
\[
	\int_{\mathbb{R}^d}\psi(u)\,\dd u=0,
	\qquad
	\|\psi\|_{L_\infty}\leq1.
\]
Set $x_0=(\ell/2,\ldots,\ell/2)^\top$, $h_0=\min\{1,\ell/4\}$, and, for $h\in(0,h_0]$, define
\[
	\psi_h(x)=\psi\left(\frac{x-x_0}{h}\right),
	\qquad
	a_h=a_0h^{s-d/2},
\]
where $a_0>0$ is chosen below. Consider the density pair
\[
	p_h=q_0+a_h\psi_h,
	\qquad
	q_h=q_0.
\]
Since $s\geq d/2$ and $h\leq1$, we have $a_h\leq a_0$. Moreover, $\supp(\psi_h)\subset Q_0$. Thus, choosing $a_0\leq M/4$ ensures
\[
	0\leq p_h\leq M,
	\qquad
	0\leq q_h\leq M,
\]
and the zero-integral property of $\psi$ ensures that both densities integrate to one.

The Fourier transform of the perturbation is
\[
	\widehat{p_h-q_h}(w)
	=
	a_hh^d e^{-ix_0^\top w}\widehat{\psi}(hw).
\]
Consequently, after the change of variables $u=hw$,
\begin{align*}
	\int_{\mathbb{R}^d}\|w\|_2^{2s}
	\bigl|\widehat{p_h-q_h}(w)\bigr|^2\,\dd w
	&=
	a_h^2h^{d-2s}
	\int_{\mathbb{R}^d}\|u\|_2^{2s}
	|\widehat{\psi}(u)|^2\,\dd u \\
	&=
	a_0^2
	\int_{\mathbb{R}^d}\|u\|_2^{2s}
	|\widehat{\psi}(u)|^2\,\dd u.
\end{align*}
Because $\psi$ is smooth and compactly supported, the last integral is finite. We can therefore choose $a_0$ small enough that the preceding display is at most $(2\pi)^dR^2$. It follows that $(P_h,Q_h)$ satisfies the smoothness and boundedness requirements of the target class.

\paragraph{Separation and coupling cost.}
The rescaling of $\psi_h$ gives
\begin{align}
	\|p_h-q_h\|_{L_2}
	&=
	a_0\|\psi\|_{L_2}h^s
	\eqqcolon \rho_h,
	\label{Eq: lower-bound L2 separation}\\
	\TV(P_h,Q_h)
	&=
	\frac{a_0}{2}\|\psi\|_{L_1}h^{s+d/2}.
	\label{Eq: lower-bound TV separation}
\end{align}
Consider the simple null and alternative
\[
	H_0:(P,Q)=(Q_h,Q_h),
	\qquad
	H_1:(P,Q)=(P_h,Q_h).
\]
Only the sample of size $n$ differs. Couple its coordinates independently by maximal coupling and couple the sample of size $m$ identically. The expected Hamming distance between the resulting datasets is
\begin{align}
	D_h
	=
	n\TV(P_h,Q_h)
	=
	\frac{a_0}{2}\|\psi\|_{L_1}
	nh^{s+d/2}.
	\label{Eq: lower-bound coupling cost}
\end{align}

\paragraph{Private testing lower bound.}
Let $\phi\in\Phi_{\alpha,\varepsilon,\delta}$. The coupling, Markov, and group-privacy argument used in the proof of \Cref{Theorem: Minimax separation in MMD} gives
\begin{align}
	\mE_{P_h,Q_h}[1-\phi]
	\geq
	e^{-10D_h\varepsilon}(0.9-\alpha)-10D_h\delta.
	\label{Eq: private Le Cam L2}
\end{align}
The condition $\alpha+\beta<0.9$ implies that there exists $c_\star>0$, depending only on $\alpha$ and $\beta$, such that
\[
	e^{-10c_\star}(0.9-\alpha)-10c_\star>\beta.
\]
Since $\varepsilon\leq\xi_{\varepsilon,\delta}$ and $\delta\leq\xi_{\varepsilon,\delta}$, \eqref{Eq: private Le Cam L2} is larger than $\beta$ whenever
\[
	D_h\xi_{\varepsilon,\delta}\leq c_\star.
\]
By \eqref{Eq: lower-bound L2 separation} and \eqref{Eq: lower-bound coupling cost},
\[
	D_h\xi_{\varepsilon,\delta}
	=
	C_\psi n\xi_{\varepsilon,\delta}
	\rho_h^{\,1+d/(2s)}
\]
for a constant $C_\psi>0$ that also depends on $a_0$ and $s$. Hence testing is impossible at type II error $\beta$ whenever
\[
	\rho_h
	\leq
	c\,(n\xi_{\varepsilon,\delta})^{-2s/(2s+d)}
\]
for a sufficiently small constant $c>0$.

It remains to make the construction cover every separation below the claimed bound. Choose $c_{d,s,R,M,\alpha,\beta}>0$ small enough that
\[
	c_{d,s,R,M,\alpha,\beta}
	\leq
	a_0\|\psi\|_{L_2}h_0^s
\]
and that the preceding privacy condition holds. For any
\[
	0<\rho
	\leq
	c_{d,s,R,M,\alpha,\beta}
	\min\Bigl\{1,\,
	(n\xi_{\varepsilon,\delta})^{-2s/(2s+d)}
	\Bigr\},
\]
take
\[
	h=
	\left(\frac{\rho}{a_0\|\psi\|_{L_2}}\right)^{1/s}.
\]
Then $h\leq h_0$, $\rho_h=\rho$, and $D_h\xi_{\varepsilon,\delta}\leq c_\star$. Thus every $(\varepsilon,\delta)$-DP level-$\alpha$ test has type II error larger than $\beta$ at an alternative in $\mathcal{P}_{L_2}^s(\rho)$. The definition in \eqref{Eq: minimax L2 separation} proves the result.
}

\subsection{Proof of \Cref{Lemma: Global sensitivity of U_MMD}} \label{Section: Proof of Lemma: Global sensitivity of U_MMD}
Let $U_{\mathrm{MMD}}'$ be the U-statistic similarly defined as $U_{\mathrm{MMD}}$ by replacing $Y_1$ with $Y_1'$. Then it can be seen that 
\begin{align*}
	|U_{\mathrm{MMD}} - U_{\mathrm{MMD}}'| ~=~& \bigg| \frac{2}{n(n-1)}\sum_{j=2}^n \big\{ k(Y_1,Y_j) - k(Y_1',Y_j) \big\} - \frac{2}{nm} \sum_{j=1}^m \big\{ k(Y_1,Z_j) - k(Y_1',Z_j) \big\} \bigg| \\[.5em]
	\leq ~ & \frac{8K}{n}.
\end{align*}
Similarly, the difference is bounded by $8K/n$ for other neighboring datasets, and thus the global sensitivity of $U_{\mathrm{MMD}}$ is bounded by $8K/n$. Let $\epsilon$ be a number between $(0, K)$. For a translation invariant kernel with non-empty level sets, we can ensure the existence of an instance where $Y_1=\ldots=Y_n$ and $Z_1=\ldots=Z_m$ such that $k(Y_1,Z_1) = \epsilon_\star$ for some $\epsilon_\star \in [0, \epsilon]$. Moreover, by taking $Y_1' = Z_1$, we have $k(Y_1,Y_1') = \epsilon_\star$ and $k(Y_1',Z_1) = K$. Under this setting, the difference between $U_{\mathrm{MMD}}$ and $U_{\mathrm{MMD}}'$ becomes
\begin{align*}
	|U_{\mathrm{MMD}} - U_{\mathrm{MMD}}'| = \frac{4(K-\epsilon_\star)}{n}. 
\end{align*}
Since $\epsilon$ (and so $\epsilon_\star$) can be arbitrarily small, we see that the global sensitivity defined with the supremum is lower bounded by $4K/n$. This concludes the claim of \Cref{Lemma: Global sensitivity of U_MMD}.

\subsection{Proof of \Cref{Theorem: Suboptimality of U-MMD}} \label{Section: Proof of Theorem: Suboptimality of U-MMD}
For notational convenience, let us write $U_{\mathrm{MMD}}(\mathcal{X}_{n+m}^{\bpi_0}) = U_{\mathrm{MMD}}$ and $U_{\mathrm{MMD}}(\mathcal{X}_{n+m}^{\bpi_i}) = U_{\bpi_i,\mathrm{MMD}}$ for $i  \in [B]$. For $\alpha > 1/(B+1)$, \Cref{Lemma: Alternative expression} yields that 
\begin{align*}
	\phi_{\mathttt{dpMMD}}^u = \mathds{1}\biggl( U_\mathrm{MMD} + \frac{2c_{m,n} K}{n\xi_{\varepsilon,\delta}} \zeta_0 > r_{1-\alpha_\star} \biggr),
\end{align*}
where $\alpha_\star = \frac{B+1}{B} \alpha - \frac{1}{B}$ and $r_{1-\alpha_\star}$ is the $1-\alpha_\star$ quantile of $\{U_{\bpi_i,\mathrm{MMD}} +\frac{2c_{m,n} K}{n\xi_{\varepsilon,\delta}} \zeta_i\}_{i=1}^B$. Let us denote the V-statistic of MMD as
\begin{align*}
	V_{\mathrm{MMD}} = \frac{1}{n^2} \sum_{i,j =1}^n k(Y_i,Y_j) + \frac{1}{m^2} \sum_{i,j=1}^m k(Z_i,Z_j) - \frac{2}{nm} \sum_{i=1}^n \sum_{j=1}^m k(Y_i,Z_j),
\end{align*}
which is greater than or equal to $U_\mathrm{MMD}$ since
\begin{align*}
	& \frac{1}{n^2} \sum_{i,j =1}^n k(Y_i,Y_j)  \geq \frac{1}{n(n-1)} \sum_{(i,j) \in \mathbf{i}_2^n} k(Y_i,Y_j) \\[.5em]
	\Longleftrightarrow ~ & \frac{1}{n} \sum_{i=1}^n k(Y_i,Y_i) = K \geq \frac{1}{n(n-1)} \sum_{(i,j) \in \mathbf{i}_2^n} k(Y_i,Y_j)
\end{align*}
and similarly
\begin{align*}
	& \frac{1}{m^2} \sum_{i,j=1}^m k(Z_i,Z_j)  \geq \frac{1}{m(m-1)} \sum_{(i,j) \in \mathbf{i}_2^m} k(Z_i,Z_j) \\[.5em]
	\Longleftrightarrow ~ & \frac{1}{m} \sum_{i=1}^m k(Z_i,Z_i) = K \geq \frac{1}{m(m-1)} \sum_{(i,j) \in \mathbf{i}_2^m} k(Z_i,Z_j).
\end{align*}
Moreover, since $V_\mathrm{MMD} \geq 0$, the U-statistic is lower bounded by $U_\mathrm{MMD} \geq U_\mathrm{MMD} - V_\mathrm{MMD}$. Consider a lower bound for $U_\mathrm{MMD}-V_\mathrm{MMD}$ as
\begin{align*}
	U_\mathrm{MMD} - V_\mathrm{MMD} ~=~&  \frac{1}{n(n-1)} \sum_{(i,j) \in \mathbf{i}_2^n} k(Y_i,Y_j) -  \frac{1}{n^2} \sum_{i,j =1}^n k(Y_i,Y_j) \\[.5em]
	+ ~ & \frac{1}{m(m-1)} \sum_{(i,j) \in \mathbf{i}_2^m} k(Z_i,Z_j) -  \frac{1}{m^2} \sum_{i,j=1}^m k(Z_i,Z_j) \\[.5em]
	= ~ & \frac{1}{n} \biggl[ \frac{1}{n(n-1)} \sum_{(i,j) \in \mathbf{i}_2^n} k(Y_i,Y_j)  - \frac{1}{n} \sum_{i=1}^n k(Y_i,Y_i) \biggr] \\[.5em]
	+ ~ & \frac{1}{m} \biggl[ \frac{1}{m(m-1)} \sum_{(i,j) \in \mathbf{i}_2^m} k(Z_i,Z_j)  - \frac{1}{m} \sum_{i=1}^m k(Z_i,Z_i) \biggr] \\[.5em]
	\geq ~ & -\frac{K}{n} - \frac{K}{m}.
\end{align*}
Thereby, $U_\mathrm{MMD} \geq - \frac{K}{n} - \frac{K}{m}$. We use this observation to lower bound the critical value $r_{1-\alpha_\star}$ as
\begin{align*}
	r_{1-\alpha_\star} \geq \frac{2c_{m,n} K}{n\xi_{\varepsilon,\delta}}  \mathrm{Quantile}_{1-\alpha_\star}\big\{\zeta_1,\ldots, \zeta_B  \big\} - \frac{K}{n} - \frac{K}{m},
\end{align*}
where $\mathrm{Quantile}_{1-\alpha_\star}\{\zeta_1,\ldots, \zeta_B\}$ denotes the $1-\alpha_\star$ quantile of $\zeta_1,\ldots, \zeta_B$. Putting pieces together yields
\begin{align*}
	& \phi_{\mathttt{dpMMD}}^u \\[.5em] 
	\leq ~ & \mathds{1} \biggl( V_\mathrm{MMD} + \frac{2c_{m,n} K}{n\xi_{\varepsilon,\delta}}  \zeta_0 > \frac{2c_{m,n} K}{n\xi_{\varepsilon,\delta}}  \mathrm{Quantile}_{1-\alpha_\star}\big\{\zeta_1,\ldots, \zeta_B  \big\} - \frac{K}{n} - \frac{K}{m} \biggr) \\[.5em]
	\leq  ~ & \mathds{1} \biggl( \frac{n \xi_{\varepsilon,\delta}}{2c_{m,n} K} V_\mathrm{MMD} + \zeta_0 >  \mathrm{Quantile}_{1-\alpha_\star}\big\{\zeta_1,\ldots, \zeta_B  \big\} -\frac{\xi_{\varepsilon,\delta}}{c_{m,n}} \biggr).
\end{align*}
Furthermore, note that the statistic $V_\mathrm{MMD}$ is the square of the empirical MMD, which satisfies
\begin{align*}
	V_\mathrm{MMD} = \big\{ \mathrm{MMD}_k(P,Q)  + R_{m,n} \}^2,
\end{align*}
where $R_{m,n} = O_P(n^{-1/2})$ due to \citet[][Theorem 7]{gretton2012kernel}, recalled in \Cref{Lemma: Concentration for MMD}.  

Let us consider $P_0,Q_0 \in \mathcal{P}_{\mathbb{S}}$ such that $\mathrm{MMD}_k(P_0,Q_0) = \varrho_0$ in the theorem statement. Let $P = P_0$, and let $Q$ be a mixture distribution $Q = w P_0 + (1-w) Q_0$. Since $\mathcal{P}_{\mathbb{S}}$ is a convex set, $Q$ belongs to $\mathcal{P}_{\mathbb{S}}$, and it can be seen that 
\begin{align*}
	\mathrm{MMD}_k(P,Q) = (1-w) \varrho_0. 
\end{align*} 
Now for $\gamma \in (1/2,1)$, we set $w = 1 - (n\xi_{\varepsilon,\delta})^{-\gamma}$. Since $\xi_{\varepsilon,\delta} \asymp n^{-1/2-r}$ with $r \in (0,1/2)$, we can ensure that $w \in (0,1)$ and $\mathrm{MMD}_k(P,Q) = (n\xi_{\varepsilon,\delta})^{-\gamma} \varrho_0$ for sufficiently large $n$. Moreover, this pair of distributions $(P,Q)$ belongs to $\mathcal{P}_{\mathrm{MMD}_k}\!(\rho)$ for $\rho$ as in \eqref{Eq: separation MMD rate for U-statistics} since 
\begin{align*}
	\mathrm{MMD}_k(P,Q) = \frac{\varrho_0}{\big(n\xi_{\varepsilon,\delta}\big)^{\gamma}} ~\geq~ \rho ~\asymp~ \frac{\log(n)}{n\xi_{\varepsilon,\delta}} \quad &\Longleftrightarrow \quad  \big(n\xi_{\varepsilon,\delta}\big)^{1-\gamma} ~\gtrsim~ \log(n) \\[.5em]
	&\Longleftrightarrow \quad   n^{\frac{(1-2r)(1-\gamma)}{2}} ~\gtrsim~ \log(n).
\end{align*}
Moreover, as $n\xi_{\varepsilon,\delta} \rightarrow \infty$ and $\xi_{\varepsilon,\delta} \rightarrow 0$, 
\begin{align*}
	\frac{n\xi_{\varepsilon,\delta}}{2c_{m,n} K} V_\mathrm{MMD}  = \frac{1}{2c_{m,n} K}  \big\{ \sqrt{n\xi_{\varepsilon,\delta}} \mathrm{MMD}_k(P,Q)  + \sqrt{n\xi_{\varepsilon,\delta}} R_{m,n} \}^2 = o_P(1),
\end{align*}
where we use the fact that $c_{m,n} \in [4,8]$. Having this observation in place, for any fixed $t>0$, we can take $N_t>0$ such that for all $n \geq N_t$, it holds that 
\begin{align*}
	\mE[\phi_{\mathttt{dpMMD}}^u] ~\leq~& \mP \Bigl(\zeta_0 > \mathrm{Quantile}_{1-\alpha_\star}\big\{\zeta_1,\ldots, \zeta_B  \big\}  - \frac{\xi_{\varepsilon,\delta}}{c_{m,n}} - \frac{n\xi_{\varepsilon,\delta}}{2c_{m,n} K} V_\mathrm{MMD},  \\[.5em]
	& \hskip 5em \Big| \frac{\xi_{\varepsilon,\delta}}{c_{m,n}} + \frac{n\xi_{\varepsilon,\delta}}{2c_{m,n} K} V_\mathrm{MMD} \Big| < t \Bigr) + \mP\Bigl( \Big| \frac{\xi_{\varepsilon,\delta}}{c_{m,n}} + \frac{n\xi_{\varepsilon,\delta}}{2c_{m,n} K} V_\mathrm{MMD} \Big| \geq t \Bigr) \\[.5em]
	\leq ~ & \mP \Bigl(\zeta_0 > \mathrm{Quantile}_{1-\alpha_\star}\big\{\zeta_1,\ldots, \zeta_B  \big\}  - t \Bigr) + t \\[.5em]
	= ~ & \mE \Bigl[ \mP \Bigl(\zeta_0 > \mathrm{Quantile}_{1-\alpha_\star}\big\{\zeta_1,\ldots, \zeta_B  \big\}  - t \,\Big|\, \zeta_1,\ldots,\zeta_B \Bigr)  \Bigr]+ t \\[.5em]
	\leq ~ & \mP \Bigl(\zeta_0 > \mathrm{Quantile}_{1-\alpha_\star}\big\{\zeta_1,\ldots, \zeta_B  \big\} \Bigr) + t \|f_{\zeta}\|_{L_\infty} + t, 
\end{align*}
where $f_\zeta$ denotes the density function of $\mathsf{Laplace}(0,1)$ and the last inequality uses the fact that
\begin{align*}
	\big|\mP(\zeta > a + b) - \mP(\zeta > b)\big| \leq \|f_{\zeta}\|_{L_\infty} |a| \quad \text{for all $a,b \in \mathbb{R}$.}
\end{align*}
Since $\|f_{\zeta}\|_{L_\infty} \leq \frac{1}{2}$ and
\begin{align*}
	\mP \Bigl(\zeta_0 > \mathrm{Quantile}_{1-\alpha_\star}\big\{\zeta_1,\ldots, \zeta_B  \big\} \Bigr)  = \mP \Bigl(\zeta_0 > \mathrm{Quantile}_{1-\alpha}\big\{\zeta_0, \zeta_1,\ldots, \zeta_B \big\} \Bigr) \leq \alpha, 
\end{align*}
we have 
\begin{align*}
	\limsup_{n \rightarrow \infty} \inf_{(P,Q) \in \mathcal{P}_{\mathrm{MMD}_k}\!(\rho)} \mE_{P,Q}[\phi_{\mathttt{dpMMD}}^u]  \leq \alpha + \frac{3t}{2}.
\end{align*}
The result follows as $t$ can be made arbitrarily small.

\subsection{Proof of \Cref{Theorem: Minimum separation of U-stat over L2}} \label{Section: Proof of Theorem: Minimum separation of U-stat over L2}
Let us denote by $C_1,C_2,\ldots$ constants that may depend on $\tau,\beta,s,R,M,d$. Following the proofs of \Cref{Theorem: General uniform power condition} and \Cref{Theorem: Minimax Separation over L2} along with the sensitivity result for $U_{\mathrm{MMD}}$ in \Cref{Lemma: Global sensitivity of U_MMD}, we may arrive at the point where the type II error of $\phi_{\mathttt{dpMMD}}^u$ is upper bounded as 
\begin{align*}
	\mE[1 - \phi_{\mathttt{dpMMD}}^u] \leq \mP \biggl( U_{\mathrm{MMD}} \leq \frac{C_1 \log(1/\alpha)}{n\sqrt{\lambda_1 \cdots \lambda_d}} + \frac{C_2\log(1/\alpha)}{n \lambda_1 \cdots \lambda_d \xi_{\varepsilon,\delta}} \biggr) + 7\beta/8.
\end{align*}
We then use the proofs of \citet[][Theorem 5]{schrab2021mmd} and \citet[][Theorem 6]{schrab2021mmd}, and show that the probability term in the above display is less than or equal to $\beta/8$ once
\begin{align*}
	\|p-q\|_{L_2}^2 ~ \geq ~  C_3 \Biggl\{ \sum_{i=1}^d \lambda_i^{2s} & +  \frac{\log(1/\alpha)}{n\sqrt{\lambda_1 \cdots \lambda_d}} + \frac{\log(1/\alpha)}{n\lambda_1 \cdots \lambda_d \xi_{\varepsilon,\delta}} \Biggr\}.
\end{align*}
This proves \Cref{Theorem: Minimum separation of U-stat over L2}.

\section{Proofs for \Cref{Section: Additional Results}} \label{Section: Proofs of the additional results}
This section collects the proofs of the results in \Cref{Section: Additional Results}.

\subsection{Proof of \Cref{Proposition: Asymptotic null distributions}} \label{Section: Proof of Proposition: Asymptotic null distributions}
Focusing on MMD, the square of the empirical MMD satisfies 
\begin{align*}
	(n+m) \widehat{\mathrm{MMD}}^2(\mathcal{X}_{n+m}) \convD  \sum_{i=1}^\infty \lambda_i Z_i^2,
\end{align*}
which can be seen by the standard asymptotic theory of V-statistics. See \citet[][Proposition 9]{fernandez2022general}. Moreover, note that the empirical MMD and the Laplace noise are independent. Therefore, the continuous mapping theorem along with Slutsky's theorem proves the result when $\sqrt{n+m} \sigma \rightarrow \eta \in [0,\infty)$. On the other hand, when $\sqrt{n+m} \sigma \rightarrow  \infty$, it holds that 
\begin{align*}
	\sigma^{-1} M_{\mathrm{MMD}} = \underbrace{\frac{1}{\sqrt{n+m}\, \sigma}}_{=\,o_P(1)} \times \underbrace{\sqrt{n+m}\,\widehat{\mathrm{MMD}}(\mathcal{X}_{n+m})}_{=\,O_P(1)} +\, \zeta = o_P(1) + \zeta.
\end{align*}
This proves the results for MMD. The proof for HSIC is completely analogous and thereby omitted, which can be derived by using \citet[][Theorem 1]{zhang2018large} in place of \citet[][Proposition 9]{fernandez2022general}.

\subsection{Proof of \Cref{Lemma: General conditions for consistency}} \label{Section: Proof of Lemma: General conditions for consistency}
We first prove the claim by considering two scenarios: (i) $B_n$ is fixed in $n$ and (ii) $B_n$ increases with $n$, and prove the consistency result. We then turn to the general case of $B_n$ and complete the proof building on the prior result. Throughout this proof, we often omit the dependence of $n$ on $B_n$ and write it as $B$ for simplicity.

\subsubsection{Simple cases} \label{Section: Preliminaries}
\paragraph{Fixed $B$.} Assume that $B$ is fixed and $B+1 > \alpha^{-1}$ or equivalently $\alpha(B+1) \geq 1$. As mentioned in the main text, we simply exploit the union bound to prove the claim for fixed $B$. Concretely,
\begin{align*}
	\mP \biggl[ \frac{1}{B+1} \biggl( \sum_{i=1}^{B} \mathds{1}(W_{0,n} \leq W_{i,n}) + 1 \biggr) \leq \alpha \biggr] ~=~&  \mP \biggl[  \sum_{i=1}^{B} \mathds{1}(W_{0,n} \leq W_{i,n})  \leq \alpha (B+1) - 1 \biggr] \\[.5em]
	\overset{\mathrm{(i)}}{\geq} ~ & \mP \biggl[  \sum_{i=1}^{B} \mathds{1}(W_{0,n} \leq W_{i,n})  \leq  0 \biggr]  \\[.5em]
	= ~ & \mP \Bigl[W_{0,n} - \max_{i \in [B]} W_{i,n}  > 0 \Bigr] \\[.5em]
	\overset{\mathrm{(ii)}}{\geq} ~ & 1 - \sum_{i=1}^B \mP(W_{0,n} \leq W_{i,n}),
\end{align*}
where step~(i) uses condition~$\alpha(B+1) > 1$ and step~(ii) uses de Morgan's law and the union bound. By the condition that $\lim_{n \rightarrow \infty} \mP(W_{0,n} \leq W_{1,n}) = 0$ and since $\{W_{0,n} - W_{i,n}\}_{i=2}^B$ are identically distributed as $W_{0,n} - W_{1,n}$, the lower bound converges to one for fixed $B$ and thus the consistency result follows. 

\smallskip 

\paragraph{Increasing $B$.} Next we consider the case where $B$ increases to infinity with $n$. For some $\eta > 0$ (specified later), define an event $\mathcal{A}$ as
\begin{align*}
	\mathcal{A} = \bigg\{ \bigg| \frac{1}{B} \sum_{i=1}^B \mathds{1}(W_{0,n} \leq W_{i,n}) - \mP(W_{0,n} \leq W_{1,n} \given \mathcal{G} ) \bigg| \leq \sqrt{\frac{1}{2B} \log \left(\frac{2}{\eta}\right)}\bigg\}.
\end{align*}
Since $\mathds{1}(W_{0,n} \leq W_{1,n}), \ldots,  \mathds{1}(W_0 \leq W_{B,n})$ are i.i.d.~random variables conditional on the sigma field $\mathcal{G}$ under the condition of \Cref{Lemma: General conditions for consistency}, Hoeffding's inequality yields that
\begin{align*}
	\mP(\mathcal{A}^c \given \mathcal{G}) \leq \eta. 
\end{align*}
By taking the expectation over $\mathcal{G}$ on both sides, it also holds marginally that $\mP(\mathcal{A}^c) \leq \eta$. Now
\begin{equation}
	\begin{aligned} \label{Eq: step 1 for increasing B}
		& \mP \biggl[ \frac{1}{B+1} \biggl( \sum_{i=1}^{B} \mathds{1}(W_{0,n} \leq W_{i,n}) + 1 \biggr) > \alpha \biggr] \\[.5em]
		=~ & \mP\biggl[ \frac{1}{B} \sum_{i=1}^B \mathds{1}(W_{0,n} \leq W_{i,n}) > \frac{B+1}{B}\biggl(\alpha - \frac{1}{B+1}\biggr) \biggr] \\[.5em]
		\leq ~ & \mP\biggl[ \frac{1}{B} \sum_{i=1}^B \mathds{1}(W_{0,n} \leq W_{i,n}) > \frac{B+1}{B}\biggl(\alpha - \frac{1}{B+1}\biggr), \, \mathcal{A} \biggr] + \mP(\mathcal{A}^c) \\[.5em]
		\leq ~ & \mP\biggl[ \mP(W_{0,n} \leq W_{1,n} \given \mathcal{G}) +  \sqrt{\frac{1}{2B} \log \left(\frac{2}{\eta}\right)} >  \frac{B+1}{B}\biggl(\alpha - \frac{1}{B+1}\biggr) \biggr] + \eta.
	\end{aligned}
\end{equation}
By taking $\eta = B^{-1}$, we know $\eta \rightarrow 0$ (since $B \rightarrow \infty$) and 
\begin{align*}
	\frac{B+1}{B}\biggl(\alpha - \frac{1}{B+1}\biggr) -  \sqrt{\frac{1}{2B} \log \left(\frac{2}{\eta}\right)} \rightarrow \alpha \quad \text{as $n \rightarrow \infty$.}
\end{align*}
Moreover our condition guarantees that $\lim_{n \rightarrow \infty} \mP(W_{0,n} \leq W_{1,n}) = 0$. Hence for any $\epsilon > 0$, we can take $N$ such that  statement~(i) $\eta < \epsilon/2$ and statement~(ii)
\begin{equation}
	\begin{aligned} \label{Eq: step 2 for increasing B}
		& \mP\biggl[ \mP(W_{0,n} \leq W_{1,n} \given \mathcal{G}) +  \sqrt{\frac{1}{2B} \log \left(\frac{2}{\eta}\right)} >  \frac{B+1}{B}\biggl(\alpha - \frac{1}{B+1}\biggr)\bigg] \\[.5em]
		\leq ~ & \mP\biggl[ \mP(W_{0,n} \leq W_{1,n} \given \mathcal{G}) > \frac{\alpha}{2} \biggr] \\[.5em]
		\overset{(\dagger)}{\leq} ~ & \frac{2}{\alpha}  \mP(W_{0,n} \leq W_{1,n})  \leq  \frac{\epsilon}{2},
	\end{aligned}
\end{equation}
hold for all $n \geq N$. Here step~($\dagger$) holds by Markov's inequality. Thus for all $n \geq N$ we have
\begin{align*}
	\mP \biggl[ \frac{1}{B_n+1} \biggl( \sum_{i=1}^{B_n} \mathds{1}(W_{0,n} \leq W_{i,n}) + 1 \biggr) > \alpha \biggr]  < \epsilon. 
\end{align*}
Since $\epsilon$ is arbitrary, the above argument proves the claim that 
\begin{align*}
	\lim_{n \rightarrow \infty} \mP \biggl[ \frac{1}{B_n+1} \biggl( \sum_{i=1}^{B_n} \mathds{1}(W_{0,n} \leq W_{i,n}) + 1 \biggr) \leq \alpha \biggr] = 1.
\end{align*}

\subsubsection{Main proof} \label{Section: arbitrary B}

We now deal with an arbitrary sequence of $B_n$ such that $\inf_{n \geq 1} B_n + 1 > \alpha^{-1}$ and prove the claim. First note that 
\begin{align*}
	& \mP \biggl[ \frac{1}{B+1} \biggl( \sum_{i=1}^{B} \mathds{1}(W_{0,n} \leq W_{i,n}) + 1 \biggr) > \alpha \biggr] = \mP\biggl[\frac{1}{B+1} \sum_{i=1}^B \mathds{1}(W_{0,n} \leq W_{i,n}) >  \biggl(\alpha - \frac{1}{B+1}\biggr) \biggr] \\[.5em] 
	\overset{(\mathrm{i})}{\leq} ~ &  \mP\biggl[ \sum_{i=1}^B \mathds{1}(W_{0,n} \leq W_{i,n}) > 0 \biggr] 	=  1 - \mP\biggl[ \sum_{i=1}^B \mathds{1}(W_{0,n} \leq W_{i,n}) = 0 \biggr] \\[.5em]
	= ~ & 1 - \mE \biggl\{ \mP\biggl[  \sum_{i=1}^B \mathds{1}(W_{0,n} \leq W_{i,n}) = 0  \, \bigg| \, \mathcal{G} \biggr] \biggr\} = 1 - \mE\Bigl[ \{\mP(W_{0,n} > W_{1,n} \given \mathcal{G}) \}^B \Bigr] \\[.5em]
	\overset{(\mathrm{ii})}{\leq} ~ & 1 - \{\mP(W_{0,n} > W_{1,n})\}^B,
\end{align*}
where step~(i) uses the condition $\alpha > \frac{1}{B+1}$ and step~(ii) uses Jensen's inequality. On the other hand, from the previous results~\eqref{Eq: step 1 for increasing B} and Markov's inequality with $\eta= B^{-1}$, we have 
\begin{align*}
	\mP \biggl[ \frac{1}{B+1} \biggl( \sum_{i=1}^{B} \mathds{1}(W_{0,n} \leq W_{i,n}) + 1 \biggr) > \alpha \biggr]  \leq \frac{1}{B} + \frac{\mP(W_{0,n} \leq W_{1,n}) + \sqrt{\frac{1}{2B} \log(2B)}}{\frac{B+1}{B} (\alpha - \frac{1}{B+1})}.
\end{align*}
Thus combining the two bounds yields
\begin{align*}
	& \mP \biggl[ \frac{1}{B+1} \biggl( \sum_{i=1}^{B} \mathds{1}(W_{0,n} \leq W_{i,n}) + 1 \biggr) > \alpha \biggr] \\[.5em]
	\leq ~ & \min \Bigg\{ 1 - \{\mP(W_{0,n} > W_{1,n})\}^B, \, \frac{1}{B} + \frac{\mP(W_{0,n} \leq W_{1,n}) + \sqrt{\frac{1}{2B} \log(2B)}}{\frac{B+1}{B} (\alpha - \frac{1}{B+1})} \Bigg\}.
\end{align*}
Interestingly, a careful argument shows that the above upper bound can be made independent of $B$. Concretely, write $p_n \coloneqq  \mP(W_{0,n} > W_{1,n})$ for simplicity, which by assumption tends to $1$ as $n$ tends to infinity, and let $b_n$ be a positive sequence that goes to infinity and satisfies $p_n^{b_n} \rightarrow 1$. For instance, one can take $b_n = \sqrt{-\log (1 - p_n)}$. Using this notation in place, observe 
\begin{align*}
	&\min \Bigg\{ 1 - \{\mP(W_{0,n} > W_{1,n})\}^B, \, \frac{1}{B} + \frac{\mP(W_{0,n} \leq W_{1,n}) + \sqrt{\frac{1}{2B} \log(2B)}}{\frac{B+1}{B} (\alpha - \frac{1}{B+1})} \Bigg\} \\[.5em]
	\leq ~ & \mathds{1}(B < b_n) \big\{ 1 -   p_n^B \big\} + \mathds{1}(B \geq b_n) \Bigg\{ \frac{1}{B} + \frac{1-p_n + \sqrt{\frac{1}{2B} \log(2B)}}{\frac{B+1}{B} (\alpha - \frac{1}{B+1})}  \Bigg\} \\[.5em]
	\leq ~ & 1 - p_n^{b_n} + \Bigg\{ \frac{1}{b_n} + \frac{1-p_n + \sqrt{\frac{1}{2b_n} \log(2b_n)}}{\frac{b_n+1}{b_n} (\alpha - \frac{1}{b_n+1})}  \Bigg\}.
\end{align*}
Now the upper bound does not depend on $B$ and converges to zero, which is ensured by the condition $\lim_{n\rightarrow \infty} p_n = 1$. As a result, we have 
\begin{align*}
	\lim_{n \rightarrow \infty} \mP \biggl[ \frac{1}{B_n+1} \biggl( \sum_{i=1}^{B_n} \mathds{1}(W_{0,n} \leq W_{i,n}) + 1 \biggr) \leq \alpha \biggr] = 1,
\end{align*}
and complete the proof of \Cref{Lemma: General conditions for consistency}.

\subsection{Proof of \Cref{Theorem: Uniform separation for HSIC}}  \label{Section: Proof of Theorem: Uniform separation for HSIC}
The proof of \Cref{Theorem: Uniform separation for HSIC} can follow a similar approach to that of \Cref{Theorem: Uniform separation for MMD} by replacing concentration inequalities for MMD statistics~(\Cref{Lemma: Bobkovs inequality} and \Cref{Lemma: Concentration for MMD}) with the corresponding ones for HSIC statistics~(\Cref{Lemma: Concentration inequality for permuted HSIC} and \Cref{Lemma: Exponential inequality for the empirical HSIC}). Throughout this proof, we denote positive constants that only depend on $K$ and $L$ by $C_1,C_2,C_3,\ldots$ 

\vskip 1em 

\noindent \textbf{Step 1 (Bounding the quantile).} Let us remind that $M_i$ is defined as $\widehat{\mathrm{HSIC}}(\mathcal{X}_{n}^{\bpi_i}) + 2\Delta_T \xi_{\varepsilon,\delta}^{-1} \zeta_i$ for $i \in \{0\} \cup [B]$. We first examine the $1-\alpha$ quantile of $M_0,M_1,\ldots,M_B$ denoted as $q_{1-\alpha,B}$ and establish a high-probability upper bound for $q_{1-\alpha,B}$. Again, the proof strategy is similar to that of \Cref{Theorem: Uniform separation for MMD}. With an abuse of notation,  denote by $q_{1-\alpha/12,\infty}^a$ and $q_{1-\alpha/12,\infty}^b$, the $1-\alpha/12$ quantile of the conditional distribution of $\widehat{\mathrm{HSIC}}(\mathcal{X}_{n}^{\bpi})$ given $\mathcal{X}_{n+m}$ and that of $2\Delta_T \xi_{\varepsilon,\delta}^{-1} \zeta$, respectively. Following the same lines of the proof of \Cref{Theorem: Uniform separation for MMD}, we see that $q_{1-\alpha,B} \leq  q_{1-\alpha/12,\infty}^a + q_{1-\alpha/12,\infty}^b$ with probability at least $1-\beta/2$ under our condition for $B$. Moreover applying \Cref{Lemma: Concentration inequality for permuted HSIC} yields 
\begin{align*}
	q_{1-\alpha/12,\infty}^a \leq C_1 \sqrt{\frac{1}{n} \max \bigg\{ \log \biggl( \frac{12}{\alpha}\biggr), \, \log^{1/2} \biggl( \frac{12}{\alpha}\biggr),\, 1\bigg\} },
\end{align*}
whereas $q_{1-\alpha/12,\infty}^b = 2\Delta_T \xi_{\varepsilon,\delta}^{-1} F_{\zeta}^{-1}(1-\alpha/12)$. Therefore, with probability at least $1-\beta/2$, it holds that 
\begin{align} \label{Eq: upper bound for q2}
	q_{1-\alpha,B} \leq  C_1 \sqrt{\frac{1}{n} \max \bigg\{ \log \biggl( \frac{12}{\alpha}\biggr), \, \log^{1/2} \biggl( \frac{12}{\alpha}\biggr),\, 1\bigg\} } + 2\Delta_T \xi_{\varepsilon,\delta}^{-1} F_{\zeta}^{-1}(1-\alpha/12).
\end{align}

\vskip 1em 

\noindent \textbf{Step 2 (Bounding the type II error).} For the empirical HSIC, \Cref{Lemma: Exponential inequality for the empirical HSIC} shows that the following event
\begin{align*}
	E_1' \coloneqq   \biggl\{ \bigg|  \mathrm{HSIC}_{k \otimes \ell}(P_{YZ}) - \widehat{\mathrm{HSIC}}(\mathcal{X}_{n}) \bigg| \leq  C_2 \sqrt{\frac{\log(8/\beta)}{n}} \biggr\}
\end{align*}
holds with probability at least $1 - \beta/4$. On the other hand, the condition $\zeta_0 \sim \mathsf{Laplace}(0,1)$ ensures that the probability of the event 
\begin{align*}
	E_2' \coloneqq  \big\{ 2\Delta_T \xi_{\varepsilon,\delta}^{-1} \zeta_0 >   2\Delta_T \xi_{\varepsilon,\delta}^{-1}F^{-1}_{\zeta}(\beta/4)\big\}
\end{align*}
is equal to $1 - \beta/4$ where $F_{\zeta}^{-1}$ is the inverse cumulative distribution function of $\zeta$. Using these results as well as inequality~\eqref{Eq: upper bound for q2} from Step~1, it can be seen in a way similar to the proof of \Cref{Theorem: Uniform separation for MMD} that 
\begin{align*}
	\mE[1 - \phi_{\mathttt{dpHSIC}}] ~\leq~ & \mP\biggl(\mathrm{HSIC}_{k \otimes \ell}(P_{YZ}) \leq q_{1-\alpha,B} + C_2 \sqrt{\frac{\log(8/\beta)}{n}} - 2\Delta_T \xi_{\varepsilon,\delta}^{-1}F^{-1}_{\zeta}(\beta/4) \biggr)  + \frac{\beta}{2} \\[.5em]
	\leq~ & \mP\biggl(\mathrm{HSIC}_{k \otimes \ell}(P_{YZ}) \leq C_1 \sqrt{\frac{1}{n} \max \bigg\{ \log \biggl( \frac{12}{\alpha}\biggr), \, \log^{1/2} \biggl( \frac{12}{\alpha}\biggr),\, 1\bigg\} } \\[.5em]
	& ~~~~~  + 2\Delta_T \xi_{\varepsilon,\delta}^{-1} F_{\zeta}^{-1}(1-\alpha/12) + C_2 \sqrt{\frac{\log(8/\beta)}{n}} - 2\Delta_T \xi_{\varepsilon,\delta}^{-1}F^{-1}_{\zeta}(\beta/4) \biggr) + \beta.
\end{align*}
Next, recall that we assume $\alpha \in (0,1)$ and $\beta \in (0,1-\alpha)$ under which it holds
\begin{align} \label{Eq: alpha beta inequality}
	\max \bigg\{ \log \biggl( \frac{12}{\alpha}\biggr), \, \log^{1/2} \biggl( \frac{12}{\alpha}\biggr),\, 1, \, \log \biggl(\frac{8}{\beta}\biggr) \bigg\} ~\leq ~ C_3 \max \bigg\{ \log \biggl( \frac{1}{\alpha}\biggr), \, \log \biggl(\frac{1}{\beta}\biggr) \bigg\}.
\end{align}
To see this inequality, first assume that $\alpha \in (0,1/2)$. Then 
\begin{align*}
	\max \bigg\{ \log \biggl( \frac{12}{\alpha}\biggr), \,  \log^{1/2} \biggl( \frac{12}{\alpha}\biggr), \, 1  \bigg\} ~\leq~ &  C_4 \log \biggl( \frac{1}{\alpha}\biggr) \leq C_4 \max \biggl\{ \log \biggl( \frac{1}{\alpha}\biggr), \, \log \biggl( \frac{8}{\beta}\biggr) \biggr\} \\[.5em]
	\leq~ & C_5 \max \biggl\{ \log \biggl( \frac{1}{\alpha}\biggr), \, \log \biggl( \frac{1}{\beta}\biggr) \biggr\} 
\end{align*} 
as $\beta\in(0,1-\alpha)$.
On the other hand, if $\alpha \in [1/2,1)$, $\beta$ should be less than $1/2$, which implies that 
\begin{align*}
	\max \bigg\{ \log \biggl( \frac{12}{\alpha}\biggr), \,  \log^{1/2} \biggl( \frac{12}{\alpha}\biggr), \, 1  \bigg\} ~\leq~ &  C_6  \log \biggl( \frac{8}{\beta}\biggr) \leq C_7  \log \biggl( \frac{1}{\beta}\biggr) \\[.5em]
	\leq ~ & C_7  \max \biggl\{ \log \biggl( \frac{1}{\alpha}\biggr), \, \log \biggl( \frac{1}{\beta}\biggr) \biggr\}. 
\end{align*}
Therefore inequality~\eqref{Eq: alpha beta inequality} follows. 

Having these ingredients in place, we follow the same lines as in the proof of \Cref{Theorem: Uniform separation for MMD} and see that
\begin{align*}
	& \mE[1 - \phi_{\mathttt{dpHSIC}}]  \\[.5em] 
	\leq ~ & \mP\biggl(\mathrm{HSIC}_{k \otimes \ell}(P_{YZ}) \leq C_8 \sqrt{\frac{\max\bigl\{\log(1/\alpha),\log(1/\beta)\}}{n}} + C_9 \frac{\max\bigl\{\log(1/\alpha),\log(1/\beta)\}}{n\xi_{\varepsilon,\delta}} \biggr) +  \beta  \\[.5em]
	\leq ~ & \beta,
\end{align*}
where the last inequality holds by taking $C_{K,L}$ to be larger than, for instance, $2 \max \{C_8, C_9\} + 1$ in the theorem statement. We conclude the proof by noting that the upper bound is independent of $P_{YZ}$ and taking the supremum on both sides over $\mathcal{P}_{\mathrm{HSIC}_{k \otimes \ell}}\!(\rho)$.

\subsection{Proof of \Cref{Theorem: Minimax separation in HSIC}} \label{Section: Proof of Theorem: Minimax separation in HSIC}
The proof of \Cref{Theorem: Minimax separation in HSIC} is based on the same idea as that of \Cref{Theorem: Minimax separation in MMD}, which uses Le Cam's two point method and coupling method. Therefore, we omit the details explained in the proof of \Cref{Theorem: Minimax separation in MMD}, and instead focus on key differences. As in the proof of \Cref{Theorem: Minimax separation in MMD}, we simply let $\rho^\star_{\mathrm{HSIC}} = \rho^\star_{\mathrm{HSIC}}(\alpha,\beta,\varepsilon,\delta,n) $ and establish the minimax separation by examining the non-privacy regime and privacy regime in order.

\subsubsection{Non-privacy regime}

We begin by proving that $\rho^\star_{\mathrm{HSIC}} \geq C_{\eta_Y,\eta_Z} \min\{\sqrt{\log(1/(\alpha+\beta))/n},1\}$ in the non-privacy regime. We pick one distribution $P_{YZ,0}$ from $\mathcal{P}_{\mathrm{HSIC}_{k \otimes \ell}}(\widetilde{\rho})$ with 
\begin{align*}
	\widetilde{\rho} = C_{\eta_Y,\eta_Z} \min \bigg\{ \sqrt{\frac{\log(1/(\alpha+\beta))}{n}}, \, 1 \bigg\}
\end{align*}
and denote the product of its marginals as $P_{Y,0}P_{Z,0}$. Then, as in \Cref{Section: Non-privacy regime MMD}, an application of Le Cam's two point method~\citep{lecam1973convergence,le2012asymptotic} and Bretagnolle--Huber inequality~\citep[][Lemma B.4]{canonne2022topics} yields
\begin{align*}
	\inf_{\phi \in \Phi_{\alpha,\varepsilon,\delta}} \sup_{P_{YZ} \in \mathcal{P}_{\mathrm{HSIC}_{k \otimes \ell}}(\widetilde{\rho})} \mE_{P_{YZ}^n} [1 - \phi] ~ \geq ~ & \inf_{\phi \in \Phi_{\alpha,\infty}} \sup_{P_{YZ} \in \mathcal{P}_{\mathrm{HSIC}_{k \otimes \ell}}(\widetilde{\rho})} \mE_{P_{YZ}} [1 - \phi]  \\[.5em]
	\geq ~ & \inf_{\phi \in \Phi_{\alpha,\infty}}  \mE_{P_{YZ,0}} [1 - \phi] = 1 -  \sup_{\phi \in \Phi_{\alpha,\infty}}  \mE_{P_{YZ,0}} [\phi] \\[.5em]
	\geq ~ & 1 - \alpha - d_{\mathrm{TV}}(P_{YZ,0}^{\otimes n}, P_{Y,0}^{\otimes n} P_{Z,0}^{\otimes n}) \\[.5em]
	\geq ~ & \frac{1}{2} e^{-n \times d_{\mathrm{KL}}(P_{YZ,0} \| P_{Y,0}P_{Z,0} )} - \alpha.
\end{align*}
Therefore the minimax type II error is at least $\beta$ if $\alpha + \beta < 0.4$ as well as 
\begin{align*}
	d_{\mathrm{KL}}(P_{YZ,0} \| P_{Y,0}P_{Z,0}) \leq \frac{1}{n} \log \biggl( \frac{1}{2(\alpha + \beta)} \biggr).
\end{align*}
Hence $\rho^\star_{\mathrm{HSIC}} \geq C_{\eta_Y,\eta_Z} \min\{\sqrt{\log(1/(\alpha+\beta))/n},1\}$ follows if we find $P_{YZ,0}$ such that 
\begin{subequations}
	\begin{align} \label{Eq: sufficient condition in non-private regime for HSIC(a)}
		& \mathrm{HSIC}_{k \otimes \ell}(P_{YZ,0}, P_{Y,0}P_{Z,0}) \geq C_{\eta_Y,\eta_Z} \min \Bigg\{ \sqrt{\frac{\log(1/(\alpha+\beta))}{n}}, \, 1 \Bigg\} \quad \text{and} \\[.5em]  \label{Eq: sufficient condition in non-private regime for HSIC(b)}
		& d_{\mathrm{KL}}(P_{YZ,0} \| P_{Y,0} P_{Z,0}) \leq \frac{1}{n} \log\biggl( \frac{1}{2(\alpha + \beta)} \biggr). 
	\end{align}
\end{subequations}
To this end, consider (discrete) random vectors $Y \in \{y_1,y_2\}$ and $Z \in \{z_1,z_2\}$ where $y_1,y_2 \in \mathbb{R}^{d_Y}$ and $z_1,z_2 \in \mathbb{R}^{d_Z}$. Further assume that $\mP(Y=y_1) = \mP(Y=y_2) = 1/2$ and $\mP(Z=z_1) = \mP(Z=z_2) = 1/2$. Suppose that the joint probabilities of $(Y,Z)$ are given as
\begin{align*}
	& \mP(Y=y_1, Z = z_1) = \mP(Y=y_2, Z=z_2) = 1/4 + \nu \quad \text{and} \\[.5em]
	& \mP(Y=y_1, Z=z_2) = \mP(Y = y_2, Z= z_1)  = 1/4 - \nu,
\end{align*}
where $\nu \in (0, 1/4]$. In this setting, it can be seen that $P_{YZ,0} \neq P_{Y,0} P_{Z,0}$ and thus $Y$ and $Z$ are trivially dependent. For such $P_{YZ,0}$ and translation invariant kernels, we have
\begin{align*}
	\mathrm{HSIC}_{k \otimes \ell} (P_{YZ,0}) ~\overset{\mathrm{(i)}}{=}~&  \sqrt{\mE[k(Y_1,Y_2) \ell(Z_1,Z_2)] + \mE[k(Y_1,Y_2)\ell(Z_3,Z_4)] - 2 \mE[k(Y_1,Y_2)\ell(Z_1,Z_3)]} \\[.5em]
	\overset{\mathrm{(ii)}}{=}~ & \sqrt{4 \nu^2\{\kappa_Y(0) - \kappa_Y(y_1-y_2)\} \{\kappa_Z(0) - \kappa_Z(z_1-z_2)\}} \\[.5em]
	\overset{\mathrm{(iii)}}{\geq} ~ &  \sqrt{4 \nu^2 \eta_Y \eta_Z} = 2 \nu \sqrt{\eta_Y \eta_Z},
\end{align*}
where step~(i) follows by \citet[][Lemma 1]{gretton2005measuring} with $\{(Y_i,Z_i)\}_{i=1}^4 \iid P_{YZ,0}$, step~(ii) can be verified through algebra and the last inequality~(iii) holds by choosing $(y_1,y_2)$ and $(z_1,z_2)$ such that $y_1 - y_2 = y_0$ and $z_1 - z_2 = z_0$. 

On the other hand, we verify condition~\eqref{Eq: sufficient condition in non-private regime for HSIC(b)} based on the well-known fact \citep[\emph{e.g.},][Lemma 2.7]{tsybakov2009} that the Kullback--Leibler divergence is upper bounded by $\chi^2$ divergence denoted by $d_{\chi^2} (P_{YZ,0} \| P_{Y,0} P_{Z,0})$. This gives
\begin{align*}
	d_{\mathrm{KL}}(P_{YZ,0} \| P_{Y,0} P_{Z,0})  ~\leq~&  d_{\chi^2} (P_{YZ,0} \| P_{Y,0} P_{Z,0}) = 8 \biggl(\frac{1}{4} + \nu \biggr)^2 + 8 \biggl(\frac{1}{4} - \nu \biggr)^2 - 1 \\[.5em]
	= ~ & 16 \nu^2.
\end{align*}
Using the above inequality, condition~\eqref{Eq: sufficient condition in non-private regime for HSIC(b)} is fulfilled by taking 
\begin{align*}
	\nu = \min \Bigg\{ \sqrt{\frac{1}{16n} \log \biggl( \frac{1}{2(\alpha + \beta)} \biggr)}, \, \frac{1}{4} \Bigg\}
\end{align*}
for which condition~\eqref{Eq: sufficient condition in non-private regime for HSIC(a)} is also satisfied as
\begin{align*}
	\mathrm{HSIC}_{k \otimes \ell} (P_{YZ,0})  ~\geq~ & 2 \nu \sqrt{\eta_Y \eta_Z} =  \frac{\sqrt{\eta_Y \eta_Z}}{2}  \min \Bigg\{ \sqrt{\frac{1}{n} \log \biggl( \frac{1}{2(\alpha + \beta)} \biggr)}, \, 1 \Bigg\} \\[.5em]
	\geq ~ & C_{\eta_Y,\eta_Z} \min \Bigg\{ \sqrt{\frac{\log(1/(\alpha+\beta))}{n}}, \, 1 \Bigg\}
\end{align*}
as $\alpha+\beta<0.4$.
Therefore it holds that $\rho^\star_{\mathrm{HSIC}} \geq C_{\eta_Y,\eta_Z} \min\{\sqrt{\log(1/(\alpha+\beta))/n},1\}$ as desired.

\subsubsection{Privacy regime} 
The proof of the separation rate under the privacy regime is essentially the same as that for the $\mathttt{dpMMD}$ test in Appendix~\ref{Section: Privacy regime MMD}. All we need is to find an instance of $P_{YZ,0}$ such that 
\begin{subequations}
	\begin{align} \label{Eq: sufficient condition in private regime for HSIC(a)}
		& \mathrm{HSIC}_{k \otimes \ell}(P_{YZ,0}) \geq C_{\eta_Y,\eta_Z} \min \Bigg\{ \frac{\log(1/\beta)}{n \big(\varepsilon+\delta\big)}, \, 1 \Bigg\} \quad \text{and} \\[.5em]  \label{Eq: sufficient condition in private regime for HSIC(b)}
		& 	\|P_{YZ,0} - P_{Y,0}P_{Z,0} \|_1 \leq \frac{1}{20n\big(\varepsilon+\delta\big)} \log \biggl( \frac{1}{4\beta} \biggr). 
	\end{align}
\end{subequations}
To this end, choose the distribution of $(Y,Z)$ as in the case of the non-privacy regime. For such $P_{YZ,0}$, a direct calculation gives $\|P_{YZ,0} - P_{Y,0}P_{Z,0} \|_1 = 4 \nu $ and thus condition~\eqref{Eq: sufficient condition in private regime for HSIC(b)} is satisfied if we take
\begin{align*}
	\nu = \min \Bigg\{ \frac{1}{80n \big(\varepsilon+\delta\big)} \log \biggl( \frac{1}{4\beta} \biggr), \, \frac{1}{4} \Bigg\}.
\end{align*}
On the other hand, the previous calculation of the lower bound for $\mathrm{HSIC}_{k \otimes \ell} (P_{YZ,0})$ yields 
\begin{align*}
	\mathrm{HSIC}_{k \otimes \ell} (P_{YZ,0}) ~ \geq ~ & 2 \nu \sqrt{\eta_Y \eta_Z} = 2 \min \Bigg\{ \frac{1}{80n \big(\varepsilon+\delta\big)} \log \biggl( \frac{1}{4\beta} \biggr), \, \frac{1}{4} \Bigg\} \sqrt{\eta_Y \eta_Z}  \\[.5em]
	\geq ~ & C_{\eta_Y,\eta_Z} \min \Bigg\{ \frac{\log(1/\beta)}{n \big(\varepsilon+\delta\big)}, \, 1 \Bigg\}
\end{align*}
which in turn verifies condition~\eqref{Eq: sufficient condition in private regime for HSIC(a)} under $\beta \in (0,1/5)$. Therefore it holds that the minimax separation under privacy regime $\rho^\star_{\mathrm{HSIC}} \geq  C_{\eta_Y,\eta_Z}  \min\!\big\{ \log(1/\beta) /(n(\varepsilon+\delta)), 1 \big\}$. 

\subsubsection{Combining bounds}
Combining the lower bounds for the non-privacy and privacy regimes, for all $\varepsilon>0$ and $\delta\in[0,1)$, we obtain
\begin{align*}
	\rho^\star_{\mathrm{HSIC}} ~ \geq ~ & C_{\eta_Y,\eta_Z}  \max \Biggl\{ \min \bigg\{  \sqrt{\frac{\log(1/(\alpha+\beta))}{n}}, \, 1 \bigg\}, \, \min \biggl\{ \frac{\log(1/\beta)}{n(\varepsilon+\delta)}, \, 1 \biggr\} \Biggr\} \\[.5em]
	\geq ~ & C_{\eta_Y,\eta_Z}  \max \Biggl\{ \min \bigg\{  \sqrt{\frac{\log(1/(\alpha+\beta))}{n}}, \, 1 \bigg\}, \, \min \biggl\{ \frac{\log(1/\beta)}{n\big(\varepsilon+\log(1/(1-\delta))\big)}, \, 1 \biggr\} \Biggr\}.
\end{align*}
The last inequality holds since $\log(1/(1-\delta))\geq\delta$ for all $\delta\in[0,1)$, and is actually tight when $\alpha\asymp\beta$ as explained in \Cref{Section: Equivalence MMD}. This concludes the proof of \Cref{Theorem: Uniform separation for HSIC}.

\subsection{Proof of \Cref{Theorem: Minimax Separation over L2 for HSIC}} \label{Section: Proof of Theorem: Minimax Separation over L2 for HSIC}
The proof of \Cref{Theorem: Minimax Separation over L2 for HSIC} follows the same structure as the proof of \Cref{Theorem: Minimax Separation over L2} in \Cref{Section: Proof of Theorem: Minimax Separation over L2}. For simplicity, write $\hat{\mathrm{HSIC}}^2$ given in \eqref{Eq: closed form HSIC} with the Gaussian kernels as $V_{\mathrm{HSIC}}$ and similarly the U-statistic given in \eqref{Eq: U-HSIC} with the Gaussian kernel as $U_\mathrm{HSIC}$. We also denote the corresponding statistics based on permuted data $\mathcal{X}_n^{\bpi}$ as $V_{\bpi,\mathrm{HSIC}}$ and $U_{\bpi,\mathrm{HSIC}}$, respectively. Letting
\begin{align*}
	K =  \prod_{i=1}^{d_Y} \frac{1}{\sqrt{2\pi} \lambda_i}  \quad \text{and} \quad L =  \prod_{i=1}^{d_Z} \frac{1}{\sqrt{2\pi} \mu_i},
\end{align*}
\Cref{Lemma: Difference between V and U} ensures that the difference between $V_\mathrm{HSIC}$ and $U_\mathrm{HSIC}$ can be written as 
\begin{align*}
	V_\mathrm{HSIC} - U_\mathrm{HSIC} = D_1 + D_2,
\end{align*}
where
\begin{align*}
	D_1 = & \frac{n-1}{n^2}KL -\frac{L}{n^3} \sum_{(i,j) \in \mathbf{i}_2^n} k_{\blambda}(Y_i,Y_j) - \frac{K}{n^3} \sum_{(i,j) \in \mathbf{i}_2^n} \ell_{\bmu}(Z_i,Z_j), \\[.5em]
	D_2= & - \frac{3n^2-4n+2}{(n-1)n^4}\sum_{(i,j) \in \mathbf{i}_2^n} k_{\blambda}(Y_i,Y_j) \ell_{\bmu}(Z_i,Z_j)+ \frac{2(5n^2-8n+4)}{n^4(n-1)(n-2) }\sum_{(i,j_1,j_2) \in \mathbf{i}_3^n} k_{\blambda}(Y_i,Y_{j_1}) \ell_{\bmu}(Z_i,Z_{j_2})\\[.5em]
	& - \frac{6n^2-11n+6}{n^4(n-1)(n-2)(n-3)} \sum_{(i_1,i_2,j_1,j_2) \in \mathbf{i}_4^n} k_{\blambda}(Y_{i_1},Y_{i_2}) \ell_{\bmu}(Z_{j_1}, Z_{j_2}).
\end{align*}
The term $D_1$ is invariant to any permutation of $Z$ values, and thus the permutation distribution of $V_{\bpi,\mathrm{HSIC}}$ is equivalent to the permutation distribution of $U_{\bpi,\mathrm{HSIC}} + D_1 + D_{2,\bpi}$ where $D_{2,\bpi}$ has the same form of $D_2$ but computed based on permuted data $\mathcal{X}_n^{\bpi}$. Writing the sensitivity of $\sqrt{V_\mathrm{HSIC}}$ as $\Delta_{V^{1/2}}$, \Cref{Lemma: Quantile representation} yields that the $\mathttt{dpHSIC}$ test rejects the null if and only if 
\begin{align*}
	\sqrt{V_\mathrm{HSIC}} + \frac{2\Delta_{V^{1/2}}}{\xi_{\varepsilon,\delta}} \zeta_0 > q_{1-\alpha,B},
\end{align*}
where $q_{1-\alpha,B}$ is the $1-\alpha$ quantile of $\{ \sqrt{V_{\bpi_i,\mathrm{HSIC}}} + 2 \Delta_{V^{1/2}}\xi_{\varepsilon,\delta}^{-1} \zeta_i \}_{i=0}^B$. As in the proof of \Cref{Theorem: Minimax Separation over L2} in \Cref{Section: Proof of Theorem: Minimax Separation over L2}, we let $q_{1-\alpha,\infty}^a$ denote the $1-\alpha$ quantile of the conditional distribution of $V_{\bpi,\mathrm{HSIC}}^{1/2}$ given $\mathcal{X}_n$. Then under the condition on $B$, the analysis given in the proof of \Cref{Theorem: Minimax Separation over L2} shows that the type II error of the $\mathttt{dpHSIC}$ test can be bounded as
\begin{align} \nonumber
	& \mP\bigl(\sqrt{V_\mathrm{HSIC}} + 2 \Delta_{V^{1/2}} \xi_{\varepsilon,\delta}^{-1} \zeta_0 \leq q_{1-\alpha,B}\bigr) \\[.5em]
	\leq ~ & \mP\bigl( \sqrt{V_\mathrm{HSIC}} \leq q_{1-\alpha/12,\infty}^a + 14 \Delta_{V^{1/2}} \xi_{\varepsilon,\delta}^{-1} \max\{\log(1/\alpha),\log(1/\beta)\} \bigr) + 5\beta/ 8. \label{Eq: type II error bound for HSIC}
\end{align}
We next delve into $q_{1-\alpha/12,\infty}^a$, and then continue to upper bound the type II error. In what follows, we use the notation $C_1, C_2, C_3, \ldots$ to represent constants, which could be dependent on $\alpha, \beta, s, R, M, d_Y, d_Z$. The values of these constants may vary in different places. 

\vskip 1em 

\noindent \textbf{Bounding $q_{1-\alpha/12,\infty}^a$.} Using the relationship between $V_\mathrm{HSIC}$ and $U_\mathrm{HSIC}$ in \Cref{Lemma: Difference between V and U} and the quantile inequality in \Cref{Lemma: Quantile inequality}, we have 
\begin{align*}
	q_{1-\alpha/12,\infty}^a ~=~&  \sqrt{\mathrm{Quantile}_{1-\alpha/12} \bigl(\big\{ U_{\bpi_i,\mathrm{HSIC}} + D_{2,\bpi_i} \big\}_{i=0}^\infty + D_1 \bigr)} \\[.5em]
	\leq ~ & \sqrt{\mathrm{Quantile}_{1-\alpha/24}\bigl( \big\{ U_{\bpi_i,\mathrm{HSIC}} \big\}_{i=0}^\infty \bigr) + \mathrm{Quantile}_{1-\alpha/24} \bigl(  \big\{ D_{2,\bpi_i} \big\}_{i=0}^\infty \bigr) +  D_1}.
\end{align*}
The quantile of the permuted U-statistic $U_{\bpi,\mathrm{HSIC}}$ has been studied by \cite{kim2020minimax}. In particular, the proof of \citet[][Theorem 5.1]{kim2020minimax} along with the proof of \citet[][Proposition 2]{albert2019adaptive} shows that  
\begin{align*}
	\mE_{\bpi}[U_{\bpi,\mathrm{HSIC}} \given \mathcal{X}_n] = 0 \quad \text{and} \quad \mE \bigl[ \mV_{\bpi} \bigl( U_{\bpi,\mathrm{HSIC}} \given \mathcal{X}_n \bigr) \bigr] \leq \frac{C_1}{n^2 \lambda_1 \cdots \lambda_{d_Y} \mu_1 \cdots \mu_{d_Z}}.
\end{align*}
Thus by Chebyshev's inequality, it can be seen that the following inequality holds with probability at least $1-\beta/8$:
\begin{align} \label{Eq: quantile bound 1}
	\mathrm{Quantile}_{1-\alpha/24}\bigl( \big\{ U_{\bpi_i,\mathrm{HSIC}} \big\}_{i=0}^\infty \bigr)  \leq \frac{C_2}{n \sqrt{\lambda_1 \cdots \lambda_{d_Y} \mu_1 \cdots \mu_{d_Z}}},
\end{align}
For the $1-\alpha/24$ quantile of $\big\{ D_{2,\bpi_i} \big\}_{i=0}^\infty$, note that 
\begin{align*}
	D_{2,\bpi} \leq  \frac{C_3}{n^4}\sum_{(i,j_1,j_2) \in \mathbf{i}_3^n} k_{\blambda}(Y_i,Y_{j_1}) \ell_{\bmu}(Z_{\pi_i},Z_{\pi_{j_2}})
\end{align*}
due to the non-negativity of $k_{\blambda}$ and $\ell_{\bmu}$, and also note that 
\begin{align*}
	\mE_{\bpi} \Biggl[  \sum_{(i,j_1,j_2) \in \mathbf{i}_3^n} k_{\blambda}(Y_i,Y_{j_1}) \ell_{\bmu}(Z_{\pi_i},Z_{\pi_{j_2}}) \,\Big| \, \mathcal{X}_n \Biggr] \leq \frac{C_4}{n} \Biggl[ \sum_{(i_1,i_2) \in \mathbf{i}_2^n} k_{\blambda}(Y_{i_1},Y_{i_2}) \Biggr] \Biggl[ \sum_{(j_1,j_2) \in \mathbf{i}_2^n}  \ell_{\bmu}(Z_{j_1},Z_{j_2}) \Biggr].
\end{align*}
Then Markov's inequality yields that the $1-\alpha/24$ quantile of $\big\{  D_{2,\bpi_i}  \big\}_{i=0}^\infty$ is bounded as
\begin{align*}
	\mathrm{Quantile}_{1-\alpha/8} \bigl(  \big\{  D_{2,\bpi_i}  \big\}_{i=0}^\infty \bigr)  \leq \frac{C_5}{n^5} \Biggl[ \sum_{(i_1,i_2) \in \mathbf{i}_2^n} k_{\blambda}(Y_{i_1},Y_{i_2}) \Biggr] \Biggl[ \sum_{(j_1,j_2) \in \mathbf{i}_2^n}  \ell_{\bmu}(Z_{j_1},Z_{j_2}) \Biggr].
\end{align*}
Moreover, since we assume $\|p_{YZ}\|_{L_\infty} \leq M$ and $\|p_Yp_Z\|_{L_\infty} \leq M$, 
\begin{align*}
	\max \Big\{\mE[k_{\blambda}(Y_1,Y_2)\ell_{\bmu}(Z_1,Z_2)], \ \mE[k_{\blambda}(Y_1,Y_2)\ell_{\bmu}(Z_1,Z_3)], \ \mE[k_{\blambda}(Y_1,Y_2)\ell_{\bmu}(Z_3,Z_4)] \Big\} \leq C_6.
\end{align*}
For instance, we see that 
\begin{align*}
	\mE[k_{\blambda}(Y_1,Y_2)\ell_{\bmu}(Z_1,Z_2)] ~=~ & \int \cdots \int k_{\blambda}(y_1,y_2)\ell_{\bmu}(z_1,z_2) p_{YZ}(y_1,z_1) p_{YZ}(y_2,z_2) \dd y_1 \dd y_2 \dd z_1 \dd z_2 \\[.5em]
	\leq ~ & \|p_{YZ}\|_{L_\infty} \int \int \underbrace{\int k_{\blambda}(y_1,y_2)\dd y_1}_{=1} \underbrace{\int\ell_{\bmu}(z_1,z_2)\dd z_1}_{=1} p_{YZ}(y_2,z_2) \dd y_2 \dd z_2 \\[.5em]
	\leq ~ & M,
\end{align*}
and the other terms can be similarly analyzed. Given this ingredient, we have
\begin{align*}
	\mE \Biggl\{ \Biggl[ \sum_{(i_1,i_2) \in \mathbf{i}_2^n} k_{\blambda}(Y_{i_1},Y_{i_2}) \Biggr] \Biggl[ \sum_{(j_1,j_2) \in \mathbf{i}_2^n}  \ell_{\bmu}(Z_{j_1},Z_{j_2}) \Biggr] \Biggr\} \leq C_7 n^4.
\end{align*}
Therefore, another application of Markov's inequality yields 
\begin{align} \label{Eq: quantile bound 2}
	\mathrm{Quantile}_{1-\alpha/24} \bigl(  \big\{D_{2,\bpi_i}  \big\}_{i=0}^\infty \bigr)  \leq  \frac{C_8}{n},
\end{align}
with probability at least $1-\beta/8$, and observe that 
$D_1 \leq \frac{n-1}{n^2} KL \leq C_9/\{n\lambda_1\cdots \lambda_{d_Y} \mu_1 \cdots \mu_{d_Z}\}$. 

In summary, with probability at least $1-\beta/4$, 
\begin{align} \label{Eq: quantile bound 3}
	q_{1-\alpha/12,\infty}^a \leq  \frac{C_{10}}{\sqrt{n \lambda_1 \cdots \lambda_{d_Y} \mu_1 \cdots \mu_{d_Z}}}.
\end{align}
given that $D_1 \leq KL/n$ and $ \lambda_1 \cdots \lambda_{d_Y} \mu_1 \cdots \mu_{d_Z} \leq 1$.

\vskip 1em 

\noindent \textbf{Bounding the type II error.} We continue from \eqref{Eq: type II error bound for HSIC}, and note that 
\begin{align*}
	& \mP \Bigl( \sqrt{V_\mathrm{HSIC}} \leq q_{1-\alpha/12,\infty}^a + 14 \Delta_{V^{1/2}} \xi_{\varepsilon,\delta}^{-1} \max\{\log(1/\alpha),\log(1/\beta)\} \Bigr) \\[.5em]
	= ~ & \mP\Bigl(V_\mathrm{HSIC} \leq (q_{1-\alpha/12,\infty}^a)^2 + 14^2 \Delta_{V^{1/2}}^2 \xi_{\varepsilon,\delta}^{-2} \max\{\log^2(1/\alpha),\log^2(1/\beta)\} \\[.5em]
	& ~~~ + 28 q_{1-\alpha/12,\infty}^a  \Delta_{V^{1/2}} \xi_{\varepsilon,\delta}^{-1} \max\{\log(1/\alpha),\log(1/\beta)\} \Bigr) \\[.5em]
	= ~ &  \mP\Bigl(U_\mathrm{HSIC}  \leq \mathrm{Quantile}_{1-\alpha/24}\bigl( \big\{ U_{\bpi_i,\mathrm{HSIC}} \big\}_{i=0}^\infty \bigr) + \mathrm{Quantile}_{1-\alpha/24} \bigl(  \big\{ D_{2,\bpi_i} \big\}_{i=0}^\infty \bigr) - D_2 \\[.5em]
	& ~~~ + 14^2 \Delta_{V^{1/2}}^2 \xi_{\varepsilon,\delta}^{-2} \max\{\log^2(1/\alpha),\log^2(1/\beta)\} + 28 q_{1-\alpha/12,\infty}^a  \Delta_{V^{1/2}} \xi_{\varepsilon,\delta}^{-1} \max\{\log(1/\alpha),\log(1/\beta)\} \Bigr), 
\end{align*}
where the first equality simply follows by taking the square on both sides, and the second equality uses the identity $V_\mathrm{HSIC} = U_\mathrm{HSIC} + D_1 + D_2$ from \Cref{Lemma: Difference between V and U}. Next, by \Cref{Lemma: Sensitivity of HSIC}, the global sensitivity of $\sqrt{V_\mathrm{HSIC}}$ is bounded as 
\begin{align*}
	\Delta_{V^{1/2}} \leq \frac{C_1}{n\sqrt{\lambda_1 \cdots \lambda_{d_Y} \mu_1 \cdots \mu_{d_Z}}}.
\end{align*}
In addition, observe that $|D_2| \leq C_2n^{-1}$ with probability at least $1-\beta/16$, which can be proven by Markov's inequality. Therefore, equipped with the previous ingredients~\eqref{Eq: type II error bound for HSIC}, \eqref{Eq: quantile bound 1}, \eqref{Eq: quantile bound 2} and \eqref{Eq: quantile bound 3}, the type II error bound can be bounded as 
\begin{align*} 
	& \mP\bigl(\sqrt{V_\mathrm{HSIC}} + 2 \Delta_{V^{1/2}} \xi_{\varepsilon,\delta}^{-1} \zeta_0 \leq q_{1-\alpha,B}\bigr) \\[.5em]
	\leq  ~& \mP\biggl( U_\mathrm{HSIC} \leq \frac{C_3}{n \sqrt{\lambda_1 \cdots \lambda_{d_Y} \mu_1 \cdots \mu_{d_Z}}} + \frac{C_4}{n^2 \lambda_1 \cdots \lambda_{d_Y} \mu_1 \cdots \mu_{d_Z} \xi_{\varepsilon,\delta}^2} \\[.5em]
	&~~~~~~~ + \frac{C_5}{n^{3/2}  \lambda_1 \cdots \lambda_{d_Y} \mu_1 \cdots \mu_{d_Z} \xi_{\varepsilon,\delta}} \biggr) +\frac{15}{16}\beta.
\end{align*}

\vskip 1em

\noindent \textbf{Condition in terms of $L_2$ distance.} A slight modification of \citet[][Lemma 1]{albert2019adaptive} yields that a sufficient condition for the first probability in the above display to be less than $\beta/16$ is 
\begin{align*}
	\mathrm{HSIC}^2_{k_{\blambda} \otimes \ell_{\bmu}} ~\geq~ & \sqrt{\mV[U_\mathrm{HSIC}]} + \frac{C_3}{n \sqrt{\lambda_1 \cdots \lambda_{d_Y} \mu_1 \cdots \mu_{d_Z}}} + \frac{C_4}{n^2 \lambda_1 \cdots \lambda_{d_Y} \mu_1 \cdots \mu_{d_Z} \xi_{\varepsilon,\delta}^2}  \\[.5em]
		& +  \frac{C_5}{n^{3/2}  \lambda_1 \cdots \lambda_{d_Y} \mu_1 \cdots \mu_{d_Z} \xi_{\varepsilon,\delta}}. 
\end{align*}
In addition, by writing $\psi = p_{YZ} - p_Yp_Z$ and the convolution of $\psi$ and $k_{\blambda} \otimes \ell_{\bmu}$ as $\psi \ast (k_{\blambda} \otimes \ell_{\bmu})$, \citet[][Theorem 1 and Proposition 4]{albert2019adaptive} show that the previous condition is implied by
\begin{align*}
	\|\psi\|_{L_2}^2 ~ \geq ~ \|\psi - \psi \ast (k_{\blambda} \otimes \ell_{\bmu})\|_{L_2}^2 + C_6 &  \Bigg\{ \frac{1}{n \sqrt{\lambda_1 \cdots \lambda_{d_Y} \mu_1 \cdots \mu_{d_Z}}}  + \frac{1}{n^2 \lambda_1 \cdots \lambda_{d_Y} \mu_1 \cdots \mu_{d_Z} \xi_{\varepsilon,\delta}^2}  \\[.5em] 
	&+ \frac{1}{n^{3/2}  \lambda_1 \cdots \lambda_{d_Y} \mu_1 \cdots \mu_{d_Z} \xi_{\varepsilon,\delta}} \Bigg\}.
\end{align*}
Finally the proof of \cite[][Theorem 2]{albert2019adaptive} yields that over the Sobolev ball, a sufficient condition for the previous inequality is 
\begin{align*}
	\|\psi\|_{L_2}^2  ~\geq~ C_7\Bigg\{& \sum_{i=1}^{d_Y} \lambda_i^{2s} + \sum_{i=1}^{d_Z} \mu_i^{2s}  + \frac{1}{n\sqrt{\lambda_1 \cdots \lambda_{d_Y} \mu_1 \cdots \mu_{d_Z}}} \\[.5em]
	& + \frac{1}{n^2 \lambda_1 \cdots \lambda_{d_Y} \mu_1 \cdots \mu_{d_Z} \xi_{\varepsilon,\delta}^2}  + \frac{1}{n^{3/2} \lambda_1 \cdots \lambda_{d_Y} \mu_1 \cdots \mu_{d_Z} \xi_{\varepsilon,\delta}}  \Bigg\}
\end{align*}
as claimed.

\subsection{Proof of \Cref{Lemma: Global sensitivity of U_HSIC}} \label{Section: Proof of Lemma: Global sensitivity of U_HSIC}
\paragraph{Sensitivity upper bound.}
For simplicity, write $U_{\mathrm{HSIC}} = U_{\mathrm{HSIC},a} + U_{\mathrm{HSIC},b}  - 2 U_{\mathrm{HSIC},c}$ where
\begin{align*}
	U_{\mathrm{HSIC},a} ~:=~ & \frac{1}{n(n-1)} \sum_{(i,j) \in \mathbf{i}_2^n} k(Y_i,Y_j) \ell(Z_i,Z_j),  \\[.5em]
	U_{\mathrm{HSIC},b} ~:=~ &  \frac{(n-4)!}{n!} \sum_{(i_1,i_2,j_1,j_2) \in \mathbf{i}_4^n} k(Y_{i_1},Y_{j_1}) \ell(Z_{i_2},Z_{j_2}),  \\[.5em]
	U_{\mathrm{HSIC},c} ~:=~ & \frac{1}{n(n-1)(n-2)} \sum_{(i,j_1,j_2) \in \mathbf{i}_3^n} k(Y_{i},Y_{j_1}) \ell(Z_{i},Z_{j_2}).
\end{align*}
To establish an upper bound for the global sensitivity, we first consider a neighboring dataset $\tilde{\mathcal{X}}_n = \{(Y_1',Z_1'), (Y_2,Z_2),\ldots,(Y_n,Z_n)\}$ and denote the U-statistic of HSIC based on $\tilde{\mathcal{X}}_n$ as $U_{\mathrm{HSIC}}':=  U_{\mathrm{HSIC},a}' + U_{\mathrm{HSIC},b}'  - 2 U_{\mathrm{HSIC},c}'$. By the triangle inequality, the absolute deviation between $U_{\mathrm{HSIC}}$ and $U_{\mathrm{HSIC}}'$ is then bounded as 
\begin{align*}
	|U_{\mathrm{HSIC}} - U_{\mathrm{HSIC}}'| ~\leq~ & | U_{\mathrm{HSIC},a} - U_{\mathrm{HSIC},a}'| + |U_{\mathrm{HSIC},b} - U_{\mathrm{HSIC},b}'| + 2 |U_{\mathrm{HSIC},c} - U_{\mathrm{HSIC},c}'| \\[.5em]
	\leq ~ & \frac{C_1 KL}{n},
\end{align*}
where $C_1,C_2,\ldots$ denote some universal positive constants throughout.

For another neighboring dataset $\tilde{\mathcal{X}}_n = \{(Y_1',Z_1), (Y_2,Z_2'),\ldots,(Y_n,Z_n)\}$, a similar analysis shows that $|U_{\mathrm{HSIC}} - U_{\mathrm{HSIC}}'| \leq \frac{C_2 KL}{n}$. Since $U_{\mathrm{HSIC}}$ is invariant to the permutation of the paired indices, we may conclude that 
\begin{align} \label{Eq: upper bound for the sensitivity of U-HSIC}
	\sup_{\bpi \in \boldsymbol{\Pi}_n} \sup_{\substack{\mathcal{X}_{n},\tilde{\mathcal{X}}_{n}:\\d_{\mathrm{ham}}(\mathcal{X}_{n},\tilde{\mathcal{X}}_{n}) \leq 1}} \bigl| U_{\mathrm{HSIC}}(\mathcal{X}_{n}^{\bpi}) - U_{\mathrm{HSIC}}(\tilde{\mathcal{X}}_{n}^{\bpi}) \bigr| \leq \frac{C_3 KL}{n}.
\end{align}
While it is loose, the inequality holds with $C_3 = 24$ for $n \geq 4$. 

\paragraph{Sensitivity lower bound.}
We now show that this upper bound is tight up to a constant factor. We treat the cases of $n$ being even or odd separately, and follow a similar reasoning to the one in \Cref{Section: Proof of Lemma: Sensitivity of HSIC}. Recall that kernels $k$ and $\ell$ are assumed to have non-empty level sets on $\mathbb{Y}$ and $\mathbb{Z}$. Hence, for a given $\epsilon \in (0, \min\{K,L\})$, we may assume that there exist $y_a,y_b \in \mathbb{Y}$ and $z_a,z_b \in \mathbb{Z}$ such that $k(y_a,y_b) = \epsilon_{1}$ and $\ell(z_a,z_b) = \epsilon_{2}$ where $0 \leq \epsilon_1, \epsilon_2 \leq \epsilon$. 

\paragraph{Sensitivity lower bound for $n=2k$ even.}
We begin by considering the case where $n$ is even, and compute the difference between U-statistics of the HSIC based on two specific neighboring datasets given below:
\begin{align*}
	\mathcal{X}_n = \begin{bmatrix}
		y_a & z_a \\
		\vdots & \vdots \\
		y_a & z_a \\
		y_a & z_a \\
		y_a & z_a \\
		y_b & z_b \\
		\vdots & \vdots \\
		y_b & z_b
	\end{bmatrix} \quad \text{and} \quad 
	\tilde{\mathcal{X}}_n = \begin{bmatrix}
		y_a & z_a \\
		\vdots & \vdots \\
		y_a & z_a \\
		y_a & z_a \\
		y_b & z_b \\
		y_b & z_b \\
		\vdots & \vdots \\
		y_b & z_b
	\end{bmatrix}.
\end{align*}
The indices displayed above are $1,\dots,k-2,k-1,k,k+1,\dots 2k$, and it is clear that $d_{\mathrm{ham}}(\mathcal{X}_n, \tilde{\mathcal{X}}_n) = 1$. We further consider the permutation $\bpi$ which permutes only the $k-1$ and $k$ entries, \emph{i.e.}, $\bpi = (1,2,\ldots, k-2,k,k-1,k+1,\ldots,n)$, so that
\begin{align*}
	\mathcal{X}_n^{\bpi} = \begin{bmatrix}
		y_a & z_a \\
		\vdots & \vdots \\
		y_a & z_a \\
		y_a & z_a \\
		y_a & z_a \\
		y_b & z_b \\
		\vdots & \vdots \\
		y_b & z_b
	\end{bmatrix} \quad \text{and} \quad 
	\tilde{\mathcal{X}}_n^{\bpi} = \begin{bmatrix}
		y_a & z_a \\
		\vdots & \vdots \\
		y_a & z_a \\
		y_a & z_b \\
		y_b & z_a \\
		y_b & z_b \\
		\vdots & \vdots \\
		y_b & z_b
	\end{bmatrix}.
\end{align*}
We now compute the U-statistic of the HSIC based on these two permuted datasets. Unfortunately, a direct computation of $U_{\mathrm{HSIC}}$ requires a complicated case-by-case analysis. Instead, we consider the following trick to simplify calculations. First of all, we pretend that $(y_a,z_a)$ and $(y_b,z_b)$ take some specific values, say $(y_a,z_a) = (1,1)$ and $(y_b,z_b) = (0,0)$, and consider indicator kernels $k(y,y') = K \times \mathds{1}(y=y')$ and $\ell(z,z') = L \times \mathds{1}(z = z')$. Under this simplified setting,  $U_{\mathrm{HSIC}}$ is essentially the statistic considered in \citet[][Proposition 1]{kim2023conditional} for multinomial data, which can be computed in a straightforward manner. In addition, we observe that $U_{\mathrm{HSIC}}$ based on the indicator kernels remains the same as $U_{\mathrm{HSIC}}$ based on generic kernels $k$ and $\ell$ up to some quantities tending to zero as $\epsilon_1,\epsilon_2 \rightarrow 0$. Using this trick along with the observation that $k(y,y) = K$ and $\ell(z,z) = L$ for all $y\in \mathbb{Y}$ and $z\in \mathbb{Z}$, it can be seen that for all $k \geq 2$
\begin{align*}
	U_{\mathrm{HSIC}}(\mathcal{X}_{n}^{\bpi}) = \frac{k(k-1)}{(2k-1)(2k-3)} KL +  C_1 \epsilon_1 + C_2 \epsilon_2 + C_3 \epsilon_1 \epsilon_2,
\end{align*}
where $C_1,C_2,C_3$ are constants that only depend on $K,L,n$.  A similar calculation yields 
\begin{align*}
	U_{\mathrm{HSIC}}(\tilde{\mathcal{X}}_n^{\bpi}) ~=~ \frac{(k-2)(k^2-4k+1)}{(k-1)(2k-1)(2k-3)} KL + C_1' \epsilon_1 + C_2' \epsilon_2 + C_3' \epsilon_1 \epsilon_2, 
\end{align*}
where $C_1',C_2',C_3'$ are constants that only depend on $K,L,n$. From these results, we deduce that
\begin{align*}
        & U_{\mathrm{HSIC}}(\mathcal{X}_{n}^{\bpi}) - U_{\mathrm{HSIC}}(\tilde{\mathcal{X}}_n^{\bpi}) \\[.5em]
	= ~ &  \left( \frac{2}{k-1} - \frac{1}{2k-1} - \frac{1}{2k-3}  \right) KL  + (C_1-C_1') \epsilon_1 + (C_2-C_2') \epsilon_2 + (C_3-C_3') \epsilon_1 \epsilon_2 \\[.5em]
        \overset{(\star)}{\to} ~ &  \left( \frac{2}{k-1} - \frac{1}{2k-1} - \frac{1}{2k-3}  \right) KL \geq \frac{KL}{k} 
\end{align*}
where convergence~($\star$) holds as $\epsilon_1,\epsilon_2 \rightarrow 0$ for each fixed $K,L,n$, and the last inequality holds for all $k \geq 2$. Therefore, by letting $k = n/2 \geq 2$, it holds that 
\begin{align*}
	\sup_{\bpi \in \boldsymbol{\Pi}_n} \sup_{\substack{\mathcal{X}_{n},\tilde{\mathcal{X}}_{n}:\\d_{\mathrm{ham}}(\mathcal{X}_{n},\tilde{\mathcal{X}}_{n}) \leq 1}} \bigl| U_{\mathrm{HSIC}}(\mathcal{X}_{n}^{\bpi}) - U_{\mathrm{HSIC}}(\tilde{\mathcal{X}}_{n}^{\bpi}) \bigr| \geq \frac{2KL}{n}.
\end{align*}

\paragraph{Sensitivity lower bound for $n=2k+1$ odd.}
We now consider the same two datasets as in the even case but with an additional row  consisting of $(y_b,z_b)$: 
\begin{align*}
	\mathcal{X}_n = \begin{bmatrix}
		y_a & z_a \\
		\vdots & \vdots \\
		y_a & z_a \\
		y_a & z_a \\
		y_a & z_a \\
		y_b & z_b \\
		\vdots & \vdots \\
		y_b & z_b \\
		y_b & z_b
	\end{bmatrix} \quad \text{and} \quad 
	\tilde{\mathcal{X}}_n = \begin{bmatrix}
		y_a & z_a \\
		\vdots & \vdots \\
		y_a & z_a \\
		y_a & z_a \\
		y_b & z_b \\
		y_b & z_b \\
		\vdots & \vdots \\
		y_b & z_b \\
		y_b & z_b
	\end{bmatrix},
\end{align*}
where the indices displayed are $1,\dots,k-2,k-1,k,k+1,\dots 2k, 2k+1$.
As before, we indeed have $d_{\mathrm{ham}}(\mathcal{X}_n, \tilde{\mathcal{X}}_n) = 1$. Once again, we consider the permutation $\bpi$ which permutes only the $k-1$ and $k$ entries, which leads to
\begin{align*}
	\mathcal{X}_n^{\bpi} = \begin{bmatrix}
		y_a & z_a \\
		\vdots & \vdots \\
		y_a & z_a \\
		y_a & z_a \\
		y_a & z_a \\
		y_b & z_b \\
		\vdots & \vdots \\
		y_b & z_b \\
		y_b & z_b
	\end{bmatrix} \quad \text{and} \quad 
	\tilde{\mathcal{X}}_n^{\bpi} = \begin{bmatrix}
		y_a & z_a \\
		\vdots & \vdots \\
		y_a & z_a \\
		y_a & z_b \\
		y_b & z_a \\
		y_b & z_b \\
		\vdots & \vdots \\
		y_b & z_b\\
		y_b & z_b
	\end{bmatrix}.
\end{align*}
To compute $U_{\mathrm{HSIC}}(\mathcal{X}_{n}^{\bpi})$ and $U_{\mathrm{HSIC}}(\tilde{\mathcal{X}}_{n}^{\bpi})$, we use the trick considered in the even case, and find that 
\begin{align*}
	U_{\mathrm{HSIC}}(\mathcal{X}_{n}^{\bpi}) ~=~ \frac{k(k+1)}{4k^2-1} KL + C_1 \epsilon_1 + C_2 \epsilon_2 + C_3 \epsilon_1 \epsilon_2 
\end{align*}
and
\begin{align*}
	U_{\mathrm{HSIC}}(\tilde{\mathcal{X}}_n^{\bpi}) ~=~ \frac{(k+1)(k-2)(k^2-3k-2)}{k(k-1)(2k-1)(2k+1)} KL + C_1' \epsilon_1 + C_2' \epsilon_2 + C_3' \epsilon_1 \epsilon_2, 
\end{align*}
where $C_1,\ldots,C_3'$ are constants that only depend on $K,L,n$. As a result, we have  
\begin{align*}
	U_{\mathrm{HSIC}}(\mathcal{X}_{n}^{\bpi}) - U_{\mathrm{HSIC}}(\tilde{\mathcal{X}}_n^{\bpi}) &~=~ \frac{4(k^3-2k+1)}{k(k-1)(4k^2-1)}KL + (C_1-C_1') \epsilon_1 + (C_2-C_2') \epsilon_2 + (C_3 - C_3') \epsilon_1 \epsilon_2 \\[.5em]
	& ~\overset{(\star)}{\to}~ \frac{4(k^3-2k+1)}{k(k-1)(4k^2-1)}KL \geq \frac{2}{2k+1} KL,
\end{align*}
where convergence ($\star$) holds as $\epsilon_1,\epsilon_2 \rightarrow 0$ for each fixed $K,L,n$, and the last inequality holds for all $k \geq 2$.  Therefore we conclude that
\begin{align*}
	\sup_{\bpi \in \boldsymbol{\Pi}_n} \sup_{\substack{\mathcal{X}_{n},\tilde{\mathcal{X}}_{n}:\\d_{\mathrm{ham}}(\mathcal{X}_{n},\tilde{\mathcal{X}}_{n}) \leq 1}} \bigl| U_{\mathrm{HSIC}}(\mathcal{X}_{n}^{\bpi}) - U_{\mathrm{HSIC}}(\tilde{\mathcal{X}}_{n}^{\bpi}) \bigr| \geq \frac{2KL}{n},
\end{align*}
for any $n \geq 4$ (either even or odd). This completes the derivation of the lower bound on the sensitivity of the HSIC U-statistic. Together with the upper bound~\eqref{Eq: upper bound for the sensitivity of U-HSIC}, this verifies the claim of \Cref{Lemma: Global sensitivity of U_HSIC}.

\subsection{Proof of \Cref{Theorem: Suboptimality of U_HSIC}}  \label{Section: Proof of Theorem: Suboptimality of U_HSIC}
Let us write the V-statistic of HSIC as
\begin{align*}
	V_{\mathrm{HSIC}} ~=~ & \frac{1}{n^2} \sum_{i,j=1}^n k(Y_i,Y_j) \ell(Z_i,Z_j) + \frac{1}{n^4} \sum_{i_1,i_2,j_1,j_2=1}^n k(Y_{i_1},Y_{j_1}) \ell(Z_{i_2},Z_{j_2})\\[.5em]
	-& \frac{2}{n^3} \sum_{i,j_1,j_2=1}^n k(Y_{i},Y_{j_1}) \ell(Z_{i},Z_{j_2}) \\[.5em]
	=~ & V_{\mathrm{HSIC},a} + V_{\mathrm{HSIC},b} - 2 V_{\mathrm{HSIC},c},
\end{align*}
which is equivalent to the square of $\widehat{\mathrm{HSIC}}$. Using the bounded kernel property, it can be seen that 
\begin{align*}
	\max \Big\{ |U_{\mathrm{HSIC},a} - V_{\mathrm{HSIC},a}|, \ |U_{\mathrm{HSIC},b} - V_{\mathrm{HSIC},b}|, \ |U_{\mathrm{HSIC},c} - V_{\mathrm{HSIC},c}| \Big\} \leq C_1 \frac{KL}{n}
\end{align*}
for some constant $C_1>0$, which leads to 
\begin{align*}
	|U_{\mathrm{HSIC}} - V_{\mathrm{HSIC}}| \leq C_2 \frac{KL}{n}. 
\end{align*}
Again, $C_2$ is some positive number. By treating $K$ and $L$ as some fixed numbers, we also note that \Cref{Lemma: Exponential inequality for the empirical HSIC} yields
\begin{align*}
	V_{\mathrm{HSIC}} = \big\{ \mathrm{HSIC}_{k \otimes \ell}(P_{YZ}) + R_n \big\}^2 \quad \text{where $R_n = O_P(n^{-1/2})$.}
\end{align*}
Having these ingredients in place, we can essentially follow the same proof strategy as for \Cref{Theorem: Suboptimality of U-MMD} and conclude that 
\begin{align*}
	\limsup_{n \rightarrow \infty} \inf_{P \in \mathcal{P}_{\mathrm{HSIC}_{k \otimes \ell}}\!(\rho)} \mE_{P_{YZ}}[\phi_{\mathttt{dpHSIC}}^u]  \leq \alpha.
\end{align*}
This completes the proof of \Cref{Theorem: Suboptimality of U_HSIC}.

\subsection{Proof of \Cref{Theorem: Minimum separation of HSIC U-stat over L2}} \label{Section: Proof of Theorem: Minimum separation of HSIC U-stat over L2}
Let us denote by $C_1,C_2,\ldots$ constants that may depend on $\alpha,\beta,R,M,d_Y,d_Z$. Following the proofs of \Cref{Theorem: General uniform power condition} and \Cref{Theorem: Minimax Separation over L2 for HSIC} along with the sensitivity result for $U_{\mathrm{HSIC}}$ in \Cref{Lemma: Global sensitivity of U_HSIC}, we may arrive at the point where the type II error of $\phi_{\mathttt{dpHSIC}}^u$ is upper bounded as 
\begin{align*}
	\mE[1 - \phi_{\mathttt{dpHSIC}}^u] \leq \mP \biggl( U_{\mathrm{HSIC}} \leq \frac{C_1}{n\sqrt{\lambda_1 \cdots \lambda_{d_Y} \mu_1 \cdots \mu_{d_Z}}} + \frac{C_2}{n \lambda_1 \cdots \lambda_{d_Y} \mu_1 \cdots \mu_{d_Z} \xi_{\varepsilon,\delta}} \biggr) + \frac{15}{16}\beta
\end{align*}
We then use the proofs of \citet[][Theorem 1]{albert2019adaptive} and \citet[][Theorem 2]{albert2019adaptive}, and show that the probability term in the above display is less than or equal to $\beta/16$ once
\begin{align*}
	\|\psi\|_{L_2}^2 ~ \geq ~  C_3 \Biggl\{ \sum_{i=1}^{d_Y} \lambda_i^{2s} + \sum_{i=1}^{d_Z} \mu_i^{2s} & + \frac{1}{n\sqrt{\lambda_1 \cdots \lambda_{d_Y} \mu_1 \cdots \mu_{d_Z}}} + \frac{1}{n \lambda_1 \cdots \lambda_{d_Y} \mu_1 \cdots \mu_{d_Z} \xi_{\varepsilon,\delta}} \Bigg\}.
\end{align*}
This proves \Cref{Theorem: Minimum separation of U-stat over L2}.

\section{Technical Lemmas} \label{Section: Technical Lemmas}
This section collects some technical lemmas used in the main proofs. The first lemma from \cite{kim2021comparing} provides an exponential concentration inequality for the permuted MMD statistic with kernel $k$. 
\begin{lemma}[Theorem 5.1 of \citealt{kim2021comparing}] \label{Lemma: Bobkovs inequality}
	Consider two samples $Y_1,\ldots,Y_{n}$ and $Z_1,\ldots,Z_{m}$ and write $\{X_1,\ldots,X_{N}\} = \{Y_1,\ldots,Y_{n},Z_1,\ldots,Z_{m}\}$. Let us further write $\widetilde{k}(x,y) = k(x,x) + k(y,y) - 2k(x,y)$ and $\gamma = nm/ (n+m)^2$. Define 
	\begin{align*}
		\widehat{\sigma}^2 = \frac{1}{N(N-1)} \sum_{(i,j) \in \boldsymbol{i}_2^N} \widetilde{k}(X_i,X_j).
	\end{align*}
	Then for all $t > 0$, 
	\begin{align*}
		\mP_{\boldsymbol{\pi}} \Biggl[ \sup_{f \in \mathcal{F}_k} \biggl( \frac{1}{n} \sum_{i=1}^{n} f(X_{\pi_i}) - \frac{1}{m} \sum_{i=1}^{m} f(X_{\pi_{i+n}}) \biggr) \geq t + \sqrt{\frac{\widehat{\sigma}^2}{2N \gamma}} \ \Bigg|\ X_1,\ldots,X_N \Biggr] \leq \exp \biggl( - \frac{N\gamma^2t^2}{2\widehat{\sigma}^2} \biggr).
	\end{align*}
\end{lemma}
An exponential tail decay in the above comes at a price --- it requires that $n$ and $m$ are well-balanced in the sense that $N \gamma^2$ is large. Thus it can be used to show that $\widehat{\mathrm{MMD}}(\mathcal{X}_{n+m}^{\bpi}) = o_P(1)$ under the limited regime where $N \gamma^2 \rightarrow \infty$. The lemma below removes this unnecessary constraint at the expense of having a polynomial tail decay.
\begin{lemma}[Markov for the Permuted MMD]  \label{Lemma: Markov for permuted MMD}
	Under the setting of \Cref{Lemma: Bobkovs inequality} with $N = n+m$, we have
	\begin{align*}
		\mP \Biggl[ \sup_{f \in \mathcal{F}_k} \biggl( \frac{1}{n} \sum_{i=1}^{n} f(X_{\pi_i}) - \frac{1}{m} \sum_{i=1}^{m} f(X_{\pi_{i+n}}) \biggr) \geq t \Biggr] \leq \frac{2K}{nt^2} \quad \text{for all $t>0$.}
	\end{align*}
	Hence $\widehat{\mathrm{MMD}}(\mathcal{X}_{n+m}^{\boldsymbol{\pi}}) = o_P(1)$, provided that $Kn^{-1} \rightarrow 0$.
	\begin{proof}
		Markov's inequality together with the V-statistic representation of the empirical MMD yields that for any $t > 0$
		\begin{align*}
			 & \mP \Biggl[ \sup_{f \in \mathcal{F}_k} \biggl( \frac{1}{n} \sum_{i=1}^{n} f(X_{\pi_i}) - \frac{1}{m} \sum_{i=1}^{m} f(X_{\pi_{i+n}}) \biggr) \geq t \Biggr] \\[.5em]			 
			 \leq~ &  \frac{1}{t^2} \biggl\{ \underbrace{\frac{1}{n^2} \sum_{i,j=1}^n \mE[k(X_{\pi_i},X_{\pi_j})] }_{\mathrm{(I)}} + \underbrace{\frac{1}{m^2} \sum_{i,j=1}^m \mE[k(X_{\pi_{i+n}},X_{\pi_{j+n}})]}_{\mathrm{(II)}} - \underbrace{\frac{2}{nm} \sum_{i=1}^n \sum_{j=1}^m \mE[k(X_{\pi_i},X_{\pi_{j+n}})]}_{\mathrm{(III)}} \biggr\}. 
		\end{align*}
		Exploiting the uniformity of permutation $\boldsymbol{\pi}$, we analyze the terms (I), (II) and (III), separately. The analyses of the first two terms (I) and (II) are similar, and we note that 
		\begin{align*}
			\mE[k(X_{\pi_{i}},X_{\pi_{i}})] = \mE[k(X_{\pi_{i+n}},X_{\pi_{i+n}})] ~=~ & \frac{1}{N} \sum_{i=1}^n \mE[k(Y_i,Y_i)] + \frac{1}{N} \sum_{j=1}^m \mE[k(Z_j,Z_j)] \\[.5em]
			= ~& \frac{n}{N} \mE[k(Y_1,Y_1)] + \frac{m}{N} \mE[k(Z_1,Z_1)]
		\end{align*}
		and for $i \neq j$,
		\begin{align*}
			& \mE[k(X_{\pi_{i}},X_{\pi_{j}})] = \mE[k(X_{\pi_{i+n}},X_{\pi_{j+n}})] \\[.5em]
			= ~  & \frac{1}{N(N-1)} \biggl\{ \sum_{(i,j) \in \boldsymbol{i}_2^n} \mE[k(Y_i,Y_j)]  + \sum_{(i,j) \in \boldsymbol{i}_2^m} \mE[k(Z_i,Z_j)]  + 2\sum_{i=1}^n \sum_{j=1}^m \mE[k(Y_i,Z_j)] \bigg\} \\[.5em]
			= ~ & \frac{n(n-1)}{N(N-1)} \mE[k(Y_1,Y_2)]  + \frac{m(m-1)}{N(N-1)} \mE[k(Z_1,Z_2)] +  \frac{2nm}{N(N-1)} \mE[k(Y_1,Z_1)]. 
		\end{align*}
		Therefore 
		\begin{align*}
			\mathrm{(I)} ~=~ & \frac{n}{nN} \mE[k(Y_1,Y_1)] + \frac{m}{nN} \mE[k(Z_1,Z_1)] \\[.5em]
			+ ~ &  \frac{n-1}{n} \biggl\{ \frac{n(n-1)}{N(N-1)} \mE[k(Y_1,Y_2)] + \frac{m(m-1)}{N(N-1)} \mE[k(Z_1,Z_2)] + \frac{2nm}{N(N-1)}\mE[k(Y_1,Z_1)]\biggr\}
		\end{align*}
		and 
		\begin{align*}
			\mathrm{(II)} ~=~ & \frac{n}{mN} \mE[k(Y_1,Y_1)] + \frac{m}{mN} \mE[k(Z_1,Z_1)] \\[.5em]
			+ ~ & \frac{m-1}{m} \biggl\{ \frac{n(n-1)}{N(N-1)} \mE[k(Y_1,Y_2)] + \frac{m(m-1)}{N(N-1)} \mE[k(Z_1,Z_2)] + \frac{2nm}{N(N-1)}\mE[k(Y_1,Z_1)] \biggr\}.
		\end{align*}
		On the other hand, the last term (III) is
		\begin{align*}
			\mathrm{(III)} ~=~ & \frac{2n(n-1)}{N(N-1)} \mE[k(Y_1,Y_2)]  + \frac{2m(m-1)}{N(N-1)} \mE[k(Z_1,Z_2)] +  \frac{4nm}{N(N-1)} \mE[k(Y_1,Z_1)].
		\end{align*}
		Given that kernel $k$ is non-negative and bounded by $K$, we can upper bound $\mathrm{(I)} + \mathrm{(II)} - \mathrm{(III)}$ as 
		\begin{align*}
			\mathrm{(I)} + \mathrm{(II)} - \mathrm{(III)} ~\leq~ & \biggl( \frac{n}{nN} + \frac{m}{nN} + \frac{n}{mN}  + \frac{m}{mN} \biggr) K = \frac{n+m}{nm} K \\[.5em]
			\leq ~ & \frac{2K}{n}.
		\end{align*}
		Hence the result follows.
	\end{proof}
\end{lemma}

The lemma below presents a concentration bound for the permuted HSIC statistic. It is worth noting that obtaining the logarithmic dependence on $\alpha$ is non-trivial, and thus we highlight it as our contribution.

\begin{lemma}[Concentration Inequality for Permuted HSIC] \label{Lemma: Concentration inequality for permuted HSIC}
	Assume that the kernels $k$ and $\ell$ are bounded as $0 \leq k(y,y') \leq K$ and $0 \leq \ell(z,z') \leq L$ for all $y,y' \in \mathbb{Y}$ and $z,z' \in \mathbb{Z}$. Then for any $\alpha \in (0,1)$,
	\begin{align*}
		\mP_{\bpi} \biggl( \widehat{\mathrm{HSIC}}(\mathcal{X}_n^{\bpi})  \geq C \sqrt{\frac{KL}{n}} \max \bigg\{ \log^{1/4} \biggl( \frac{1}{\alpha} \biggr), \log^{1/2}\biggl( \frac{1}{\alpha} \biggr), 1 \bigg\} \,\bigg| \, \mathcal{X}_n \biggr) \leq \alpha,
	\end{align*}
	where $C$ is some positive constant.
	\begin{proof}
		The proof is based on \citet[][Theorem 6.2]{kim2020minimax} that establishes an exponential tail bound for a permuted U-statistic of HSIC. Given $\mathcal{X}_n^{\bpi}$, let us denote the permuted U-statistic by
		\begin{align*}
			U_{\bpi,\mathrm{HSIC}} ~=~ & \frac{(n-2)!}{n!} \sum_{(i,j) \in \boldsymbol{i}_2^n} k(Y_i,Y_j) \ell(Z_{\pi_i},Z_{\pi_j}) + \frac{(n-4)!}{n!} \sum_{(i_1,i_2,j_1,j_2) \in \boldsymbol{i}_4^n} k(Y_{i_1},Y_{j_1}) \ell(Z_{\pi_{i_2}},Z_{\pi_{j_2}}) \\[.5em]
			-& \frac{2(n-3)!}{n!} \sum_{(i,j_1,j_2) \in \boldsymbol{i}_3^n} k(Y_{i},Y_{j_1}) \ell(Z_{\pi_i},Z_{\pi_{j_2}}),
		\end{align*}
		where $\boldsymbol{i}_m^n$ stands for the set of all $m$-tuples drawn from $[n]$ without replacement. \citet[][Theorem 6.2]{kim2020minimax} shows that $U_{\bpi,\mathrm{HSIC}}$ satisfies
		\begin{align*}
			\mP_{\bpi} \bigl( U_{\bpi,\mathrm{HSIC}} \geq t \given \mathcal{X}_n \bigr) \leq  \exp \bigg\{ \! - C_1 \min \biggl( \frac{t^2}{\Sigma_n^2}, \ \frac{t}{\Sigma_n} \biggr)  \bigg\} \quad \text{for all $t>0$,}
		\end{align*}
		where 
		\begin{align*}
			\Sigma_n^2 = \frac{1}{n^2(n-1)^2} \sup_{\bpi \in \boldsymbol{\Pi}_n} \Biggl\{ \sum_{(i,j) \in \boldsymbol{i}_2^n} k^2(Y_i,Y_j) \ell^2(Z_{\pi_i},Z_{\pi_j})  \Bigg\} \leq \frac{1}{n(n-1)} K^2L^2.
		\end{align*}
		The inequality above holds under the assumption that $k$ and $\ell$ are uniformly bounded by $K$ and $L$. Therefore, under the assumption of \Cref{Lemma: Concentration inequality for permuted HSIC}, 
		\begin{align} \label{Eq: concentration inequality for U-stat}
			\mP_{\bpi} \bigl( U_{\bpi,\mathrm{HSIC}} \geq t \given \mathcal{X}_n \bigr) \leq  \exp \bigg\{ \! - C_2 \min \biggl( \frac{n^2t^2}{K^2L^2}, \ \frac{nt}{KL} \biggr) \bigg\} \quad \text{for all $t>0$.}
		\end{align}
		On the other hand, the squared empirical HSIC can be written in the form of a V-statistic given as
		\begin{align*}
			\widehat{\mathrm{HSIC}}^2(\mathcal{X}_n^{\bpi}) ~=~ & \frac{1}{n^2} \sum_{i,j=1}^n k(Y_i,Y_j) \ell(Z_{\pi_i},Z_{\pi_j}) + \frac{1}{n^4} \sum_{i_1,i_2,j_1,j_2=1}^n k(Y_{i_1},Y_{j_1}) \ell(Z_{\pi_{i_2}},Z_{\pi_{j_2}})\\
			-& \frac{2}{n^3} \sum_{i,j_1,j_2=1}^n k(Y_{i},Y_{j_1}) \ell(Z_{\pi_i},Z_{\pi_{j_2}}),
		\end{align*}
		and we note that the difference between the U-statistic and the V-statistic is bounded as
		\begin{align} \nonumber
			& \big| U_{\bpi,\mathrm{HSIC}} -  \widehat{\mathrm{HSIC}}^2(\mathcal{X}_n^{\bpi}) \big| \\[.5em] \nonumber
			\leq~& \bigg| \frac{(n-2)!}{n!} \sum_{(i,j) \in \boldsymbol{i}_2^n} k(Y_i,Y_j) \ell(Z_{\pi_i},Z_{\pi_j}) - \frac{1}{n^2} \sum_{i,j=1}^n k(Y_i,Y_j) \ell(Z_{\pi_i},Z_{\pi_j})   \bigg| \\[.5em] \nonumber
			+~& \bigg| \frac{(n-4)!}{n!} \sum_{(i_1,i_2,j_1,j_2) \in \boldsymbol{i}_4^n} k(Y_{i_1},Y_{j_1}) \ell(Z_{\pi_{i_2}},Z_{\pi_{j_2}})  - \frac{1}{n^4} \sum_{i_1,i_2,j_1,j_2=1}^n k(Y_{i_1},Y_{j_1}) \ell(Z_{\pi_{i_2}},Z_{\pi_{j_2}})  \bigg| \\[.5em] \nonumber
			+~& \bigg| \frac{2(n-3)!}{n!} \sum_{(i,j_1,j_2) \in \boldsymbol{i}_3^n} k(Y_{i},Y_{j_1}) \ell(Z_{\pi_i},Z_{\pi_{j_2}}) - \frac{2}{n^3} \sum_{i,j_1,j_2=1}^n k(Y_{i},Y_{j_1}) \ell(Z_{\pi_i},Z_{\pi_{j_2}})  \bigg| \\[.5em] \label{Eq: temp for the difference between HSIC nad U}
			\leq ~ & C_3 \frac{KL}{n},
		\end{align}
		for any $\bpi$ and $\mathcal{X}_n$. Using the above bound~\eqref{Eq: temp for the difference between HSIC nad U} and letting $t>0$, we convert the concentration inequality for the permuted empirical HSIC into that for the permuted U-statistic as
		\begin{align*}
			\mP_{\bpi} \bigl( \widehat{\mathrm{HSIC}}(\mathcal{X}_n^{\bpi})  \geq t \given \mathcal{X}_n \bigr) ~=~&  \mP_{\bpi} \bigl( \widehat{\mathrm{HSIC}}^2(\mathcal{X}_n^{\bpi})  \geq t^2 \given \mathcal{X}_n \bigr) \\[.5em]
			= ~ &  \mP_{\bpi} \biggl( \widehat{\mathrm{HSIC}}^2(\mathcal{X}_n^{\bpi}) - U_{\bpi,\mathrm{HSIC}}  + U_{\bpi,\mathrm{HSIC}}  \geq t^2 \given \mathcal{X}_n \biggr) \\[.5em] 
			\leq ~ & \mP_{\bpi} \bigl( \big| \widehat{\mathrm{HSIC}}^2(\mathcal{X}_n^{\bpi}) - U_{\bpi,\mathrm{HSIC}} \big|  + U_{\bpi,\mathrm{HSIC}}  \geq t^2 \given \mathcal{X}_n \bigr) \\[.5em]
			\overset{\mathrm{(i)}}{\leq} ~ &   \mP_{\bpi} \biggl( U_{\bpi,\mathrm{HSIC}}  \geq t^2 - C_3 \frac{KL}{n} \given \mathcal{X}_n \biggr) \\[.5em]
			\overset{\mathrm{(ii)}}{\leq} ~ & \mP_{\bpi} \biggl( U_{\bpi,\mathrm{HSIC}}  \geq \frac{t^2}{2} \given \mathcal{X}_n \biggr) \\[.5em]
			\overset{\mathrm{(iii)}}{\leq} ~ & \exp \bigg\{ \! - C_4 \min \biggl( \frac{n^2t^4}{K^2L^2}, \ \frac{nt^2}{KL} \biggr) \bigg\}
		\end{align*}
		where step~(i) uses the bound~\eqref{Eq: temp for the difference between HSIC nad U}, step~(ii) assumes that $t^2 \geq 2 C_3 \frac{KL}{n}$ and step~(iii) follows by concentration inequality~\eqref{Eq: concentration inequality for U-stat}. Setting the last exponential bound to $\alpha$ and solving for $t$ yield the desired result.
	\end{proof}
\end{lemma}

We recall the concentration inequality of the empirical MMD presented by \cite{gretton2012kernel}, which has been used in various places throughout the paper. 

\begin{lemma}[Theorem 7 of \citealt{gretton2012kernel}]  \label{Lemma: Concentration for MMD}
	Assume that the kernel $k$ is bounded as $0 \leq k(x,y) \leq K$ for all $x,y \in \mathbb{S}$. Then for any $t>0$
	\begin{align*}
		\mP\biggl\{ \bigg|  \widehat{\mathrm{MMD}}(\mathcal{X}_{n+m}) - \mathrm{MMD}_k(P,Q) \bigg| > 2 \biggl( \sqrt{\frac{K}{m}} + \sqrt{\frac{K}{n}} \biggr) +t \biggr\} \leq 2\exp \biggl( -\frac{t^2mn}{2K(m+n)} \biggr).
	\end{align*}
\end{lemma}

\bigskip

The following is a counterpart to \Cref{Lemma: Concentration for MMD} for the empirical HSIC. While similar results exist in the literature~\citep[e.g.,][Theorem 3]{gretton2005measuring}, none precisely align with our specific requirements. Hence we opt to present a detailed proof of the following result. 

\begin{lemma}[Exponential Inequality for the Empirical HSIC] \label{Lemma: Exponential inequality for the empirical HSIC}
	Assume that the kernels $k$ and $\ell$ are bounded as $0 \leq k(y,y') \leq K$ and $0 \leq \ell(z,z') \leq L$ for all $y,y' \in \mathbb{Y}$ and $z,z' \in \mathbb{Z}$. Then for any $t \geq 0$,
	\begin{align*}
		\mP \biggl\{ \big| \widehat{\mathrm{HSIC}}(\mathcal{X}_n) - \mathrm{HSIC}_{k \otimes \ell}(P_{YZ}) \big| \geq C_1 \sqrt{\frac{KL}{n}} + t \bigg\} \leq 2 \exp \biggl( - \frac{C_2 t^2 n}{KL} \biggr),
	\end{align*}
	where $C_1,C_2$ are positive constants.
	\begin{proof}
		Throughout the proof, we denote by $C_1,C_2,\ldots$ some generic positive constants. By \Cref{Lemma: Sensitivity of HSIC}, the empirical HSIC has global sensitivity at most $C_1\sqrt{KL}/n$. Thus an application of McDiarmid's inequality~\citep[][Theorem 29]{gretton2012kernel} yields
		\begin{align} \label{Eq: McDiarmid for HSIC}
			\mP \bigl\{ \big| \widehat{\mathrm{HSIC}}(\mathcal{X}_n) - \mE\bigl[ \widehat{\mathrm{HSIC}}(\mathcal{X}_n)\bigr] \big| \geq t \big\} \leq 2 \exp \biggl( - \frac{C_2 t^2 n}{KL} \biggr) \quad \text{for all $t \geq 0$.}
		\end{align}
	    Now we bound the difference between the expectation of $\widehat{\mathrm{HSIC}}(\mathcal{X}_n)$ and the population HSIC as
	    \begin{align*}
	    	& \big| \mE\bigl[ \widehat{\mathrm{HSIC}}(\mathcal{X}_n)\bigr] - \mathrm{HSIC}_{k \otimes \ell}(P_{YZ}) \big| \\[.5em] 
	    	=~ & \bigg|  \mE \biggl[ \sup_{f \in \mathcal{F}_{k \otimes \ell}} \bigg\{ \frac{1}{n} \sum_{i=1}^n f(Y_i,Z_i) - \frac{1}{n^2}\sum_{i,j=1}^n f(Y_i,Z_j) \bigg\} \biggr] -  \sup_{f \in \mathcal{F}_{k \otimes \ell}}  \big\{ \mE_{P_{YZ}}[f(Y,Z)] - \mE_{P_YP_Z}[f(Y,Z)] \big\} \bigg| \\[.5em]
	    	\overset{\mathrm{(i)}}{\leq} ~ &\mE  \bigg|  \sup_{f \in \mathcal{F}_{k \otimes \ell}} \bigg\{ \frac{1}{n} \sum_{i=1}^n f(Y_i,Z_i) - \frac{1}{n^2}\sum_{i,j=1}^n f(Y_i,Z_j) \bigg\}  -  \sup_{f \in \mathcal{F}_{k \otimes \ell}}  \big\{ \mE_{P_{YZ}}[f(Y,Z)] - \mE_{P_YP_Z}[f(Y,Z)] \big\} \bigg|  \\[.5em]
	    	\overset{\mathrm{(ii)}}{\leq} ~ & \mE \bigg| \sup_{f \in \mathcal{F}_{k \otimes \ell}} \bigg\{ \frac{1}{n} \sum_{i=1}^n f(Y_i,Z_i) - \mE_{P_{YZ}}[f(Y,Z)]  - \frac{1}{n^2}\sum_{i,j=1}^n f(Y_i,Z_j) + \mE_{P_YP_Z}[f(Y,Z)]  \bigg\} \bigg| \\[.5em]
	    	\overset{\mathrm{(iii)}}{\leq} ~ & \underbrace{\mE  \sup_{f \in \mathcal{F}_{k \otimes \ell}} \bigg| \frac{1}{n} \sum_{i=1}^n f(Y_i,Z_i) - \mE_{P_{YZ}}[f(Y,Z)] \bigg|}_{\coloneqq \mathrm{(I)}}  + \underbrace{\mE  \sup_{f \in \mathcal{F}_{k \otimes \ell}} \bigg| \frac{1}{n^2}\sum_{i,j=1}^n f(Y_i,Z_j) - \mE_{P_YP_Z}[f(Y,Z)] \bigg|}_{\coloneqq \mathrm{(II)}},
	    \end{align*}
		where step~(i) uses Jensen's inequality, step~(ii) uses the reverse triangle inequality and step~(iii) follows by the triangle inequality. Let $\{(\tilde{Y}_i,\tilde{Z}_i)\}_{i=1}^n$ be i.i.d.~copies of $(Y_1,Z_1)$, and $\{\epsilon_i\}_{i=1}^n$ be i.i.d.~Rademacher random variables. Then Jensen's inequality in conjunction with symmetrization yields
		\begin{align} \nonumber
			\mathrm{(I)} ~\leq~ & \mE  \sup_{f \in \mathcal{F}_{k \otimes \ell}} \bigg| \frac{1}{n} \sum_{i=1}^n \epsilon_i \big\{ f(Y_i,Z_i) - f(\tilde{Y}_i, \tilde{Z}_i) \big\} \bigg| \\[.5em]
			\leq ~ & 2 \mE \sup_{f \in \mathcal{F}_{k \otimes \ell}} \bigg| \frac{1}{n} \sum_{i=1}^n \epsilon_i  f(Y_i,Z_i)  \bigg|	~\leq ~ 2 \sqrt{\frac{KL}{n}}, \label{Eq: Eq: inequality 1}
		\end{align}
		where the last inequality is due to \citet[][Lemma 22]{bartlett2002rademacher}. For the second term, 
		\begin{align*}
			\mathrm{(II)} ~ \overset{(\mathrm{i})}{\leq} ~ & \frac{1}{n} \sum_{j=1}^n \mE  \sup_{f \in \mathcal{F}_{k \otimes \ell}} \bigg| \frac{1}{n} \sum_{i=1}^n  f(Y_i,Z_j) - \mE_{P_YP_Z}[f(Y,Z)] \bigg| \\[.5em] 
			= ~ & \mE  \sup_{f \in \mathcal{F}_{k \otimes \ell}} \bigg| \frac{1}{n} \sum_{i=2}^n  f(Y_i,Z_1) + \frac{f(Y_1,Z_1)}{n} - \mE_{P_YP_Z}[f(Y,Z)] \bigg| \\[.5em] 
			\overset{(\mathrm{ii})}{\leq} ~ & \mE  \sup_{f \in \mathcal{F}_{k \otimes \ell}} \bigg| \frac{1}{n} \sum_{i=2}^n  f(Y_i,Z_1) - \mE_{P_YP_Z}[f(Y,Z)] \bigg|  + \frac{1}{n}  \mE  \sup_{f \in \mathcal{F}_{k \otimes \ell}} |f(Y_1,Z_1)| \\[.5em]
			\overset{(\mathrm{iii})}{\leq} ~ & \mE  \sup_{f \in \mathcal{F}_{k \otimes \ell}} \bigg| \frac{1}{n-1} \sum_{i=2}^n  f(Y_i,Z_1) - \mE_{P_YP_Z}[f(Y,Z)] \bigg|  \\[.5em]
			& + \frac{1}{n(n-1)} \sum_{i=2}^n \mE  \sup_{f \in \mathcal{F}_{k \otimes \ell}} |f(Y_i,Z_1)|  + \frac{1}{n} \mE \sup_{f \in \mathcal{F}_{k \otimes \ell}} |f(Y_1,Z_1)|,
		\end{align*}
		where step~(i) follows by Jensen's inequality and step~(ii) and step~(iii) use the triangle inequality. Let $\mathcal{H}_{k \otimes \ell}$ be the reproducing kernel Hilbert space with the product kernel $k \otimes \ell$. Then by the reproducing property and the Cauchy--Schwarz inequality, we have 
		\begin{align*}
			\sup_{f \in \mathcal{F}_{k \otimes \ell}} |f(y,z)| = \sup_{f \in \mathcal{F}_{k \otimes \ell}} |\langle k(\cdot,y) \ell(\cdot,z), f \rangle_{\mathcal{H}_{k \otimes \ell}}| \leq \sqrt{k(y,y) \ell(z,z)} \sup_{f \in \mathcal{F}_{k \otimes \ell}} \|f\|_{\mathcal{H}_{k \otimes \ell}} \leq \sqrt{KL}
		\end{align*}
		for all $y \in \mathbb{Y},z \in \mathbb{Z}$. Using this inequality and letting $\{(\tilde{Y}_i,\tilde{Z}_1)\}_{i=2}^n$ be i.i.d.~copies of $\{(Y_i,Z_1)\}_{i=2}^n$ and recalling that $\{\epsilon_i\}_{i=1}^n$ are i.i.d.~Rademacher random variables, the second term can be further bounded by
		\begin{align} \nonumber
			\mathrm{(II)} ~ \leq ~ & \mE  \sup_{f \in \mathcal{F}_{k \otimes \ell}} \bigg| \frac{1}{n-1} \sum_{i=2}^n  f(Y_i,Z_1) - \mE_{P_YP_Z}[f(Y,Z)] \bigg|  + \frac{2\sqrt{KL}}{n} \\[.5em]
			= ~ &  \mE  \sup_{f \in \mathcal{F}_{k \otimes \ell}} \bigg| \frac{1}{n-1} \sum_{i=2}^n  f(Y_i,Z_1) - \frac{1}{n-1} \sum_{i=2}^n  \mE[f(\tilde{Y}_i, \tilde{Z}_1)] \bigg|  + \frac{2\sqrt{KL}}{n} \\[.5em] \nonumber
			\leq ~ &  \mE  \sup_{f \in \mathcal{F}_{k \otimes \ell}} \bigg| \frac{1}{n-1} \sum_{i=2}^n  f(Y_i,Z_1) - f(\tilde{Y}_i, \tilde{Z}_1) \bigg|  + \frac{2\sqrt{KL}}{n} \\[.5em]
			= ~ & \mE\biggl\{ \mE \biggl[ \sup_{f \in \mathcal{F}_{k \otimes \ell}} \bigg| \frac{1}{n-1} \sum_{i=2}^n \epsilon_i \{f(Y_i,Z_1) - f(\tilde{Y}_i, \tilde{Z}_1)\} \bigg| \,\Big|\, Z_1,\tilde{Z}_1 \biggr]\biggr\} + \frac{2\sqrt{KL}}{n} \\[.5em] \nonumber
			\leq~ & 2 \mE\biggl\{ \mE \biggl[ \sup_{f \in \mathcal{F}_{k \otimes \ell}}\bigg| \frac{1}{n-1} \sum_{i=2}^n \epsilon_i f(Y_i,Z_1) \bigg| \,\Big|\, Z_1 \biggr] \biggr\} + \frac{2\sqrt{KL}}{n} \\[.5em] \label{Eq: inequality 2}
			\leq ~ & 2 \sqrt{\frac{KL}{n-1}} + \frac{2\sqrt{KL}}{n},
		\end{align}
		where the last inequality is due to \citet[][Lemma 22]{bartlett2002rademacher}. Putting the inequalities~\eqref{Eq: Eq: inequality 1} and \eqref{Eq: inequality 2} together, we have 
		\begin{align*}
			\big| \mE\bigl[ \widehat{\mathrm{HSIC}}(\mathcal{X}_n)\bigr] - \mathrm{HSIC}_{k \otimes \ell}(P_{YZ}) \big|  \leq C_3 \sqrt{\frac{KL}{n}}.  
		\end{align*}
		The above inequality along with McDiarmid's inequality~\eqref{Eq: McDiarmid for HSIC} yields
		\begin{align*}
			\mP \biggl\{ \big| \widehat{\mathrm{HSIC}}(\mathcal{X}_n) - \mathrm{HSIC}_{k \otimes \ell}(P_{YZ})  \big| \geq t + C_3 \sqrt{\frac{KL}{n}} \bigg\} \leq 2 \exp \biggl( - \frac{C_2 t^2 n}{KL} \biggr) \quad \text{for all $t \geq 0$.}
		\end{align*}
		Hence we complete the proof of \Cref{Lemma: Exponential inequality for the empirical HSIC}. 
	\end{proof}
\end{lemma}

The following lemma is folklore in the literature~\citep[\emph{e.g.},][Lemma 1]{romano2005exact}, where we recall that $\floor{x}$ denotes the largest integer smaller than or equal to $x$. We provide a proof for completeness. 

\begin{lemma}[Permutation $p$-value] \label{Lemma: permutation p-value}
	Suppose that $X_1,\ldots,X_n,X_{n+1}$ are exchangeable random variables. Then for any $\alpha \in [0,1]$, it holds that 
	\begin{align*}
		\mathbb{P} \biggl( \frac{1}{n+1} \bigg\{\sum_{i=1}^n \mathds{1}(X_{n+1} \leq X_i)  + 1 \bigg\} \leq \alpha \biggr) \leq \frac{\floor{(n+1)\alpha}}{n+1} \leq \alpha.
	\end{align*}
 	Suppose further that $X_1,\ldots,X_{n+1}$ are all distinct with probability one. Then
	\begin{align*}
		\mathbb{P} \biggl( \frac{1}{n+1} \bigg\{\sum_{i=1}^n \mathds{1}(X_{n+1} \leq X_i)  + 1 \bigg\} \leq \alpha \biggr) = \frac{\floor{(n+1)\alpha}}{n+1}.
	\end{align*}
	\begin{proof}
		Let $X_{(1)} \leq X_{(2)} \leq \ldots \leq X_{(n+1)}$ be the order statistics of $X_1,\ldots,X_{n+1}$. Then we have a series of the identities:
		\begin{align*}
			& \frac{1}{n+1} \bigg\{\sum_{i=1}^n \mathds{1}(X_{n+1} \leq X_i)  + 1 \bigg\} = \frac{1}{n+1} \sum_{i=1}^{n+1} \mathds{1}(X_{n+1} \leq X_i)  \leq \alpha \\[.5em]
			\Longleftrightarrow ~ & \sum_{i=1}^{n+1} \mathds{1}(X_{n+1} \leq X_i)  \leq \floor{(n+1)\alpha} \\[.5em]
			\Longleftrightarrow ~  & \sum_{i=1}^{n+1} \mathds{1}(X_i < X_{n+1}) \geq (n+1) - \floor{(n+1)\alpha}\coloneqq  k \\[.5em]
			\Longleftrightarrow ~ & X_{n+1} > X_{(k)}.
		\end{align*}
		Now by the exchangeability condition, we have
		\begin{align*}
			\mathbb{P}(X_{n+1} > X_{(k)}) ~=~& \mathbb{E}\bigg[ \frac{1}{n+1} \sum_{i=1}^{n+1} \mathds{1}(X_{i} > X_{(k)}) \bigg].
		\end{align*}
		On the other hand, by the definition of $X_{(k)}$,
		\begin{align*}
			\frac{1}{n+1} \sum_{i=1}^{n+1} \mathds{1}(X_{i} > X_{(k)}) \leq \frac{n+1 - k}{n+1} = \frac{\floor{(n+1)\alpha}}{n+1}.
		\end{align*}
		Hence the first result follows. When $X_1,\ldots,X_{n+1}$ are all distinct, observe
		\begin{align*}
			\sum_{i=1}^{n+1} \mathds{1}(X_{i} > X_{(k)}) = n+1-k.
		\end{align*}
		Thus
		\begin{align*}
			\mathbb{E}\bigg[ \frac{1}{n+1} \sum_{i=1}^{n+1} \mathds{1}(X_{i} > X_{(k)}) \bigg] = \frac{n+1-k}{n+1} =  \frac{\floor{(n+1)\alpha}}{n+1}.
		\end{align*}
	\end{proof} 
\end{lemma} 

The permutation test is often slightly conservative due to its discrete nature. The following randomization trick is well-known in the literature~\citep[e.g.,][Chapter 15]{lehmann2005testing}, which modifies the permutation test to have type I error exactly equal to $\alpha$.

\begin{lemma}[Randomized Tests] \label{Lemma: Randomized tests}
	For given $\alpha \in (0,1)$ and $\gamma \in (0,\alpha]$, consider a test function $\phi$ such that $\mP_{H_0}(\phi = 1) = \gamma \leq \alpha$. Letting $U$ be a uniform random variable on $[0,1]$ independent of $\phi$, define a randomized test $\phi_{\mathrm{rand}}$ as
	\begin{align*}
		\phi_{\mathrm{rand}} = \phi + (1 - \phi) \times \mathds{1}\biggl(U \leq \frac{\alpha - \gamma}{1-\gamma} \biggr).
	\end{align*}
	The randomized test $\phi_{\mathrm{rand}}$ has type I error exactly equal to $\alpha$, and power never worse than that of $\phi$, \emph{i.e.}, $\mP_{H_0}(\phi_{\mathrm{rand}} = 1) = \alpha$ and $\mP_{H_1}(\phi_{\mathrm{rand}} = 1) \geq \mP_{H_1}(\phi = 1)$.
	\begin{proof}
		Type I error control is immediate given that 
		\begin{align*}
			\mE_{H_0}[\phi_{\mathrm{rand}}] ~=~ & \mE_{H_0}[\phi] + \mE_{H_0}[1 - \phi] \times \mP_{H_0}\biggl(U \leq \frac{\alpha - \gamma}{1-\gamma} \biggr) \\[.5em]
			= ~ & \gamma + (1-\gamma) \times \frac{\alpha - \gamma}{1 - \gamma} = \alpha. 
		\end{align*}
		The second claim about the power is also immediate given that $\phi_{\mathrm{rand}} \geq \phi$.
	\end{proof}
\end{lemma}

The permutation test can be expressed in terms of the permutation $p$-value as well as the quantile of the permutation distribution, which can be justified by the following lemma. 

\begin{lemma}[Quantile Representation] \label{Lemma: Quantile representation}
	For any $\alpha \in [0,1]$ and dataset $\{X_1,\ldots,X_{n+1}\}$, we have the identity 
	\begin{align*}
		\mathds{1}\biggl( \frac{1}{n+1} \bigg\{\sum_{i=1}^n \mathds{1}(X_{n+1} \leq X_i)  + 1 \bigg\} \leq \alpha  \biggr) = \mathds{1} \bigl(X_{n+1} > q_{1-\alpha}\bigr),
	\end{align*}
	where $q_{1-\alpha}$ is the $1-\alpha$ quantile of $\{X_1,\ldots,X_{n+1}\}$ given as
	\begin{align*}
		q_{1-\alpha}   = & \inf \biggl\{ t \in \mathbb{R} : \frac{1}{n+1} \sum_{i=1}^{n+1} \mathds{1}(X_i \leq t) \geq 1- \alpha  \biggr\}. 
	\end{align*}
	Moreover, by letting $X_{(\ceil{(1-\alpha)(n+1)})}$ denote the $\ceil{(1-\alpha)(n+1)}$th order statistic of $X_1,\ldots,X_{n+1}$ and $X_{(0)} = -\infty$, we have $q_{1-\alpha} = X_{(\ceil{(1-\alpha)(n+1)})}$. 	
	\begin{proof}
		The second claim follows by noting that 
		\begin{align} \nonumber
			q_{1-\alpha}  ~=~ &  \inf \biggl\{ t \in \mathbb{R} : \frac{1}{n+1} \sum_{i=1}^{n+1} \mathds{1}(X_i \leq t) \geq 1- \alpha  \biggr\} \\[.5em]  \label{eq: quantile definition}
			 = ~ & \inf \biggl\{ t \in \mathbb{R} : \frac{1}{n+1} \sum_{i=1}^{n+1} \mathds{1}(X_i <  t) \geq 1- \alpha \biggr\} \\[.5em] \nonumber
			  = ~ & \inf \biggl\{ t \in \mathbb{R} : \sum_{i=1}^{n+1} \mathds{1}(X_i <  t) \geq (1- \alpha)(n+1) \biggr\} \\[.5em] 
			 = ~ & X_{(\ceil{(1-\alpha)(n+1)})}. \nonumber
		\end{align}
		For the first claim, denote $G(x) = \sum_{i=1}^{n+1} \mathds{1}(X_i < x)$, which is a left-continuous step function. We then have
		\begin{align*}
			\mathds{1} \biggl(\frac{1}{n+1} \bigg\{\sum_{i=1}^n \mathds{1}(X_{n+1} \leq X_i)  + 1 \bigg\} \leq \alpha \biggr)  ~=~&  \mathds{1} \bigl( G(X_{n+1}) \geq (1- \alpha)(n+1) \bigr) \\[.5em]
		 	\leq ~ & \mathds{1}(X_{n+1} > q_{1-\alpha})
		\end{align*}  
		where the last step follows since $G$ is a left-continuous step function. In more detail, suppose that $G(X_{n+1}) \geq (1 - \alpha)(n+1)$. Since $G$ is a left-continuous step function, there exists a small constant $\epsilon >0$ such that $G(X_{n+1} - \epsilon) = G(X_{n+1}) \geq (1-\alpha)(n+1)$. Therefore, $X_{n+1}$ cannot be the $1-\alpha$ quantile and should be greater than $q_{1-\alpha}$. 
		
		Moreover, the event $X_{n+1} > q_{1-\alpha}$ implies that $G(X_{n+1}) \geq (1-\alpha)(n+1)$ by the definition of $q_{1-\alpha}$ as in expression~\eqref{eq: quantile definition}. Hence we conclude that 
		\begin{align*}
			\mathds{1} \bigl( G(X_{n+1}) \geq (1- \alpha)(n+1) \bigr) ~ = ~ \mathds{1}(X_{n+1} > q_{1-\alpha}),
		\end{align*}
		and the first claim follows. 
	\end{proof}
\end{lemma}

The next lemma plays a crucial role in proving the validity of the differentially private permutation test in \Cref{Algorithm: DP permutation test}.

\begin{lemma}[Alternative Expression] \label{Lemma: Alternative expression}
	Given $\alpha \in (0,1)$ and $n \geq 1$, set
	\begin{align*}
		\alpha_\star = \max\bigg\{ \biggl(\frac{n+1}{n} \alpha - \frac{1}{n}\biggr), \, 0 \bigg\}.
	\end{align*}
	Then for dataset $\{X_1,\ldots,X_{n+1}\}$, we have the identity
	\begin{align*}
		\mathds{1}(X_{n+1} > q_{1-\alpha}) =  \mathds{1}(X_{n+1} > r_{1-\alpha_\star})\mathds{1}\biggl(\alpha \geq \frac{1}{n+1}\biggr),
	\end{align*}
	where $q_{1-\alpha}$ and $r_{1-\alpha_\star}$ are the $1-\alpha$ quantile of $\{X_i\}_{i=1}^{n+1}$ and the $1-\alpha_\star$ quantile of $\{X_i\}_{i=1}^{n}$, respectively.   
	\begin{proof}
		For $\alpha \in (0,\frac{1}{n+1})$, $q_{1-\alpha}$ becomes the maximum of $X_1,\ldots,X_{n+1}$ for which $\mathds{1}(X_{n+1} > q_{1-\alpha}) = 0$. Similarly, it becomes  
		\begin{align*}
			\mathds{1}(X_{n+1} > r_{1-\alpha_\star})\mathds{1}\biggl(\alpha \geq \frac{1}{n+1}\biggr) = 0.
		\end{align*}
		 Hence we only need to verify the identity under $\alpha \geq \frac{1}{n+1}$. In what follows, we assume $\alpha \geq \frac{1}{n+1}$ and show that $\mathds{1}(X_{n+1} > q_{1-\alpha}) = \mathds{1}(X_{n+1} > r_{1-\alpha_\star})$.
		 
		 Remark that the $1-\alpha$ quantile of $\{X_1 + c,\ldots,X_{n+1} + c\}$ is the same as the $1-\alpha$ quantile of $\{X_1,\ldots,X_{n+1}\}$ plus $c$ for any $c \in \mathbb{R}$. Using this location-shift property of quantiles, observe that
		\begin{align*}
			\mathds{1}(X_{n+1} > q_{1-\alpha}) ~ = ~ & \mathds{1} \biggl( 0 > \inf \biggl\{ x \in \mathbb{R}: \frac{1}{n+1} \sum_{i=1}^{n+1} \mathds{1}(X_i - X_{n+1} \leq x) \geq 1 - \alpha \biggr\} \biggr) \\[.5em]
			= ~ & \mathds{1} \biggl( 0 > \inf \biggl\{ x \in \mathbb{R}: \frac{1}{n} \sum_{i=1}^n \mathds{1}(X_i - X_{n+1} \leq x) \geq \frac{n+1}{n}(1 - \alpha) - \frac{\mathds{1}(0 \leq x)}{n} \biggr\} \biggr) \\[.5em] 
			\geq ~ & \mathds{1} \biggl( 0 > \inf \biggl\{ x \in \mathbb{R} : \frac{1}{n} \sum_{i=1}^n \mathds{1}(X_i - X_{n+1} \leq x) \geq \frac{n+1}{n}(1 - \alpha) \biggr\} \biggr) \\[.5em] 
			\overset{(\dagger)}{=} ~ & \mathds{1} \biggl( 0 > \inf \biggl\{ x \in \mathbb{R} : \frac{1}{n} \sum_{i=1}^n \mathds{1}(X_i - X_{n+1} \leq x) \geq \min \biggl[ \frac{n+1}{n}(1 - \alpha), \, 1 \biggr] \biggr\} \biggr) \\[.5em]
			= ~ & \mathds{1}(X_{n+1} > r_{1-\alpha_\star}).
		\end{align*}
		where the equality~($\dagger$) holds under $\alpha \geq \frac{1}{n+1}$. Therefore it holds that $\mathds{1}(X_{n+1} > q_{1-\alpha}) \geq \mathds{1}(X_{n+1} > r_{1-\alpha_\star})$.
		
		Next we prove the other direction $\mathds{1}(X_{n+1} > q_{1-\alpha}) \leq \mathds{1}(X_{n+1} > r_{1-\alpha_\star})$. Note that the infimum in the definition of a quantile can be replaced by the minimum so that 
		\begin{align*}
			\frac{1}{n+1} \sum_{i=1}^{n+1} \mathds{1}(X_i - X_{n+1} \leq \underbrace{q_{1-\alpha} - X_{n+1}}_{\coloneqq  \tilde{q}_{1-\alpha}}) \geq 1 - \alpha.
		\end{align*}
		Having this in mind, for $\alpha \geq \frac{1}{n+1}$, assume $X_{n+1} > q_{1-\alpha}$ (equivalently $0 > \tilde{q}_{1-\alpha}$) under which it holds that 
		\begin{align*}
			& \frac{1}{n} \sum_{i=1}^n \mathds{1}(X_i - X_{n+1} \leq \tilde{q}_{1-\alpha}) \geq \frac{n+1}{n}(1 - \alpha) - \frac{\mathds{1}(0 \leq \tilde{q}_{1-\alpha})}{n}  = \frac{n+1}{n}(1 - \alpha) \\[.5em]
			\Longleftrightarrow ~ &  \frac{1}{n} \sum_{i=1}^n \mathds{1}(X_i - X_{n+1} \leq \tilde{q}_{1-\alpha}) \geq  \min \biggl\{ \frac{n+1}{n}(1 - \alpha), \, 1 \biggr\}.
		\end{align*}
		This implies that $\tilde{q}_{1-\alpha} \geq r_{1-\alpha_\star} - X_{n+1} \coloneqq  \tilde{r}_{1-\alpha_\star}$ by the definition of $r_{1-\alpha_\star}$. Consequently, $\mathds{1}(X_{n+1} > q_{1-\alpha}) \leq \mathds{1}(\tilde{q}_{1-\alpha} \geq \tilde{r}_{1-\alpha_\star})$, which further implies that 
		\begin{align*}
			\mathds{1}(X_{n+1} > q_{1-\alpha}) \leq \mathds{1}(\tilde{q}_{1-\alpha} \geq \tilde{r}_{1-\alpha_\star}) \mathds{1}(X_{n+1} > q_{1-\alpha}) \leq \mathds{1}(0 > \tilde{r}_{1-\alpha_\star})= \mathds{1}(X_{n+1} > r_{1-\alpha_\star}).
		\end{align*}
		Thus we conclude that $\mathds{1}(X_{n+1} > q_{1-\alpha}) = \mathds{1}(X_{n+1} > r_{1-\alpha_\star})$ for $\alpha \geq \frac{1}{n+1}$ as well. This completes the proof of \Cref{Lemma: Alternative expression}.
	\end{proof}
\end{lemma}

The next result is concerned with the global sensitivity of the quantiles.

\begin{lemma}[Sensitivity of Quantiles] \label{Lemma: Sensitivity of quantiles}
	Suppose that the test statistic $T$ has the global sensitivity at most $\Delta_T$ as in~\eqref{Eq: global sensitivity}. Let us denote by 
	$r_{1-\alpha} \bigl(\mathcal{X}_n;\{\bpi_i, \zeta_i \}_{i=1}^B \bigr)$ the $1-\alpha$ quantile of $\{M_i\}_{i=1}^B$ where $M_i = T(\mathcal{X}_n^{\bpi_i}) + 2 \Delta_T \xi_{\varepsilon,\delta}^{-1} \zeta_i$ with $\xi_{\varepsilon,\delta}= \varepsilon+\log(1/(1-\delta))$ and $\zeta_i \iid  \mathsf{Laplace}(0,1)$ for $i \in [B]$. Then for any $\alpha \in [0,1)$, the global sensitivity of the $1-\alpha$ quantile satisfies
	\begin{align*}
		\sup_{\substack{\mathcal{X}_n,\tilde{\mathcal{X}}_n:\\d_{\mathrm{ham}}(\mathcal{X}_n,\tilde{\mathcal{X}}_n) \leq 1}} \big| r_{1-\alpha} \bigl(\mathcal{X}_n;\{\bpi_i, \zeta_i \}_{i=1}^B \bigr) - r_{1-\alpha} \bigl(\tilde{\mathcal{X}}_n; \{ \bpi_i, \zeta_i \}_{i=1}^B \bigr) \big| \leq \Delta_T,
	\end{align*}
	for any permutations $\bpi_1,\ldots,\bpi_B$ and any $\zeta_1,\ldots,\zeta_B\iid  \mathsf{Laplace}(0,1)$.
	\begin{proof}
		Let $\tilde{\mathcal{X}}_n$ be a neighboring dataset of $\mathcal{X}_n$ where $\tilde{\mathcal{X}}_n$ and $\mathcal{X}_n $ differ only in their $k$th component for some $k \in [n]$. Denote the permuted test statistics computed on $\tilde{\mathcal{X}}_n^{\bpi_1},\ldots,\tilde{\mathcal{X}}_n^{\bpi_B}$ by $\tilde{T}_1,\ldots,\tilde{T}_B$. For simplicity, we write $r_{1-\alpha} = r_{1-\alpha} \bigl(\mathcal{X}_n; \{ \bpi_i, \zeta_i \}_{i=1}^B \bigr)$ and $\tilde{r}_{1-\alpha} = r_{1-\alpha} \bigl(\tilde{\mathcal{X}}_n; \{ \bpi_i, \zeta_i \}_{i=1}^B \bigr)$. Having this notation, first note that $T_i \geq \tilde{T}_i - \Delta_T$ for all $i \in [B]$ under the assumption of \eqref{Eq: global sensitivity} and thus  
		\begin{align*}
			{1-\alpha} \leq \frac{1}{B} \sum_{i=1}^B \mathds{1}(T_i + 2 \Delta_T \xi_{\varepsilon,\delta}^{-1} \zeta_i \leq r_{1-\alpha}) \leq \frac{1}{B} \sum_{i=1}^B \mathds{1}(\tilde{T}_i + 2 \Delta_T \xi_{\varepsilon,\delta}^{-1} \zeta_i \leq r_{1-\alpha} + \Delta_T),
		\end{align*}
		which implies that $\tilde{r}_{1-\alpha} \leq r_{1-\alpha} + \Delta_T$. 
		
		Next we argue that $\tilde{r}_{1-\alpha} \geq r_{1-\alpha} - \Delta_T$. For this direction, let $\epsilon > 0$ be an arbitrary constant. Then by the definition of $r_{1-\alpha}$ and $T_i \leq \tilde{T}_i + \Delta_T$ for all $i \in [B]$,
		\begin{align*}
			{1-\alpha} > \frac{1}{B} \sum_{i=1}^B \mathds{1}(T_i + 2 \Delta_T \xi_{\varepsilon,\delta}^{-1} \zeta_i \leq r_{1-\alpha} - \epsilon) \geq \frac{1}{B} \sum_{i=1}^B \mathds{1}(\tilde{T}_i + 2 \Delta_T \xi_{\varepsilon,\delta}^{-1} \zeta_i \leq r_{1-\alpha} - \epsilon - \Delta_T). 
		\end{align*}
		Hence $\tilde{r}_{1-\alpha} > r_{1-\alpha} - \epsilon - \Delta_T$. Since $\epsilon$ is arbitrary, we conclude $\tilde{r}_{1-\alpha} \geq r_{1-\alpha} - \Delta_T$. In summary, we have established that $|r_{1-\alpha} - \tilde{r}_{1-\alpha}| \leq \Delta_T$, which holds for any $k \in [n]$, any permutations $\bpi_1,\ldots,\bpi_B$ and any $\zeta_1,\ldots,\zeta_B\iid  \mathsf{Laplace}(0,1)$. Therefore, the desired claim follows.
	\end{proof}
\end{lemma}

The Laplace mechanism introduces a perturbed statistic by adding Laplace noise. Instead of directly studying quantiles of this perturbed statistic, it is often easier to analyze the quantiles of the original statistic and the quantiles of the added Laplace noise separately. These results can then be combined using the following quantile inequality.

\begin{lemma}[Quantile Inequality] \label{Lemma: Quantile inequality}
	Let $q_{1 - \alpha}^{X+Y}$ be the $1-\alpha$ quantile of $X+Y$. Similarly, let $q_{1-\alpha/2}^X$ and $q_{1-\alpha/2}^Y$ be the $1-\alpha/2$ quantile of $X$ and $Y$, respectively. Then $q_{1- \alpha}^{X+Y} \leq q_{1-\alpha/2}^X + q_{1-\alpha/2}^Y$. 
	\begin{proof}
		Note that if $\{X+Y > q_{1-\alpha/2}^X + q_{1-\alpha/2}^Y\}$, then at least one of the events $\{X > q_{1-\alpha/2}^X\}$ and $\{Y > q_{1-\alpha/2}^Y\}$ should occur (otherwise contradiction). This together with the union bound yields
		\begin{align*}
			\mP\bigl(X+Y > q_{1-\alpha/2}^X + q_{1-\alpha/2}^Y \bigr) ~\leq~ & \mP\bigl( X > q_{1-\alpha/2}^X \bigr) + \mP\bigl( Y > q_{1-\alpha/2}^Y \bigr)  \\[.5em]
			\leq ~ & \alpha,
		\end{align*}
		where the last inequality follows by the definition of quantile. The above implies that the $1-\alpha$ quantile of $X+Y$ is less than or equal to $q_{1-\alpha/2}^X + q_{1-\alpha/2}^Y$ and therefore the claim follows. 
	\end{proof}
\end{lemma}

The following lemma provides a high-probability upper bound for the quantile of an empirical distribution in relation to a theoretical distribution. In the context of the permutation test, this result enables us to eliminate the randomness from Monte Carlo simulations during the power analysis. 
\begin{lemma}[Quantile Approximation] \label{Lemma: Quantile Approximation}
	Let $\{X,X_1,\ldots,X_B\}$ be i.i.d.~random variables conditional on a sigma field $\mathcal{G}$, and let $X_0$ be another random variable defined on the same probability space as $X$. For $\alpha \in (0,1)$, define $\hat{q}_{1-\alpha}$ and $q_{1-\alpha}$ as 
	\begin{align*}
		& \hat{q}_{1-\alpha}  := \inf\bigg\{ x \in \mathbb{R}:  \frac{1}{B+1} \sum_{i=0}^B \mathds{1}(X_i \leq x) \geq 1-\alpha \bigg\} \quad \text{and} \\
		& q_{1-\alpha} :=  \inf\bigg\{ x \in \mathbb{R}:  \mP(X \leq x \given \mathcal{G}) \geq 1-\alpha \bigg\}.
	\end{align*}
	Assume that the conditional distribution of $X$ given $\mathcal{G}$ has no atom. Then for any $\beta \in (0,1)$, $\delta >0$, $\gamma \in \bigl(0, \alpha/(1+\delta)\bigr)$, and $B$ satisfying 
	\begin{align*}
		B \geq \max \bigg\{ \frac{2+\delta}{\delta^2 \gamma}\log\biggl(\frac{1}{\beta}\biggr), \ \frac{1-\alpha}{\alpha - \gamma(1+\delta)} \bigg\},
	\end{align*}
	it holds that $\hat{q}_{1-\alpha} \leq q_{1-\gamma}$ with probability at least $1-\beta$. 
	\begin{remark} \label{Remark: quantile approximation} 
		\leavevmode \normalfont 
		\begin{itemize}
			\item For simplicity, if we set $\gamma = \alpha/6$ and $\delta = 2$, then \Cref{Lemma: Quantile Approximation} ensures $\hat{q}_{1-\alpha} \leq q_{1-\alpha/6}$ with probability at least $1-\beta$ provided that 
			\begin{align} \label{Eq: condition for B}
				B \geq 2\alpha^{-1} \max \big\{3 \log(1/\beta), \, 1- \alpha \big\},
			\end{align}
			or more simply $B \geq 6 \alpha^{-1} \log(1/\beta)$ if $\beta < e^{-1/3} \approx 0.716$. The constant factors in the above condition can be further improved by choosing a smaller value of $\gamma$ and a larger value of $\delta$. This improvement, on the other hand, comes at the cost of inflating the constant factor in the minimum separation of a test as we are ultimately concerned with a smaller value of the significance level. 
			\item We also remark that this condition offers an improvement over the previous requirements for $B$ imposed, such as, in \cite{kim2020minimax,schrab2021mmd}. It particularly improves the dependence of $\alpha$ from $\alpha^{-2}$ to $\alpha^{-1}$. The previous approach relies on the additive form of concentration inequalities such as the Dvoretzky--Kiefer--Wolfowitz (DKW) inequality~\citep{massart1990}. On the other hand, we use a multiplicative Chernoff bound, which is the key to this improvement. Another technique to obtain the $\alpha^{-1}$ factor in the number of $B$ is 
			\item \citet[][Lemma 6]{domingo2023compress} derive a high-probability upper bound for $B$ using the property of order statistics. This technique, completely different from ours, also yields a factor of $\alpha^{-1}$ for the number of $B$. In particular, their condition requires that $B \geq \alpha^{-1} - 1$, and moreover, the coreset count in their Theorem 1 must be larger than some function of $B,\alpha,\beta$. 
			\item \Cref{Lemma: Quantile Approximation} assumes that the conditional distribution of $W$ given $\mathcal{G}$ has no atom. In our context, this assumption is easily met as we add continuous noise to a test statistic through the Laplace mechanism. In a more general context, such a guarantee can be achieved, for example, through jittering.
		\end{itemize}
	\end{remark}
	\begin{proof}
		 As mentioned in \Cref{Remark: quantile approximation}, a main tool for \Cref{Lemma: Quantile Approximation} is a multiplicative Chernoff bound\footnote{\url{https://en.wikipedia.org/wiki/Chernoff_bound}}, which we recall here: \emph{Suppose that $Y_1,\ldots,Y_n$ are independent random variables taking values in $\{0,1\}$. Let $Y = \sum_{i=1}^n Y_i$ denote their sum and let $\mu = \mE[Y]$ denote the expected value of the sum. Then for any $\delta \geq 0$,}
		\begin{align} \label{Eq: multiplicative chernoff}
			\mP\bigl(Y \geq (1+\delta) \mu)\bigr) \leq e^{-\delta^2 \mu / (2+\delta)}.
		\end{align}
		Given this multiplicative bound, note that $\hat{q}_{1-\alpha}$ is equivalent to 
		\begin{align*}
		\hat{q}_{1-\alpha}  = \inf\bigg\{ x \in \mathbb{R}:  \frac{1}{B+1} \sum_{i=0}^B \mathds{1}(X_i > x) \leq \alpha \bigg\}.
		\end{align*}
		For a fixed $x$ (which will be specified later on) and $\delta >0$, observe that
		\begin{align} \nonumber
		\frac{1}{B+1} \sum_{i=0}^B \mathds{1}(X_i > x) ~\leq~ & \frac{1}{B+1} + \frac{1}{B+1} \sum_{i=1}^B \mathds{1}(X_i > x) \\[.5em] \label{Equation: inequality 1}
		= ~ & \frac{1}{B+1}  + \frac{B(1+\delta)}{B+1} \biggl[ \frac{1}{1+\delta} \times \frac{1}{B} \sum_{i=1}^B \mathds{1}(X_i > x) \biggr].
		\end{align}
		Note that the multiplicative Chernoff~\eqref{Eq: multiplicative chernoff} yields
		\begin{align*}
		\mP\biggl(\frac{1}{B} \sum_{i=1}^B \mathds{1}(X_i > x) \geq (1+\delta) \mP(X > x \given \mathcal{G})  \,\bigg|\, \mathcal{G} \biggr) \leq \exp \biggl(-\frac{\delta^2}{2+\delta} B \mP(X > x \given \mathcal{G}) \biggr) \overset{\text{set}}{\leq} \beta.
		\end{align*}
		Therefore with probability at least $1-\beta$, the following event holds
		\begin{align*}
		\frac{1}{B} \sum_{i=1}^B \mathds{1}(X_i > x) < (1+\delta) \mP(X > x \given \mathcal{G})  ~\Longleftrightarrow~ \frac{1}{1+\delta} \times \frac{1}{B} \sum_{i=1}^B \mathds{1}(X_i > x) <  \mP(X > x\given \mathcal{G}). 
		\end{align*}
		In particular, this condition holds when 
		\begin{align}  \label{Eq: condition 1}
		B \geq \frac{2+\delta}{\delta^2} \frac{1}{\mP(X > x \given \mathcal{G})} \log\biggl(\frac{1}{\beta}\biggr).
		\end{align}
		Now we take $x=q_{1-\gamma}$, which gives that $ \mP(X > x \given \mathcal{G}) = \gamma$ since the conditional distribution of $X$ given $\mathcal{G}$ has no atom. Then, continuing from inequality~\eqref{Equation: inequality 1}, we have under the event stated above that 
		\begin{align*}
		\frac{1}{B+1} \sum_{i=0}^B \mathds{1}(X_i > q_{1-\gamma}) ~\leq~ &   \frac{1}{B+1}  + \frac{B(1+\delta)}{B+1} \mP(X > q_{1-\gamma} \given \mathcal{G}) \\[.5em]
		= ~ &  \frac{1}{B+1}  + \frac{B(1+\delta)}{B+1} \gamma.
		\end{align*}
		Now assume that 
		\begin{align} \label{Eq: condition 2}
		\frac{1}{B+1}  + \frac{B(1+\delta)}{B+1} \gamma \leq \alpha.
		\end{align}
		Then with probability at least $1-\beta$, the definition of $\hat{q}_{1-\alpha}$ yields
		\begin{align*}
		\hat{q}_{1-\alpha} \leq q_{1-\gamma}.
		\end{align*}
		There are two conditions that ensures the above inequality, namely condition~\eqref{Eq: condition 1} and condition~\eqref{Eq: condition 2}. Under the condition that $\alpha > \gamma(1+\delta)$, condition~\eqref{Eq: condition 2} becomes equivalent to
		\begin{align*}
		B \geq \frac{1-\alpha}{\alpha - \gamma(1+\delta)}.
		\end{align*}
		Hence, a sufficient condition for $B$ that ensures $\hat{q}_{1-\alpha} \leq q_{1-\gamma}$ is given as
		\begin{align*}
		B \geq \max \bigg\{ \frac{2+\delta}{\delta^2 \gamma}\log\biggl(\frac{1}{\beta}\biggr), \ \frac{1-\alpha}{\alpha - \gamma(1+\delta)} \bigg\}.
		\end{align*}
		This completes the proof of \Cref{Lemma: Quantile Approximation}. 
	\end{proof}
\end{lemma}

The following lemma computes the difference between the V-statistic and U-statistic of the squared HSIC, which is useful in the proof of \Cref{Theorem: Minimax Separation over L2 for HSIC}.
\begin{lemma}[Difference between $V_{\mathrm{HSIC}}$ and $U_{\mathrm{HSIC}}$] \label{Lemma: Difference between V and U}
	Recall $\hat{\mathrm{HSIC}}^2$ given in \eqref{Eq: closed form HSIC} and $U_{\mathrm{HSIC}}$ given in \eqref{Eq: U-HSIC}. Suppose that kernels satisfy $k(y,y) = K$ and $\ell(z,z) = L$ for all $y \in \mathbb{Y}$ and $z \in \mathbb{Z}$. Then the difference between $\hat{\mathrm{HSIC}}^2$ and $U_{\mathrm{HSIC}}$ is computed as follows:
	\begin{align*}
		\hat{\mathrm{HSIC}}^2 - U_{\mathrm{HSIC}} = D_1 + D_2,
	\end{align*}
	where
	\begin{align*}
		D_1 ~=~ & \frac{n-1}{n^2}KL -\frac{L}{n^3} \sum_{(i,j) \in \mathbf{i}_2^n} k(Y_i,Y_j) - \frac{K}{n^3} \sum_{(i,j) \in \mathbf{i}_2^n} \ell(Z_i,Z_j) \\[.5em]
		D_2 ~=~ & - \frac{3n^2-4n+2}{(n-1)n^4}\sum_{(i,j) \in \mathbf{i}_2^n} k(Y_i,Y_j) \ell(Z_i,Z_j)+ \frac{2(5n^2-8n+4)}{n^4(n-1)(n-2) }\sum_{(i,j_1,j_2) \in \mathbf{i}_3^n} k(Y_i,Y_{j_1}) \ell(Z_i,Z_{j_2})\\[.5em]
		& - \frac{6n^2-11n+6}{n^4(n-1)(n-2)(n-3)} \sum_{(i_1,i_2,j_1,j_2) \in \mathbf{i}_4^n} k(Y_{i_1},Y_{i_2}) \ell(Z_{j_1}, Z_{j_2}).
	\end{align*}
	\emph{\textbf{Remark}.} \emph{The term $D_1$ is invariant to any permutation of either $Y$ values or $Z$ values.}
	\begin{proof}
		Write the squared empirical HSIC in \eqref{Eq: closed form HSIC} as $V_{\mathrm{HSIC}}$ and the U-statistic of the squared HSIC as $U_{\mathrm{HSIC}}$. Note that $V_{\mathrm{HSIC}} = V_a + V_b - 2V_c$ where
		\begin{align*}
			& V_a ~=~ \frac{1}{n^2}\sum_{i,j=1}^n k(Y_i,Y_j) \ell(Z_i,Z_j), \quad  V_b ~=~ \frac{1}{n^4}\sum_{i_1,i_2,j_1,j_2=1}^n k(Y_{i_1},Y_{i_2}) \ell(Z_{j_1}, Z_{j_2}), \\[.5em]
			& V_c ~=~ \frac{1}{n^3}\sum_{i,j_1,j_2=1}^n k(Y_i,Y_{j_1}) \ell(Z_i,Z_{j_2}),
		\end{align*}
		and each term can be expressed as
		\begin{align*}
			V_a ~=~ & \frac{KL}{n} + \frac{1}{n^2} \sum_{(i,j) \in \mathbf{i}_2^n} k(Y_i,Y_j) \ell(Z_i,Z_j), \\[.5em]
			V_b ~=~ & \frac{KL}{n^3}  + \frac{2K}{n^4} \sum_{(i,j) \in \mathbf{i}_2^n} \ell(Z_i,Z_j) + \frac{2L}{n^4} \sum_{(i,j) \in \mathbf{i}_2^n} k(Y_i,Y_j)  + \frac{KLn(n-1)}{n^4}  \\[.5em]
			+ ~ & \frac{K(n-2)}{n^4} \sum_{(i,j) \in \mathbf{i}_2^n} \ell(Z_i,Z_j) + \frac{L(n-2)}{n^4} \sum_{(i,j) \in \mathbf{i}_2^n} k(Y_i,Y_j) \\[.5em]
			+ ~ & \frac{2}{n^4} \sum_{(i,j) \in \mathbf{i}_2^n} k(Y_i,Y_j) \ell(Z_i,Z_j) +  \frac{4}{n^4} \sum_{(i,j_1,j_2) \in \mathbf{i}_3^n} k(Y_i,Y_{j_1}) \ell(Z_i,Z_{j_2}) \\[.5em] 
			+ ~ &  \frac{1}{n^4} \sum_{(i_1,i_2,j_1,j_2) \in \mathbf{i}_3^n} k(Y_{i_1},Y_{i_2}) \ell(Z_{j_1}, Z_{j_2}) \quad \text{and} \\[.5em]
			V_c ~=~ &\frac{KL}{n^2} + \frac{K}{n^3} \sum_{(i,j) \in \mathbf{i}_2^n} \ell(Z_i,Z_j) + \frac{L}{n^3} \sum_{(i,j) \in \mathbf{i}_2^n} k(Y_i,Y_j) \\[.5em]
			+ ~ & \frac{1}{n^3} \sum_{(i,j) \in \mathbf{i}_2^n} k(Y_i,Y_j) \ell(Z_i,Z_j) + \frac{1}{n^3} \sum_{(i,j_1,j_2) \in \mathbf{i}_3^n} k(Y_i,Y_{j_1}) \ell(Z_i,Z_{j_2}). 
		\end{align*}
		Note also that $U_{\mathrm{HSIC}}$ can be written as $U_{\mathrm{HSIC}} = U_a + U_b - 2U_c$ where
		\begin{align*}
			U_{a} ~=~ & \frac{1}{n(n-1)} \sum_{(i,j) \in \mathbf{i}_2^n} k(Y_i,Y_j) \ell(Z_i,Z_j),  \\[.5em]
			U_{b} ~=~ &  \frac{(n-4)!}{n!} \sum_{(i_1,i_2,j_1,j_2) \in \mathbf{i}_4^n} k(Y_{i_1},Y_{j_1}) \ell(Z_{i_2},Z_{j_2}) \quad \text{and} \\[.5em]
			U_{c} ~=~ & \frac{1}{n(n-1)(n-2)} \sum_{(i,j_1,j_2) \in \mathbf{i}_3^n} k(Y_{i},Y_{j_1}) \ell(Z_{i},Z_{j_2}).
		\end{align*}
		The result follows by calculating $V_a - U_a$, $V_b - U_b$ and $-2 V_c +2 U_c$, and simplifying their sum via algebra. 
	\end{proof}
\end{lemma}

\end{document}